\documentclass[11pt]{report}

\usepackage{amsmath}
\usepackage{tikz}
\usetikzlibrary{automata,arrows.meta,positioning,}
\usepackage{pgfplots}
\usepackage{pgf}
\usepackage{mathrsfs}
\usepackage{hyperref}
\usepackage{graphicx}
\usepackage[linesnumbered,ruled,vlined]{algorithm2e}
\usepackage[toc,page]{appendix}
\usepackage{amssymb,amsthm}
\usepackage{authblk}
\usepackage{blkarray}
\usepackage{appendix}
\usepackage{wasysym}
\usepackage{float}
\usepackage{mathtools}
\usepackage{pythonhighlight}
\usepackage{setspace}

\usepackage[left=1.5in,right=1in,top=1in,bottom=1in,]{geometry}
\usepackage{enumitem}
\usepackage{kantlipsum}
\usepackage{bm}
\usepackage{mathtools}
\usepackage{caption}
\usepackage{amsbsy}
\usepackage{mdframed}

\newmdenv[
topline=false,
bottomline=false,
rightline=false,
skipabove=\topsep,
skipbelow=\topsep,
leftmargin=-10pt,
rightmargin=-10pt,
innertopmargin=0pt,
innerbottommargin=0pt
]{sidelines}

\theoremstyle{definition}
\newtheorem{theorem}{Theorem}[chapter]

\newtheorem{example}[theorem]{Example}
\newtheorem{proposition}[theorem]{Proposition}
\newtheorem{corollary}[theorem]{Corollary}
\newtheorem{lemma}[theorem]{Lemma}
\theoremstyle{definition}
\newtheorem{definition}[theorem]{Definition}

\setlist[enumerate]{label*=\arabic*.}

\title{Algorithms for Reachability Problems on Stochastic Markov Reward Models}
\author{Irfan Muhammad}
\affil{School of Computer Science, University of Birmingham}
\begin{document}
    
    \begin{titlepage}
   	 \huge
   		\begin{flushleft}
   	   	 
   	   	 ALGORITHMS FOR REACHABILITY PROBLEMS ON STOCHASTIC MARKOV REWARD MODELS    
   		\end{flushleft}
      		 \begin{flushleft}

      			 by
      		 \end{flushleft}

  			 \begin{flushleft}
  				 
  				 IRFAN MUHAMMAD
  			 \end{flushleft}
  		 \begin{center}
			   			 
  			 \vspace{3cm}
  			 A thesis submitted to the University of Birmingham for the degree of\\
  			 \vspace{0.5cm}
  			 DOCTOR OF PHILOSOPHY
  		 \end{center}
  		 \begin{flushright}
  			 
  			 \vspace{3cm}
  			 School of Computer Science\\
  			 College of Engineering and Physical Sciences\\
  			 University of Birmingham\\
  			 September 2020\\
  		 \end{flushright}

    \end{titlepage}

    \begin{abstract}
   	 
   	 \paragraph{} Probabilistic model-checking is a field which seeks to automate the formal analysis of probabilistic models such as Markov chains. In this thesis, we study and develop the stochastic Markov reward model (sMRM) which extends the Markov chain with rewards as random variables. The model recently being introduced, does not have much in the way of techniques and algorithms for their analysis. The purpose of this study is to derive such algorithms that are both scalable and accurate.
   	 
   	 \paragraph{} Additionally, we derive the necessary theory for probabilistic model-checking of sMRMs against existing temporal logics such as PRCTL. We present the equations for computing \textit{first-passage reward densities}, \textit{expected value problems}, and other \textit{reachability problems}. Our focus however is on finding strictly numerical solutions for \textit{first-passage reward densities}. We solve for these by firstly adapting known direct linear algebra algorithms such as Gaussian elimination, and iterative methods such as the power method, Jacobi and Gauss-Seidel. We provide solutions for both discrete-reward sMRMs, where all rewards discrete (lattice) random variables. And also for continuous-reward sMRMs, where all rewards are strictly continuous random variables, but not necessarily having continuous probability density functions (pdfs). Our solutions involve the use of fast Fourier transform (FFT) for faster computation, and we adapted existing quadrature rules for convolution to gain more accurate solutions, rules such as the trapezoid rule, Simpson's rule or Romberg's method.
   	 
   	 \paragraph{} In the discrete-reward setting, existing solutions are either derived by hands, or a combination of graph-reduction algorithms and symbolically solving them via computer algebra systems. The symbolic approach is not scalable, and for this we present strictly numerical but relatively more scalable algorithms. We found each - direct and iterative - capable of solving problems with larger state spaces. The best performer was the power method, owed partially to its simplicity, leading to easier vectorization of its implementation. Whilst, the Gauss-Seidel method was shown to converge with fewer iterations, it was slower due to costs of deconvolution. The Gaussian Elimination algorithm performed poorly relative to these.
   	 
   	 \paragraph{} In the continuous-reward setting, existing solutions are adaptable from literature on semi-Markov processes. However, it appears that other algorithms should still be researched for the cases where rewards have discontinuous pdfs. The algorithm we have developed has the ability to resolve such a case, albeit the solution does not appear as scalable as the discrete-reward setting.
   	 
    \end{abstract}
    
    \renewcommand{\abstractname}{Acknowledgements}
    \begin{abstract}
   	 
   	 \paragraph{} Much thanks goes to my supervisor Professor David Parker, who opened for me the doors to academia. Again thanks goes to him for his guidance, and his reviewing of my work. I also thank my RSMG group members, Professor Peter Hancox and Dr. Shan He, both of whom motivated me and trusted that I could accomplish this. And also thanks are to be given to Associate Professor Nick Hawes, and Professor Jeremy Wyatt, my previous co-supervisors who still benefited me significantly in the little of what they said.
   	 
   	 \paragraph{} I also would like to gratefully acknowledge my family and friends who helped me pull through during times difficult.
    
   	 \paragraph{} Finally, I am gratefully to DARPA for funding me this opportunity through the PRINCESS project (contract FA8750-16-C-0045), via the DARPA BRASS programme.
    \end{abstract}

    \newcommand{\deconv}{\ \pentagon\ }
    \tableofcontents
    
    
    \chapter{Introduction} Probabilistic models have been utilized effectively in many areas of research and have numerous industrial applications. By abstracting real world systems, and making some sound assumptions of how they behave, we can use such models to represent a simplification of a natural phenomenon, or an artificial or engineered one. Doing so grants us a less expensive representation, one that we can study for insights as to how the true system actually behaves. The understanding we obtain is useful for many reasons. If the system being modelled is artificial, we can use the knowledge to determine whether the true system meets safety specifications. If it is a real world phenomena, we can use it to make predictions and forecasts, or perhaps to develop or prove scientific theories.

\paragraph{} A \textit{probabilistic} model is one which models systems predictable only through relative ratios of outcomes, i.e. systems best characterized as \textit{random}. These models can then be used to determine the likelihood of certain critical events, such as failure rates of an industrial facility \cite{warr2014numerical}, the likelihood of cancer re-emission or death \cite{weiss1965semi} or the survivability of coronary patients \cite{warr2012introduction}. Other applications of these models could be for weather-forecasting, speech-recognition, in the field of computational biology and robotics. A sufficiently flexible model that is well developed has the capacity to model complex real world problems. A flexible model is one that imposes fewer requirements on real-world systems, to be representable.

\paragraph{} The field of \textbf{probabilistic verification} comprises various formal techniques to investigate behaviours of probabilistic models. These models range from simple structures such as Markov chains and probabilistic Petri nets to complicated models such as those written with probabilistic programming languages. A contribution to this field would be for example designing methods for investigating behaviours of models not previously analysable. Or alternatively, presenting new models. One new class of model in literature is the stochastic Markov reward model (sMRM) \cite{bacci2019model}, a model which captures the accumulation of rewards over a random (Markov) process. These rewards are random themselves, and are accumulated as independent increments. Its development is recent, and we do not find much work directly on them. In this thesis, we develop algorithms for investigating a set of questions concerning them. We consider both the discrete-reward and continuous-reward variants of these models, and we focus mostly on a class of problems known as \textit{reachability problems}.

\paragraph{}  The techniques investigated in probabilistic verification are various, however one popular family of methods is \textit{probabilistic model checking}. The method allows us to investigate both \textit{qualitative} and \textit{quantitative } properties of particular probabilistic models, via a combination of \textit{specification languages} of which are usually \textit{temporal logics} - to write what we intend to investigate of a model, and a set of model checking algorithms that consists of logic translators, model transformers, and numerical algorithms such as Gaussian elimination, or the Gauss-Seidel methods and their variants. We focus on designing scalable practical numerical algorithms for sMRMs, implementing and experimenting with them.

\section{Motivation} Let us consider for the time being application in robotics. The ability to automate robots such that they succeed at their tasks without risk of failure opens the door for many an industrial application. For example, they may be used to investigate the sea floors for undetonated mines left from previous wars or used in nuclear facilities where humans are at risk from radiation poisoning. Or mundane tasks within factories or farming.

\paragraph{} We motivate our research into a class of probabilistic models - stochastic Markov reward models - with the following problem:  A robot leaves a charging station to perform a set of tasks in its surrounding environment. The tasks consist of performing inspections on physical structures. The order in which it performs its tasks is determined by some planner; it visits a set of objects and inspects them before returning to recharge. Each movement or inspection consumes energy and time of the robot, both being random variables.

\paragraph{} Concerning such a robot, we can ask some questions about its behaviour: 1) What is the probability that the robot runs out of energy during any of its missions? 2) Is it possible for the robot to complete the mission under some time bound $t$? 3) What is the likelihood that it fails to complete more than $x$ number of tasks due to running out of energy?

\paragraph{} Problems like these can be modelled with discrete-time Markov chains (DTMCs), a model which allows us to approximate all possible behaviour of the robot with a finite set of configurations known as \textit{states}, and which captures its behaviour evolution through time by using \textit{transitions} between these states. For example, the mission can be modelled as a DTMC where the state of a robot is two variables: location, and whether the robot is inspecting or not. The state transitions represent movements of the robot from one location to another, and they also represent the starting, continuation and ending of an inspection task. To represent the energy and time costs, we have the choice of integrating these costs into the state space of the DTMC, that is the state of the robot will also include the amount of energy the robot currently has and how much time it has spent away from the recharging station. For complicated models however, this approach would blow up the memory requirements exponentially, thus leading us to considering better techniques.

\paragraph{} Another solution is to use a related model known as the Markov reward model (MRM). In these models, \textit{rewards} (or synonymously \textit{costs}) can be represented in a manner whereby they are not integrated into the state space. Thus, the memory requirements are generally lesser than using DTMCs. An MRM is essentially a DTMC with a \textit{reward structure} connected to it. Both the DTMC and the reward structure can be represented as individual matrices of equal size. We can represent time/energy costs more efficiently with these models. However, the MRM is more restrictive than a regular DTMC allowing only non-negative rewards/costs over transitions. For example with MRMs, we cannot transition to states which would reduce time spent on the mission or increase the robot's current energy amount (by recharging for example). This is not a problem we focus on at this time. Note that modelling both time/energy costs leads to a bi-dimensional (bi-variate) reward structure, thus two matrices (of equal size) are required for it instead of one.   

\paragraph{} A drawback with MRMs is that it only deals with deterministic rewards, i.e. it only allows us to model the energy between movements as fixed amounts. A more recent model has been introduced into literature at the time of this work \cite{bacci2019model}, which allows these rewards to be random variables. These models are referred to as stochastic Markov reward models (sMRM). It is this class of model that is the focus of this thesis  and what it builds upon.

\paragraph{} However, a challenge with sMRMs is that multivariate random variable rewards (or random vector rewards) are not generally composed of mutually independent univariate random variables, unlike MRMs with deterministic rewards (i.e. degenerate random variables), where they are always independent of each other. Unfortunately, this blow up in complexity whilst important to us, we have no general solution for. The focus therefore is solely on univariate reward structures, or random vectors with mutually independent components.

\paragraph{} One might propose another set of models instead of sMRMs, known as (transition-based) hidden Markov models (HMMs). This model extends DTMCs, by allowing there to be a further set of states called \textit{hidden states}. Hence the original set of states of the underlying DTMC are known as \textit{observable} states.  The complete state space of the model is the Cartesian product $\{\textit{hidden states} \times \textit{observable states}\}$. A required property of the model is that the likelihood of being in a hidden state at any particular time step is solely dependent on the transition made at the previous time step (between two observable states) by the HMM, and not on any other (previous) step.\footnote{In a state-based HMM, the likelihood of being in a hidden state at a particular time instance solely depends on the state the HMM is in at that time instance.} One might try to model question 1 above using an HMM by first using the observable states to represent the location of the robot and whether the robot is inspecting or not, just like with the DTMC. Then, introducing two hidden states for the model: the robot still has energy, and the robot has run out of energy. When the robot transits between observable states, e.g. moves location, it has a likelihood of running out of energy, or still having energy, i.e. it has a likelihood of ending up in one of the two hidden states. If the process is to be modelled as a HMM, the likelihood of being in a hidden state (e.g. having no energy) at any time must be modelled as a function of the most recent transition the robot took. However, this is theoretically unsound as the likelihood of having no energy is dependent on all previous transitions  and not just the most recent. As for a formal explanation, we delay this until after we have defined sMRMs properly (see Section \ref{subsec:alt_models}).

\paragraph{} If the robot does not yet have a plan of action (or mission), and we have a set of actions we can program the robot to do at any given state, if we were to ask a fourth question - for any state the robot is in, what is the best action the robot can take, in terms of preventing the robot from running out of energy and increasing the likelihood of completing a set of inspections? Then, a more general model or framework is required to answer this. One way to handle this is by extending previous models (DTMCs, MRMs, or sMRMs) with an \textit{action set} which is used to annotate each state with \textit{actions} we want the robot to be able to make. State transitions are now dependent on the action chosen by the robot at a given state. These modified models can then be studied to find the best action to make at any state given a particular problem. Extending DTMCs this way gives us the class of models known as Markov decision processes (MDPs). In the case of MRMs or sMRMs, we can call their extensions Markov reward decision processes (MRDPs) or stochastic MRDPs (sMRDPs). However, we have explained that the HMM is not a valid choice of model for our problem, and therefore its extension, called the partially observable Markov decision process (POMDP) cannot be considered.  However, an alternative manner of solving this is by using a black-box (or off-the-shelf) optimization algorithm that produces a plan of action. The black-box algorithm chooses a set of actions for the robot, one for each state. This set of actions then translates to a set of behaviours for the robot which can be modelled as previously via a DTMC, MRM or sMRM. The models are then investigated to see if for example the robot will run out of energy. If it is proven that the plan is sufficient (e.g. the robot will not run out of energy with high probability), the robot can be configured to run the plan. If not, the black-box algorithm will attempt to find a better plan through a new set of actions. Therefore an iteration between plan generation (via the black box algorithm), and plan verification (via the sMRM or other) is needed to determine a good set of actions that would be sufficient. In this thesis, we will not focus on planning problems. Rather the focus is strictly on sMRMs. However, as just explained, sMRMs can still be used as part of an optimization procedure to generate usable plans for a robot.

\paragraph{} The robot problem serves to introduce the class of models being investigated in this thesis. It is not however our main focus, but it does clarify how these models can be used.

\section{Objectives and contributions}

\paragraph{} The stochastic Markov reward models (sMRMs) were recently introduced by \cite{bacci2019model}.  The authors would analyse their behaviour via simulation techniques. We seek to derive numerical algorithms for reachability problems defined over sMRMs, that are scalable and accurate. We avoid simulation techniques which offer generally statistical guarantees on the accuracy of the result, in favour of numerical approaches which can give formal guarantees. Additionally, simulation techniques are considered slow when accuracy is of concern.

\paragraph{} Our contributions are as follows:
\begin{enumerate}
	\item We tie together model-checking with sMRMs using a temporal logic called PRCTL - \textit{probabilistic reward control tree logic}. \cite{andova2003discrete}. We lay a foundation for sMRM theory and present theoretical solutions for several problems.
	\item We present algorithms for our main problem, the computation of the \textit{the passage-time reward mass functions} (described later), over the discrete-reward sMRMs that are fast and scalable. The algorithms involve adapting existing solutions for solving systems of linear equations: The Gaussian elimination algorithm, the power, Jacobi and Gauss Seidel methods.
	\item Algorithms for the continuous-reward sMRMs were also derived, that are fast but only slightly scalable. We discovered new quadrature rules for convolution that are amenable for use within the sMRMs that can provide more accurate answers. The solution we provide is perhaps a more general solution relative to some algorithms in literature which are adaptable for use, allowing us to solve problems that other algorithms cannot do without further development.
\end{enumerate}

\section{Thesis layout}
The remainder of this thesis is organized as follows. The following chapter discusses relevant literature to the sMRM models. It discusses probabilistic model checking in the context of such models, as well as existing numerical algorithms utilizable for solving some problems of these models. Chapter 3 presents the necessary background for the thesis. It discusses Markov models, temporal logics and the algorithms for resolving properties over Markov models. Chapter 4 presents the theoretical foundations for sMRMs, introducing the system of convolution equations and provides proofs for the basis of our work. Chapter 5 begins with deriving an exact solution for discrete-reward sMRM problems via the Gaussian elimination algorithm. Chapter 6 introduces iterative methods for solving discrete-reward sMRMs via the power, Jacobi and Gauss-Seidel method. These methods are more scalable relative to the direct Gaussian elimination method. Chapter 7 concerns an attempt to resolve continuous-reward sMRMs. Finally, Chapter 8 concludes this thesis and presents direction for future work.

\section{Computer details}
\label{sec:computer details}
For the experiments found in this thesis, two different computers were used. Additionally, the algorithms we present are implemented in python. We made heavy use of the following python libraries: 1) numpy \cite{oliphant2006guide,van2011numpy} for their numerical algorithms and matrix operations. 2) fftw, a python wrapper for the FFTW library \cite{FFTW} for FFT operations. 3) PaCal \cite{korzen2014pacal} - a probabilistic arithmetic calculator. This was used heavily as a benchmarking tool for the algorithms we develop.

\paragraph{} The specs. of the two computers are as follows.
\begin{enumerate}
	\item \textit{Computer 1:} A laptop with 7.7GB of RAM, and an Intel Core i7-5500U CPU (@ 2.40GHz x 4). Here, we are using the OPENBLAS \cite{openblas} package as a back-end for numpy, but multi-threading was turned off.
	\item \textit{Computer 2:} A desktop PC  with 31.9GB of RAM, and an AMD Ryzen 5 3600 6-Core CPU (@ 3.593GHz x 6).
\end{enumerate}

\paragraph{} \textit{Computer 1} is our default computer, and unless mentioned otherwise, can be assumed to be the computer used for a particular experiment in this thesis.

    \chapter{Related work}

\section{Introduction} In this chapter, we discuss models related to the sMRM such as the Markov reward model (MRM) with deterministic rewards, and semi-Markov processes (SMP). We also delve into probabilistic model checking with temporal logics, an area concerned with solving general problems over probabilistic models such as DTMCs, MRMs, and others. Doing so will give us a reference as to what questions may be asked of sMRM behaviour. Then, we also discuss our focus and compare the algorithms derived here with possible algorithms that exist in literature.

\section{Probabilistic model checking} Probabilistic model checking is a wide field covering formal methods used to verify, prove or investigate the behaviour of probabilistic models. Models can be as specific as software code, or an abstraction of real-word systems. These models can be developed by hand, or generated automatically by software. When it is the latter, work to prevent memory overflow include using symbolic representations \cite{bryant1986graph}, and SAT methods \cite{biere1999symbolic}. Once these models have been generated, there are a host of algorithms to analyse them, with results from fields including logic theory, automata theory, numerical methods and graph theory. The set of techniques we focus on are discussed in the book \cite{principles_book}, involving state spaces and temporal logics (or more generally specification languages).

\paragraph{} \textit{Probabilistic model checking with temporal logics} consists of a combination of three separate components, the first two being: A probabilistic model, and a specification language (e.g. a temporal logic). The model is used to represent the phenomena or system at hand, whilst the language allows us to write questions (synonymous to \textit{statements} or \textit{properties} or \textit{problems}) we want to investigate concerning the model. We have seen briefly what the models and questions could be in the introduction. For those questions that are \textit{decidable} and solutions exist for them which are unique, the algorithms that resolve them are an integral part of model checking, and form the third component.

\paragraph{} A more formal understanding is that these logics present a way to express particular behaviours of the model by allowing us to write \textit{properties} it may potentially exhibit. Then, the goal is to determine whether the model has these properties, or formally that the model \textit{satisfies} these properties. The action of determining so, is termed \textit{model-checking}.

\paragraph{}  We summarize below previously investigated problems on probabilistic specific models related to sMRMs, and also detail some existing temporal logics defined over them. The relevant logics are detailed in the following chapter and specific algorithms related to our work for proving satisfaction of their properties on models are presented.

\subsection{DTMCs} The discrete time Markov chain (DTMC) forms the foundation for sMRMs and has been studied extensively. Problems for sMRMs that are independent of notions of rewards reduces to problems for DTMCs. Such problems can be partitioned into two: the long term behaviour of a Markov chain, and the short-term (or bounded) behaviour. The former includes investigating (i) \textit{expected first passage/arrival times} - the expected number of steps to reach a state $j$ beginning from another state $i$, (ii) \textit{equilibrium/steady state distribution} - the  probability distribution over being in each state of a Markov chain as time tends to infinity, (iii) and the \textit{mean recurrence time} - the expected amount of steps it takes for a Markov chain to return to state $i$ beginning from the same state. See the textbook \cite{kirkwood2015markov} for further details. For bounded behaviour, we have for example (iv) \textit{transient state probabilities} - the probability distribution over being in a particular state at some given time $n$.

\paragraph{} Probabilistic model-checking is a field which also studies DTMCs. The authors Hansson and Jonsson \cite{hansson1994logic} introduced a temporal logic PCTL - \textit{probabilistic computation tree logic}, which allows the expression of particular properties for DTMCs. All properties written in PCTL are entirely computable, and the algorithms developed to solve PCTL problems cover the entire class. The logic allows expression for a wide range of properties, which in turns allows us to solve problems including determining (i) \textit{first passage (or reachability) probabilities} - the probability of first entering a state $j$ from each state of the process, (ii) \textit{step-bounded reachability probabilities} - the probability of first entering a state $j$ from each state under $n$ steps, transient state probabilities, and (iii) \textit{repeated reachability and persistence probabilities} - the probability that the DTMC repeatedly enters a set of states, and the probability of a DTMC transitioning only within a particular set of states respectively. See the textbook \cite{principles_book} for an exposition to the subject.    

\paragraph{}A property written in such specification languages is to be resolved over its respective DTMCs. If a DTMC \textit{satisfies} a property, then we mean by this that the DTMC is guaranteed to exhibit such a behaviour. The algorithms involved for resolving properties over DTMCs can be categorized into three groups: 1) Translating the logical property into an equivalent form amenable to simpler computation. This would involve logic theory. 2) Transforming the DTMC with respect to the new property, preparing it for computation. This yields a simpler property as well. This is a combination of results and algorithms from graph theory and automata theory. 3) Solving the transformed DTMC with respect to the remaining property. As for algorithms for solving PCTL statements, they include those that range from being numerical or symbolic, global or local (on-the-fly), and deterministic or statistical. See for example \cite{principles_book,daws2004symbolic,latella2014fly,younes2002probabilistic}. If numerical approaches are used, then part of the solution may involve solving a system of linear equations, or repeated matrix-vector multiplications depending on the statement.

\paragraph{} Whilst PCTL is our language of focus, there are other languages which allow resolving of other properties over DTMCs. For example there is PCTL* - \textit{probabilistic computation tree logic star} by \cite{aziz1995usually}. It is a logic that includes as subsets PCTL, and LTL - \textit{linear temporal logic} introduced by \cite{pnueli1977temporal}, which allows resolving of $\omega$-regular properties over DTMCs.

\paragraph{} Additionally from model checking literature, is the problem of \textit{parametric} model checking, the problem where models are not completely described and have parameters instead, and the goal is to determine if such models satisfy particular logical properties for a range of different values of these parameters. For DTMCs, see \cite{daws2004symbolic}. Additionally the problem of \textit{repairing} Markov models has also been studied for \textit{controllable} DTMCs \cite{bartocci2011model}, where if such a DTMC does not satisfy a logical property, another closely related DTMC is sought to ensure the satisfaction of that property. PCTL has also been studied for DTMCs with continuous state spaces, where transition matrices are replaced with kernels instead \cite{ramponi2010connections}. Solutions may involve analogues of existing numerical solutions of finite-space DTMCs, see \cite{townsend2015continuous} as a guide. Available general surveys on probabilistic model checking include \cite{katoensurvey,kwiatkowska2010advances}.

\subsection{Discrete-time MRMs and sMRMs} Markov reward models (MRMs) are essentially (discrete-time) Markov chains extended with a reward structure, allowing us to weight the occurrence of events of the process with a reward (or dually, cost) function. Then, we can ask not only for probabilities of events as with DTMCs but also the expected reward accumulated for these events. These rewards are either attached to states \cite{principles_book}, or attached to transitions of the process \cite{cloth2005model}, and the process accumulates these rewards when being in a state or transitioning between them respectively. If the rewards are deterministic, then we have regular MRMs. If they are random variables, then we have stochastic MRMs (sMRMs).

\paragraph{} The study of Markov reward models similarly includes determining long-term and short-term behaviour. As for the former, it includes (i) the \textit{expected cumulated reward} for reaching a set of states $B$ beginning in any particular state, (ii) and a conditional variant where we have the expected reward to reach $B$ under the condition that $B$ is eventually reached, called \textit{conditional expected cumulated reward}. (iii) Also, we have \textit{long running averages} of states - the expected cumulated reward earned when beginning in a state $s$, and (iv) \textit{quantile probabilities} - the minimum reward bound such that the probability of reaching a set of states $B$ from each state $s$ whilst earning a reward less than the bound, is greater than a pre-specified probability $p$ \cite{baier2014energy,ummels2013computing}, and this is generalized for multivariate rewards by \cite{haase2017computing}.

\paragraph{} One existing specification language for MRM is PRCTL - \textit{probabilistic reward computation tree logic}, introduced by \cite{andova2003discrete}. It is an extension of PCTL, but allows the computation of various expected value problems over MRMs that includes those above (except for quantile probabilities) but provides algorithms to other properties also. See the paper for details.

\paragraph{} As for existing work directly on sMRMs, the \textit{variance of the cumulated reward} has already been studied, as has the \textit{covariance of cumulated reward between two sMRMs} with algorithms for both found in \cite{verhoeff2004reward}. As for those who introduced the sMRM and gave it its name \cite{bacci2019model}, they presented Monte-Carlo algorithms for model-checking a class of (dependent) multivariate-reward  sMRMs with the logic PRCTL. They presented an example problem where they solved for the probability of reaching a set of states $B$ from a particular state $i$, with the mean cumulated reward being less than or equal to $r$. We however chose to focus on deriving numerical algorithms instead and focused on multivariate-rewards that are mutually independent as a first. Additionally, it would appear that expected value problems for sMRMs can draw ideas from regular MRMs, as we will show in Chapter 4, two expected value problems including the one above can use solutions similar to that for regular MRMs.

\subsection{Continuous-time Markov models} Other related models include continuous-time Markov chains (CTMCs) and semi-Markov processes (SMPs).

\paragraph{} As for CTMCs, this is an extension of DTMCs, where transition times are now random and distributed by exponential distributions. CTMCs have been studied to solve for the steady state distribution, and for the \textit{Kolmogorov forward and backward equations}, with the latter being considered a major goal for CTMCs \cite{kirkwood2015markov}. Roughly speaking, the forward equation describes the probability of being in a state $j$ at some time $t+h$, given that it was in state $i$ at time zero. The backwards equation is the reverse of this, it gives the probability of being in state $i$ at time zero, given that the process is in state $j$ at time $t+h$.  Within CTMCs are other popular models, such as birth-processes and birth-death processes, both of which forward and backward equations are studied for.

\paragraph{} In probabilistic model checking, the temporal logic CSL - \textit{continuous stochastic logic} has been introduced for CTMCs by \cite{aziz1996verifying}, which is a continuous-time variant of PCTL, hence similar properties can be model-checked. It is extended by \cite{baier2003model} who also presented approximate model checking algorithms for the logic. The logic allows writing properties for (i) determining the probability of reaching a set of states (from a particular initial state) within a specific time interval, (ii)  finding the probability of remaining within a set of states within a specific time interval, and others.

\paragraph{} If CTMCs are extended with (deterministic) rewards, then the new model is called a continuous MRM (CMRM). CMRMs are a subset of univariate sMRMs, where each reward is distributed as an exponential distribution, and is related to time. Here, CSRL - \textit{continuous stochastic reward logic}  has been introduced for model-checking CMRMs by \cite{baier2000logical}. It includes as sub-logics, both CSL and CRL - \textit{continuous reward logic}. CSRL allows writing properties similar to CSL, they can be used additionally for (i) determining the probability of reaching a set of states (from a particular initial state) within a specified time interval and within a specified reward interval, (ii) finding the probability of remaining within a set of states within a specific time interval and within a specified reward interval.

\paragraph{} If we remove the model's restriction on transition time distributions being exponential distributions, we obtain the semi-Markov process (SMP) model. Hence, CTMCs are subsets of SMPs. Further, SMPs are syntactically univariate sMRMs (with non-negative rewards [i.e. time]) and are one of the closest models to it. Therefore much of their theory can be borrowed. Problems for SMPs include computing (i) \textit{first-passage time densities} - the probability distribution function over the reward accumulated when starting from a state $s$ and reaching a set of states $B$ \cite{warr2014numerical}, (ii) their moments \cite{harrison2002passage}, (iii) the \textit{cumulative distributions} for said densities and (iv) the \textit{hazard functions} derived from these densities \cite{warr2012introduction}. Additionally, the (v) mean recurrence time (as defined earlier for DTMCs) and
(vi) \textit{asymptotic state probabilities} - the probability an SMP is in state $j$ given that it began in state $i$, when time tends to infinity  \cite{warr2014numerical}.

\paragraph{} Summarily, these various models presented help lay the foundation for sMRM theory. Not just their theory can be borrowed but also the practical algorithms developed for them for problem solving.

\paragraph{} More generalized models can be studied that are beyond the scope of this thesis. For example by adding \textit{actions} to DTMCs, we arrive at an important model, the Markov decision process (MDP). This can be annotated with rewards, yielding Markov reward decision processes (MRDPs or sMRDPs). SMPs have a parametrized variant where covariates have been introduced to them \cite{huzurbazar2010incorporating}. There is also work on stochastic hybrid systems \cite{abate2010approximate,dhople2014stochastic,cauchi2019}  and probabilistic programming languages \cite{gordon2014probabilistic}, both being quite general models. There are other probabilistic models which exist that have been extended with reward structures. For example, Markov automata \cite{guck2014modelling} and stochastic Petri nets \cite{ciardo1989spnp}.

\section{Algorithms for sMRMs} Our main focus is in resolving \textit{first passage reward densities} or \textit{reward reachability densities} - the distribution function over cumulated reward for first reaching a set of states $B$, having begun from any state of the sMRM. The reason for this focus is that reachability problems are cornerstones of probabilistic model checking, and solving them enables us to solve a large class of problems. An understanding of this can be grasped after reading the next two chapters.

\paragraph{} We intend to determine these densities for two types of sMRMs: Continuous-reward and discrete-reward sMRMs. Their definitions will be given formally in Chapter 4, however the distinction is that continuous-reward sMRMs have their rewards all characterized by continuous random variables, whilst discrete-reward sMRMs have rewards characterized by discrete (lattice) random variables. These two classes are chosen as we have found them to have forms amenable to faster computations.

\subsection{Continuous-Reward sMRM and SMPs}
\label{subsec:algorithms_for_SMPs}

\paragraph{} Syntactically, a semi-Markov process (SMP) is a stochastic Markov reward model, where the rewards are univariate random variables and represent time. Therefore many problems for semi-Markov processes are analogous to problems for sMRMs. Likewise are we able to adapt their solutions for use.

\paragraph{} When solving for \textit{first passage time densities} in an SMP or equivalently, \textit{first passage reward densities} in an sMRM, we are generally confronted with a system of equations to solve \cite{warr2014numerical}. Then, we have found that traditional algorithms apply: Direct numerical methods such as Gaussian elimination can be used, iterative numerical methods such as the power, Jacobi or Gauss-Seidel methods, or symbolic approaches for small problems.

\paragraph{} In the original form, the system to solve is a \textit{system of convolution equations} (described later in Chapter 4). This system is transformable into a particular system of linear equations, where each term is a function. A solution to the system is usually computed at samples of these functions. Thus, a numerical approach generally requires solving a system of equations for each sample. This space-complexity blow up means that symbolic approaches for small problems are useful, as solving the system avoids sampling the functions.

\paragraph{} One algorithm for solving first passage time densities consists of three sub-algorithms: 1) Transforming the system into a set of linear equations. 2) Solving the system. 3) Inverse transforming the solution to obtain the densities.

\paragraph{} Transforming the system can be done via the continuous Fourier transform, discrete Fourier, or Laplace. This is done exactly either algebraically by hand, or via a computer algebra system (CAS). If neither are possible, it can be approximated numerically. Note that the Laplace transform does not exist for every random variable, whilst the continuous Fourier does. The discrete Fourier is used as an approximation to the problem and is always done numerically.

\paragraph{} Solving the system can be done as previously mentioned, either numerically, or algebraically (i.e. symbolically). However, combining numerical and symbolic approaches is also possible.

\paragraph{} Inverse transforming the solution is less straightforward than the initial transformation of the system. This is because inversion problems can be ill-conditioned. The inverse Laplace transform is considered one of them \cite{epstein2008bad}. However, for SMPs the inversion procedure appears unaffected since the transform is applied to a special class of functions - non-negative random variables, that are absolutely continuous with respect to the Lebesgue measure, this stated in \cite{warr2012introduction}. However, in general sMRMs may contain random variables that are not non-negative. Nevertheless, there appears to be quite a few different algorithms for inverse transforms of the continuous Fourier, for example \cite{abate1992fourier,witkovsky2016numerical} and Laplace, where we have a multitude \cite{hassanzadeh2007comparison,cohen2007numerical}.

\paragraph{} As for existing solutions to SMPs, we find \cite{bradley2004hypergraph} using the Laplace transform, the power method, and two inversions transforms: Laplace-Euler \cite{abate1995numerical} or Laplace-Laguerre \cite{abate1996laguerre}, as their three sub-algorithms. The paper details solving a system with more than a million states (under 10 minutes), hence showing the scalability of their approach. Another solution is by \cite{warr2014numerical} who also uses the Laplace transform, a perhaps direct numerical approach for solving the system, and the Laplace-Euler inversion. Their problems included one with 9 states, and is mentioned to be resolvable under a second. An earlier paper by the same authors \cite{huzurbazarstatistical}, used the Laplace transform, a graph reduction algorithm which they have described for solving the system and experimented with two inversion algorithms: the Laplace-Euler and a saddlepoint approximation algorithm. For each of the inversion procedures, the paper describes a case where they would perform poorly. It would appear that the Laplace-Euler does not work for empirical distributions naively. They did propose a fix for this, however it is specifically for empirical Laplace transforms (ELTs) \cite{elt_fast_inversion}\footnote{This paper may not be visible to the public, and we could not find access to it.}. Additionally from paper \cite{huzurbazarstatistical}, the Laplace transforms of the empirical distributions are numerically derived, and for a problem with around three states (and therefore up to nine transitions/pmfs) would take several minutes, which is slow. Another approach is \cite{huzurbazar2010incorporating} which uses the continuous Fourier transform, and experimented with two inverse transforms: The algorithm by \cite{abate1992fourier}, and a saddlepoint approximation via \cite{strawderman2004computing}. The paper showed that the former inversion to be better relative to the saddlepoint approximation. The models they presented were small and algebraic solutions were presented for the first passage time densities.

\paragraph{} Another algorithm from SMP literature available for finding first passage reward densities involves using the moments of these densities to infer the density itself. This is done using a vector of (non-negative integer moments) via the method of moments algorithm presented by \cite{au2004efficient}. It appears that whilst these algorithms may be useful as an approximation if using a few moments, may lead to \textit{representation explosion} if high precision is required, this stated in  \cite{bradley2004hypergraph}. As a vector is required for each state of the system to store the moments, these vectors will grow as long as the required accuracy of the density has not been achieved.

\paragraph{} Yet the most recent work \cite{bacci2019model}, that of which introduced the sMRM model into literature, resorted to sampling. This choice is perhaps due to sMRMs having generally $n$-dimensional reward random vectors, and therefore the complexity of such algorithms would be too high to resolve $n$-dimensional pdfs numerically. However, if the reward random vector consists of mutually independent random variables, then this reduces to solving 1-dimensional sMRMs $n$ times, which is significantly more tractable. This was the case with regular MRMs, where the reward random vector is solely composed of degenerate random variables (or constants), and hence are always independent of each other. Therefore independence between non-degenerate random variables is still a step up, and not sidewards. Sampling is generally regarded as slow, when high precision is required.

\paragraph{} As for general n-dimensional reward random vectors. This is presumably future work. Whatever future technique is to be considered, one has to keep in mind the dangers of \textit{representation explosion} or time-complexity growth when moving away from lattice representations.

\paragraph{} In our thesis, we focus on univariate sMRMs. When rewards are continuous, we approximate the density via the discrete Fourier transform (DFT), which benefits from some quadrature rules we have developed. The system is solved via the power method only.  The transform is numerical and fast, and the inverse appears to work over discontinuous distributions, which was a problem with the popular Laplace-Euler technique for the Laplace inversion transform. In fact, \cite{warr2012introduction} suggested to consider the DFT for the case of discontinuous distributions, or perhaps more generally, for problems where the Laplace inversion is generally ill-posed for. The power method was a good choice, as \cite{bradley2004hypergraph} showed how scalable it could be in solving problems with millions of states, although it is not sure if our algorithm scales as well with the DFT. Additionally the DFT requires just the pdf of these random variables, whether analytical or empirical. This is unlike the Laplace transform of a pdf which does not always exist. And the DFT is strictly numerical and can be computed quickly without requiring algebraic derivations or precise numerical integrations.

\subsection{Discrete-reward sMRMs and SMPs} As for work on solving \textit{first passage time mass functions} for discrete-time SMPs numerically, the only work we are familiar with is that by \cite{warr2014numerical}, which used the continuous Fourier transform of a function when known, and the discrete Fourier when not known. The inversion was done by the inverse discrete Fourier transform. They presented a problem with 3 states, which was solvable algebraically by hand. In this case, the discrete Fourier transform is used as an approximation to the solution. We will reproduce their problem later in Section \ref{problem:waste_treatment}. The authors stated that they were deterred from the Laplace transform (perhaps due to the difficulty of inversion with discrete random variables).

\paragraph{} For the case where rewards are discrete, we develop solutions for solving first passage reward mass functions, using the discrete Fourier transform and inverse transform. We are able to develop algorithms obtaining machine precision, via iterative methods such as the power method, Jacobi and Gauss-Seidel. Secondly, we also present a direct (exact) solution via the Gaussian elimination algorithm adapted for the system of convolution equations, which may prove useful for slowly converging problems that are not too big (due to its time complexity), if they occur. Thus, we move away from algebraic solutions found by hand (or computers) towards numerical algorithms.

    \chapter{Preliminaries}

\section{Introduction} We introduce here the theoretical foundations for probabilistic model checking for sMRMs with univariate rewards. To do so, we first define DTMCs, and use it to introduce the language PCTL - \textit{probabilistic tree control logic}. Then MRMs are introduced as well as PRCTL - \textit{probabilistic reward control tree logic}.

\paragraph{} Furthermore, theory is presented on the topics of summations of random variables, and characteristic functions due to their relevance to sMRMs. Within the topic of summation of random variables, is the topic of convolution (as is seen later). We find it necessary to also introduce the inverse operation of (discrete) convolution; deconvolution, which will be used for subsequent chapters.

\paragraph{} This chapter partially summarizes the book \cite{principles_book} on the topics of probabilistic model checking. Thus further results and more in-depth explanations can be found there. We will also adopt their notation quite considerably. Additionally, some main derivations presented in this thesis will be adaptations of a strategy found in \cite{verhoeff2004reward}. The technique they presented is helpful in deriving solutions. Their notation will also be in this thesis.

\section{DTMCs and PCTL} Discrete time Markov chains are a class of probabilistic models that is used to represent a system consisting of a (finite) number of states, and where the system can only be in one state at any given time. As time passes, the system can transition between states randomly. Randomness is represented quantitatively as probabilities. Additionally, the model being Markovian assumes that if the system is in a particular state, the probability of transiting into another state must not be dependent on the previous states that the model was in.

\begin{definition}
	\label{definition:DTMC}
	Formally, A DTMC $\mathcal{M}$ is a tuple $(S, \textbf{P}, i_{init})$, where:
	\begin{itemize}
		\item $S$ is the state space, with the size being finite, i.e. $|S| < \infty$,
		\item $\textbf{P}: S \times S \rightarrow [0,1]$ is a map for the transition probabilities between any two states of the Markov chain. We set $\textbf{P}$ to be constant with respect to time, hence it is stationary.
		\item $i_{init}: S \rightarrow [0,1]$ is the initial distribution function of the chain, with $\sum_{s \in S}i_{init}(s) = 1$. This is to specify in which state the Markov chain (the system) is likely to have started in.
	\end{itemize}
	
\end{definition}

\subsubsection{Sample space, events, and a probability measure}

Let us define a sequence of states as a \textit{path}, e.g. a path could be $s_0,s_1,s_2,\cdots s_n$ with each $s_i \in S$. A path can then be used to denote a possible \textit{outcome} of the DTMC starting from state $s_0$, transiting consecutively to states $s_1,s_2,\cdots,s_{n-1}$, and ending in $s_n$. Then the set of all unique infinite length paths a DTMC exhibits defines a \textit{sample space}. Any \textit{measurable} subset of the sample space is generally called an \textit{event}.

\paragraph{} The manner in which probabilities are assigned to events must satisfy Kolmogorov's axioms, to be a valid \textit{classical probability theory}. However, we can define the probability of a finite path $s_0,s_1,s_2,\cdots s_n$ without much complication,  written as $Pr(s_0,s_1,s_2,\cdots s_n)$, which is equal to the probability that the system began in $s_0$, and performs transitions until it reaches $s_n$. Since the process is Markovian, we have
\begin{flalign*}
Pr(s_0,s_1,s_2,\cdots s_n)
&= Pr(s_0)Pr(s_1|s_0)Pr(s_2|s_1,s_0)\cdots \cdot Pr(s_n|s_{n-1}\cdots,s_0) &\\
&= Pr(s_0)Pr(s_1|s_0)Pr(s_2|s_1)Pr(s_3|s_2)\cdots\cdot Pr(s_n|s_{n-1})&\\
&= i_{init}(s_0)\textbf{P}(s_0,s_1)\textbf{P}(s_1,s_2)\cdots\textbf{P}(s_{n-1},s_n)
\end{flalign*}

\paragraph{} Let $\pi$ be an arbitrary finite path, and $\Omega_{\pi}$ be the set of all infinite paths that begin with $\pi$. Then it must be the case that $$Pr(\Omega_{\pi}) = Pr(\pi)$$
Therefore, we are now able to talk about the probability of certain \textit{events}.

\paragraph{} In literature, $\Omega_{\pi}$ is called the \textit{cylinder set} of $\pi$.  Let $\Pi$ be the set of all finite paths possible through a DTMC. Then the set $\{\Omega_{\pi}| \pi \in \Pi \}$ characterises the set of all basis events, and the smallest $\sigma$-algebra derived from it represents the \textit{event space}.

\paragraph{} Then let $\mathscr{S}$ be the sample space, and $\mathscr{E}$ be the event space (or $\sigma$-algebra). We can define a complete \textit{probability space} for DTMCs as $(\mathscr{S},\mathscr{E},Pr)$.

\paragraph{} By definition, the empty path fragment $\pi = \{\} = \mathscr{S}$ (the sample space), occurs with probability 1, i.e. $Pr(\pi) = 1$.

\subsubsection{PCTL notation for events and probability} We can define a grammar to specify events using probabilistic computation tree logic (PCTL) \cite{hansson1994logic}. The logic introduces the following operators:$\ \Diamond, \Square, \Circle, \text{U}$, which translates loosely as \textit{eventually}, \textit{always}, \textit{next}, and \textit{until}.

\paragraph{} Define $B,C$, both to be arbitrary subsets of $S$. Then for example:

\begin{enumerate}
	\item$\Diamond B$ - means the set of paths of a DTMC which \textit{eventually} reaches $B$, i.e. all outcomes where a transition into a state of $B$ occurred,
	\item $C\ \text{U}\ B$ -  denotes the event, the collection of all outcomes, where the DTMC was in $C$ (\textit{until}) before directly transiting into $B$. This excludes all outcomes not beginning in $C$, or not eventually transiting into $B$.
	\item $C\ \text{U}^{\leq n}\ B$ - denotes an event similar to above, however the set of outcomes is restricted only to those that transit into $B$ under $n$ steps, where $n \in \mathbb{N}$.
	\item $\Square B$ - the set of outcomes that never transition out of the states of $B$, i.e. \textit{always} remaining within these states.
	\item $\Circle B$ - the set of outcomes that enter $\next B$ after the \textit{next} (immediate) transition.
\end{enumerate}

\paragraph{} Let $s_i \in S$ be some state. Additionally from PCTL, the notation $$s_i \vDash \Diamond B$$
means that the event $\Diamond B$ \textit{holds} for $s_i$. In terms of events, this is \textit{conditioning} the sample space to the set of paths beginning strictly from $s_i$, and restricting the set of paths to those satisfying $\Diamond B$. Thus in \cite{principles_book}, this is creating a new DTMC with $s_i$ as a deterministic initial state, which leads us to a new sample space $\mathscr{S}_s$ and then determining the event $\Diamond B$. An alternative representation could be $$\Diamond B\ |\ (\text{initial state is $s_i$})$$

\paragraph{} Finally, PCTL is specifically used to ask for the probability of these events, therefore we can ask for $Pr(s_i \vDash \Diamond B)$. Then
$$Pr(s_i \vDash \Diamond B) = Pr(\Diamond B | \ \text{initial state is $s_i$})$$

\subsubsection{PCTL grammar} Valid statements of PCTL are determined by the following grammar. For the probability of the event $\psi$ beginning in a state $s$, i.e. $Pr(s \vDash \psi)$, $\psi$ can be defined from any of:

\begin{align*}
\psi \Coloneqq \Circle\Phi\ |\ \Phi_1\ \text{U}\ \Phi_2 \ |\  \Phi_1\ \text{U}^{\leq n}\ \Phi_2
\end{align*}
\begin{align*}
\Phi \Coloneqq \text{true}\ |\ C\ |\ \Phi_1 \cap \Phi_2\ |\ \neg\Phi\ |\ \mathbb{P}_{[a,b]}(\psi)
\end{align*}
where $C \subseteq S$ and $$\mathbb{P}_{[a,b]}(\psi) = Pr(\psi) \in [a,b]$$ denoting whether or not the probability of the event $\psi$ is within some interval $[a,b] \subseteq [0,1]$ (excluding the empty-set $\{\}$), and where $n \in \mathbb{N}$. In the literature, $\Phi$ is termed a \textit{state formula}, whilst $\psi$ is termed a \textit{path formula}. A state $s$ satisfies a state formula $\Phi$ when $s \in \Phi$, and a path $\pi$ satisfies a path formula $\psi$ when $\pi \vDash \psi$, or  $\pi \in \psi$ if we interpret $\psi$ as an event in a probability space.

\paragraph{} The grammar is written quite concisely, hiding the inclusion of the other operators previously mentioned. For example, the remaining \textit{logical} operators: OR $\vee$, and implication $\rightarrow$, can be derived using $\neg,\wedge$ alone. Additionally, the \textit{temporal-logic} operators $\Square$, $\Diamond$, are both derivable using the until operator $\text{U}$. For example, (perhaps intuitively) we have $\Diamond B \equiv \text{true U}\ B$, and $\square B \equiv \neg \Diamond \neg B$. The operator $\text{U}^{\leq n}$ is called \textit{constrained}-until or \textit{step-bounded}-until. Other constrained operators can be defined, e.g. $\Diamond^{\leq n}$, or $\Square^{\leq n}$.

\paragraph{Reachability problems} The event $\Diamond \Phi$ characterizes the class of \textit{reachability}  problems. The phrase $\Diamond \Phi$ is synonymous to the set of paths where the set of states satisfying the state formula $\Phi$ is \textit{eventually} \textit{reached}. This is important since by logic theory, all unconstrained statements above reduces to mostly resolving strictly $\Diamond B_1$ or $\Circle B_2$ statements, i.e. the complexity of all non-constrained PCTL statements is not significantly much more than being able to resolve reachability statements and \textit{next} statements.

\paragraph{} Notice that the state formula $\Phi$ always reduces to a set of states $B \subseteq S$. In total therefore, there are only three main algorithms required for resolving PCTL statements. An algorithm for next statements, e.g. $\Circle B$, until statements $B_2\ \text{U}\ B_2$, and constrained-until statements e.g. $B_2\ \text{U}^{\leq n} \ B_2$, where $B,B_1,B_2$ are arbitrary subsets of $S$. Repeating from above, it is provable that until, and constrained-until statements for a DTMC $\mathcal{M}$ reduces to reachability and constrained reachability of a new DTMC $\mathcal{M}'$ \cite{principles_book}. Aside from the algorithm needed for this transformation, we mostly only need to solve statements of the form $\Circle B, \Diamond B_1$ and $\Diamond^{\leq r} B_2$.

\paragraph{Thesis focus} In this thesis, we focus solely on solving unconstrained reachability problems. Hence, we will not delve much into $\Circle \Phi$ properties, nor constrained properties such as $\Phi_1\ \text{U}^{\leq n}\ \Phi_2$. For model checking algorithms for these, please refer to \cite{principles_book}. Reachability problems for sMRMs will be defined in the next chapter.  We detail below the means of computing reachability probabilities - $Pr(s \vDash \Diamond B)$.

\subsubsection{Solution to reachability statements: $Pr(s \vDash \Diamond B)$}
\label{markov_chain_reachability_problems}

\paragraph{} Given the event $\psi = \Diamond B$, then the set of all paths in $\psi$ can be uniquely defined as the set consisting of all cylinder sets of all unique finite paths $\pi= s_0,s_1,\cdots,s_n$, where $s_0 = s$, with $s_1,\cdots,s_{n-1} \notin B$ and $s_n \in B$. Then let this set of finite paths be $\Pi_\psi$, and the set $\Omega_{\psi} = \{\Omega_{\pi} | \pi \in \Pi_\psi \}$. In this case, the elements of $\Omega_{\pi}$ do not overlap , i.e. there are no two cylinder sets $\Omega_{\pi_1}, \Omega_{\pi_2}$ such that $\Omega_{\pi_1} \in \Omega_{\pi_2}$. Then,

$$Pr(\psi)  = \sum_{\Omega_{\pi} \in \Omega_{\psi} } Pr(\Omega_{\pi}) = \sum_{\pi \in \Pi_\psi} Pr(\pi)$$

\paragraph{Notation} Before we proceed further, for the derivations to come, we will typically reserve letters $s,t$ to denote states, and symbols $\hat{\pi}, \hat{\psi}$ to denote paths. Additionally, when we write $\Pi.t$, we mean the set of paths of $\Pi$ that begin from state $t$, i.e. $\{\hat{\pi} \in \Pi | \ \hat{\pi}[0] = t \}$. We can concatenate states and paths via dot notation, for example $P(s.\hat{\pi})$ denotes the probability of a path event, one that begins with state $s$ and is followed by path $\hat{\pi}$.

\paragraph{}Given a DTMC $\mathcal{M}$, let $B$ be a set of target states we are interested in. Given an initial state $s$, the probability of the event that $B$ is eventually reached is defined to be the collection of finite paths $\Pi_{\Diamond B}$, where each path $s_0,s_1,\cdots ,s_{n-1},s_n$ has $s_n \in B$, and $s_0,\cdots s_{n-1} \notin B$. Let $\Pi_{\Diamond B}.s$ denote the set of paths satisfying $\Diamond B$ beginning from $s$. Then,
\begin{flalign*}
Pr(s \vDash \Diamond B)
&= \sum_{\pi \in \Pi_{\Diamond B}.s} Pr(\pi) &\\
&= \sum_{t \in S} \sum_{\pi \in \Pi_{\Diamond B}.t} Pr(s.\pi) &\\
&= \sum_{t \in S} \sum_{\pi \in \Pi_{\Diamond B}.t} \textbf{P}(s,t)Pr(\pi) &\\
&= \sum_{t \in S} \textbf{P}(s,t) \sum_{\pi \in \Pi_{\Diamond B}.t}  Pr(\pi)   &\\    
&= \sum_{t \in S} \textbf{P}(s,t) Pr(t \vDash \Diamond B)
\end{flalign*}

\paragraph{} We also have $Pr(s \vDash \Diamond B) = 1$ for all $s \in B$. Thus, the above is equivalent to
\begin{flalign*}
Pr(s \vDash \Diamond B)  
&= \sum_{q \in S_?} \textbf{P}(s,q) Pr(q \vDash \Diamond B) + \sum_{u \in B}P(s,u)
\end{flalign*}
with $S_?$ being the set of states that reach $B$ with non-zero probability but are not within it.

\paragraph{} One way to solve this for DTMCs with finite state spaces, is to represent the problem as a system of linear equations. Let $x_s \triangleq Pr(s \vDash \Diamond B)$ for all $s \in S$, then $$x_s = \sum_{t \in S_?}\textbf{P}(s,t)x_t + \sum_{u \in B}\textbf{P}(s,u)$$

\paragraph{} Then let \textbf{x} be the vector $(x)_{s \in S_?}$, and $\textbf{b} = (\sum_{u \in B}\textbf{P}(s,u))_{s \in S_?}$, and we define a matrix \textbf{A} where $A_{i,j} = \textbf{P}(i,j)$. Now we can write the system as $$\textbf{x} = \textbf{Ax} + \textbf{b}$$ with the solution being $$\textbf{x} = (\textbf{I}-\textbf{A})^{-1}\textbf{b}$$

\paragraph{} Methods available for solving the system could be Gaussian-elimination, the power method, Jacobi or Gauss-Seidel methods.

\section{MRMs and PRCTL} Markov reward models (MRMs) extend DTMCs, by allowing transitions to be detailed with rewards or costs.\footnote{Rewards and costs are synonymous in this thesis. However, in the context of a problem, the use of one may be more suitable than the other.} Thus, this process not only transitions from state to state as time progresses (like a DTMC), but also accumulates reward along these transitions. MRMs can then be used to study the expected rewards of reaching a particular state beginning in another, for example.

\begin{definition}
	An MRM is a tuple ($\mathcal{M}$, \textit{rew}), where:
	\begin{itemize}
		\item $\mathcal{M}$ is a DTMC, i.e. $\mathcal{M} = (S, \textbf{P}, i_{init})$.
		\item \textit{rew}: $S \times S \rightarrow \mathbb{R_{+}}$ is a reward function, that assigns a reward to each transition $s \rightarrow t$ in $\mathcal{M}$ for all $s,t \in S$. Note that the rewards are strictly non-negative for each transition.
	\end{itemize}
	
\end{definition}

\paragraph{} Let $\pi = s_0,s_1,\cdots,s_n$ be a finite path, then the \textit{cumulated reward} over $\pi$, denoted $Rew(\pi)$ is defined as $$ Rew(\pi) = rew(s_0,s_1)+rew(s_1,s_2)+rew(s_2,s_3)+\cdots +rew(s_{n-1},s_n) $$

\paragraph{} The probability of this path occurring, i.e. $Pr(\pi)$, is defined as it was for DTMCs. The same probability space as DTMCs can be used for MRMs.

\subsubsection{Contingent expected rewards} One important measure for MRMs is to compute the expected (or mean) cumulated reward earned for a particular event, beginning in some state $s$. Let us denote this as $ExpRew(\psi)$, for the event $\psi$. Recall that the definition of the expected value of a discrete random variable defined over the real line is typically:
$$\mathbb{E}[X] = \sum_x x\cdot f(x)$$
\paragraph{} In the context of MRMs, the expected (cumulated) reward with respect to the sample space $\mathscr{S}$ has to be defined in a manner that avoids $\infty$ being introduced unnecessarily. One consideration is that like DTMCs we can make use of cylinder sets. Let $\Omega_{\pi}$ be an event, where $\pi$ is a finite path. Then define the cumulated reward earned by the event $\Omega_{\pi}$ as $Rew(\pi)$. Thus all basis events $\Omega_{\pi}$ are well defined and less than $\infty$.  

\paragraph{} An event of interest can be decomposed into a set of non-overlapping cylinder sets (or basis events). Let this event be $\psi$, and $\Omega_{\psi}$ be the set of non-overlapping cylinder sets that describes all paths of the event. Define also $\Pi_\psi = \{\pi |\ \Omega_{\pi} \in \Omega_{\psi} \}$ . Then we define the expected reward \textit{contingent} on this event as:

$$\mathbb{E}[\psi] = \sum_{\Omega_{\pi} \in \Omega_{\psi}} Rew(\pi)\cdot Pr(\Omega_{\pi})
=  \sum_{\pi \in \Pi_\psi} Rew(\pi)\cdot Pr(\pi) $$

\paragraph{} However, $\psi$ being decomposable into cylinder sets is not enough to guarantee that $\mathbb{E}[\psi] \neq \infty$. This is since an event can characterize a single infinite path $\pi \in \mathscr{S}$, of which generally accumulates infinite reward. Thus, a possible consideration (in certain applications) is to assign zero-reward to such events.

\paragraph{Expected rewards of reachability problems} The class of expected rewards contingent on the event $s \vDash \Diamond \Phi$ with $\Phi$ being a PCTL state formula and satisfying the condition that $Pr(s \vDash \Diamond \Phi) = 1$ for all $s \in S$, is known to always be finite. I.e. $\mathbb{E}[s \vDash \Diamond \Phi] < \infty$ for any $\Phi$. This class of reward problems is added to PCTL in an extension known as PRCTL - \textit{probabilistic reward computation tree logic} \cite{andova2003discrete}.

\subsubsection{PRCTL grammar} The grammar for PRCTL extends PCTL with two terms:

\begin{enumerate}
	\item $\mathbb{E}(\Phi)$ - the expected reward for reaching a set of states determined by the state formula $\Phi$. Its definition is $$\mathbb{E}_{[c,d]}(\Phi) = ExpRew(s \vDash \Diamond \Phi)$$
	\item $\Phi_1\ \text{U}_{\leq r}\ \Phi_2$ - reachability probabilities constrained on bounded rewards. Elaborating, this is the event where $\Phi_2$ is reached without accumulating more than $r \in \mathbb{N}$ in reward, with the DTMC being in states $\Phi_1$ up till that point. From this formula, we can derive \textit{reward-bounded} reachability problems of the form $\Diamond_{\leq r} \Phi$ -  the event where $\Phi$ is reached without accumulating more than $r$ reward. These formulas are used to compute the probabilities of such events, i.e. $Pr(s \vDash \Diamond_{\leq r} \Phi)$.
\end{enumerate}

\paragraph{} The full PRCTL grammar is as follows:
\label{grammar_PRCTL}
\begin{align*}
\psi \Coloneqq \Circle\Phi\ |\ \Phi_1\ \text{U}\ \Phi_2 \ |\  \Phi_1\ \text{U}^{\leq n}\ \Phi_2  |\ \Phi_1\ \text{U}_{\leq r}\ \Phi_2 \
\end{align*}
\begin{align*}
\Phi \Coloneqq \text{true}\ |\ B\ |\ \Phi_1 \cap \Phi_2\ |\ \neg\Phi\ |\ \mathbb{P}_{[a,b]}(\psi)|\ \mathbb{E}_{[c,d]}(\Phi)
\end{align*}

where $[c,d] \subseteq [0,\infty)$, and $ [a,b] \subseteq [0,1]$ excluding the empty-set, and both $n,r \in \mathbb{N}$ and
$$ \mathbb{E}_{[c,d]}(\Phi) = ExpRew(s \vDash \Diamond B) \in [c,d]$$

\paragraph{} The grammar above can be made to include the term $\Phi_1\ \text{U}_{[a,b]}\ \Phi_2$ - the event where $\Phi_2$ is reached in accumulating reward only between $a$ and $b$ (both in $\mathbb{N}$), but the DTMC also remaining within $\Phi_1$ until transitioning into $\Phi_2$. Then the computation $Pr(\Phi_1\ \text{U}_{[a,b]}\ \Phi_2)$ is equivalent to $$Pr(\Phi_1\ \text{U}_{\leq b}\ \Phi_2) - Pr(\Phi_1\ \text{U}_{\leq a}\ \Phi_2)$$
Hence knowing how to solve for the reward-bounded reachability is all that is required.

\paragraph{Reachability problems} Note that $ExpRew(s \vDash \Diamond \Phi)$  reduces to $ExpRew(s \vDash \Diamond C)$, for some $C \subseteq S$. Also, determining the probability of $Pr(s \vDash \Phi_1\ \text{U}_{\leq r}\ \Phi_2)$ is not much harder than determining $Pr(s \vDash \Diamond_{\leq r}\ B)$ for some $B \subseteq S$ \cite{principles_book}. We proceed to present algorithms for these two reachability reward problems.

\subsubsection{Solution to expected reward of reachability statements: $\mathbb{E}(B) = ExpRew(s \vDash \Diamond B)$}  Given an MRM $\mathcal{R}$, let $B \subseteq S$. Let us denote $\Pi_{\Diamond B}.s$ as the set of paths of the event $s \vDash \Diamond B$. Then

\begin{flalign}
\label{eqn:expected_rew}
ExpRew(s \vDash \Diamond B)
&= \sum_{\pi \in \Pi_{\Diamond B}.s} Rew(\pi) Pr(\pi) \notag &\\
&= \sum_{t \in S} \sum_{\pi \in \Pi_{\Diamond B}.t} Rew(s.\pi) Pr(s.\pi) \notag &\\
&= \sum_{t \in S} \sum_{\pi \in \Pi_{\Diamond B}.t} (rew(s,t) + Rew(\pi))\cdot \textbf{P}(s,t)Pr(\pi) \notag &\\
&= \sum_{t \in S} \textbf{P}(s,t) \sum_{\pi \in \Pi_{\Diamond B}.t} rew(s,t) Pr(\pi) + Rew(\pi) Pr(\pi) \notag &\\    
&= \sum_{t \in S} \textbf{P}(s,t) \Big(\sum_{\pi \in \Pi_{\Diamond B}.t} rew(s,t) Pr(\pi) + \sum_{\pi \in \Pi_{\Diamond B}.t}Rew(\pi) Pr(\pi) \Big) \notag &\\
&= \sum_{t \in S} \textbf{P}(s,t) \Big(rew(s,t) \sum_{\pi \in \Pi_{\Diamond B}.t}  Pr(\pi) + ExpRew(t \vDash \Diamond B) \Big) \notag &\\    
&= \sum_{t \in S} \textbf{P}(s,t)\Big(rew(s,t)Pr(t \vDash \Diamond B) + ExpRew(t \vDash \Diamond B)  \Big)
\end{flalign}

where we have also defined $ExpRew(t \vDash \Diamond B) = 0$, for all $t \in B$. I.e. the expected cumulated reward for states already in $B$ is zero. Note that the result above is strange in that typically we do not find the term $Pr(t \vDash \Diamond B)$ within the derivation. For example \cite{verhoeff2004reward,principles_book} both do not present such a term. However, in both, their derivations assumed that $Pr(s \vDash \Diamond B)=1$.

\paragraph{} If $Pr(s \vDash \Diamond B)=1 $ for all $s \in S$, then the above is equivalent to:
\begin{flalign*}
ExpRew(s \vDash \Diamond B)
&= \sum_{t \in S} \textbf{P}(s,t)\Big(rew(s,t) + ExpRew(t \vDash \Diamond B)  \Big) &\\
&= \sum_{t \in S_?} \textbf{P}(s,t)\Big(rew(s,t) + ExpRew(t \vDash \Diamond B)  \Big) +  \sum_{u \in B} \textbf{P}(s,u)rew(s,u)
\end{flalign*}
with $S_?$ being the set of states that can reach $B$ with non-zero probability, but not in $B$, i.e. $S \setminus B$ since $Pr(s \vDash \Diamond B)=1 $ for all $s \in S$.

\paragraph{} Let $r_s \triangleq ExpRew(s \vDash \Diamond B)$, and $x_s \triangleq Pr(s \vDash \Diamond B)$ for all $s \in S$. Then the general form \eqref{eqn:expected_rew} can be written as
$$r_s = \sum_{t \in S}\textbf{P}(s,t)(rew(s,t)\cdot x_t + r_t)$$
Then define
\begin{align*}
&\textbf{r} \triangleq      (r_s)_{s \in S_?} &\\
&\textbf{b} \triangleq (\sum_{t \in S_?}P(s,t)rew(s,t)x_s + \sum_{u \in B}\textbf{P}(s,u)rew(s,u)_{s \in S_?} &\\
&\textbf{A} \triangleq (A_{i,j})_{i,j \in S_?^2} \triangleq  (\textbf{P}(i,j)))_{i,j \in S_?^2}
\end{align*}

\paragraph{}Now we can write a system of equations in the form $\textbf{r} = \textbf{Ar} + \textbf{b}$. And the solution is $$\textbf{r} = (\textbf{I}-\textbf{A})^{-1}\textbf{b}$$

\subsubsection{Solution to reward-bounded reachability probabilities: $Pr(s \vDash \Diamond_{\leq r} B)$}  
\label{proof:reward_bounded reachability}

\paragraph{} The approach in \cite{principles_book} solves this problem as follows:

\begin{flalign*}
Pr(s \vDash \Diamond_{\leq r} B)
&= \sum_{\pi \in \Pi_{\Diamond B}.s}Pr(Rew \leq r| \pi)Pr(\pi) &\\
&= \sum_{t \in S}\sum_{\pi \in \Pi_{\Diamond B}.t}Pr(Rew \leq r| s.\pi)Pr(s.\pi) &\\
&= \sum_{t \in S}\sum_{\pi \in \Pi_{\Diamond B}.t}Pr(Rew \leq r| s.\pi)\textbf{P}(s,t)Pr(\pi) &\\
&= \sum_{t \in S}\textbf{P}(s,t)\sum_{\pi \in \Pi_{\Diamond B}.t}Pr(Rew \leq r - rew(s,t)| \pi)Pr(\pi)&\\
&= \sum_{t \in S}\textbf{P}(s,t)Pr(t \vDash \Diamond_{\leq (r - rew(s,t))}\ B)
\end{flalign*}

\paragraph{} Also, we have $Pr(s \vDash \Diamond_{\leq r} B)  = 1$ if $s \in B$. This is the same for $Pr(s \vDash \Diamond_{\leq q} B)$ for any $q \in \mathbb{N}$, when $s \in B$. This is since no reward is accumulated for these states, thus they satisfy the inequality $(Rew \leq r)$ immediately. If $s \notin B$, then $Pr(Rew \leq r - rew(s,t)| \pi) = 0$ if $rew(s,t) > r$, for any $t \in S$, since it is impossible to reach $B$ through $t$ satisfying $(Rew \leq r)$.  We also assume that $Pr(s \vDash \Diamond B) = 1$ for all $s \in S$.

\paragraph{} Let $x_{s,r} \triangleq Pr(s \vDash \Diamond_{\leq r} B)$ for all $s \in S$. Then the results above can be written as
$$x_{s,r} = \sum_{t \in S}\textbf{P}(s,t)x_{t,(r-rew(s,t))}
$$

\paragraph{} Notice that $x_{s,p}$ can be computed independently from all $x_{t,q}$ where $q > p$, and $t \in S$. Thus, $x_{s,r}$ is solved by computing successively, the terms $$(x_{s,0})_{s \in S_?}, (x_{s,1})_{s \in S_?},\cdots,(x_{s,r-1})_{s \in S_?},(x_{s,r})_{s \in S_?}$$ where each $(x_{s,p})_{s \in S_?}$ uses the previous terms as is seen in the equation above.

\paragraph{} The above solution requires recursive computations, but yields a system of equations to solve when zero-rewards exist, i.e. when there exists $rew(s,t) = 0$, for some $s,t \in S^2$. For a state $t$, let $S_{t,0}$ be the set of states such that $rew(t,u) = 0$ for all $u \in S$.  Then, we can write the solution above as:
$$x_{s,r} = \sum_{t \in S\setminus S_{t,0}}\textbf{P}(s,t)x_{t,(r-rew(s,t))}    
+ \sum_{q \in S_{t,0}}\textbf{P}(s,t)x_{t,r}
$$   

\paragraph{} Define
\begin{align*}
&\textbf{x}_p \triangleq      (x_{s,p})_{s,p \in S \times \{0,1,\cdots,r\}} &\\
&\textbf{b} \triangleq (\sum_{t \in S\setminus S_{t,0}}\textbf{P}(s,t)x_{t,(r-rew(s,t))})_{s \in S}     &\\
&\textbf{A} \triangleq (A_{i,j})_{i,j \in S^2} \triangleq  \Bigg(\begin{cases}
\textbf{P}(i,j)) & rew(i,j) = 0\\
0 & otherwise
\end{cases}
\Bigg)_{i,j \in S^2}
\end{align*}
\paragraph{}Now we can write a system of equations in the form $\textbf{x}_p = \textbf{A}\textbf{x}_{p} + \textbf{b}$ with the solution being just $\textbf{x}_p = (\textbf{I}-\textbf{A})^{-1}\textbf{b}$. In this way we can compute $$\textbf{x}_0,\textbf{x}_1,\cdots,\textbf{x}_r$$ consecutively, to solve the problem.

\section{SMPs}  A semi-Markov process is a DTMC, except that transitions do not occur deterministically with respect to time, but rather by random. This randomness is characterized by transition time distributions.
Semi-Markov processes define a rich class of probabilistic models, including the DTMC and CTMC.

\begin{definition}
	A semi-Markov process is a tuple ($\mathcal{M},\textbf{G}$) where:
	
	\begin{itemize}
		\item $\mathcal{M}$ is a DTMC, i.e. $\mathcal{M} = (S, \textbf{P}, i_{init})$.    
		\item $\textbf{G}$: $(S \times S) \rightarrow f$ is a map between every state transition $t \in S \times S$ of the DTMC, and a probability distribution $f$ that determines the transition time distribution.
	\end{itemize}
\end{definition}

\paragraph{} The difference between an SMP and a univariate sMRM is that $\textbf{G}$ is generalized to not be related to time. It is just a map over state transitions to reward distributions. Due to the closeness of these models, we choose not to present reachability problems and their theoretical solutions for SMPs here to prevent overlapping results. Instead, the details are presented only for sMRMs in the next chapter.

\section{Random variables and characteristic functions}

\subsection{Representation of random variables}

Given two random variables \textit{X, Y}, they can be represented in a variety of ways, e.g. using their densities (pdf) $f_X(x), f_Y(y)$, or their cumulative distributions $F_X(x), F_Y(y)$, or their characteristic functions $\phi_X(\tau), \phi_Y(\tau)$. For any random variable $X$, its probability density function $f_X(x)$ can be transformed to a \textit{characteristic function} $\phi_X(\tau)$, and then inverse transformed back into $f_X(x)$.

\begin{definition} A \textit{characteristic function} of a random variable $X$, written as $\phi_X(\tau)$ can be derived via the formula: $\phi_X(\tau) = E[e^{i \tau X}]$. Therefore if $X$ is continuous, then
	$E[e^{i \tau X}]= \int_{- \infty}^{\infty}e^{\tau ix}f_X(x)dx$
	and if $X$ is a discrete (lattice) random variable, for example if it has $\mathbb{N}$ as a support, then
	$E[e^{i \tau X}]= \sum_{x=- \infty}^{\infty}e^{\tau ix}f_X(x)$
	
	\paragraph{} The characteristic functions can equivalently be derived by applying the correct Fourier transforms to either the pdf or pmf of $X$.
	
\end{definition}    

\begin{definition}
	\label{def:fourier_transforms}
	The \textit{continuous Fourier transform (FT)} applied to a function $f$ is written as $\mathcal{F}\{f\}$ and defined as
	$\mathcal{F}\{f\}(\tau) = \int_{- \infty}^{\infty}e^{\tau ix}f(x)dx$. The \textit{inverse continuous FT} applied to a function $\phi$ is written as $\mathcal{F}^{-1}\{\phi\}$ with definition
	$\mathcal{F}^{-1}\{\phi\}(x) = \int_{-\infty}^{\infty}\phi(\tau)e^{2\pi ix\tau}d\tau$.
	
	\paragraph{} The \textit{discrete-time FT} applied to a function $g$ is written $\mathcal{F}_d\{g\}$ and defined as
	$\mathcal{F}_d\{g\}(\tau) = \sum_{x=- \infty}^{\infty}e^{\tau ix}g(x)$. The \textit{inverse discrete-time FT} applied to a function $\psi$ is written as $\mathcal{F}_d^{-1}\{\psi\}$ with definition
	$\mathcal{F}_d^{-1}\{\psi\}(x) = \sum_{-\infty}^{\infty}\psi(\tau)e^{2\pi ix\tau}$.
	
	\paragraph{} Let $X$ be a random variable. If $X$ is continuous, then $\phi_X(\tau) = E[e^{i \tau X}] = \mathcal{F}\{f_X\}(\tau)$. If $X$ is a discrete lattice defined over $\mathbb{N}$, then $E[e^{i \tau X}] = \mathcal{F}_d\{f_X\}(\tau)$. Each of these characteristic functions can then be transformed back to the pdf or pmf with the respective inverse Fourier transform. Note that in the future we may drop the subscript $d$ from $F_d$, and therefore the type of Fourier transform is to be inferred from the context.
	
\end{definition}

\paragraph{} Note however that there are random variables that have characteristic functions but no analytical expressions are known for their probability density function (e.g. the Stable distribution). Nevertheless, every random variable  has a characteristic function. This is due to the proposition below.

\begin{proposition}
	\label{proposition:bounded cf}
	For any random variable $X$, we have that $|\phi_X(\tau)| \leq 1$ for all $\tau$. Hence, the integral (or summation) above in $\mathcal{F}\{f_X\}(\tau)$ absolutely converges and always exists. See \cite[page 97]{kendall1977} for proof and details (of the continuous case).
\end{proposition}

\paragraph{}A random variable $Z$ defined as a result of a sum of two other (independent) random variables $X,Y$, e.g $Z=X+Y$, can be represented as a density function, derived from a convolution operation on the other two respective probability density functions: $f_Z(z) = \int_{-\infty}^{\infty}f_X(y)f_Y(y-z)dy$. We will denote this operation as $f_Z(z) = f_X\ast f_Y$ where we use $\ast$ as the symbol for the \textit{linear convolution} operator.

\paragraph{} If we use characteristic functions instead, then the summation of these random variables can be performed via multiplication instead: $\phi_Z(\tau) = \phi_X(\tau)\phi_Y(\tau)$ .

\subsubsection{Convolution and deconvolution}
\label{subsubsection:conv_and_deconv}
Convolution and deconvolution are two operators, each the inverse of the other. In this work we will denote convolution as $\ast$ and deconvolution as $\deconv$. The properties of convolution and deconvolution are similar to that of multiplication and division.  Let $f,g,h$ denote probability density functions or discrete lattice functions (i.e. arrays or vectors). Then, it is known that convolution is:
\begin{itemize}
	\item  Commutative: $f \ast g = g \ast f$
	\item Associative: $f \ast (g \ast h) = (f \ast g) \ast h$
	\item Distributive: $f \ast (g + h) = (f\ast g) + (f \ast h)$
\end{itemize}

Additionally, we have that $f \ast \delta_{x,0} = f$ and that $a(f \ast g) = (af)\ast g$, where $a$ is a constant, and $\delta_{x,0}$ is the Dirac delta if $f,g,h$ are continuous, or the Kronecker delta if discrete (lattices).

\paragraph{} With respect to equations involving deconvolutions, then deconvolution has the properties:
\begin{itemize}
	\item Right distributive: $((f + g) \deconv h) = (f \deconv h) + (g \deconv h)$, but not left distributive: $h \deconv (f+g) \neq  (f \deconv h) + (g \deconv h)$. This is like division.
	\item Yields the identity: $f \deconv f = \delta_{k,0}$, where $\delta_{k,0}$ is the Dirac or Kronecker delta, depending whether $f$ is continuous or discrete respectively.
\end{itemize}
Therefore for example we have that $((f + g) \pentagon h) - ((f + k) \pentagon h) \equiv (g-k) \pentagon h$.

\paragraph{} Additionally, the following properties hold between convolution and deconvolution:
\begin{itemize}
	\item $(f \ast g) \pentagon h = (f \pentagon h) \ast g$. \paragraph{}\textit{Proof: } Let $\mathcal{F}\{f\} = F$ be the Fourier transform of the pdf $f$. Then $\mathcal{F}\{(f \ast g) \pentagon h\} = \frac{F\cdot G}{H} = \frac{F}{G}H = \mathcal{F}\{(f \pentagon h) \ast g\}$. Applying the inverse Fourier transform both sides of this yields the result.\qed
	
	\item $|f\ast g| \leq |f|\ast|g|$ and $|f \pentagon g| \geq |f|\pentagon|g|$.
	\label{proof:conv_deconv_inequality}
	\paragraph{} \textit{Proof: } For convolution we have $|f \ast g|(t) = |\int_{x}dx\cdot f(t-x)g(x)| \leq \int_{x}dx\cdot|f(t-x)||g(x)| = (|f| \ast |g|)(t)$. For deconvolution, then since we have $|f| = |g \ast h| \leq |g| \ast |h|$, deconvolving both sides by $|h|$ yields, $|f| \deconv |h| \leq |g|$. Also we have $|g| = |f \deconv h|$. Therefore $|g| = |f \deconv h| \geq |f| \deconv |h|$.     \qed
	
\end{itemize}    

\section{Summary} In this chapter we defined DTMCs and MRMs, and presented an introduction to model checking with temporal logics, more specifically PCTL and PRCTL. For particular problems, we explained their solutions, and by doing so we introduced much of the notation we will be using in this thesis.

\paragraph{} More importantly, we highlighted that reachability problems are one of the main problems for model-checking DTMCs as the event $\Phi_1 \ \text{U}\ \Phi_2$ for a DTMC $\mathcal{M}_1$ can be reduced to some an event $\Diamond B$ in a transformed DTMC $\mathcal{M}_2$. This holds true for sMRMs too.

\paragraph{} In the next chapter, we lay the theoretical foundations of sMRMs and define several problems of interest over them.

\newcommand{\mrmpdf}{Pr(r\ \cap\ s \vDash \Diamond B)}
\newcommand{\mrmpdfsub}{Pr}
\newcommand{\mrmmixeddist}{\sum_{\hat{\pi} \in \Pi}Pr(\hat{\pi})\mrmpdfsub(r | \hat{\pi})}
\newcommand{\mrmfouriers}{\phi_{s}(\tau)}
\newcommand{\mrmfourierrews}{\phi_{Rew(\hat{\pi})}(\tau)}
\newcommand{\mrmpdfcostbounded}{Pr(s \vDash \Diamond_{\leq r} B)}
\newcommand{\mrmexpected}{ExpRew(s \vDash \Diamond B)}

    \chapter{Stochastic Markov Reward Models (sMRMs)}

\section{Introduction}

\paragraph{} In this chapter, we introduce the theory on stochastic MRMs (sMRMs), an extension of the traditional MRMs which allows rewards to be random variables (or random vectors). The theory on sMRM was introduced into probabilistic model checking recently by \cite{bacci2019model}. In the literature, sMRMs may be known previously as Markov processes with \textit{random rewards} \cite{benito1982calculating}, \textit{statistical flowgraphs} \cite{huzurbazarstatistical}, simply rewards defined over Markov chains \cite{verhoeff2004reward}. Whilst we present here derivations and theory of our own, literature on Markov chains with rewards and SMPs exists that share a similar theory, for example see the previously cited articles  and \cite{bradley2004hypergraph}.

\begin{figure}\centering

	\begin{tikzpicture}[shorten >=2pt,node distance=3cm,auto]
	\begin{scope}[]
	\node[state,initial] (l0) at (0,0) {$l_0$};
	\node[state, right of=l0] (l1)  {$l_1$};
	\node[state,accepting, above of=l1]  (l2) {$l_2$};
	\node[state, right of=l1] (l3)  {$l_3$};
	\node[state,accepting, right of=l3] (l4)  {$l_4$};
	
	\end{scope}
	
	\begin{scope}[]
	\path [->] (l0) edge node {\tt $0.2$} (l1);
	\path [->] (l0) [bend left] edge node {$0.8$} (l2);
	\path [->] (l1) edge node {$0.9$} (l3);
	\path [->] (l1) [bend right] edge node {$0.1$} (l2);
	\path [->] (l2) edge [loop above] node {$1$} (l2);
	\path [->] (l3) edge[bend left] node {$0.4$} (l1);
	\path [->] (l3) edge node {$0.6$} (l4);
	\path [->] (l4) edge [loop above] node {$1$} (l4);
	\end{scope}
	\end{tikzpicture}
	\caption{A Markov chain with two absorbing states $l_2,l_4$.}
	\label{fig:example1}
\end{figure}
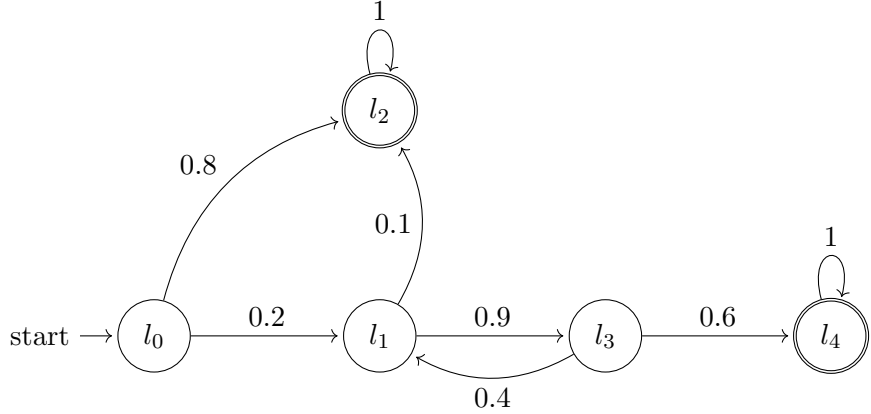

\label{study1:initial problem}
\paragraph{} The Markov chain in Figure \ref{fig:example1} captures the movement of a robot within a building. The robot begins in position $l_0$ of the building, and by moving probabilistically between places, it ends up eventually in $l_2$ or $l_4$. These states are its final destination, and we want to determine if it is capable of reaching such positions without running out of energy.

\paragraph{} The energy cost of each transition is a random variable. Assuming we know these random variables in advance, we are able to incorporate them into our Markov chain to form stochastic Markov reward models.   

\begin{definition}
	\label{definition:sMRM}
	A stochastic Markov reward model (sMRM) is a tuple (\textit{M}, \textit{rew}) where:
	\begin{itemize}
		\item $\mathcal{M}$ is a DTMC, a tuple $(S, \textbf{P}, i_{init})$.
		\item \textit{rew}: $(S \times S) \rightarrow f$ is a map between every state transition $t \in S \times S$ of the DTMC, and a probability distribution $f$.
	\end{itemize}
\end{definition}

\paragraph{} If the rewards are discrete random variables, then $f$ is a probability mass function (pmf), defined over lattices i.e. $f: \mathbb{L} \rightarrow [0,1]$, where $ \mathbb{L}$ is a countable lattice subset of  $\mathbb{R}$. $f$ is a pdf if the rewards are continuous instead.

\paragraph{} $f$ may be a joint distribution in the case where the rewards are random vectors (or multivariate). Whilst we focus however on the univariate case, mutually independent piecewise multivariate rewards can be solved by considering each variable separately.

\paragraph{} The \textit{sample space} of the sMRM is represented as the product $\Pi \times \mathbb{R}$, where $\Pi$ is the set of all paths that the underlying DTMC can generate and $\mathbb{R}$ is the set of values for the reward. The basis events for $\Pi$ are the cylinder sets, whereas the basis events for rewards is $\mathcal{E}= \{(-\infty,x] | x \in \mathbb{R} \}$. The event space for rewards is the Borel $\sigma$-algebra, whilst the event space for paths is the smallest $\sigma$-algebra over the set of all cylinder sets. Let us denote them as $\Sigma_1$ and $\Sigma_2$ respectively. Then, we can define the product $\sigma$-algebra as $\Sigma_1 \times \Sigma_2$, or the product event space.

\paragraph{} The probability space can be defined as  ($\mathscr{S} = \Pi \times \mathbb{R}, \mathscr{E} = \Sigma_1 \times \Sigma_2, Pr$).

\paragraph{} For the majority of this work, we enforce that for all states in $s \in S$, the probability of reaching $B$ is one, i.e. $$Pr(s \vDash \Diamond B) = 1$$ A manner to remove this restriction will be detailed later (see Sec. \ref{subsection:partial_densities}). Also, for the remainder of the thesis, we only solve for the case where each reward random variable $rew(s,t)$ \textbf{is strictly non-negative}, for all pairs $s,t \in S^2$.

\section{Reachability problems}
\label{sec:reachability_problems}

\paragraph{} Given an sMRM $R_M$, there are several reachability problems we can try to solve:
\begin{enumerate}[wide, labelwidth=!, labelindent=5pt]
	\item $Pr(s \vDash \Diamond_{=r} B)$ \label{question1}- The probability of accumulating $r$ reward, and \textit{eventually reaching} $B \subseteq S$ starting from a state $s \in S$. However, we will use the notation $Pr(Rew = r\ \cap\ s \vDash \Diamond B)$ instead, and it will frequently be seen shortened as $$Pr(r\ \cap\ s \vDash \Diamond B)$$
	
	\paragraph{} Let $\Pi$ be the set of paths starting from $s$ and ending in $B$, (i.e. those that satisfy $s \vDash \Diamond B$) then we have that
	\begin{flalign}
	\label{eqn:pdf_convex_sums}
	\mrmpdf = Pr(r\ \cap\ \Pi) = \sum_{\hat{\pi} \in \Pi}Pr(r  \cap \hat{\pi}) = \sum_{\hat{\pi} \in \Pi}Pr(\hat{\pi})Pr(r | \hat{\pi})
	\end{flalign}
	where $Pr(r | \hat{\pi})$ is the probability density (or mass) function of the accumulated reward given a particular path $\hat{\pi} = s_0,\cdots,s_n$, i.e. it is the pdf of the random variable $Rew(r;\hat{\pi}) =  rew(s_0,s_1)+rew(s_1,s_2)+rew(s_2,s_3)+\cdots +rew(s_{n-1},s_n)$.
	
	\paragraph{} If $\sum_{\hat{\pi} \in \Pi}Pr(\hat{\pi}) = 1$, then $\mrmpdf$ is a \textit{convex combination} of probability density functions; a mixture distribution. This is what we call the \textit{first-passage reward density (or mass function)}, the focus of our work. In the model checking literature, this may be better termed as \textit{reachability reward density}. This may be a pmf or pdf depending on whether it is a discrete-reward sMRM or continuous-reward. \textbf{The semantics is slightly different} to PRCTL for MRMs as now we can interpret $Pr(Rew = r\ \cap\ s \vDash \Diamond B)$  as a function where $r$ is the variable of the function, i.e. we can denote it as $f_{s \vDash \Diamond B}(r)$ instead, or $f_s(r)$ if the property is obvious.
	
	\item $\mrmpdfcostbounded$ - The probability of reaching a set $B$ from $s$ with the cumulated reward being less than or equal to $r$. This is equivalent to \textit{reward-bounded reachability probability} from PRCTL. Alternative notation would be $Pr(Rew \leq r\ \cap\ s \vDash \Diamond B)$. If interpreted as a function, this is the cumulative distribution function (cdf) of $Pr(Rew = r\ \cap\ s \vDash \Diamond B)$, and is denoted as $F_{s \vDash \Diamond B}(r)$ or $F_{s}(r)$.
	
	\item $Pr(s \vDash \Diamond_{\mathbb{E}[\leq r]} B)$ - the probability of reaching a set $B$ from $s$, with the mean cumulated reward being less than or equal to $r$. An example of this appeared in \cite{bacci2019model} for sMRMs.
	
	\item $\mrmexpected$ \label{question4}- The expected amount of reward accrued when starting from a state $s$ before reaching the target set $B$. This value is a scalar in $\mathbb{R}$. This is the expected reward for reachability from PRCTL.
	
	\item $Pr(s \vDash \Diamond_{\leq ?} B) > p$ - asks for the \textit{minimal} reward bound $r$ such that the probability of reaching $B$ from $s$ is greater than $p$. This is known as a \textit{quantile query} \cite{baier2014energy}.
	\item $Pr(s \vDash \Diamond^{\leq t}_{=r} B)$ - The \textit{constrained} (or step-bounded) variant of $Pr(s \vDash \Diamond_{=r} B)$. Alternative notation for this is $Pr(r\ \cap\ s \vDash \Diamond^{\leq t} B)$, which is a pdf/pmf. A cdf variant of this can be constructed.
	\item $Pr(s \vDash \Circle_{=r} B)$ - The probability of reaching $B$ in the next step, beginning from $s$, whilst accumulating $r$ reward. The alternative notation for this is $Pr(r\ \cap\ s \vDash \Circle B)$, again a pdf/pmf, and we can have a cdf variant of this as well.
\end{enumerate}

\paragraph{} Hence, after learning how to solve for these properties, we can define PRCTL for sMRMs. The grammar is identical as for regular MRMs (see Section \ref{grammar_PRCTL}). The PCTL subset of PRCTL can be resolved with DTMC algorithms. Resolution of the expected value problems  (3, 4) for sMRMs is almost identical for those of MRMs. Only, (1, 2, 5, 6, 7) are differently computed, with (1) being able to borrow solution ideas from SMP theory.

\paragraph{} The focus of this thesis will be on problem (1) above. As for (2), then this is just the numerical integration of (1).  As for (5), the quantile query is the smallest $r$ of the cdf $\mrmpdfcostbounded$, such that the probability is greater than $p$. Hence having solved (1), we have the ability to derive this. The solution for (6, 7) can be derived indirectly from the power method algorithm for solving (1). This is shown in Theorem \ref{proof:time_bounded_reachability_sMRMs}.

\paragraph{} We proceed to derive the solutions for the reachability problems of (3,4,2,5,1) in that respective order. We have left (1) for last as this is the main focus and will be elaborated.

\subsection{$\mrmexpected$ properties} This is computed as
\begin{flalign*}
ExpRew(s \vDash \Diamond B)
&= \mathbb{E}[\sum_{\pi \in \Pi_{\Diamond B}.s} Rew(\pi) Pr(\pi)] &\\
&= \sum_{\pi \in \Pi_{\Diamond B}.s} \mathbb{E}[Rew(\pi)] Pr(\pi) &\\
&= \sum_{t \in S} \sum_{\pi \in \Pi_{\Diamond B}.t} \mathbb{E}[Rew(s.\pi)] Pr(s.\pi) &\\
&= \sum_{t \in S} \sum_{\pi \in \Pi_{\Diamond B}.t} \mathbb{E}[(rew(s,t) + Rew(\pi))]\cdot \textbf{P}(s,t)Pr(\pi) &\\
&= \sum_{t \in S} \textbf{P}(s,t) \sum_{\pi \in \Pi_{\Diamond B}.t} \mathbb{E}[rew(s,t)] Pr(\pi) + \mathbb{E}[Rew(\pi)] Pr(\pi) &\\    
&= \sum_{t \in S} \textbf{P}(s,t) \Big(\sum_{\pi \in \Pi_{\Diamond B}.t} \mathbb{E}[rew(s,t)] Pr(\pi) + \sum_{\pi \in \Pi_{\Diamond B}.t}\mathbb{E}[Rew(\pi)] Pr(\pi) \Big) &\\
&= \sum_{t \in S} \textbf{P}(s,t) \Big(\mathbb{E}[rew(s,t)] \sum_{\pi \in \Pi_{\Diamond B}.t}  Pr(\pi) + ExpRew(t \vDash \Diamond B) \Big) &\\    
&= \sum_{t \in S} \textbf{P}(s,t)\Big(\mathbb{E}[rew(s,t)]Pr(t \vDash \Diamond B) + ExpRew(t \vDash \Diamond B)  \Big)
\end{flalign*}
or equivalently since $Pr(s \vDash \Diamond B) = 1$ for all $s \in S$,
\begin{flalign*}
ExpRew(s \vDash \Diamond B)
&= \sum_{t \in S} \textbf{P}(s,t)\Big(\mathbb{E}[rew(s,t)] + ExpRew(t \vDash \Diamond B)  \Big)
\end{flalign*}
which is just a system of linear equations as before. \qed

\paragraph{} Therefore, all that is needed is to be able to compute the expected value of all rewards, e.g. $\mathbb{E}[rew(s,t)]$, then the computation is identical to that of regular MRMs.

\subsection{$Pr(s \vDash \Diamond_{\mathbb{E}[\leq r]} B)$ properties} It can be shown that
\begin{flalign*}
Pr(s \vDash \Diamond_{\mathbb{E}[\leq r]} B)
&= \sum_{t \in S}\textbf{P}(s,t)Pr(t \vDash \Diamond_{\mathbb{E}[\leq (r - rew(s,t))]}\ B)
\end{flalign*}
using a similar derivation for MRMs described in Section \eqref{proof:reward_bounded reachability}. Therefore the algorithm for solving this property can borrow the algorithm from Section \eqref{proof:reward_bounded reachability}. Like the previous property, we are required first to compute the expected value of all rewards $\mathbb{E}[rew(s,t)]$. \qed

\subsection{$Pr(s \vDash \Diamond_{\leq r} B)$ properties}  Assume for now we can already compute $Pr(r \cap s \vDash \Diamond B)$. Let us denote this pdf/pmf as $f_s(x)$ for short, and let $\Pi$ be the set of finite paths starting from $s$ and ending in $B$, (i.e. those that satisfy $s \vDash \Diamond B$). Then we have
$$
Pr(s \vDash \Diamond_{\leq r} B) = \int_{x=0}^{r}f_{s}(x)dx
$$
since
\begin{flalign*}
Pr(s \vDash \Diamond_{\leq r} B) &= \sum_{\hat{\pi}}Pr(\hat{\pi} \cap (Rew \leq r))= \sum_{\hat{\pi}}Pr(\hat{\pi})Pr(Rew \leq r | \hat{\pi})    &&\\
&= \sum_{\hat{\pi}}Pr(\hat{\pi})\int_{x=0}^{r}f_{Rew(\hat{\pi})}(x)dx= \int_{x=0}^{r}\sum_{\hat{\pi}}Pr(\hat{\pi})f_{Rew(\hat{\pi})}(x)dx  &&\\
&= \int_{x=0}^{r}f_{s}(x)dx
\end{flalign*}

\paragraph{} Hence, if we have already computed $f_s(x)$, we can compute reward-bounded reachability probabilities by integration.

\paragraph{} Note that if we compute the cdf
$$F_s(r) = \int_{x=0}^{r}f_{s}(x)dx$$
then we can obtain
$$Pr(s \vDash \Diamond_{[a,b]} B) = F_s(b) - F_s(a)$$
the probability that $B$ is reached from $s$ with reward accumulated only within the interval $[a,b]$.

\subsubsection{Multivariate mutually-independent rewards} Consider an sMRM problem where the rewards are random vectors of dimension $n$. Then let \textbf{Rew} be the random vector, denoting the multivariate accumulated reward. Firstly, we have  $Pr(s \vDash \Diamond_{\leq r} B) = Pr(Rew \leq r \cap s \vDash \Diamond B)$. Then define
$$Pr(\textbf{Rew} \leq \vec{r}, s_0 \vDash B)$$ to be the probability that we can reach the set of states B from $s_0$ with reward accumulated under or equal to $\vec{r}$, where $\textbf{Rew} \leq \vec{r} = \cap_{i =1}^{n} Rew_i \leq r_i$. If the components of the random vectors are mutually independent of each other, then
$$Pr(\textbf{Rew} \leq \vec{r}, s_0 \vDash B)  = Pr(\cap_{i=1}^{n}Rew_i \leq r_i, s_0 \vDash B) = \prod_{i=1}^{n}Pr(Rew_i \leq r_i, s_0 \vDash B)$$
where for each $i$, $Pr(Rew_i \leq r_i, s_0 \vDash B)$ can be solved independently. \qed

\paragraph{} The solution above also gives us a means to compute quantile queries for multivariate independent rewards, which was studied for regular MRMs in \cite{haase2017computing}.

\subsection{$Pr(s \vDash \Diamond_{\leq ?} B) > p$ quantile queries} The quantile query can be solved via interpolation of $F_s(x)$. Given $F_s(x)$, we are to find the smallest $r$ such that $F_s(r) > p$. A naive algorithm based on trial-and-error is to sample the distribution and compare it to $p$ and then move towards the direction of $p$ knowing that $F$ is monotonic. The algorithms for computing $F_s(x)$ in our thesis only computes $F_s(x)$ for values $x \in [0,k]$. If $ F_s(k)< p$, we have to recompute the problem with larger $k$, otherwise an algorithm can be used to find the point. However, whilst typically $\lim_{r \rightarrow \infty}F_s(r) = 1$, if $Pr(s \vDash \Diamond B) < 1$, then this no longer holds. See Section \ref{subsection:partial_densities}.

\paragraph{} However, in the setting where this is true, i.e. when $Pr(s \vDash \Diamond B) = 1$, a guide to arrive at a sufficiently large interval $[0,k]$ at which $F_s(k) = 1$ is to use the expected value and variances of the first passage-reward density, for example $$k = ExpRew(s \vDash \Diamond B) + \gamma VarRew(s \vDash \Diamond B)$$ where $\gamma$ is a parameter used to increase the range involved. The computation for $VarRew(s \vDash \Diamond B)$ can be found in \cite{verhoeff2004reward}.

\subsection{$\mrmpdf$ properties}
\paragraph{System of convolution equations} We present a computation of $\mrmpdf$ for each state of a finite sMRM with a common graph-based technique. Let us first state the solution, and leave the derivation till later.

\paragraph{} Firstly for any states in $B$, then their cumulated rewards are assigned to zero. We also make them entirely self-absorbent, and we fix the reward transition distributions (for these self-loops) as the \textit{zero distribution}, expressed by the Dirac delta (or Kronecker in the discrete case) $\delta_{x,0}(x)$. Consequently then, for any path absorbed by $B$, further reward is not cumulated in the sMRM process after entering $B$.

\paragraph{} We can shorten $\mrmpdf$ to $f_{s,B}(r)$. However, if we write $f_s(r)$ it is assumed that the set of states to reach is $B$. Then  define $f_{rew(s,t)}(r)$ as the reward pdf for the transition $s \rightarrow t$. Let $Pre^*(B)$ denotes the set of states that reach $B$ with non-zero probability. Then define $S_? = Pre^*(B)\char`\\B$, to be the set of states that reach $B$ but are not in it. Firstly, since for each state $u \in B$ its reward  is zero, we have
$$f_u(r) = \delta_{x,0}(r)$$
which is the zero distribution (returns zero with probability one). Then for the remaining states, we have
\begin{flalign}
\label{eqn:linear_sys_eqn_pdf}
f_s(r) &= \sum_{t \in S}(\textbf{P}(s,t)f_{rew(s,t)} \ast f_{t})(r)  \notag &\\
&= \sum_{t \in S_?}(\textbf{P}(s,t)f_{rew(s,t)} \ast f_{t})(r) + \sum_{u \in B}\textbf{P}(s,u)f_{rew(s,u)}(r) \ast f_{u}(r)  \notag &\\
&= \sum_{t \in S_?}(\textbf{P}(s,t)f_{rew(s,t)} \ast f_{t})(r) + \sum_{u \in B}\textbf{P}(s,u)f_{rew(s,u)}(r) \ast \delta_{x,0}(r)  \notag &\\    
&= \sum_{t \in S_?}(\textbf{P}(s,t)f_{rew(s,t)} \ast f_{t})(r) + \sum_{u \in B}\textbf{P}(s,u)f_{rew(s,u)}(r)
\end{flalign}
for all $s \in S_?$. This then yields what we call a \textit{system of convolution equations}.

\paragraph{} Then since $S_? \cup B = S$ (due to all states reaching $B$ with probability 1), we have a solution for all states. Intuitively, the expression means that the reward accumulated from a state $s$ before arriving at $B$, is the convex combination of the rewards accumulated by all states that $s$ can immediately transition into (the combination weighted by the probability of entering these states).

\paragraph{} Instead of solving a system of convolution equations, we can convert it to a set of systems of linear equations via the Fourier transform. Let the characteristic function transform operator (or Fourier transform) be represented as $\mathcal{F}$, and the transform applied to a function $f(x)$ be denoted as $\mathcal{F}\{f(x)\}$. Then, the Fourier transform of $\mrmpdf$ (which is $f_s$ from \eqref{eqn:linear_sys_eqn_pdf}) is
\begin{flalign}
\label{eqn:char_representation}
\mathcal{F}\{\mrmpdf\} = \mathcal{F}\{\mrmmixeddist\}
= \sum_{\hat{\pi} \in \Pi}Pr(\hat{\pi})\mathcal{F}\{\mrmpdfsub(r | \hat{\pi})\}
\end{flalign}
which always exists, since Fourier transforms of pdfs and pmfs always exist.

\paragraph{}From now on, we will denote $\mathcal{F}\{\mrmpdf\}$ as $\phi_{s}(\tau)$. And $\mathcal{F}\{\mrmpdfsub(r | \hat{\pi})\}$ as $\phi_{Rew(\hat{\pi})}(\tau)$. We also assume that the temporal logic property is always $\Diamond B$. Now we rewrite the equation \eqref{eqn:char_representation} above as
$$ \mrmfouriers =  \sum_{\hat{\pi} \in \Pi} Pr(\hat{\pi})\mrmfourierrews  $$

\paragraph{} The results above can be derived knowing the properties for Fourier transforms over functions; its linearity with respect to constants and addition. Alternatively, using the law of total expectation: using traditional probability notation, then for two dependent random variables $X,Y$ we have $\mathcal{F}\{X\} = E[e^{itX}] = \mathbb{E}[\mathbb{E}[e^{itX}|Y] = \sum_{\hat{y} \in Y}P(\hat{y})*\mathbb{E}[e^{itX}| \hat{y}] = \sum_{\hat{y} \in Y}P(\hat{y})*\mathcal{F}\{X|y\}$. If this is not clear, let $X \sim Pr(Rew=r)$ with the domain being the real line, and $Y \sim Pr(\hat{\pi})$ with the domain being $\Pi$, then the above holds.    

\paragraph{} If using characteristic functions instead, i.e. using $\mrmfouriers =     \mathcal{F}\{\mrmpdf\}$, we have
\begin{flalign}
\label{eqn:linear_sys_eqn_cf}
\phi_s(\tau) = \sum_{t \in S_?}\textbf{P}(s,t)\phi_{rew(s,t)}(\tau)\phi_t(\tau) + \sum_{u \in B}\textbf{P}(s,u)\phi_{rew(s,u)}(\tau)
\end{flalign}
providing us with a system of linear equations, for each $\tau$.  And for each $u \in B$,
$$\phi_u(r) = \mathcal{F}\{\delta_{r,0}(r)\}=1$$

\paragraph{} There is a one-to-one correspondence between the two formulas \eqref{eqn:linear_sys_eqn_cf} and \eqref{eqn:linear_sys_eqn_pdf}. This is due to the Fourier transform being bijective.

\begin{example}
	\label{example:study1_manual_solve}
	We now solve the problem presented in Sec. \ref{study1:initial problem}. The characteristic functions will be used instead of pdfs.
	
	\paragraph{} For each transition $s\rightarrow t$  of the system we assign the symbolic reward $rew(s,t)(r)$, of which the Fourier transform is $\mathcal{F}\{rew(s,t)(r)\} = \phi_{rew(s,t)}(\tau)$. The goal is to compute $\mrmfouriers$, for the temporal property $\Diamond B$, where $B =\{l_2, l_4 \}$. In the following, we may drop $(\tau)$ from our notation for a simpler representation.
	
	\paragraph{} Firstly, since $B = \{l_2,l_4\}$, then we make them self-looping and assign the loop reward transition the characteristic function of the Dirac delta. Then, using equation \eqref{eqn:linear_sys_eqn_cf}, we arrive at the following linear system of equations for this sMRM:
	
	\begin{figure}[H]
		\centering

		\begin{tikzpicture}[shorten >=2pt,node distance=3cm,auto]
		\begin{scope}[]
		\node[state,initial] (l0) at (0,0) {$l_0$};
		\node[state, right of=l0] (l1)  {$l_1$};
		\node[state,accepting, above of=l1]  (l2) {$l_2$};
		\node[state, right of=l1] (l3)  {$l_3$};
		\node[state,accepting, right of=l3] (l4)  {$l_4$};
		
		\end{scope}
		
		\begin{scope}[]
		\path [->] (l0) edge node {\tt $0.2\phi_{rew(0,1)}$} (l1);
		\path [->] (l0) edge[bend left]node {$0.8\phi_{rew(0,2)}$} (l2);
		\path [->] (l1) edge node {$0.9\phi_{rew(1,3)}$} (l3);
		\path [->] (l1) edge [bend right] node {$0.1\phi_{rew(1,2)}$} (l2);
		\path [->] (l2) edge [loop above] node {$\phi_{\delta}$} (l2);
		\path [->] (l3) edge[bend left] node {$0.4\phi_{rew(3,1)}$} (l1);
		\path [->] (l3) edge node {$0.6\phi_{rew(3,4)}$} (l4);
		\path [->] (l4) edge [loop above] node {$\phi_{\delta}$} (l4);
		\end{scope}
		\end{tikzpicture}
		\caption{A stochastic Markov reward model (sMRM) with two absorbing states $l_2,l_4$. Each transition has been annotated with the probabilities multiplied by the Fourier transform of the reward/cost distribution (which we have left as variables rather than define them explicitly).}
	\end{figure}
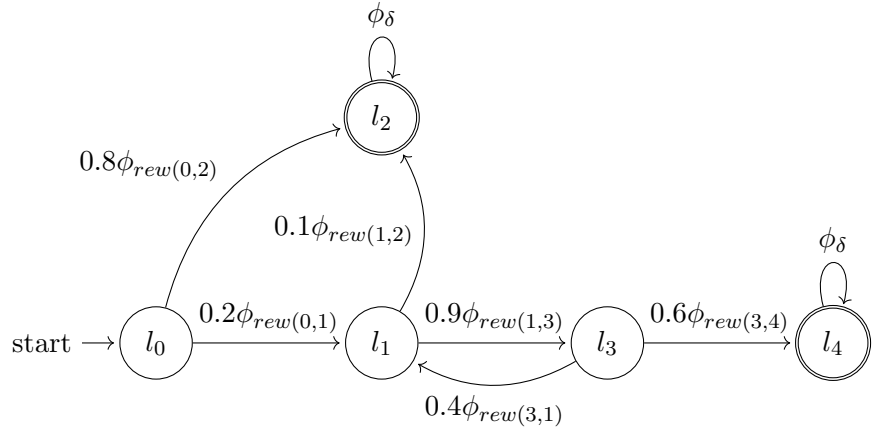
	
	\begin{flalign*}
	\phi_0 &= 0.2\phi_{rew(0,1)}\phi_1 + 0.8\phi_{rew(0,2)}\phi_2 \\
	\phi_1 &= 0.1\phi_{rew(1,2)}\phi_2 + 0.9\phi_{rew(1,3)}\phi_3 \\
	\phi_2 &= \phi_{\delta} \\
	\phi_3 &= 0.4\phi_{rew(3,1)}\phi_1 + 0.6\phi_{rew(3,4)}\phi_4 \\
	\phi_4 &= \phi_{\delta}  
	\end{flalign*}
	One can then derive manually that the solution for $\phi_0,\phi_1,\phi_3$ (i.e. characteristic functions for state $s_0,s_1,s_3$) is
	\begin{flalign*}
	\phi_0 &= 0.2\phi_{rew(0,1)}\phi_1 +  0.8\phi_{rew(0,2)}\phi_{\delta} \\
	\phi_3 &= 0.4\phi_{rew(3,1)}\phi_1 + 0.6\phi_{rew(3,4)}\phi_{\delta} \\
	\phi_1 &= \frac{0.1\phi_{rew(1,2)} + 0.9\phi_{rew(1,3)}0.6\phi_{rew(3,4)}}{1 - (0.9\phi_{rew(1,3)}0.4\phi_{rew(3,1)})}
	\end{flalign*}
\end{example}

\paragraph{} The manual solution here can be automated via a symbolic solver such as Sympy, but as we shall see later, we will need to move away from solving the solution symbolically due to symbolic solvers not being very scalable. The reason is simply the intractability that comes with symbolic solvers.

\paragraph{} We proceed to detail how equation \eqref{eqn:linear_sys_eqn_cf} was derived. It derives indirectly the equation \eqref{eqn:linear_sys_eqn_pdf}.

\begin{theorem}[Derivation of the set of systems of linear equations]
	\label{theorem:derivation_of_system}
	
	We are given an sMRM  with state space $S$, with a set of goal states $B \subset S$. Let it be the case that every state in $S$ can reach $B$. Then we assign to each state $u \in B$ zero rewards, such that the first-passage reward density is the (Dirac or Kronecker) delta $\delta(r)$. Let $S_? = S\char`\\B$ be the set of states not in $B$ (but can reach $B$).  For each $s \in S$, let $\mrmfouriers$ be the Fourier transform of $\mrmpdf$, the first-passage reward density, i.e. $\mrmfouriers = \mathcal{F}\{\mrmpdf\}$. Thus for each state $u \in B$, we immediately have $\phi_u{(\tau)} =\mathcal{F}\{\delta\}(\tau) = \phi_{\delta}{(\tau)} = 1 $.
	\paragraph{} Then, for each state $s \in S_?$, we have the equivalence
	$$\phi_s(\tau) = \sum_{t \in S_?}\textbf{P}(s,t)\phi_{rew(s,t)}(\tau)\phi_t(\tau) + \sum_{u \in B}\textbf{P}(s,u)\phi_{rew(s,u)}(\tau)$$
	
	\paragraph{} \textit{Proof:} Our derivation methodology is analogous to that found in \cite{verhoeff2004reward}, except there the author derived solutions for the expected reachability reward (see Question \ref{question4}), and the variance of the reachability reward. This type of derivation allows us to obtain the set equations via \textit{paths}.
	
	\paragraph{} Let us define $\Pi.s$ to be the set of (finite) paths starting in $s$ ending at a state in $B$. This is a shortened variant of the earlier notation which would require us to write $\Pi_{\Diamond B}.s$ instead.

	\paragraph{} Then for each $s \in S_?$:       
	\begin{flalign*}
	\phi_s(\tau) && \\
	&= \{  \text{definition of $\phi_s(\tau)$} \}&& \\
	&\sum_{\hat{\pi} \in \Pi.s} Pr(\hat{\pi}) \phi_{Rew(\hat{\pi})}(\tau) && \\
	&= \{ \text{letting $\hat{\pi} = s.\hat{\psi}$. Then since $s \notin B$, $\hat{\psi} \neq \emptyset$.}  \} && \\
	& \sum_{t \in S} \ \sum_{\hat{\psi} \in \Pi.t} Pr(s.\hat{\psi}) \phi_{rew(s.\hat{\psi})}(\tau) && \\
	&= \{ \text{definition of path probability: $Pr(s.\hat{\psi}) = \textbf{P}(s,t)Pr(\hat{\psi})$ } \} && \\
	& \sum_{t \in S} \ \sum_{\hat{\psi} \in \Pi.t} \textbf{P}(s,t)Pr(\hat{\psi}) \phi_{rew(s.\hat{\psi})}(\tau) && \\
	&= \{ \text{sum of 2 rvs using their cfs: } \phi_{rew(s. \hat{\psi})}(\tau) =\phi_{rew(s,t)}(\tau)\phi_{rew(\hat{\psi})}(\tau) \} && \\
	& \sum_{t \in S} \ \sum_{\hat{\psi} \in \Pi.t} \textbf{P}(s,t)Pr(\hat{\psi})\phi_{rew(s,t)}(\tau)\phi_{rew(\hat{\psi})}(\tau) && \\
	&= \{ \text{distribute \textbf{P}(s,t) and $\phi_{rew(s,t)}(\tau)$ outside since they are independent of $\hat{\psi}$}  \} && \\
	&\sum_{t \in S}\textbf{P}(s,t)\phi_{rew(s,t)}(\tau)\sum_{\hat{\psi} \in \Pi.t} Pr(\hat{\psi})\phi_{rew(\hat{\psi})}(\tau) && \\
	&= \{ \text{simplify since by definition$\sum_{\hat{\psi} \in \Pi.t} Pr(\hat{\psi})\phi_{rew(\hat{\psi})}(\tau) = \phi_t(\tau)$} \} && \\
	&\sum_{t \in S}\textbf{P}(s,t)\phi_{rew(s,t)}(\tau)\phi_{t}(\tau)
	\end{flalign*}      
	However, since we have that $\phi_u(\tau) \triangleq \phi_\delta(\tau)$ for all states $u \in B$, we  now have that
	\begin{flalign*}
	\phi_s(\tau) &= \sum_{t \in S_?}\textbf{P}(s,t)\phi_{rew(s,t)}(\tau)\phi_{t}(\tau) +  
	\sum_{u \in B}\textbf{P}(s,u)\phi_{rew(s,u)}(\tau)\phi_u(\tau) & \\
	&= \sum_{t \in S_?}\textbf{P}(s,t)\phi_{rew(s,t)}(\tau)\phi_{t}(\tau) +  
	\sum_{u \in B}\textbf{P}(s,u)\phi_{rew(s,u)}(\tau)      
	\end{flalign*}                  
	which completes the proof. \qed

	\paragraph{} \textit{Remark: } If we apply the inverse Fourier transform to the formula above, we arrive at $$f_s(r) = \sum_{t \in S_?}((\textbf{P}(s,t)f_{rew(s,t)}) \ast f_{t})(r) + \sum_{u \in B}\textbf{P}(s,u)f_{rew(s,u)}(r)$$
	which provides us with a system of convolution equations. \qed
\end{theorem}

\subsubsection{Matrix notation}
\label{subsubsec:definitions_of_terms}
We now present the solution to $\mrmpdf$ in matrix form. Define
\begin{align*}
\textbf{f} &\triangleq (f_s)_{s \in S_?} &\\
A &\triangleq (\textbf{P}(s,t))_{s,t \in S_?}&\\
\textbf{G} &\triangleq (G_{s,t})_{s,t \in S_?} \triangleq (f_{(rew(s,t))})_{s,t \in S_?}&\\
\textbf{h} &\triangleq  (h_s)_{s \in S_?} \triangleq (\sum_{u \in B}\textbf{P}(s,u)f_{rew(s,u)})_{s \in S_?}
\end{align*}

where each pdf of $f_s$, $f_{rew(s,i)}$ with $s \in S_?,\ i \in S$ can be represented as an analytical function, denoted symbolically (or algebraically). Alternatively, we can use vectors to represent them, such that
\begin{align*}
f_s &= (f_s(r))_{r \in \mathbb{R}} &\\
f_{rew(s,i)} &= (f_{rew(s,i)}(r))_{r \in \mathbb{R}}
\end{align*}
\paragraph{} This implies that $\textbf{G}$ is a `three-dimensional matrix', i.e. a \textit{hypermatrix}, and both \textbf{f},\textbf{h} are two-dimensional vectors, or \textit{hypervectors}. \textbf{We will use the latter representation for the majority of our work} however we do avoid it in the examples of this chapter. Using this representation, in the continuous case, \textbf{G} has dimensions $(\mathbb{R} \times S_? \times S_?)$, whilst \textbf{f,h} have dimensions $(\mathbb{R} \times S_?)$.

\paragraph{Notation} From now on, we will write regular vectors and matrices in the ordinary manner e.g. $f,g$ for vectors and $A,B$ as matrices. Hypervectors and hypermatrices are similar except bold instead, e.g. \textbf{h},\textbf{z} for hypervectors and \textbf{C},\textbf{D} for hypermatrices. This also applies to vectors and matrices of functions. There will be some exceptions to this, for example \textbf{P} is the probability matrix, or map over Markov chain transitions. Additionally, when we index hypermatrices or hypervectors, we generally make their letters light again, e.g. $\textbf{G}$, but $G_{s,t}$, or $G_{s,t}(r)$. We do not change cases for subscripted terms, e.g. not $a_{i,j}$, but $A_{i,j}$ for a matrix $A$.    

\paragraph{} Now we can rewrite \eqref{eqn:linear_sys_eqn_pdf} in matrix form as
\begin{flalign}
\label{eqn:lineq_matrix_trans_pdf}
\textbf{f} = (A\circ\textbf{G})\textcircled{$\ast$} \textbf{f} + \textbf{h}
\end{flalign}  

where \textcircled{$\ast$} denotes \textit{hypermatrix convolution}, analogous to matrix multiplication, and $A\circ\textbf{G}$ is the \textit{element-wise hypermatrix-matrix multiplication} between the matrix $A$ and hypermatrix $\textbf{G}$. Let us now explicate the definition of the operators on hypermatrices.

\begin{definition} Let  $\textbf{A},\textbf{B}, \textbf{C}$ be three hypermatrices, with dimensions $(\mathbb{R} \times N \times M), (\mathbb{R} \times M \times L)$ and $(\mathbb{R} \times M \times L)$ respectively. Then,
	
	\begin{enumerate}[wide, labelwidth=!, labelindent=5pt]
		\item     The \textit{hypermatrix convolution} product $\textbf{D} = \textbf{A} \textcircled{$\ast$} \textbf{B}$, is defined with dimensions $(\mathbb{R} \times N \times L)$, of which values $D_{i,j}(r)$, for all $i \in N, j \in L, r \in \mathbb{R}$, are determined by    
		$$D_{i,j}(r) =  \sum_{q=0}^{N-1}(A_{i,q} \ast B_{q,j})(r) $$
		
		\item The \textit{hypermatrix multiplication} product $\textbf{U} = \textbf{AB}$, is defined with dimensions $(\mathbb{R} \times N \times L)$ of which values $U_{i,j}(r)$, for all $i \in N, j \in L, r \in \mathbb{R}$, are determined by    
		$$U_{i,j}(r) = \sum_{q=0}^{N-1}A_{i,q}(r)B_{q,j}(r)$$
		or $U(r) = A(r)B(r)$ where the right-hand side is a matrix multiplication.
		
		\item The resulting sum/subtraction $\textbf{S} = \textbf{B} \pm \textbf{C}$, is defined with dimensions equal to $\textbf{B}$ of which values $S_{i,j}(r)$, for all $i \in M, j \in L, r \in \mathbb{R}$, are determined by    
		$$S_{i,j}(r) =  A_{i,j}(r) \pm B_{i,j}(r) $$
		or $S(r) = A(r) \pm B(r)$.
		
	\end{enumerate}
	
\end{definition}

\paragraph{} We now define operators between hypermatrices and matrices.

\begin{definition}
	Let $D$ be a matrix of dimensions $(N \times M)$ and $\textbf{A}$ the hypermatrix with dimensions $(\mathbb{R} \times N \times M)$. Then,
	\begin{enumerate}[wide, labelwidth=!, labelindent=5pt]
		\item The \textit{hypermatrix-matrix Hadamard} product $\textbf{H} = D\circ \textbf{A} = \textbf{A} \circ D$ is a hypermatrix with dimensions equal to \textbf{A}, with values $H_{i,j}(r)$ determined by
		$$ H_{i,j}(r) =  D_{i,j}A_{i,j}(r) = A_{i,j}(r)D_{i,j}$$
		for $i,j \in (N \times M)$ and $r \in \mathbb{R}$, or alternatively  $ H(r) = D \circ A(r) =  A(r) \circ D $
		where  $ D \circ A(r)$ is the Hadamard product (or element-wise product) between the two matrices.
		\item Let $E$ be a matrix with dimensions $(O \times N)$. The \textit{hypermatrix-matrix multiplication} product $\textbf{Z} = E\textbf{A}$ is a hypermatrix with dimensions equal to ($\mathbb{R} \times O \times M$), with values $Z_{i,j}(r)$ determined by
		$$Z_{i,j}(r) = \sum_{q=0}^{N-1}E_{i,q}A_{q,j}(r)$$
		or $Z(r) = EA(r)$.
		
		\item The \textit{hypermatrix-matrix sum/subtraction} $\textbf{F} = \textbf{A} \pm D = \mp D + \textbf{A}$, is defined with size equal to $\textbf{A}$, with values $F_{i,j}(r)$ determined by
		$$F_{i,j}(r) =  A_{i,j}(r) \pm  D_{i,j} =  \mp D_{i,j} + A_{i,j}(r)$$
		or $F(r) =  A(r) \pm  D =  \mp D + A(r)$.
	\end{enumerate}

\end{definition}

\paragraph{} Moving forward, if we use characteristic functions to represent our problem, then firstly define
\begin{align}
\label{eqn:definitions_terms_DFT}
\textbf{x} &\triangleq (x_s)_{s \in S_?} \triangleq (\phi_s)_{s \in S_?} &\\
\textbf{C} &\triangleq  (C_{s,t})_{s,t \in S_?} \triangleq (\phi_{(rew(s,t))})_{s,t \in S_?} \notag &\\
\textbf{d} &\triangleq  (d_s)_{s \in S_?} \triangleq (\sum_{u \in B}\textbf{P}(s,u)\phi_{rew(s,u)})_{s \in S_?} \notag
\end{align}

\paragraph{} Then the set of equations \eqref{eqn:linear_sys_eqn_cf} has the hypermatrix form,
\begin{flalign}
\label{eqn:lineqn_matrix_trans}
\textbf{x}(\tau) =(A\circ \textbf{C}(\tau))\textbf{x}(\tau) + \textbf{d}(\tau)
\end{flalign}

\paragraph{} The solution of which is just
\begin{flalign}
\label{eqn:matrix_solution}
\textbf{x} = (I - (A\circ\textbf{C}))^{-1}\textbf{d}
\end{flalign}

with $I$ being the identity matrix of size $|S_?| \times |S_?|$. The inverse of a square hypermatrix $\textbf{B}^{-1}$ (having an equal number of rows and columns) is defined such that $\textbf{B}^{-1}\textbf{B} = I\circ \textbf{1}$ with \textbf{1} having size equal to \textbf{B}. If we use the analytical (or symbolic) form of these characteristic functions, the general solution to this system will give us rational characteristic functions as we have seen in an earlier example and was shown earlier by for example \cite{vidal2004rational,warr2014numerical}.

\begin{proposition}[Weak and irreducibly diagonal dominance]
	\label{prop:wirr_diagonal_dom}     
	The system in \eqref{eqn:lineqn_matrix_trans} can be rewritten in the form
	\begin{flalign}
	\label{eqn:eqn_matrix_FT_re_arranged}
	(I - (A\circ\textbf{C}))\textbf{x} = \textbf{d}
	\end{flalign}
	\paragraph{} Then, the matrix $Z(\tau) = (I - A\circ \textbf{C}(\tau))$ is \textit{weak and irreducibly diagonally dominant}, i.e. there exists a set of states $q \in Q \subseteq S_?$, such that ${Z(\tau)}_{q,q}(\tau)$ exhibits \textit{strict} diagonal dominance, i.e.
	\begin{flalign}
	\sum_{t \in S_?\setminus  q}|Z(\tau)_{q,t}(\tau)| < |Z(\tau)_{q,q}(\tau)|
	\end{flalign}
	and for any remaining states $s \in S_?\setminus Q$, we have \textit{weak} diagonal dominance, i.e.  
	\begin{flalign}
	\sum_{t \in S_?\slash s}|Z(\tau)_{s,t}(\tau)| \leq |Z(\tau)_{s,s}(\tau)|
	\end{flalign}
	\paragraph{} \textit{Proof: } Making use of the fact that for all $\tau$, $|C_{s,t}(\tau)| \leq 1$ (see Definition \ref{def:fourier_transforms}), and that $A_{s,t} = \textbf{P}(s,t) \leq 1$, for all $s,t \in S_?^2$, we can easily deduce that $Z(\tau)$ is weakly diagonally dominant. To show that it is also irreducibly diagonally dominant, then for some row $i$ of $Z$, we have that $\sum_{t \in S_?} A_{i,t} < 1$. This is since there must be some state that reaches $B$ with non-zero probability. Knowing this, it is trivial to show that for any row such as $i$, $Z$ is strictly diagonally dominant.
\end{proposition}

\begin{theorem}[Unique solution]
	\label{theorem:unique sol}
	The system of equations \eqref{eqn:lineqn_matrix_trans} has a unique solution, and therefore \eqref{eqn:matrix_solution} is resolvable.
	
	\paragraph{} \textit{Proof:} From Proposition \ref{prop:wirr_diagonal_dom}, we know that $(I - A\circ \textbf{C}(\tau))$ is weak and irreducibly diagonally dominant. Then, this is a sufficient condition to prove that a unique solution exists, a proof of which can be found in \cite[Theorem 6.2.27]{horn_johnson_1985}. Hence, we have a unique solution for all $\tau$.
	
	\paragraph{} Corollary, this theorem proves that reward random variables positive and negative can be used together within a single sMRM and a unique solution will still exist. As another corollary, this implies that a unique solution exists in the time domain (i.e. without having transformed the equations), i.e. since the Fourier transform is bijective and always exists for random variables. \qed
	
	\begin{example}[Matrix (symbolic) approach]
		\label{example:study1_automatic}
		Continuing with the study problem in Example \ref{example:study1_manual_solve}, instead of solving the system manually, we can now solve it automatically via symbolic matrix solvers. One such solver is \textit{Sympy} \cite{sympy}.  Note that in this example, we use characteristic functions once more and also forgo the vector representation and use symbolic terms, as previously done in Example  \ref{example:study1_manual_solve}.
		
		\paragraph{} Firstly, let us denote terms in the form $p_{i,j}\phi_{rew(i,j)}$ as $\sigma_{i,j}$ for each pair $i,j \in S_?^2$. Then we rewrite the set of linear equations in matrix form  \eqref{eqn:lineqn_matrix_trans} for states in $S_? = {l_0,l_1,l_3}$  giving:
		$$\begin{pmatrix}
		\phi_0 \\
		\phi_1 \\
		\phi_3
		\end{pmatrix} =  \begin{pmatrix}
		0 & \sigma_{0,1} & 0 \\
		0 & 0 & \sigma_{1,3}\\
		0 & \sigma_{3,1} & 0
		\end{pmatrix}
		\begin{pmatrix}
		\phi_0 \\
		\phi_1 \\
		\phi_3
		\end{pmatrix}    +
		\begin{pmatrix}
		\sigma_{0,2} \\
		\sigma_{1,2} \\
		\sigma_{3,4}
		\end{pmatrix}
		$$
		Then by \eqref{eqn:matrix_solution},
		$$    \begin{pmatrix}
		\phi_0 \\
		\phi_1 \\
		\phi_3
		\end{pmatrix} =
		\begin{pmatrix}
		1 & -\sigma_{0,1} & 0 \\
		0 & 1 & -\sigma_{1,3}\\
		0 & -\sigma_{3,1} & 1
		\end{pmatrix}^{-1}
		\begin{pmatrix}
		\sigma_{0,2} \\
		\sigma_{1,2} \\
		\sigma_{3,4}
		\end{pmatrix}
		$$
		\paragraph{}Solving this system (symbolically) via     \eqref{eqn:matrix_solution}, gives:
		$$\begin{pmatrix}
		\phi_0 \\
		\phi_1 \\
		\phi_3
		\end{pmatrix}
		=
		\alpha
		\begin{pmatrix}
		\sigma_{0,1}(- \sigma_{3,4}\sigma_{1,3} - \sigma_{1,2}) +
		\sigma_{0,2}(\sigma_{3,1}\sigma_{1,3} - 1) \\
		\sigma_{3,4}\sigma_{1,3} + \sigma_{1,2} \\
		\sigma_{3,1}\sigma_{1,2} +         \sigma_{3,4} \\
		\end{pmatrix}    
		$$
		
		where $\alpha \triangleq \frac{1}{\sigma_{3,1}\sigma_{1,3} - 1}$. \qed
		
		\paragraph{} These equations are equal to the equations found manually earlier. This was proven using Sympy \cite{sympy} to show that these equations have the same solution set with the corresponding equations manually derived earlier.
	\end{example}

\end{theorem}   

\subsubsection{Partial passage-reward densities}
\label{subsection:partial_densities}
When presenting the solution to find the first-passage reward density  $f_s(r) = Pr(s \vDash \Diamond_{=r} B)$ (Question \ref{question1}),
we enforced earlier that all states $s \in S$ could reach $B$. And thus for the set of paths $\Pi$ that reach $B$ beginning from $s_0$, we had  $$\sum_{\hat{\pi} \in \Pi.s_0}Pr(\hat{\pi}) = 1$$
\paragraph{} Let us relax this assumption in two ways. Firstly, let there always be absorbing states within our sMRM which every state will eventually reach. Then we define this set of states as $Abs$, a subset of $S$, and let $\Pi_{\alpha}.s$ be the set of paths that end in a state $\alpha \in Abs$, beginning in $s$. We now get
$$\sum_{\alpha \in Abs}\sum_{\hat{\pi} \in \Pi_\alpha.s}Pr(\hat{\pi}) = 1 = \sum_{\alpha \in Abs} Pr(s \vDash \Diamond \alpha)$$

\paragraph{} If we set $B = Abs$, we have $Pr(s \vDash \Diamond B) = \sum_{\alpha \in Abs} Pr(s \vDash \Diamond \alpha) = 1$ and if $B \subset Abs$, then $Pr(s \vDash \Diamond B) \leq 1$. Therefore with respect to accumulated rewards, we have in the case where $B= Abs$,
$$\int dr \cdot  f_s(r) = \int dr \cdot \mrmpdf   = 1$$
\paragraph{} And if $ B \subset Abs$,    
\begin{flalign*}
\int_{-\infty}^{\infty}dr \cdot f_s(r) &= \int_{-\infty}^{\infty}dr \cdot \mrmpdf  &\\
&= \int_{-\infty}^{\infty}dr \cdot \sum_{\alpha \in B} \sum_{\hat{\pi} \in \Pi_{\alpha}.s}Pr(\hat{\pi})Pr(r | \hat{\pi}) &\\
&= \sum_{s \in B} \sum_{\hat{\pi_s} \in \Pi_{\alpha}.s}Pr(\hat{\pi}) \int_{-\infty}^{\infty}dr \cdot Pr(r | \hat{\pi}) &\\
&=  \sum_{s \in B} \sum_{\hat{\pi_s} \in \Pi_{\alpha}.s}Pr(\hat{\pi}) \int_{-\infty}^{\infty}dr \cdot f_{rew(\hat{\pi})}(r) &\\
&= \sum_{s \in B} \sum_{\hat{\pi_s} \in \Pi_{\alpha}.s}Pr(\hat{\pi}) \cdot 1 &\\
&= Pr(s \vDash \Diamond B) &\\
&\leq 1
\end{flalign*}
of which implies that $f_s(r)$ generally represents a \textit{partial} probability density function (seeing that it no longer always integrates to 1).

\paragraph{} To obtain the first-passage reward density for such a system above, it is sufficient to remove all states that do not reach $B$ first and then solve it (e.g. via matrix inversion in Example \ref{example:study1_automatic}). More generally if $B \subset Abs$, which includes the case where $B \cap Abs  =\emptyset$, then to determine the first-passage reward density  $Pr(s \vDash \Diamond_{=r} B)$, we first modify the sMRM such that states in $B$ are absorbing artificially. After this we remove any states that do not reach $B$ from the system and then solve.

\section{Experimenting with symbolic inversion}

\paragraph{} We have found that symbolic matrix inversion (using the Sympy package for python \cite{sympy}), is generally intractable for complete (square) matrices with more than six dimensions (or states in our case).

\paragraph{} For the problem $$(I - (A\circ \textbf{C}))\textbf{x} = \textbf{d}$$
where $\textbf{C},\textbf{x},\textbf{d}$ are matrices (or vectors) of functions. Let  $(I - (A\circ \textbf{C}))$ be completely symbolic (hence we are not using the vector representation of characteristic functions). For example for a $2 \times 2$ matrix, we have $$\begin{pmatrix}
a & b \\
c & d \\
\end{pmatrix}$$
where $a,b,c,d$ are strictly variables. Then our goal is to find $(I - (A\circ \textbf{C}))^{-1}$ symbolically.  From some basic experiments, we found the following results:    
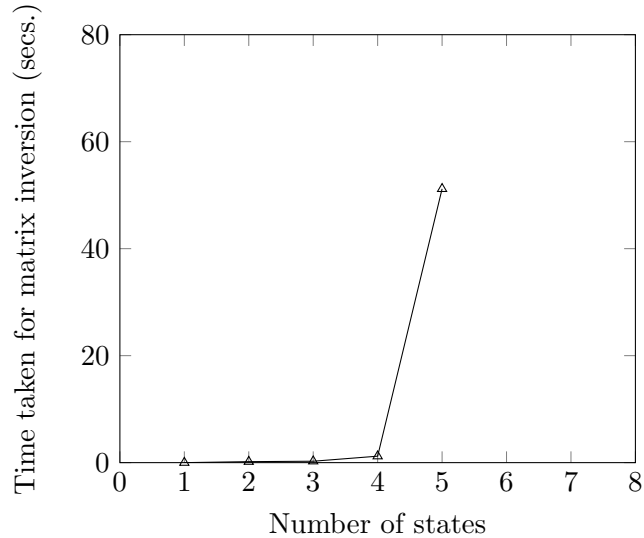
\begin{figure}[H]
	\centering
	\begin{tikzpicture}
	\begin{axis}[
	xlabel={Number of states},
	ylabel={Time taken for matrix inversion (secs.)},
	xmin=0, xmax=8,
	ymin=0, ymax=80,
	xtick={0,1,2,3,4,5,6,7,8},
	grid style=dashed,
	]
	\addplot[
	color=black,
	mark=triangle,
	]
	coordinates {
		(1,    0.0004210472107)(2,    0.1848955154)(3,0.2642114162)(4,1.207051277)(5,51.1838026)
	};
	\end{axis}
	\end{tikzpicture}
	\caption{The time taken to invert the (dense) matrix $\textbf{I} - (\textbf{A} \circ \textbf{C})$ blows up drastically.}    
	\label{fig:size-of-s-vs-time-for-solving-system1}
\end{figure}
where we have given up trying to solve problems larger than six states. Whilst the above shows a bottleneck, it is not a fair analysis alone as usually there is only one variable $x$ of a univariate-reward sMRM.

\begin{flalign}
\label{matrix:6x6}
\left(\begin{smallmatrix}
x&2+x&3x +5&2x + 10& 11x-9& 8x-20\\
2+x&x&2+x&x& 27x - 10& 15x - 20\\
x&2+x&5x- 3&30x + 1& x+50 -10& 11x + 12x^2\\
x+5 - 10x&20x+ x^2 +x^3&x+8 -10&2+x& 42& 36x\\
x+5 - 10x + 10x^3& 6x&x+8 -10 - 25*3 + x^2&10+x& -15x& 16x-4\\
x + x^2 + 3&6*x + 12&x- 10&2x&x^2&x-15
\end{smallmatrix}\right)
\end{flalign}

\paragraph{} When using a single variable $x$, we find that solving times does get easier. For example, inverting the $6\times 6$ matrix above \eqref{matrix:6x6} took around 8.5 seconds.  However, the resulting matrix from the inversion was unwieldy, being over hundreds of thousands of characters long. A simplification procedure exists in Sympy which took roughly 215 seconds to complete. The result however was still quite large being around 11000 characters long. One might argue that the matrix above is not a true representation of an sMRM problem as the functions above are arbitrary and not characteristic functions. Whilst true, certain characteristic functions have more complicated forms, such as those of the Weibull distribution, Gaussian or log-normal distribution. Thus with such functions, it may be better to use a separate variable to represent them. This however leads to the blow-up previously shown.

\paragraph{} Therefore our current bottleneck is in the inversion process (sympy uses Gaussian elimination as the default algorithm \cite{sympy_inversion} which is known to have a complexity of $O(n^3)$ for a matrix of dimensions $n\times n$)). Whilst better symbolic algorithms could be researched, to mitigate this problem, we proceed to replacing the symbolic step with a numerical one; the solutions presented in the next two chapters for discrete-reward sMRMs, and one subsequent chapter for continuous-reward sMRMs.

\section{Summary} In this chapter, we presented the theoretical foundations for sMRM model checking. We explained that many problems for sMRMs reduce to solving for the first-passage reward density $Pr(Rew = r\ \cap\ s \vDash \Diamond B)$. Expected value problems were found to reduce to algorithms similar to that for MRMs.

\paragraph{} In the subsequent chapters, we move forward to deriving numericals solutions for the first-passage reward density. We first present a direct approach for solving the system of convolution equations, and in the following chapter we present iterative algorithms.

\paragraph{} Note that whilst we have presented techniques in the following chapters for solving the system of convolution equations previously derived, there are other ways to derive the first-passage reward densities, for example see \cite{kao1974modeling,warr2012introduction}.

\chapter{Direct Methods for sMRMs}

\section{Introduction} In this chapter, direct algorithms for solving $\mrmpdf$ numerically are presented. A direct algorithm is one that terminates under a finite number of steps, and provides a full solution only upon termination. The first algorithm we present is the naive Gauss elimination algorithm, but adapted for systems of convolution equations. The second is an approximation algorithm. To summarize this chapter, our major finding is that Gaussian elimination can be adapted for the system of convolution equations by carefully replacing multiplication with convolution and division by deconvolution. Addition and subtraction amounts to piecewise addition or subtraction of vectors.

\paragraph{} We shall restrict ourselves to using probability mass functions, i.e. discrete reward random variables (defined only over $\mathbb{N}$). Continuous rewards are treated in Chapter 7. Secondly, we present a computation of $\mrmpdf$ only for $r \in \{0,1,2,\cdots,k\}$, for some $k \in \mathbb{N}$. For any sMRM, we assume that all states $s \in S$ \textbf{can reach} $B$ with probability one. We also continue to focus on \textbf{univariate} reward random variable. We first present the necessary definitions for this chapter (of which will be used in subsequent chapters). We then define the algorithms aforementioned.

\section{Motivation} As previously discussed, the reason for finding new algorithms is due to the symbolic approach presented in the previous chapter being intractable. The papers \cite{warr2014numerical,sakamoto1997probability} provide for us the basis of the following work. They suggest replacing characteristic functions of probability density functions with the discrete Fourier transform (DFT) instead via discretization of the characteristic function. This is useful as there exists a fast numerical implementation of the DFT known as the fast Fourier transform (FFT) with which we can perform linear convolutions. However in this case we restrict ourselves to probability mass functions (e.g. discrete R.V.s) by choice, and therefore any distribution to be used that is continuous should be discretized first, or a discrete analogue found for it. For example see \cite{chakraborty2015generating} for a survey of analogues.

\paragraph{} Before we present the Gaussian elimination algorithm, we will first define the convolution and deconvolution operations that will be used.

\section{Discrete convolution and deconvolution}

Given a non-negative discrete lattice random variable $X$ defined over $\mathbb{N}$, let $p_X[x]$ be its probability mass function (pmf). From now on, we represent pmfs as vectors (arrays) (with indices in $\mathbb{N}$), and are not expressed analytically (algebraically).  Additionally, for the computation of $\mrmpdf$, we will restrict ourselves to $r \in \{0,1,2,\cdots,k\} \in \mathbb{N}$ only. Therefore for example, $f_s[r]$ will be computed only for such a range.

\paragraph{} The discrete Fourier transform operator $\mathcal{D}$ applied to the pmf of a random variable $X$ is written as $\mathcal{D}\{p_X[x]\}$, with its result being denoted as $\varphi_X[\tau]$. Then when applying $\mathcal{D}$ to $\mrmpdf$, we write this transform as
$$\mathcal{D}\{\mrmpdf\} = \varphi_{s \vDash \Diamond B}[\tau]$$
or $\varphi_{s}[\tau]$ for short.

\subsubsection{Sums of random variables \& linear convolution via DFTs} Given two independent random variables $X,Y$ defined over $\mathbb{N}$ (for simplicity), with pmfs $p_X, p_Y$ respectively, then the pmf of the operation $X+Y$ is determined via a discrete linear convolution, i.e.

\begin{flalign}
\label{eqn:naive_convolution}
p_Z(k) = (X+Y)(k) = \sum_{x=0}^{k}p_X(x)p_Y(k-x)
\end{flalign}

If we intend to compute the first $k$ values of $p_Z(x)$, i.e. for all $x \in \{0,1,\cdots,k-1\}$, then using DFTs to perform the convolution, the pmf of $Z$ can be computed as
\begin{flalign}
\label{eqn:DFT_sums}
p_Z(k)= \mathcal{D}^{-1}\{\mathcal{D}\{\widetilde{p}_X[x]\}\mathcal{D}\{\widetilde{p}_Y[x]\}\}(k)
\end{flalign}
where $\mathcal{D}^{-1}$ is the inverse DFT operator, and $\widetilde{p_i}$ is defined as
\begin{flalign}
\label{eqn:pi_widetilde}
\widetilde{p}_i[x] \triangleq \begin{cases}
p_i[x] & 0 \leq x \leq k \\
0 &  k < x \leq 2k-1
\end{cases}
\end{flalign}
for each $i \in \{X,Y\}$. Note that $\widetilde{p_i}$ is just $p_i$ truncated and zero padded (with $k-1$ zeroes). The zero-padding is necessary with DFTs to prevent a known problem called \textit{time-aliasing}. Ultimately, the above states that the convolution requires two vectors of length $2k-1$, to determine the first $k$ values of $p_Z$ exactly.

\paragraph{} Let us introduce the function $$\mathtt{conv}_k(v_1,v_2) \triangleq \mathcal{D}^{-1}\{\mathcal{D}\{\widetilde{v}_1[x]\}\mathcal{D}\{\widetilde{v}_2[x]\}\}$$ as notation to specify the convolution of vector $v_1$ with $v_2$ via FFTs for up to $k$ points only. If $v_1,v_2$ are derived from analytical pmfs, only the first $k$ values are sampled, and the resulting vector padded with $k-1$ zeroes. The length of both $\widetilde{v_1},\widetilde{v_2}$ is $2k-1$.

\paragraph{} Efficient algorithms exist to compute DFTs (and their inverses) via the fast Fourier transform (FFT) with complexity $O(nlog_2n)$ where $n$ is the size of the vector being transformed. For example see \cite{FFTW}. Thus the computation of \eqref{eqn:DFT_sums} has complexity $3\cdot O(nlog_2n) + O(n)$ where $n=2k-1$. This is since $\mathcal{D}$ is used twice, $\mathcal{D}^{-1}$ once and element-wise vector multiplication is performed once. Asymptotically, this is equivalent to $O(klog_2k)$.\qed    

\subsubsection{Linear deconvolution via polynomial division}
\label{subsubsection:linear_deconv}
Given three non-negative discrete lattice random variables $X,Y,Z$ defined over $\mathbb{N}$, where $Z = X+Y$, and $X$ is independent of $Y$. Then the deconvolution operation is sought when we know for example the pmfs of $Z$ and $Y$ only but we would like to find the pmf of $X$. Algebraically, then the pmf $p_X$ is computed as
\begin{align}
p_X[k] =\frac{1}{p_Y[0]}(p_Z[k] - \sum_{x=0}^{k-1}p_X[x]p_Y[k-x])
\end{align}
and for clarity, $p_X[0] = \frac{p_Z[0]}{p_Y[0]}$.

\paragraph{} This can be validated by simply rearranging the definition of convolution. Or if starting from the result above, we obtain
\begin{align*}
p_X[k] &=\frac{1}{p_Y[0]}(p_Z[k] - \sum_{x=0}^{k-1}p_X[x]p_Y[k-x]) &\\
p_X[k]p_Y[0] &=p_Z[k] - \sum_{x=0}^{k-1}p_X[x]p_Y[k-x] &\\
p_Z[k] &= p_X[k]p_Y[0] + \sum_{x=0}^{k-1}p_X[x]p_Y[k-x]  &\\
p_Z[k] &= \sum_{x=0}^{k}p_X[x]p_Y[k-x]  &\\
p_Z[k] &= (X+Y)[k]
\end{align*}
which is the original definition of $p_Z[k]$.

\paragraph{} Note that deconvolution is computed recursively as seen above (unlike convolution). As such, computing $p_X[k]$ requires an ordered computation, i.e. computing (in order) $p_X[0],p_X[1],\cdots,p_X[k-1]$ before $p_X[k]$. In another sense, deconvolution is a recursive filter, but convolution is non-recursive.

\label{eqn:full_deconv}
\paragraph{} An important thing to note is that if $p_Y[0] = 0$, the computation above will lead to a division by zero, i.e. $1/p_Y[0]$. Generally, it is only possible to resolve this if $p_X,p_Y$ are not of fixed length, e.g. infinitely long. Starting from the third equation above, then
\begin{align*}
p_Z[k] &= p_X[k]p_Y[0] + \sum_{x=0}^{k-1}p_X[x]p_Y[k-x]  &\\
&= 0 + \sum_{x=0}^{k-1}p_X[x]p_Y[k-x]  &\\
&= p_X[k-1]p_Y[1] + \sum_{x=0}^{k-2}p_X[x]p_Y[k-x]
\end{align*}    
Then, we can rearrange the result to obtain $p_X[r]$:
\begin{align*}
p_Z[k] &= p_X[k-1]p_Y[1] + \sum_{x=0}^{k-2}p_X[x]p_Y[k-x]  &\\
p_X[k-1] &= \frac{1}{p_Y[1]}(p_Z[k] - \sum_{x=0}^{k-2}p_X[x]p_Y[k-x]) &\\
p_X[k] &= \frac{1}{p_Y[1]}(p_Z[k+1] - \sum_{x=0}^{k-1}p_X[x]p_Y[k+1-x]) &\\
p_X[k] &= \frac{1}{p_Y[1]}(p_Z[k+1] - \sum_{x=0}^{k-1}p_X[x]p_Y[k+1-x])
\end{align*}
where $p_X[0] = \frac{p_Z[1]}{p_Y[1]}$. From a computational point of view, this result is equal to computing $p_X  =p_Z[1:] \deconv p_Y[1:]$, where for example $p_Z[1:]$ is the vector $p_Z$ with each element shifted left by one place. However this latter computation requires the values of $p_Z[1:k+1]$ to compute the first $k$ values of $p_X$. If we only have access to the first k values of $p_Z$ this computation will yield therefore only $k-1$ values which is problematic for our system of convolution equations since the vectors are fixed to some finite length $k$. In the context of finitely long vectors, let us define a \textbf{full deconvolution} $f \deconv g$ to be where $g[0] \neq 0$. We will call  $g$ a \textbf{full deconvolutor}.

\paragraph{} Another area where convolution arises is with polynomial multiplication. This is since the coefficients of the polynomial resulting from multiplication are derived via convolution of the coefficients of the multiplicands. For example, we can encode a pmf $p_X[r]$ (with $r = 0,1,\cdots,k-1$) as a polynomial of degree $k-1$\footnote{Note that the degree of the polynomial is one less than the length of $p_X$.} by $$P_{k-1}(x) = p_X[0] + p_X[1]x + p_X[2]x^2 + \cdots + p_X[k]x^{k-1}$$
and we will denote $P_{k-1}(x;p_X)$ to mean that the coefficients of the polynomial is encoded via the vector $p_X$. Then for example, we have the multiplication $P_{k-1}(x;p_X)\times P_{k-1}(x;p_Y) = P_{2k-2}(x;p_Z)$ where $p_Z = p_X * p_Y$.

\paragraph{}  Likewise, the inverse operation - polynomial division - can be used to reverse the process, e.g. $ P_{2k-2}(x;p_Z)/P_{k-1}(x;p_Y) = P_{k-1}(x;p_X) + R(x)$, where $R$ is the remainder polynomial and $P_{k-1}(x;p_X)$ is termed the quotient. Hence, we can use polynomial division for deconvolution since the quotient obtains $p_X = p_Z \deconv p_Y$. Note that the remainder $R$ is not required and can be discarded. The division above assumes $p_Z$ is known up to $2k-1$ values, and $p_Y$ up to $k$ values. However, if we are given only $k$ values of $p_Z$ and $p_Y$ we can still obtain $P_{k-1}(x;p_X)$ (and therefore the first $k$ values of $p_X$) by zero-padding $p_Z$ with k-1 zeroes. This is provided that we have a \textbf{full deconvolutor}.

\paragraph{} Let $v_1,v_2$ be vectors of length $k$. Then define the polynomial
$$ P_{k-1}(x;v_3) + R = P_{2k-2}(x;\widetilde{v_1})/P_{k-1}(x;v_2)$$
where $\widetilde{v_1}$ is equal to $v_1$ zero padded with $k-1$ zeroes. Then, we introduce the function
$$\mathtt{deconv}_k(v_1,v_2) \triangleq \mathtt{coeffs}(P_{k-1}(x;v_3))$$

where $\mathtt{coeffs}(P_{k-1}(x;v_3)) = v_3$, i.e. $\mathtt{coeffs}$ returns the coefficients of the polynomial. The subscript $k$ in front of $\mathtt{deconv}$ denotes that the first $k$ points of deconvolution is obtained.

\paragraph{} Let $p_X,p_Y$ be pmfs. Then to obtain the quotient $p_Z[x] = p_X \deconv p_Y$ for $x = 0,1,\cdots,k-1$, define $v_1[x] = p_X[x], v_2[x] = p_Y[x]$ for $x = 0,1,\cdots,k-1$. Then $p_Z[x] = \mathtt{deconv}_k(v_1,v_2)$.

\paragraph{} There appears to be algorithms that compute polynomial divisions using FFTs, for example see the lecture notes \cite{poly_div,poly_div2}. They are stated to have the same worst-case complexity as the convolution case with the FFT - $O(nlog_2n)$ where $n=2k-1$ for the computation of $\mathtt{deconv}_k$. Note however that the deconvolution operation is generally numerically unstable.

\newcommand{\redechform}{\bar{{\mathcal{A}}}}

\section{Naive Gaussian elimination}
\label{sec:naive_GE}  We can rearrange \eqref{eqn:lineq_matrix_trans_pdf} into the form
\begin{flalign}
\label{eqn:matrix_sol_pdf}
((I \circ \bm{\Delta}_{x,0}) -A\circ \textbf{G})\ \textcircled{$\ast$}\ \textbf{f} = \textbf{h}
\end{flalign}
where $\bm{\Delta}_{x,0}$ is a hypermatrix where each element is the Kronecker delta $\delta_{x,0}$ with mass 1 over zero (and thus a hypermatrix overall). Thus $(I \circ \bm{\Delta}_{x,0})$ returns a diagonal hypermatrix instead, which represents the equivalent of the identity matrix in traditional linear algebra. Using vectors for pmfs (of discrete (lattice) random variables), then each of $\textbf{f}$, $(A \circ \textbf{G})$, and $\textbf{h}$ is a three dimensional hypermatrix or hypervector. Again, we will restrict ourselves to $r \in \{0,1,2,\cdots,k\} \in \mathbb{N}$ only, e.g. $(G_{s,t}[r])_{s,t \in S_?^2}$ is defined for $r=0,1,2,\cdots$. Since, we are interested in computing $\mrmpdf$ only for $r = 0,1,2,\cdots,k-1$, i.e. just $k$ values, it suffices to fix the sizes of each term as
\begin{enumerate}
	\label{cardinalities_pdf}
	\item $(1 \times {S_?| \times |S_?|})$ for $A$ and $I$.
	\item $({k \times |S_?| \times |S_?|})$ for $\textbf{G}$ and $\bm{\Delta}_{x,0}$.
	\item $({k \times |S_?| \times 1})$ for $\textbf{f}$ and $\textbf{h}$.
\end{enumerate}

\paragraph{} Let us denote $((I \circ \bm{\Delta}_{x,0}) -A\circ \textbf{G})$ as $\mathbf{\mathcal{A}}$, leaving us with $$\mathbf{\mathcal{A}}\ \textcircled{$\ast$}\ \textbf{f} = \textbf{h}$$

\paragraph{} We refer to the second and third dimension as the \textit{rows} and \textit{columns} of the hypermatrix respectively, i.e. $\mathcal{A}$ has size $({k\times \text{no. rows} \times \text{no. columns}})$. Hence $\mathcal{A}_{i,j}(r)$ refers to indexing the hypermatrix by the $i^{th}$ row, the $j^{th}$ column, and the $r^{th}$ value. The indices start from zero, i.e. $i=0,1,\cdots,|S_?|-1$. The $i^{th}$ equation  of the system is denoted as
$$(\mathbf{\mathcal{A}}\ \textcircled{$\ast$}\ \textbf{f})_i =  \sum_{j} \mathcal{A}_{i,j} \ast {f}_{j} =  h_i
$$

\paragraph{} We now present a Gaussian elimination (GE) type algorithm to solve for $\textbf{f}$ in the system above. It is helpful to realise at this stage that convolution and deconvolution are used for the same purposes as multiplication and division in the traditional algorithm (for linear equations), whilst addition and subtraction are just element-wise vector addition and subtraction.

\paragraph{Description of GE:}
Gaussian elimination solves the system in two stages: The first stage involves eliminating from each equation $i$ where $i \geq 1$, all terms $f_j$ where $j \leq i$. The procedure is done algorithmically: Starting from equation $0$, this equation is used to eliminate from equation $i$ for all $i\geq 1$, the term $f_0$. This can be done by subtracting from each equation $i$ a convolution (multiple) of equation 0. That is, we replace equation $i$ with
\begin{flalign}
\label{eqn:gauss_reduction_update}
(\sum_{j} \mathcal{A}_{i,j} \ast {f}_{j}) - \sigma_i * (\sum_{j} \mathcal{A}_{0,j} \ast {f}_{j}) =  h_i - (\sigma_i * h_0)
\end{flalign}
where $\sigma_i \triangleq \mathcal{A}_{i,0} \deconv \mathcal{A}_{0,0}$ is the convolution (multiple) term. and it is assumed that $\mathcal{A}_{0,0}$ (called a \textit{pivot} in the literature on Gaussian elimination) is a full deconvolutor (described in Section \ref{eqn:full_deconv}). This assumption is proved true in Theorem \ref{subsec:full_deconv}. The update for equation $i$ \eqref{eqn:gauss_reduction_update} can also be written as
$$ \sum_{j} \mathcal{A}^{'}_{i,j} \ast {f}_{j} =  h^{'}_i$$
where $\mathcal{A}^{'}_{i,j}  = (\mathcal{A}_{i,j} - \sigma_{i} * \mathcal{A}_{0,j})$, and $h^{'}_{i} = h_i - (\sigma_i * h_0)$. Note then that
$$\mathcal{A}^{'}_{i,0} = \mathcal{A}_{i,0} - (\mathcal{A}_{i,0} \deconv \mathcal{A}_{0,0}) * \mathcal{A}_{0,0} = \mathcal{A}_{i,0} - \mathcal{A}_{i,0} = \vec{0}$$
Thus, $f_0$ is eliminated from equation $i$.

\paragraph{} Given a system of convolution equations, substituting an equation by a sum of itself and a multiple of another does not change the solution set (see Proposition \ref{theorem:gauss_solution_set}). Hence, the system after being repeatedly updated (via \eqref{eqn:gauss_reduction_update}) will still have the same solution.

\paragraph{} The procedure above ensures $f_0$ is eliminated from equations $i$, for all  $i\geq 1$. Then, re-assign $\mathcal{A} \coloneqq \mathcal{A}'$, and $\textbf{h} \coloneqq \textbf{h}'$, and the elimination procedure repeats as before, but now starting with equation 1, and eliminating $f_1$ from equations $i,$ for all $i \geq 2$. This repetition continues, e.g. starting with equation $j=2,3,\cdots,|S_?|-2$, and eliminating $f_j$, from all equations $i$ where $i>j$. The end product of all of this is the system:
\begin{align}
\label{eqn:row_echelon_form}
\begin{pmatrix}
\mathcal{A}'_{0,0} &  \mathcal{A}'_{0,1} & \mathcal{A}'_{0,2} & \cdots &\cdots & \mathcal{A}'_{0,|S_?|-1}\\
\vec{0} &   \mathcal{A}'_{1,1} & \mathcal{A}'_{1,2} & \cdots &\cdots& \mathcal{A}'_{1,|S_?|-1} \\
\vec{0} &  \vec{0} & \mathcal{A}'_{2,2} & \cdots &\cdots& \mathcal{A}'_{2,|S_?|-1} \\
\vdots & \vdots & \ddots & \ddots & \ddots & \vdots  \\
\vec{0} & \vec{0} & \vec{0} & \cdots & \mathcal{A}'_{|S_?|-2,|S_?|-2} & \mathcal{A}'_{|S_?|-1,|S_?|-2}\\
\vec{0} & \vec{0} & \vec{0} & \cdots & \vec{0} & \mathcal{A}'_{|S_?|-1,|S_?|-1}
\end{pmatrix}\textcircled{$*$}
\begin{pmatrix}
f_0 \\
f_1 \\
f_2 \\
\vdots \\
f_{|S_?|-2} \\
f_{|S_?|-1}
\end{pmatrix}
=
\begin{pmatrix}
h'_0 \\
h'_1 \\
h'_2 \\
\vdots \\
h'_{|S_?|-2} \\
h'_{|S_?|-1}
\end{pmatrix}
\end{align}
or written in matrix form,
$$\mathcal{A}'\textcircled{$*$}\textbf{f} = \textbf{h}'$$

The resulting hypermatrix $\mathcal{A}'$ is \textit{upper-triangular}, that is the entries $\mathcal{A}'_{i,j} = \vec{0}$ for all $i >j$.

\paragraph{} Next, the second stage of Gauss elimination (known as \textit{back-substitution}) begins. Note from \eqref{eqn:row_echelon_form} that the solution hypervector \textbf{f} can be determined via the reverse order $f_{|S_?|-1}$, $f_{|S_?|-1}$, $\cdots$, $f_0$. Starting from the last equation of the system \eqref{eqn:row_echelon_form} i.e. equation $|S_?|-1$,  we have
$$\mathcal{A}^{'}_{|S_?|-1,|S_?|-1}* f_{|S_?|-1} = h^{'}_{|S_?|-1}$$
and deconvolving both sides by $\mathcal{A}^{'}_{|S_?|-1,|S_?|-1}$  gives
$$f_{|S_?|-1} = h^{'}_i \deconv \mathcal{A}^{'}_{|S_?|-1,|S_?|-1}$$
yielding solution for $f_{|S_?|-1}$. Knowing $f_{|S_?|-1}$, then $f_{|S_?|-2}$ can be determined, since by rearranging equation ${|S_?|-2}$, we obtain
$$f_{|S_?|-2} = (h^{'}_{|S_?|-1} -
\mathcal{A}^{'}_{|S_?|-2,|S_?|-1}*f_{|S_?|-1}) \deconv \mathcal{A}^{'}_{|S_?|-2,|S_?|-2}$$
of which the terms on the RHS are all known. In the same vein, we can determine $f_{i}$, in the order $i=|S_?|-3,|S_?|-4,\cdots,0,$ via the formula
\begin{flalign}
\label{eqn:deconv_GE_part2}
f_{i} = (h'_i -  \sum_{j=i}^{|S_?|-1} \mathcal{A}^{'}_{i,j}*f_{j}) \deconv \mathcal{A}^{'}_{i,i}
\end{flalign}
which results in solving for \textbf{f} completely. Note that the pivots $(\mathcal{A}_{l,l}^{'})_{l \in S_?}$ are  assumed to be full deconvolutors.

\label{par:up_to_k}
\paragraph{} Whilst we have fixed the sizes of the hypermatrices in \eqref{eqn:matrix_sol_pdf} (see Section \ref{cardinalities_pdf}), the operators $(*,\pentagon)$ will cause the system to grow in size. For example, the convolution of two vectors of length $k$ will result in a vector of length $2k-1$. Hence, for a practical implementation, we substitute these operators with $(\mathtt{conv}_k,\mathtt{deconv}_k)$ which will ensure that the results are fixed to length $k$. This does not affect the solution set as the set of all vector operators in Gaussian elimination do not require values beyond the $k^{th}$ value (assuming that deconvolutions are only performed with full deconvolutors). Precisely speaking, for any vector operator $\square \in \{/,\cdot,+,-,*,\pentagon \}$ with the first four being element-wise operators, then $(f\ \square\ g)[r]$ is always independently computed of $(f\ \square\ g)[r+l]$ for all $l>0$. Hence, values $r$ where $r \geq k-1$ are irrelevant for the computation for the first $k$ values of the solution set.

\subsection{Full deconvolutions}
\label{subsec:full_deconv}
Earlier, it was assumed that the \textit{pivots} $(\mathcal{A}^{'}_{l,l})_{l \in S_?}$ were full deconvolutors for both stages of GE (where $\mathcal{A}^{'}$ is the final hypermatrix in upper-triangular form). Firstly, note that if this assumption holds in the first stage, it immediately holds in the second. This is since the \textit{pivots} used for deconvolution in the first stage, are unchanged in the second stage. If this assumption does not hold, either divisions by zero will be introduced, or full deconvolution cannot occur.

\paragraph{} With linear systems, one way to handle pivots that lead to divisions by zero, is by utilizing \textit{pivoting strategies} and are described for example in \cite[p. 150]{dahlquist1974numerical}. However, for particular systems pivoting is not required and our system is such a case. The proof is as follows.

\begin{theorem}
	\label{theorem:full_deconv}
	Let $\mathcal{A}^{'}$ be the final updated hypermatrix (i.e. when it is upper-triangular) obtained after the first stage of GE, and $\mathcal{A}$ be the original hypermatrix. Then, $\mathcal{A}^{'}_{i,j}$  has the closed form
	\begin{align}
	\label{eqn:full_deconvolutor}
	\mathcal{A}^{'}_{i,j}
	&= \mathcal{A}_{i,j} - \sum_{j=0}^{i-1} \sigma_{j} * \mathcal{A}_{j,i}
	\end{align}
	where $ \sigma_{j} = \mathcal{A}^{'}_{i,j} \deconv \mathcal{A}^{'}_{j,j} $. Then, we have that
	\begin{flalign}
	\label{eqn:zero_index_GE}
	\mathcal{A}^{'}_{i,j}[0]
	&= \mathcal{A}_{i,j}[0] - (\sum_{j=0}^{i-1} \sigma_{j} * \mathcal{A}_{j,i}) [0] \notag &\\
	&= \mathcal{A}_{i,j}[0] - \sum_{j=0}^{i-1} \sigma_{j}[0]\mathcal{A}_{j,i}[0] \notag  &\\
	&= \mathcal{A}_{i,j}[0] - \sum_{j=0}^{i-1} (\mathcal{A}^{'}_{i,j} \deconv \mathcal{A}^{'}_{j,j})[0]\mathcal{A}_{j,i}[0] \notag  &\\
	&= \mathcal{A}_{i,j}[0] - \sum_{j=0}^{i-1} \frac{\mathcal{A}^{'}_{i,j}[0] }{\mathcal{A}^{'}_{j,j}[0]}\mathcal{A}_{j,i}[0]
	\end{flalign}
	implying that $\mathcal{A}^{'}[0]$ is only dependent on $\mathcal{A}[0]$.
	
	\paragraph{} Then, we have that $\mathcal{A}^{'}_{l,l}[0] \neq 0$ for all $l \in S_?$. This ensures that $\mathcal{A}^{'}_{l,l}$ is a full deconvolutor, as described in Section \ref{eqn:full_deconv}.

	\paragraph{} \textit{Proof: } Eq. \eqref{eqn:zero_index_GE} is the update rule for Gaussian elimination for linear systems. Hence Gaussian elimination with our convolution equations reduces to Gaussian elimination with linear equations at the zeroth value, i.e. for the system $(\mathbf{\mathcal{A}}\ \textcircled{$\ast$}\ \textbf{f})[0] = \textbf{h}[0]$. The pivots for GE with this linear system are $(\mathcal{A}^{'}_{l,l}[0])_{s \in S_?}$. We can deduce that the matrix $\mathcal{A}[0]$ is \textit{diagonally dominant}, i.e. for all $l \in S_?$, we have
	$$ \mathcal{A}_{l,l}[0] \geq \sum_{j \in S_?} \mathcal{A}_{l,j}[0]$$
	\par{}Then, Dahlquist \& Bj$\ddot{\text{o}}$rk \cite[p. 151-152]{dahlquist1974numerical} states that diagonal dominance is a sufficient condition to ensure that $\mathcal{A}^{'}_{l,l}[0]$ are \textbf{non-zero} for all $l \in S_?$, without using pivoting. Hence, full deconvolutions are guaranteed for our system of convolution equations. \qed

\end{theorem}

\subsection{Same solution set}
The first stage of Gaussian elimination (known as forward elimination or Gauss reduction) preserves the solution of the system, such that solving with the hypermatrix $\mathcal{A}'$ (being upper triangular) yields the same answer as solving directly from the original hypermatrix $\mathcal{A}$ . The hypermatrix $\mathcal{A}'$ and the corresponding hypervector $\textbf{h}'$ are achieved using at most different four operations (detailed in Proposition \ref{theorem:gauss_solution_set}), each of which do not change the solution set. We now state this formally as a theorem. It is a simple restating of a theorem for systems of linear equations found in \cite[ Theorem 1.5]{hefferon2017linear}, but now for convolution equations.

\begin{proposition}      \label{theorem:gauss_solution_set} Given a system of convolution equations $\textbf{E}\textcircled{$\ast$}\textbf{f} = \textbf{h}$, with equations $E_0,E_2,\cdots E_N$, and the unique solution $\textbf{x}$ (a hypervector), where each $E_i$ represents the equation $$\sum_j E_{i,j} \ast f_j = h_i$$
	If the system is transformed via these two updating operations:
	
	\begin{enumerate}
		\item replacing an equation by a convolution of itself with a non-zero vector.
		\item replacing an equation with the equation itself added to another equation that has been convolved with a non-zero vector.
	\end{enumerate}
	then the transformed system will still have the same solution set \textbf{x}, with the exception of the case where columns are swapped. Due to the simplicity of the statements here, we present the proof to the Appendix, see Theorem \ref{theorem:gauss_solution_set_full}.
\end{proposition}

\subsection{The algorithm in pseudocode} The algorithm for reducing the system into upper-triangular form is as follows:

\begin{algorithm}[H]
	\caption{Gauss reduction (without pivoting)}
	\KwData{$\mathcal{A},\textbf{h}$}
	\KwResult{$\mathcal{A}$ in upper-triangular form, and the corresponding hypervector $\textbf{h}$}
	\ForEach{$j = 0,1,2,\cdots,|S_?|-2$}{
		\ForEach{$i = j+1,j+2,\cdots,|S_?|-1$}{
			$\sigma_i = \mathtt{deconv}_k(\mathcal{A}_{i,j},\mathcal{A}_{j,j})$\;
			\ForEach{$l = j,j+1,\cdots,|S_?|-1$}{$\mathcal{A}_{i,l} = \mathcal{A}_{i,l} - \mathtt{conv}_k(\mathcal{A}_{j,l},\sigma_i)$
			}
			$h_i = h_i - \mathtt{conv}_{k}(h_j,\sigma_i)$
		}
	}
\end{algorithm}
\paragraph{} Once the algorithm terminates, then as previously described, $\mathcal{A}$ is now in upper-triangular form. Note that the updates for $\mathcal{A},\textbf{h}$ are done in-place in the algorithm above. Now solving for $\textbf{f}$, then

\begin{algorithm}[H]
	\caption{Back substitution}
	\KwData{$\mathcal{A}$ in upper-triangular form, and the corresponding hypervector $\textbf{h}$}
	\KwResult{\textbf{f}, the solution hypervector}
	\ForEach{$i = |S_?|-1,|S_?|-2,\cdots,1,0$}{
		$f_i = \mathtt{deconv}_k(h_i - \sum_{j=i+1}^{|S_?|-1} \mathtt{conv}_k(\mathcal{A}_{i,j},f_j),\mathcal{A}_{i,i}) $\;
	}
\end{algorithm}

\paragraph{} Theoretically, the algorithm terminates with the true $\textbf{f}$ (disregarding numerical instabilities, and issues related to \textit{pivoting}).  \qed

\subsubsection{Algorithmic complexity} For reducing the system into upper-triangular form, for each row we need to perform up to $|S_?|$ convolutions. For the reduction phase, i.e. converting the augmented matrix into upper-triangular form, setting the pivot to a row $j \in [1,2,\cdots,|S_?|-1]$, where $j=1$ is the second-last row and $j=|S_?|-1$ being the upper-most, then the complexity for each $j$ is: 1) $O(j)$ $\mathtt{deconv}_k$ operations, 2) $O(j^2)$ $\mathtt{conv}_k$ operations and 3) $O(jk)$ subtractions. Then, assuming $\mathtt{deconv}_k$ and $\mathtt{conv}_k$ have complexity $O(klog_2(k))$, this results in a \textit{row-complexity}

$$O(jklog_2(k)) + O(j^2klog_2(k)) + O(jk) = O(j^2klog_2(k))$$

with a total complexity of $\sum_{j=1}^{|S_?|-1}O(j^2klog_2(k)) = O(|S_?|^3klog_2(k))$, which is therefore cubic in $|S_?|$ and linearithmic in $k$.

\begin{example}
	
	Let us solve a toy system of two equations: $a\ast f_1 + b\ast f_2 = h_1$, $c\ast f_1 + d\ast f_2 = h_2$. Then, in matrix form, we have:    
	\[
	\begin{pmatrix}
	a & b \\
	c & d \\
	\end{pmatrix}\textcircled{$\ast$}
	\begin{pmatrix}
	f_1 \\
	f_2 \\
	\end{pmatrix}    =
	\begin{pmatrix}
	h_1 \\
	h_2
	\end{pmatrix}
	\]
	and this can be written in \textit{augmented-matrix form} as
	\[
	\left(\begin{array}{cc|c}
	a & b & h_1\\
	c & d  & h_2\\
	\end{array}\right)
	\]
	
	To reduce the above into upper-triangular form, we simply need to \textit{zero} $c$. To do so we calculate $\sigma = \mathtt{deconv}_k(c,a)$. Then we convolve the first row by $\sigma$ and then subtract the result from the second row. This giving:
	\[
	\left(\begin{array}{cc|c}
	a & b & h1\\
	0 & (d - \mathtt{conv}_k(b,\sigma))  &     h_2 - \mathtt{conv}_{k}(h1,\sigma)\\
	\end{array}\right)
	\]    
	Now having the equations in upper-triangular form, $f_2$ can now be computed. Let $\bar{h}_2 =     h_2 - \mathtt{conv}_{k}(h1,\sigma)$, and $\bar{d} = (d - \mathtt{conv}_k(b,\sigma)$. Then
	$$f_2 = \mathtt{deconv}_k(\bar{h}_2, \bar{d})$$ Finally, $f_1$ can be computed using $f_2$: $$f_1 = \mathtt{deconv}_k(h_1 - \mathtt{conv}_k(b,f_2),a)$$
	
\end{example}

\section{Approximate Gaussian elimination}
\label{sec:approx_direct}

\paragraph{} It is possible to derive an approximate solution to the system \eqref{eqn:lineq_matrix_trans_pdf} by applying the DFT transform to both sides of the equations giving
\begin{flalign}
\label{eqn:approx_eqn_DFT}
x_s[\tau] = \sum_{t \in S_?}\textbf{P}(s,t)\mathcal{D}\{\widehat{f_{rew(s,t)}}\}[\tau]x_t[\tau] + \sum_{u \in B}\textbf{P}(s,u)\mathcal{D}\{\widehat{f_{rew(s,u)}}\}[\tau]
\end{flalign}
\label{key}
where $\mathcal{D}^{-1}\{x_s\}[r] \approx f_s[r]$, and we define
\begin{flalign}
\label{eqn:pi_widehat}
\widehat{p}_i[x] \triangleq \begin{cases}
p_i[x] & 0 \leq x \leq k \\
0 &  k < x \leq n+k
\end{cases}
\end{flalign}
which introduces a parameter $n$, the amount of zero-padding we will use.  Note that $\widehat{p}_i = \widetilde{p}_i$ when $n = k-1$.

\paragraph{} In matrix notation, we arrive at a form similar to \eqref{eqn:lineqn_matrix_trans}, i.e.
\begin{flalign}
\label{eqn:approx_GE}
\textbf{x}(\tau) =(A\circ \textbf{C}_n(\tau))\textbf{x}(\tau) + \textbf{d}_n(\tau)
\end{flalign}
where
\begin{enumerate}
	\item $(A_{s,t})_{s,t \in S_?^2} = (\textbf{P}(s,t))_{s,t \in S_?^2} $.
	\item $(C_{s,t})_{s,t \in S_?^2} = \mathcal{D}\{\widehat{f_{rew(s,t)}}\}_{s,t \in S_?^2}$.
	\item  $(d_s)_{s \in S_?} = (\sum_{u \in B}\textbf{P}(s,u)\mathcal{D}\{\widehat{f_{rew(s,u)}}\})_{s \in S_?}$.
\end{enumerate}
with the sizes:
\begin{enumerate}
	\item $\textbf{x} : ({(n+k) \times |S_?| \times 1})$.
	\item $A \circ \mathbf{C} : ({(n+k) \times |S_?| \times |S_?|})$.
	\item $\textbf{d} : ({(n+k) \times |S_?| \times 1})$.
\end{enumerate}

\paragraph{} In this setting, the system has the solution
\begin{flalign}
\label{eqn:matrix_sol_DFT}
\textbf{x}(\tau) = ((I - A) \circ \mathbf{C}_n(\tau))^{-1}\textbf{d}_n(\tau)
\end{flalign}

which can be solved via mature linear algebra libraries or tools (for systems of linear equations) for each $\tau = 0,1,2,\cdots,(n+k-1)$. After solving, we use $(\mathcal{D}^{-1}\{x_s\})_{s \in S_?}$ to give us an approximation of $(f_s[r])_{s \in S_?}$, for $r=0,1,\cdots,k-1$.  It is proven later in Theorem \ref{theorem:unique_approx_DFT} that the system above has a unique solution. Additionally, Theorem \ref{theorem:increased_padding} proves that when taking the limit $n\rightarrow \infty$, then the approximate solution equals the exact solution of \eqref{eqn:lineq_matrix_trans_pdf} at least for the first $k$ values.

\paragraph{} Note that GE can be optimized via particular decompositions such as LU and Cholesky  decompositions. Later in our experiment we shall use a library for the approximate variant here, that uses the LU decomposition with GE. In the remainder of this thesis, we refer to this method as either the \textbf{LU approx. method} or the \textbf{approximate LU method}. Although, more appropriately, we should have referred to it as the approximate GE-LU method.

\subsubsection{Algorithmic complexity} For this approximate solution, the worst case complexity for solving a linear system (directly, e.g. by Gaussian elim.) for each $\tau \in [0,1,\cdots,(n+k-1)]$ is $O(|S_?|^3)$. Therefore the overall complexity is $O(|S_?|^3(n+k))$, which is cubic in $|S_?|$ and linear in $n+k$. However creating the system in the first place  has $O(|S_?|^2klog_2k)$ complexity, which is a power less than solving it. Therefore, the time-complexity is slightly better than (exact) Gaussian elimination, which is linearithmic in $k$.

\begin{theorem}[Unique approximate solution via DFTs]
	\label{theorem:unique_approx_DFT}
	
	The Fourier transform of the system of convolution equations defined in \eqref{eqn:lineqn_matrix_trans} has a unique solution  $\textbf{x}(\tau)= (x_s(\tau))_{s \in S_?}$ for all $\tau \in \{0,1,\cdots,n+k-1\}$.

	\paragraph{} \textit{Proof:}  The proof is nearly identical to Theorem \ref{theorem:unique sol}. Essentially, we can replace the characteristic functions (continuous Fourier transforms) with the DFT in the theorem, and pdfs are replaced with pmfs and this is sufficient to prove uniqueness. A property that is required of DFTs to ensure the theorem is possible is that for any pmf vector $p_X[r]$, then its DFT $\mathcal{D}\{\widehat{p_X}\}$ is absolutely bounded, i.e. $|\mathcal{D}\{\widehat{p_X}\}[\tau]| = |\sum_{r=0}^{n+k-1}\widehat{p_X}[r]e^{-j2\pi \cdot r \cdot\tau/n+k}| \leq \sum_{r=0}^{n+k-1}\widehat{p_X}[r]|e^{-j2\pi\cdot r\cdot \tau/n+k}| \leq 1$. This is similar to characteristic functions being absolutely bounded by 1 with pdfs. \qed
	
\end{theorem}

\begin{theorem}[Convergence with increased padding]
	\label{theorem:increased_padding} Let $n$ be the zero padding length used in \eqref{eqn:pi_widehat}, for defining the system of equations \eqref{eqn:approx_GE}. We have that as $n \rightarrow \infty \implies |\hat{f}_s[r] - f_s[r]| \rightarrow \vec{0}$, where $\hat{f}_s[r]$ is the approximate solution via \eqref{eqn:matrix_sol_DFT} and $f_s[r]$ is the true solution of \eqref{eqn:matrix_sol_pdf}. Additionally, the sizes of the hypermatrices in the approximate system and exact system are exactly those described in Sec. \ref{sec:approx_direct} and Sec. \ref{sec:naive_GE} respectively.
	
	\paragraph{} \textit{Proof:} For each state $s \in S_?$, it can be shown that by analogy of \eqref{eqn:pdf_convex_sums},
	\begin{flalign}
	\label{eqn:approx_pdf_convex_sums}
	\hat{f_s}[r] = \sum_{\hat{\pi} \in \Pi.s} Pr(\hat{\pi})\widehat{f_{Rew(\hat{\pi})}}[r]
	\end{flalign}
	where $\Pi.s$ is the set of paths beginning from $s$ ending in $B$. That is, the approximate pmf $\hat{f_s}[r]$ is equal to a convex combination of (approximate) pmfs each signifying the accumulated reward over a unique path in $\Pi.s$. Note that the widehat notation was introduced in \eqref{eqn:pi_widehat}. For the true solution, then from \eqref{eqn:approx_pdf_convex_sums},
	\begin{flalign*}
	f_s[r] = \sum_{\hat{\pi} \in \Pi.s} Pr(\hat{\pi})f_{Rew(\hat{\pi})}[r]
	\end{flalign*}
	\par{} For any path $\hat{\pi} \in \Pi.s$,
	\begin{flalign}
	\label{eqn:approx_computation}
	\widehat{f_{Rew(\hat{\pi})}}[r]     = \mathcal{D}^{-1}\{    \mathcal{D}\{\widehat{f_{rew(\pi[0],\pi[1])}}\}\mathcal{D}\{\widehat{f_{rew(\pi[1],\pi[2])}}\}\cdots\mathcal{D}\{\widehat{f_{rew(\pi[|\pi|-2],\pi[|\pi|-1]}}\} \}
	\end{flalign}
	where $\widehat{f_{Rew(\hat{\pi})}}[r]$ is a vector strictly of length $n+k$.  However, we have for the exact computation,
	\begin{flalign}
	\label{eqn:true_computation}
	{f_{Rew(\hat{\pi})}}[r] = f_{rew(\pi[0],\pi[1])} * f_{rew(\pi[1],\pi[0])} * \cdots * f_{rew(\pi[|\pi|-2],\pi[|\pi|-1])}
	\end{flalign}
	where we have defined each $f_{rew(i,j)}$ (for $i,j \in S^2$) to be of length $k$ (see Sec. \ref{sec:naive_GE}). Let $m = |\pi|-1$, then ${f_{Rew(\hat{\pi})}}[r]$ is of length $mk - (m-1)$ since we have $m$ convolutions (over vectors of length $k$).
	\paragraph{} For any path $\hat{\pi} \in \Pi.s$, The approximation is exact when $n = mk - (m-1) - k$. If $n$ is less, errors will be introduced, or if using terminology from signal processing, we can say that the resulting vector is \textit{time-aliased} or affected by \textit{wrap-around error}. If $n$ is more, the approximation is still exact, although relative to the case where $n= mk - (m-1) - k$, it will have an additional $(mk - (m-1))- (n+k)$ zeros appended to the end of the resulting vector.
	
	\paragraph{} Next, let us partition the set of paths $\Pi.s$ into two sets: 1) A set of paths each with length (i.e. the number of steps) less than or equal to $m$, and 2) the set of paths each having length greater than $m$. We denote these sets as $\Pi.s[0:m]$ and $\Pi.s[m:\infty]$ respectively. Trivially, $\Pi.s[0:m] \cup \Pi.s[m:\infty] = \Pi.s$. Now we obtain
	\begin{flalign}
	\label{eqn:approx_sol_split_length}
	\hat{f_s}[r] = \sum_{\hat{\pi} \in \Pi.s} Pr(\hat{\pi}) \widehat{f_{Rew(\hat{\pi})}}[r]   
	= \sum_{\hat{\pi} \in \Pi.s[0:m]} Pr(\hat{\pi}) \widehat{f_{Rew(\hat{\pi})}}[r]    
	+ \sum_{\hat{\pi} \in \Pi.s[m:\infty]} Pr(\hat{\pi})\widehat{f_{Rew(\hat{\pi})}}[r]
	\end{flalign}
	\paragraph{} Hence, for the approximate solution \eqref{eqn:approx_sol_split_length}, if we keep increasing $n$ indefinitely, paths longer than $m$ will also be computed  correctly (for their first $k$ values). Hence, as $n \rightarrow \infty$, we have $m \rightarrow \infty$.
	
	\paragraph{} We can now prove formally that the \textbf{absolute} error converges:      
	\begin{flalign*}
	\lim_{n \rightarrow \infty} \triangle\epsilon &= \lim_{n \rightarrow \infty} |\hat{f_s}[r] - f_s[r]|    & \\
	&= \lim_{m \rightarrow \infty} |\Big( \sum_{\hat{\pi} \in L_{\Pi_{\Diamond B}}.s[0:m]} Pr(\hat{\pi}) \widehat{f_{Rew(\hat{\pi})}}[r] + \sum_{\hat{\pi} \in L_{\Pi_{\Diamond B}}.s[m:\infty]} Pr(\hat{\pi}) \widehat{f_{Rew(\hat{\pi})}}[r] \Big) & \\
	&\    - \Big(  \sum_{\hat{\pi} \in L_{\Pi_{\Diamond B}}.s[0:m]} Pr(\hat{\pi}) f_{Rew(\hat{\pi})}[r] + \sum_{\hat{\pi} \in L_{\Pi_{\Diamond B}}.s[m:\infty]} Pr(\hat{\pi}) f_{Rew(\hat{\pi})}[r] \Big)|    & \\
	&= \lim_{m \rightarrow \infty} |\sum_{\hat{\pi} \in L_{\Pi_{\Diamond B}}.s[m:\infty]} Pr(\hat{\pi}) \widehat{f_{Rew(\hat{\pi})}}[r] - \sum_{\hat{\pi} \in L_{\Pi_{\Diamond B}}.s[m:\infty]} Pr(\hat{\pi}) f_{Rew(\hat{\pi})}[r]|    & \\
	&= \lim_{m \rightarrow \infty} |\sum_{\hat{\pi} \in L_{\Pi_{\Diamond B}}.s[m:\infty]} Pr(\hat{\pi})(\widehat{f_{Rew(\hat{\pi})}}[r] - f_{Rew(\hat{\pi})}[r])|    & \\
	&\leq  \lim_{m \rightarrow \infty} \sum_{\hat{\pi} \in L_{\Pi_{\Diamond B}}.s[m:\infty]} Pr(\hat{\pi})|(\widehat{f_{Rew(\hat{\pi})}}[r] - f_{Rew(\hat{\pi})}[r])|    & \\
	&\leq  \lim_{m \rightarrow \infty} \sum_{\hat{\pi} \in L_{\Pi_{\Diamond B}}.s[m:\infty]} Pr(\hat{\pi})& \\
	&= 1 - \lim_{m \rightarrow \infty} \sum_{\hat{\psi \in L_{\Pi_{\Diamond B}}.s[0:m]}}P(\hat{\pi}) & \\
	&= 1-1
	\end{flalign*}      
	\paragraph{} The third equality above is due to the fact that the first $m$ paths are correctly computed and therefore the error between the approximate solution and the exact solution is zero for $r= 0,1,\cdots,k-1$. The last inequality is obtained after proving that $|\widehat{f_{Rew(\hat{\pi})}}[r] - f_{Rew(\hat{\pi})}[r])| \leq 1$, due to both $\widehat{f_{Rew(\hat{\pi})}}[r]$ and $f_{Rew(\hat{\pi})}[r])$ being bounded in $[0,1]$. The approximate pmf $\widehat{f_{Rew(\hat{\pi})}}[r]$ can be shown to be bounded in this interval by firstly showing that its absolute value is less than or equal to one via \eqref{eqn:approx_computation}, for $r=0,1,\cdots,n+k-1$. And secondly, by proving that \eqref{eqn:approx_computation} is equivalent to performing \textit{circular convolutions} over vectors each of length $k$ and having non-negative entries. This of which results in a non-negative vector. Thus we have shown that the absolute error tends to zero for each $s \in S_?$ as $n \rightarrow \infty$.        \qed
\end{theorem}

\subsubsection{Special case for approximation} Note that for a given sMRM if all pmf reward vectors have masses concentrated within some finite interval $[0,L]$. Then if our original interval of interest is $[0,k]$, convergence is also met if we keep increasing $k$ beyond $L$ and just use the $N-1$ padding scheme (see \eqref{eqn:pi_widetilde}). This is since each pmf vector will be non-zero for length $L$ before being strictly zero up to the $k^{th}$ point. Therefore this becomes equivalent to the approximate solution, with the added benefit of not truncating non-zero values of the reward pmf vectors.

\section{Experiments}
\label{section_experiments}

We proceed to experiment with the direct algorithms developed in this section. We present firstly two problems as a sanity check for our work, one a toy problem, the other a potential real world problem from literature. We then compare the scalability of Gaussian elimination with the approximate variant. The computer used for this entire section is \textit{Computer 1}. See Section \ref{sec:computer details} for details of the computer and software used.

\paragraph{Implementation details}
\label{imp:Gauss_&_LU}
Note that we implemented the $\mathtt{conv_k}$ operations of Gaussian elimination with a time complexity of $O(klog_2k)$. The $\mathtt{deconv_k}$ operations were implemented with a complexity of $O(k^2)$, although there are algorithms with complexity $O(klog_2k)$. The LU approximation method (see Section \ref{sec:approx_direct}) was solved for each $\tau$ using \textit{numpy}'s \textit{solve} function via Gaussian elimination with LU decomposition.

\subsection{Example 1: Toy problem}

\label{problem:5statesMRM}
We are given an sMRM with five states $\{s_0, s_1, s_2, s_3,s_4\}$, with each state being able to reach the goal state $s_4$. The \textbf{property} we are interested in is $Pr(r\ \cap\ s \vDash \Diamond s_4)$, for all $s \in S_? = \{s_0, s_1, s_2, s_3\}$ and $r = 0,1,2,\cdots,N-1$, where $N= 150$.

\paragraph{} We define a matrix $A$ to represent transition probabilities between states in $S_?$ and a vector $b$ to represent the probability of transiting from a state in $S_?$ to $s_4$.

\[ A =
\begin{blockarray}{ccccc}
& s_0 & s_1 & s_2 & s_3 &  \\
\begin{block}{c(cccc)}
s_0 & 0.1288838  &0.38242891& 0.12495781 &0.13139189\\
s_1 & 0.27758284 &0.09654253& 0.15592425 &0.24690511 \\
s_2 & 0.10418887 &0.18054794& 0.1492027  &0.32815053 \\
s_3 & 0.33540355 &0.31410283& 0.16746947 &0.1316041 \\
\end{block}
\end{blockarray}
\]
\[    
b =  \begin{blockarray}{cccc}
s_0 & s_1 & s_2 & s_3   \\
\begin{block}{(cccc)}
0.23233759 &0.22304527 &0.23790995 &0.05142005\\
\end{block}
\end{blockarray}
\]

\begin{figure}[H]
	\centering

	\begin{tikzpicture}[shorten >=2pt,node distance=4cm,auto,thick]
	\begin{scope}[]
	\node[state,initial,scale=1.3] (l0) at (0,0) {$s_0$};
	\node[state, right of=l0,scale=1.3] (l1)  {$s_1$};
	\node[state, above of=l0,scale=1.3]  (l2) {$s_2$};
	\node[state, above of=l1,scale=1.3] (l3)  {$s_3$};
	\node[state,accepting, above of=l3,scale=1.3] (l4)  {$s_4$};
	
	\end{scope}
	
	\begin{scope}[]
	\path [->] (l0) edge[loop below] node {\tt $0.1288838$} (l0);    
	\path [->] (l0) edge node {\tt } (l1);
	\path [->] (l0) edge[bend left]node {} (l2);
	\path [->] (l0) edge node {\tt } (l3);
	\path [->] (l0) edge[bend left=100] node[left] {\tt $0.23233759$} (l4);
	\path [->] (l1) edge[bend left] node {\tt } (l0);    
	\path [->] (l1) edge[loop below] node {\tt $0.09654253$} (l1);
	\path [->] (l1) edge[bend left]node {} (l2);
	\path [->] (l1) edge node {\tt } (l3);
	\path [->] (l1) edge[bend right=120] node[right] {\tt $0.22304527$} (l4);
	\path [->] (l2) edge node {\tt } (l0);    
	\path [->] (l2) edge node {\tt} (l1);
	\path [->] (l2) edge[loop above]node {\tt $0.1492027$} (l2);
	\path [->] (l2) edge node {\tt } (l3);
	\path [->] (l2) edge node[right] {\tt $0.23790995$} (l4);
	\path [->] (l3) edge[bend left] node {\tt } (l0);    
	\path [->] (l3) edge[bend left] node {\tt } (l1);
	\path [->] (l3) edge[bend left]node {} (l2);
	\path [->] (l3) edge[loop above] node {\tt $0.1316041$} (l3);
	\path [->] (l3) edge[bend right=45] node[above] {\tt $0.05142005$} (l4);
	\path [->] (l4) edge[loop above] node[above] {\tt $1$} (l4);
	\end{scope}
	\end{tikzpicture}
	\caption{A DTMC with one absorbing state $s_4$. We annotated only some of the probabilities of $A,b$, to reduce clutter. The figure also represents a stochastic Markov reward model (sMRM) if we annotate each transition with the binomial distribution for the reward random variable.}
	\label{fig:toy_example_discrete_reward}
\end{figure}
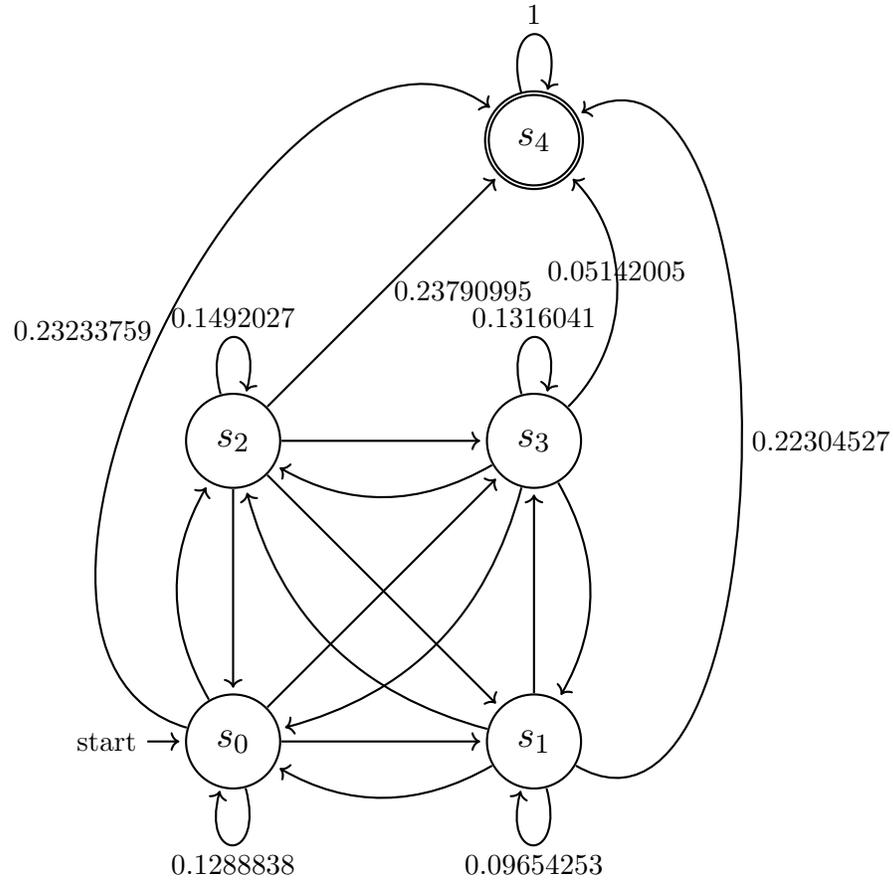

\paragraph{} Transition probabilities from $s_i$ to the three other states corresponds to the $i^{th}$ row of the matrix $A$. Likewise, its probability of entering $B$ immediately is the $i^{th}$ element of $b$. The respective (underlying) DTMC of the sMRM is shown in Figure \ref{fig:toy_example_discrete_reward}.

\paragraph{} We define the reward random variable for every transition to be equal to the binomial distribution with parameters $n = 100, p = 0.5$. We denote the pmf of this distribution as $binom_{n,p}[r]$. From this we can write a system of convolution equations $ \textbf{f} = (A \circ \textbf{G}) \textcircled{$\ast$} \textbf{f} + \textbf{h}$
where $G_{i,j}[r] = binom_{n,p}[r]$, and $h_i[r] = b_i\cdot binom_{n,p}[r]$ for
$r= 0,1,\cdots,N-1$. Figure \ref{fig:toy_example_discrete_reward} can be used to represent the sMRM if we annotate each transition with the binomial random variable.

\begin{figure}[H]
	\centerline{
		\includegraphics[width=1\linewidth]{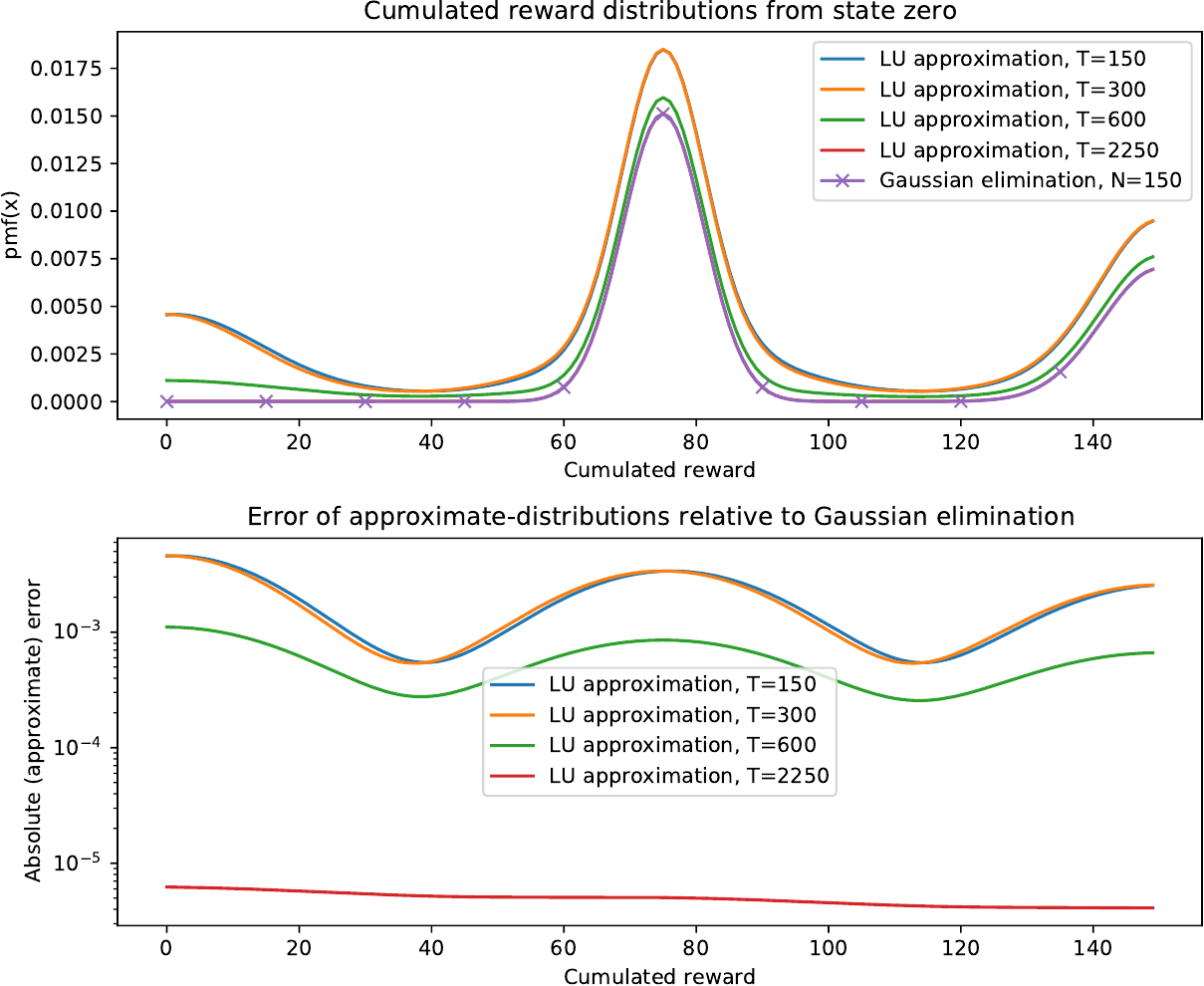}}
	\caption{Top figure: The pmf of the accumulated reward from $s_0$ until $s_4$, (truncated) up to the $N^{th}$ point. Bottom figure: the absolute error between the approximate solutions relative to Gaussian elimination. For the LU approximations, the different padding lengths can be determined by $T-N$ where $T$ stands for Total length.}
	\label{fig:gevslu}
\end{figure}

\paragraph{} In Fig. \ref{fig:gevslu}, we compare the approximate LU method (see Section \ref{sec:approx_direct}) with naive Gaussian elimination. For the LU method, we compute $Pr(r\ \cap\ s \vDash \Diamond s_4)$ for each $s \in S_?$ (with $B = s_4$) and for $r= 0,1,\cdots,T$, where $T = N + n$, and $n$ is the zero padding length (described in Section \ref{sec:approx_direct}). We also plot the absolute (approximate) error between the Gaussian elimination solution and the approximations. As previously described, the errors between the pmfs obtained via the approximation solutions and Gaussian elimination are the effects of time-aliasing (or wrap-around error) relating to the choice of padding length $T-N$. This example shows that the longer the padding length, the more accurate the result.

\paragraph{} The times taken for each method are found in Table \ref{fig:guvslutimes-crop}. The times measured include both the time taken for \textbf{solving the system of equations} (the solving phase), as well as time taken \textbf{creating the hypermatrices} required for solving the system  (i.e. the preparation phase). This is how we will time all our algorithms in the future unless mentioned otherwise. 


\begin{table}
	\centerline{
		\includegraphics[width=1\linewidth]{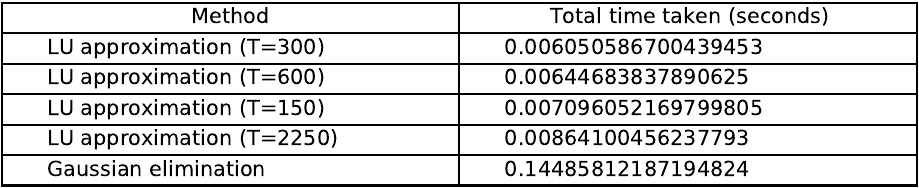}}
	\caption{Time taken to solve the problem for various direct methods, ordered from fastest to slowest. Gaussian elimination (GE) took about 18 times longer than the LU method with $T=2250$. However, the respective LU method required 15 (2250/150) times more space. }
	\label{fig:guvslutimes-crop}
\end{table}

\begin{figure}[H]
	\leavevmode
	\centerline{
		\includegraphics[width=0.8\linewidth]{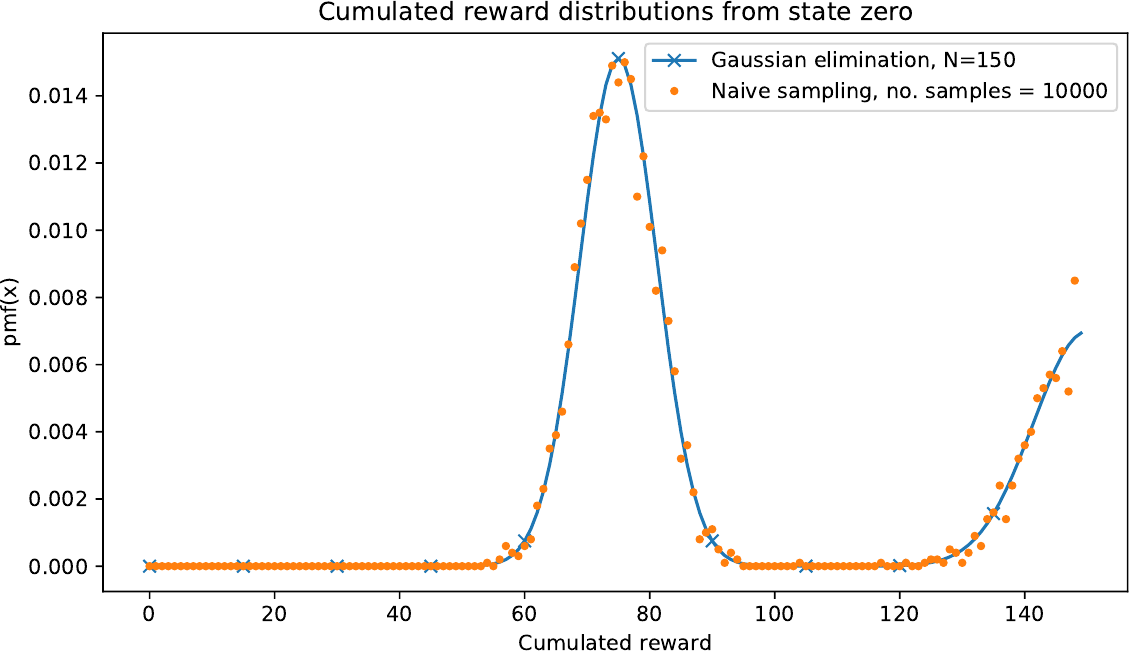}}
	\caption{The pmf obtained via sampling aligns with that obtained via Gaussian elimination. However, the naive sampling method did take around 12 times longer to achieve this result.}
	\label{fig:guvssampling}
\end{figure}

\paragraph{} Next, we sampled 10000 traces from the sMRM starting from state $s_0$, computed the cumulated reward over each trace, and then computed the relative frequencies of the cumulated rewards to  form the pmf $Pr(r\ \cap\ s \vDash \Diamond s_4)$ for $r = 0,1,\cdots,N-1$ (see Fig. \ref{fig:guvssampling}) . The length of each trace was determined by the sampling termination conditions: 1) Entering the goal state $s_4$. 2) Accumulating reward equal to or greater than $N$. We find that the results  from sampling align with those obtained from Gaussian elimination providing validation for the numerical algorithm. The total time taken for the sampling algorithm was 1.77 secs (2dp).

\subsection{Example 2: Waste treatment semi-Markov process}
\label{problem:waste_treatment}
This problem comes from \cite{warr2014numerical}, with the model originally from \cite{barbu2008semi}, who studied the computation of the passage-time density for a waste treatment model. The sMRM represents a discrete-time semi-Markov process, therefore giving us a discrete-reward sMRM, where reward represents time.

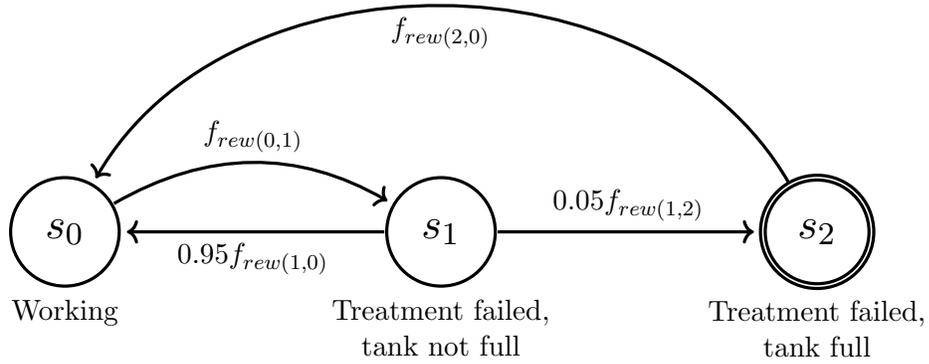
\begin{figure}[H]
	\centering
	\begin{tikzpicture}[shorten >=2pt,node distance=5cm,auto,very thick]
	\begin{scope}[]
	\node[label={[align=center]below:Working},state,scale=1.5] (l0) at (0,0) {$s_0$};
	\node[label={[align=center]below:Treatment failed,\\tank not full},state, right of=l0,scale=1.5] (l1)  {$s_1$};
	\node[label={[align=center]below:Treatment failed,\\tank full},state, accepting, right of=l1,scale=1.5]  (l2) {$s_2$};
	
	\end{scope}
	
	\begin{scope}[]
	\path [->] (l0) edge[bend left] node {$f_{rew(0,1)}$} (l1);    
	\path [->] (l1) edge node {\tt $0.95f_{rew(1,0)}$} (l0);    
	\path [->] (l1) edge node {\tt $0.05f_{rew(1,2)}$} (l2);
	\path [->] (l2) edge[bend right=60] node {$f_{rew(2,0)}$} (l0);      
	\end{scope}
	\end{tikzpicture}
	\caption{A graphical representation of the sMRM modelling a function of the textile factory.}
	\label{fig:treatment_Markovchain}
\end{figure}

\paragraph{} Figure \ref{fig:treatment_Markovchain} captures a model of the waste treatment in a textile factory. The factory has a production which generates waste and also a waste treatment facility to handle such waste. If the facility fails, untreated waste can be stored in a holding tank. However, if the facility fails for too long, for example if repairs are not carried out fast enough, the holding tank fills completely and production of the factory stops.

\paragraph{}The reward distributions $f_{rew(i,j)}$ are defined as
\begin{flalign*}
&f_{rew(0,1)} \sim geometric(0.8) &\\
&f_{rew(1,0)} \sim discreteWeibull(0.3, 0.5) &\\
&f_{rew(1,2)} \sim discreteWeibull(0.5, 0.7) &\\
&f_{rew(2,0)} \sim discreteWeibull(0.6, 0.9)
\end{flalign*}
and the discrete Weibull distribution is given as
\[
f(t;q,b) = \begin{cases}
q^{(t-1)^b} - q^{t^b} & t \in [1,2,3,\cdots] \\
0 & otherwise
\end{cases}
\]

\paragraph{} The authors of \cite{warr2014numerical} computed the probability of reaching state $s_2$, the state where production halts, starting from $s_0$ when the treatment facility is functional, i.e. $Pr(r\ \cap\ s_0 \vDash \Diamond s_2)$. We compute this distribution for $r= 0,1,\cdots,N-1$ where $N=100$ and present the pmf below for the cumulated reward (which is time in this case). We used the Gaussian elimination method to solve this problem.

\begin{figure}[H]
	\centering
	\includegraphics[width=1\linewidth]{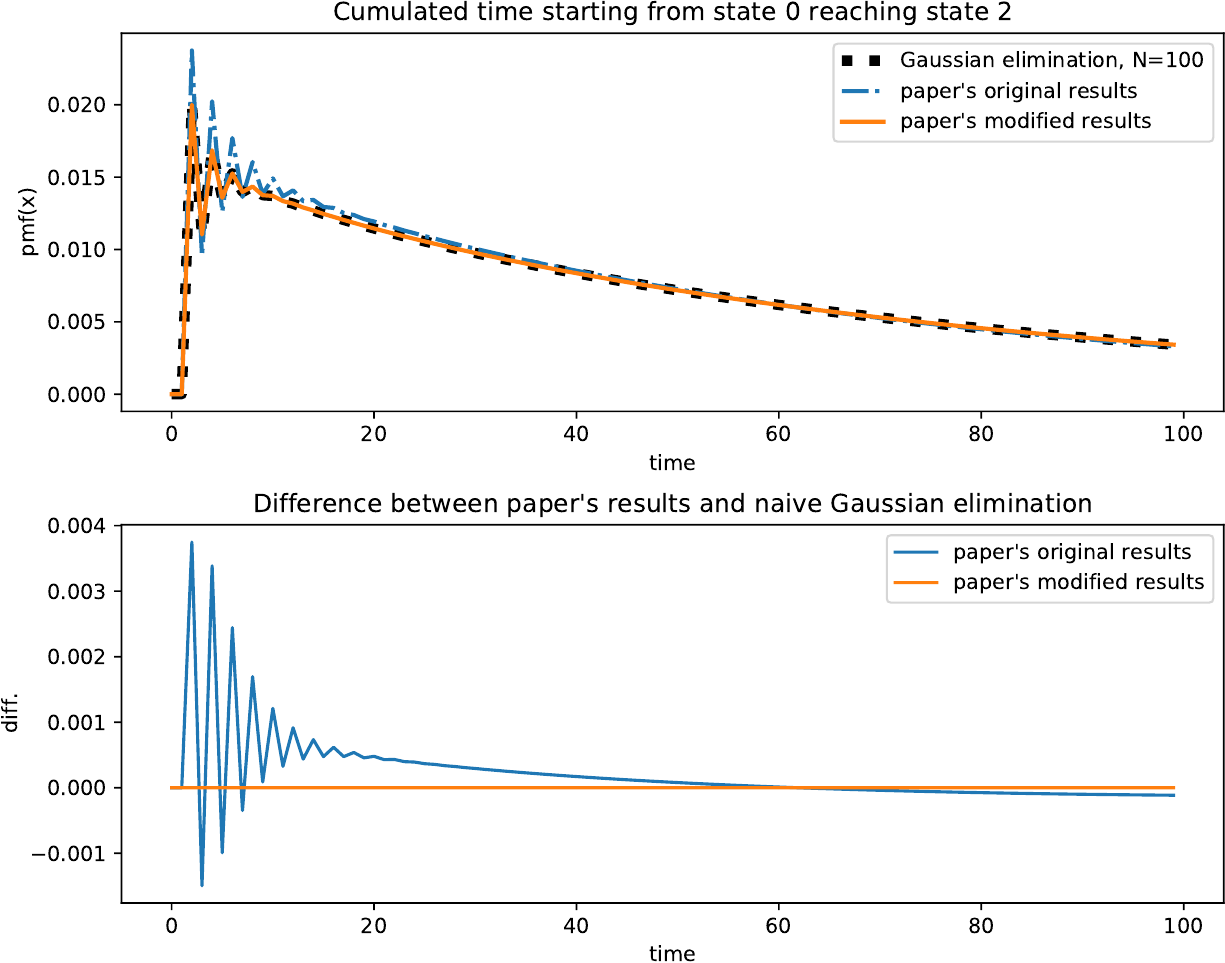}
	\caption{Top figure: The pmf over time elapsed before treatment has failed and the tank is full ($s2$), starting from $s0$. After replacing the DTFT with the DFT in the paper's implementation, it aligned better with our results. Bottom figure: The difference in pmf values relative to the naive Gaussian elimination algorithm.}
	\label{fig:wastetreatmentproblem-crop}
\end{figure}

\paragraph{} When comparing our solution to those presented in \cite{warr2014numerical} (see Fig. \ref{fig:wastetreatmentproblem-crop}), our results (via Gaussian elimination) mostly align with theirs. We experimented with their code to find that they had not used the DFT for the geometric distribution, but rather used the known algebraic DTFT (characteristic function) of the distribution (as was mentioned in their paper), thus leading to their slightly incorrect results. When we modified their implementation by using the DFT for the geometric distribution, the values of the resulting pdf is nearer to that obtained from Gaussian elimination. From our understanding, the inverse DFT cannot be used naively to invert the continuous FT.

\subsection{Scalability of GE vs LU approximation}    
\label{subsection:scalability_GE_vs_LU}

We experimented with Markov chains with random discrete reward random variables, and timed how long it took to solve for the property $\mrmpdf$. Both the preparation phase and solving phase are included in the timing. We measure the performance of the Gaussian elimination algorithm with the LU approximation method. With the latter method, we used the $numpy$ package in python to solve each linear system using the LU decomposition as part of the procedure. Therefore the test is not fair, however the comparison may still be useful.

\paragraph{} To measure the performance of the algorithms, we vary two variables (each independently of the other): 1) $|S_?|$, the size of the set of states that can eventually reach $B$ our goal states. 2) $k$, such that $\mrmpdf$ is computed for $r = 0,1,\cdots,k-1$. Performance is characterized only by time. When varying $|S_?|$ alone, $k$ was fixed to $501$. And when varying $k$ alone, $|S_?|$ was fixed to six.

\paragraph{} For any valuation of the parameters ($|S_?|,k$), 200 sMRMs were sampled with these parameters, their average and worst times recorded and plotted. For each of these 200 sMRMs, the underlying Markov chain probability matrix is generated randomly and uniformly. We have also only one goal state. The rewards are generated as follows: If the reward is on a transition between a state $s \in S_?$ and the goal state, it is distributed via a geometric distribution with its parameter $p_1$ randomly and uniformly  sampled from $[0,1]$. If it is other than that (i.e. transitions between states of $S_?$), then the reward is a binomial distribution, with its parameters ($n= k\beta$,$p=p_2$), where $k$ is the interval length above, and both $\beta$ and $p_2$ are randomly and uniformly sampled from $[0,1]$. Thus for any transition $s_i \rightarrow B$ for each $s_i \in S_?$, we have
$$f_{rew(s_i,B)}(r) = Geo(r;p_1) = (1-p_1)^{r-1}p_1$$
and for transitions $s_i \rightarrow s_j$ for each $s_i, s_j \in S_?^2$, then
$$f_{rew(s_i,s_j)}(r) = Binomial(r;k\beta,p_2) = \binom{k\beta}{r}p_2^{r}(1-p_2)^{n-r}$$
Then with the system generated, the property $\mrmpdf$ is solved with the different methods (for each experiment). This choice of distributions allows us to sample tail distributions (by the Geometric) and peak distributions (by the Binomial).

\begin{figure}[H]
	\centering{
		\includegraphics[width=1\linewidth]{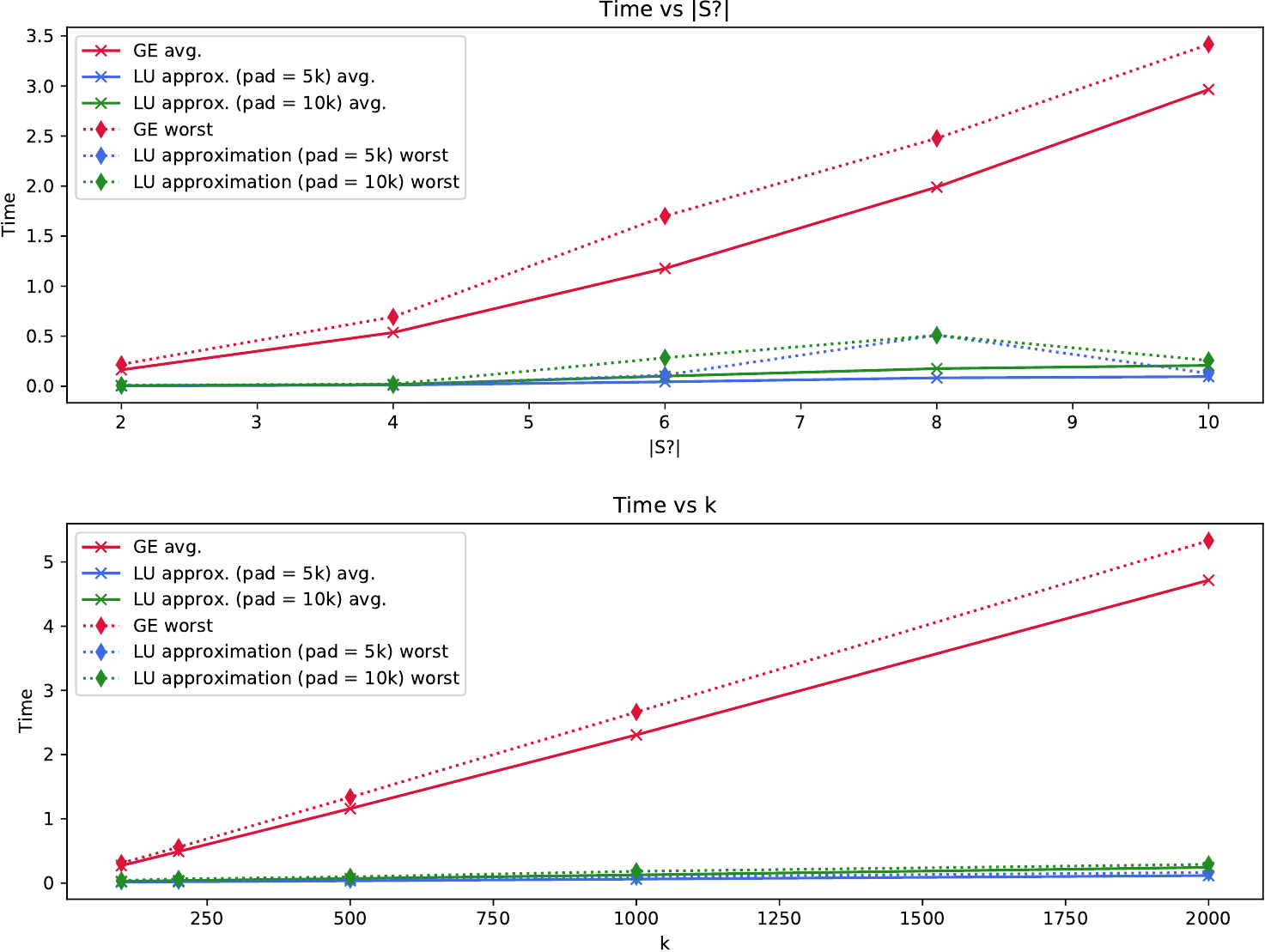}}
	\caption{Time taken to solve $\mrmpdf$ for various values of $|S_?|$ and $k$. For the top graph, $k = 501$, and for the bottom $|S_?|=6$. For each graph, a total of 1000 experiments are performed. We find that the LU approximation method scales significantly better than GE in time. In the legend, (pad = 5k) means that the LU approx. method uses a padding length (see \eqref{eqn:pi_widehat}) of $5k$.}
	\label{fig:scalabilitydirectmethods-crop}
\end{figure}

\paragraph{} The results (Fig. \ref{fig:scalabilitydirectmethods-crop}) show that the exact algorithm is significantly less scalable than the approximate algorithm in terms of time complexity. Furthermore, doubling the padding size of the approximate algorithm did not increase the time taken by much. In fact, adding together the times to compute both LU approximation (with $5k$ and $10k$ padding) still beats the average time for GE when either $|S_?|$ or $k$ are large. Additionally, by computing both of these LU approximations we are able to determine the relative absolute error (by taking the maximal difference of their pmfs) and thus determine the accuracy of the approximate solutions. Therefore an iterative algorithm involving the LU approximations may prove useful and solve faster overall. Note however that when the LU approximations are padded with $n$ zeros, then GE will typically have a lower space complexity when $n > k-1$.

\section{Summary and discussion} In this chapter, we presented direct algorithms for solving discrete-reward sMRMs. We presented the Gaussian elimination algorithm for a system of convolution equations, as well as an approximate variant of the algorithm, which has the benefit of being able to be solved directly from existing numerical algebra tools. These algorithms allow us to solve problems with more than six states, as shown from the previous examples and therefore improves upon the symbolic approach in such a manner.


\subsection{Alternative models}
\label{subsec:alt_models}

\paragraph{DTMCs} For reachability problems such as $\mrmpdf$ or $Pr(s \vDash \Diamond_{\leq k} B)$, the approach here with the sMRM is better than naively using DTMCs (defined in Definition \ref{definition:DTMC}), in that it yields better space complexity. For example, if we integrate into the state space of a DTMC our reward space, then the new state space is $ S \times \{0,1,\cdots,k\}$, which means in the worst-case the space complexity of our transition matrix \textbf{P} is $O((S \times \{0,1,\cdots,k\}) \times (S \times \{0,1,\cdots,k\}))$. Compare this with an sMRM, where the largest term is the matrix \textbf{G}, with a worst-case complexity of $S \times S \times \{0,1,\cdots,k\}$.

\paragraph{HMMs} For the property $Pr(s \vDash \Diamond_{\leq k} B)$, it is question as to whether we can use a (transition-based) Hidden Markov model (HMM) instead of an sMRM to reduce the computational costs. Firstly, let $\Pi.s$ be the set of paths that begin in state $s$ and end in $B$. Then,
\begin{flalign*}
Pr(s \vDash \Diamond_{\leq k} B)
= \sum_{\hat{\pi} \in \Pi.s}Pr(\hat{\pi})Pr(Rew \leq k | \hat{\pi})
= \sum_{\hat{\pi} \in \Pi.s}Pr(\hat{\pi})F_{Rew(\hat{\pi})}(k)
\end{flalign*}
where for each path $\hat{\pi}$, the probability of accumulating $\leq k$ reward is determined by the cdf $F_{Rew(\hat{\pi})}$. The probability $F_{Rew(\hat{\pi})}(k)$ is dependent on the whole path $\hat{\pi}$. It is not possible to define the probability of accumulating $\leq k$ reward at any step $n$ in $\hat{\pi}$ to be strictly dependent on the transition made at time $n-1$.

\paragraph{} To elaborate, to use an HMM instead of an sMRM, we need to firstly replace the reward structure (see Def. \ref{definition:sMRM}) with a set of \textit{hidden states}. Let us first give a formal definition of HMMs.
\begin{definition}        
	A HMM is a tuple $(\mathcal{M}, O, \textbf{P}_o)$, where:
	\begin{itemize}
		\item $\mathcal{M}$ is a DTMC, i.e. $\mathcal{M} = (S, \textbf{P}, i_{init})$.
		\item $O$ is the hidden state space (or observation state space).
		\item $\textbf{P}_o: S \times S \times O \rightarrow \mathbb{R}$ is a map for the emission probabilities of states in $O$, for each transition $s \rightarrow t$, with $s,t \in S^2$.
	\end{itemize}
	\paragraph{} A property of HMMs is that at time step $n$, the probability of being in a hidden state $o$ is dependent on the transition $s\rightarrow t$ that the HMM is made at time step $n-1$.
\end{definition}
\paragraph{} Thus, let us create two hidden states: \textit{less-than-equal-to-k} and \textit{more-than-k}. The probability of being in \textit{less-than-equal-to-k} at a given time is only dependent on the transition made just prior to that time. Let $\textit{p}_n$ denote the probability of being in \textit{less-than-equal-to-k} at step $n$. Thus for a path $\hat{\pi}= s_0,s_1,\cdots,s_{|\hat{\pi}|-1}$, we want
\begin{flalign*}
F_{Rew(\hat{\pi})}(k) = Pr(\ \textit{p}_1\ \cap\ \textit{p}_2\ \cap\ \cdots\ \cap\ \textit{p}_{|\hat{\pi}|-1} \ |\ \hat{\pi}))
= \prod_{i=0}^{|\hat{\pi}|-2}Pr(\ \textit{p}_{i+1}\ |\ s_{i} \rightarrow s_{i+1})
\end{flalign*}

for some definition of $Pr(\ \textit{p}_{i+1}\ |\ s_i \rightarrow s_{i+1})$. There are no general definitions available that would satisfy the equality above. One false definition would be $Pr(\ \textit{p}_{i+1}\ |\ s_i \rightarrow s_{i+1}) = F_{rew(s_i,s_{i+1})}(k)$. Therefore, we cannot define an HMM (or Markov risk model) this way.

\paragraph{} In the next chapter, we move towards deriving iterative algorithms instead: The power, Jacobi and Gauss-Seidel methods from existing literature will be adapted for sMRMs.

\chapter{Iterative Methods for sMRMs}

\section{Introduction}
In this chapter we propose adaptations of existing iterative algorithms to solve the set of equations \eqref{eqn:lineq_matrix_trans_pdf} for the discrete case. Iterative methods begin with an initial guess for the solution, and for those that converge, they continuously produce better solutions at each \textit{iteration} of the method. Iterative methods are regarded as being more scalable than direct algorithms such as Gaussian elimination, i.e. they can be more efficient when $|S_?|$ or $k$ is large. Additionally, an advantage over direct algorithms is that the (numerical) error produced by the algorithm can be controlled quite well. However, a possible disadvantage is that the time it takes to solve a problem is not generally known in advance, i.e. it is problem dependent, unlike for direct methods.

\paragraph{} Here, we will discuss the power, Jacobi and Gauss-Seidel methods. Like Gaussian elimination, these algorithms also have an approximate form.  A few theorems are also presented to show that results for systems of linear equations exist in a similar fashion for systems of convolution equations. We also empirically evaluate the performance of these methods, with the results presented at the end of the chapter.

\paragraph{} The definitions of terms from the previous chapter will be used here, for example of convolution, deconvolution, and symbols like $\tilde{p}_i$ and $\widehat{p}_i$, from equations \eqref{eqn:pi_widetilde} and \eqref{eqn:pi_widehat} respectively.

\section{The power method}
\label{sec:power_method}  Let us devise an iterative sequence using the set of equations \eqref{eqn:lineq_matrix_trans_pdf}:
\begin{flalign}
\label{eqn:lineq_power_method_pdf}
\textbf{f}^{(n+1)} = (A\circ \textbf{G})\textcircled{$\ast$} \textbf{f}^{(n)} + \textbf{h}
\end{flalign}
of which if we let \begin{flalign}\label{eqn:initial_cond_f}
\textbf{f}^{(0)} = 0
\end{flalign}    
then we have $\lim_{n \rightarrow \infty}\textbf{f}^{(n)} = \textbf{f}$, where \textbf{f} is the solution for \eqref{eqn:lineq_matrix_trans_pdf}. The proof of this is given later in Theorem \ref{theorem:fixed_point}. This method is a generalization of the \textit{power method} described in \cite{principles_book,zapreev2006safe}. As for the historical origins of the method, in discussion of it, the article \cite{golub2000eigenvalue} states that

\paragraph{} \textit{ Householder \cite{householder2013theory}  called this Simple Iteration, and attributed the  first treatment of it to M$\ddot{\text{u}}$ntz \cite{muntz1913solution}. Bodewig \cite[p. 250]{bai1999able} attributes the power method to von Mises \cite{mises1929praktische}, and acknowledges M$\ddot{\text{u}}$ntz for computing approximate eigenvalues from quotients of minors of the explicitly computed matrix $A^k$, for increasing values of k.}

\paragraph{} We can approximate $\lim_{n \rightarrow \infty}\textbf{f}^{(n)}$ by iterating through $n$ and stopping when the absolute max error between the approximation at time $n$, and $n+1$ is small enough. In effect when
\begin{flalign}
\label{eqn:stopping_criteria}
\mathtt{max} | \textbf{f}^{(n)} - \textbf{f}^{(n+1)} | \leq \epsilon    
\end{flalign}
where $\epsilon$ is the absolute error tolerance level, and where we have defined
$$\mathtt{max}|\textbf{f}^{(n+1)} - \textbf{f}^{(n)}| \triangleq \mathtt{max}_{s,r} |f_s[r]^{(n+1)} - f_s[r]^{(n)}|$$

\paragraph{} The termination or stopping criterion above however is not sound, in that it does not prove that $\mathtt{max} | \textbf{f}^{(n)} - \textbf{f}^{(n+M)}| \leq \epsilon$ for some large natural number $M$. For an example of false convergence, see \cite{haddad2018interval}. A way to extend this algorithm in a manner that ensures soundness is detailed in Section \ref{sec:summ_disc}.

\paragraph{} A practical form of \eqref{eqn:lineq_power_method_pdf} is achieved by replacing $*,\deconv$ with $\mathtt{conv_k},\mathtt{deconv_k}$, i.e. for each $s \in S_?$, we have
\begin{flalign}
\label{eqn:power_method_early_conv}
f_s[r]^{(n+1)} =  \sum_{t \in S_?}\mathtt{conv}_k(A_{s,t}G_{s,t},\ f_t^{(n)}) + h_s[r]
\end{flalign}
where $f_s[r]^{(n+1)}$ is computed only for $r=0,1,\cdots,k-1$. The size of each term is defined as in Gaussian elimination, Section \eqref{sec:naive_GE}.

\paragraph{} Since the FFT is used for convolution above, each iteration has a worst-case complexity of $O(|S_?|^2(2k-1)log_2(2k-1))$. However, we can improve upon this if we realise that the hypermatrix $A\circ\textbf{G}$ is constant and so for each iteration, we do not need to repeatedly compute their DFT transforms used in $\mathtt{conv}_k$.    Let us introduce the hypermatrix $\bm{\mathcal{G}}$ with size $(2k-1 \times |S_?| \times |S_?|)$ where $\mathcal{G}_{i,j}[\tau] = \mathcal{D}\{\widetilde{A_{i,j}G_{i,j}}\}[\tau]$. Of which we introduce the matrix notation as  $$\bm{\mathcal{G}} = \mathcal{D}\{\widetilde{A\circ \textbf{G}}\}$$
\paragraph{} Then, we can use the iteration
\begin{flalign}
\label{eqn:power_method_matrix_practical}
\textbf{f}^{(n+1)} =  \mathcal{D}^{-1}\{\bm{\mathcal{G}}(\mathcal{D}\{\widetilde{\textbf{f}}^{(n)}\})\} + \textbf{h}
\end{flalign}
where we have overloaded the definition of $\mathcal{D}$ as just previously.  The proof that this iterative system converges to the true solution $\textbf{f}[r]$ for $r=0,1,\cdots,k-1$, is presented in Theorem \ref{theorem:length_of_DFTs}.

\subsubsection{Algorithmic complexity} The DFT transform of $\widetilde{A\circ\textbf{G}}$ has complexity $O(|S_?|^2(2k-1)log_2(2k-1))$, but only has to be performed once. Next, in each iteration we perform a DFT transform and an inverse transform over a hypervector. This has total complexity $O(2|S_?|(2k-1)log_2(2k-1))$. Convolutions are now multiplications, and we require $O(|S_?|^2(2k-1))$ multiplications. The number of additions required is  $O((|S_?|^2(2k-1) + |S_?|k)$. The convergence check has $O(|S_?|(2k-1))$ subtractions and $O(|S_?|(2k-1))$ scans to find the maximum value. Thus for each iteration, asymptotically in the worst-case we have $$O(|S_?|(2k-1)log_2(2k-1)) + O(|S_?|^2(2k-1))$$
with a one-time cost of $$O(|S_?|^2(2k-1)log_2(2k-1))$$ for the preparation-phase of the solution.

\subsection{Fixed point characterization}
\label{proofs:power method}

\paragraph{} We reproduce here the three sub-theorems presented in \cite[Theorem 10.15]{principles_book}, in the context of sMRMs. These theorems were originally presented for Markov chains. Note that whilst we prove these statements for the discrete case (i.e. with pmfs), they also \textbf{hold for the continuous case (with pdfs)}.  

\begin{theorem}[Fixed point characterization (for pmfs)]
	\label{theorem:fixed_point}
	The cumulated reward hypervector $\textbf{f} = (\textbf{f}[r])_{r \in \mathbb{N}} = (f_{s \vDash \Diamond B}[r])_{r \in \mathbb{N},s \in S_?} =  ( Pr(r \cap s \vDash \Diamond B) )_{r \in \mathbb{N},s \in S_?}$ is the (unique) fixed point of the operator $\Upsilon: [0,1]^{\mathbb{N} \times  S_?\times 1} \rightarrow [0,1]^{\mathbb{N} \times S_? \times 1}$. This operator is defined as:
	\[\Upsilon(\textbf{f}) = ((A\circ\textbf{G})\ \textcircled{$\ast$} \ {\textbf{f}}) + \textbf{h} \]
	
	
	\paragraph{} Additionally, let $\textbf{f}^{(0)} = \textbf{0}$, and $\textbf{f}^{(n+1)} = \Upsilon(\textbf{f}^{(n)})$ where $n \geq 0$. Then, for any ${r \in \mathbb{N}}$ and  $s \in S_?$ the following three statements hold:
	
	\begin{enumerate}
		\item $f^{(n)}_s[r] = Pr(\ r\  \cap [s \vDash \Diamond^{\leq n} B])$, for all $n \geq 0$.
		\item $ \lim\limits_{n \rightarrow \infty} f_s^{(n)}[r] = f_s[r]$.
		\item $f_s^{(0)}[r] \leq  f_s^{(1)}[r] \leq  f_s^{(2)}[r]  \leq \cdots \leq  f_s[r]$.
	\end{enumerate}
	where statement 2. states that the solution converges to a fixed point, and 3. states that the convergence is monotonic. Due to the simplicity of these statements, we have presented the proofs in the Appendix, see Theorem \ref{theorem:fixed_point_full}.
\end{theorem}

\begin{corollary}[Uniqueness of the fixed-point solution]
	\label{corollary:unique_power}
	Note that the solution to the fixed point operator defined in Theorem \ref{theorem:fixed_point} is unique.
	\paragraph{} \textit{Proof:} This follows from Theorem \ref{theorem:unique sol}. Whilst the proof was given for the system in the Fourier domain, since there is a one-to-one relationship between pmfs and characteristic functions, then $\phi_{s}$ being unique $\implies$ $f_{s}$ is unique. \qed
	
\end{corollary}

\subsection{Convergence of the exact power method} For the iterative sequence
$$\textbf{f}^{(n+1)} = (A\circ \textbf{G})\textcircled{$\ast$} \textbf{f}^{(n)} + \textbf{h}$$
we want to make sure that the computation is correct for each iteration $n$,
up to the $k^{th}$ point, i.e. that $\mathbf{f}^{(n)}[r]$ is computed correctly for $r=0,1,\cdots,k-1$, and for any $n \geq 0$. As shown from the following lemma, this ensures that $\lim_{n \rightarrow \infty}\textbf{f}^{(n)}[0:k] = \textbf{f}[0:k]$ is computed correctly.

\begin{lemma}[Convergence of the first $k$ values]
	\label{lemma:convergence of first k values}
	Note that the pmf vector $f_s^{(n+1)}$ can be split into a concatenation of two vectors: $f_s^{(n+1)} = f_s^{(n+1)}[0:k]^\frown f_s^{(n+1)}[k:\infty]$.
	
	\paragraph{} Then $\lim_{n \rightarrow \infty}f_s^{(n)}[0:k]$ converges without requiring the computation of $\lim_{n \rightarrow \infty}f_s^{(n)}[k:\infty]$. In effect, we can compute $f_s^{(1)}[0:k], f_s^{(2)}[0:k], \cdots, f_s^{(n)}[0:k],\cdots$ and this is sufficient for the convergence of  $\lim_{n \rightarrow \infty}f_s^{(n)}[0:k]$.
	
	\paragraph{} \textit{Proof:} For any $f_s^{(n+1)}[r]$, we have that $$f_s^{(n+1)}[r] = \sum_{t \in S_?}\textbf{P}(s,t)(\sum_{i=0}^{r})f_{rew(s,t)}[r-i]f_t^{(n)}[r] + \sum_{u \in B}\textbf{P}(s,u)f_{rew(s,u)}[r]$$
	Therefore to calculate $f_s^{(n+1)}[r]$ we only need values $(f_t^{(n)}[i])_{t \in S_?}$ for all $i$ s.t. $i \leq r$. Thus $(f_s^{(n+1)}[0:k])_{s \in S_?}$ requires only the values  $(f_t^{(n)}[0:k])_{t \in S_?}$, and therefore does not require computing $(f_s^{(n+1)}[k:\infty])_{s \in S_?}$. Hence  $\lim_{n \rightarrow \infty}f_s^{(n)}[0:k]$ can be determined without $\lim_{n \rightarrow \infty}f_s^{(n)}[k:\infty]$.  \qed
	
	\paragraph{} Note that the same result also holds for continuous rewards (with pdfs).
\end{lemma}

\paragraph{} The following theorem presents a way to compute $\textbf{f}^{(n)}[0:k]$ correctly for any $n  \geq 0$ and using DFTs. Using this theorem and Lemma \ref{lemma:convergence of first k values}, we can iterate towards $\lim_{n \rightarrow \infty}f_s^{(n)}[0:k]$ correctly.

\begin{theorem}[Convergence of the (practical) power method]  \label{theorem:length_of_DFTs} The iterative sequence via the practical form of the power method \eqref{eqn:power_method_matrix_practical} converges to the true solution $f_s[r]$ for $r = 0,1,\cdots,k-1$ and for all $s \in S_?$.

	\paragraph{} \textit{Proof: } We prove that the following equivalence
	$$((A\circ \textbf{G})\textcircled{$\ast$} \textbf{f}^{(n)})[r] \equiv \mathcal{D}^{-1}\{\mathcal{D}\{\widetilde{A \circ \textbf{G}}\}(\mathcal{D}\{\widetilde{\textbf{f}}^{(n)}\})\}[r]$$
	holds for $r = 0,1,\cdots,k-1$. This will then ensure that $\textbf{f}^{(n+1)}[0:k]$ is computed correctly for all $n$. First define the hypervectors $\textbf{g}= (A\circ \textbf{G})\textcircled{$\ast$} \textbf{f}^{(n)}$, and $\textbf{l} = \mathcal{D}^{-1}\{\widetilde{A \circ \textbf{G}}\}(\mathcal{D}\{\widetilde{\textbf{f}}^{(n)}\})\}$. Then $g_s[r] = \sum_{t \in S_?}(A_{s,t}G_{s,t} \ast f_t^{(n)})[r]$ and
	$l_s[r] = \sum_{t \in S_?} \mathcal{D}^{-1}\{\mathcal{D}\{\widetilde{A_{s,t}G_{s,t}}\}   \mathcal{D}\{\widetilde{f}_t^{(n)}\} \}[r] $. By \eqref{eqn:DFT_sums}, for any $s,t \in S_?^2$, we have
	$$(A_{s,t}G_{s,t} \ast f_t^{(n)})[r] =  \mathcal{D}^{-1}\{\mathcal{D}\{\widetilde{A_{s,t}G_{s,t}}\}   \mathcal{D}\{\widetilde{f}_t^{(n)}\}\}[r]$$ for $r=0,1,\cdots,k-1$.
	This is sufficient to deduce that $g_s[0:k] \equiv l_s[0:k]$ proving the equivalence. Then, $\textbf{f}^{(n+1)}[0:k] = \textbf{g}[0:k] + \textbf{h}[0:k] = \textbf{l}[0:k] + \textbf{h}[0:k]$. Hence, we can compute $\textbf{f}^{(n)}[0:k]$ correctly for any $n$ via the practical form of the power method. By Lemma \ref{lemma:convergence of first k values}, we only need to correctly compute $\textbf{f}^{(n)}[0:k]$ for each $n \geq 0$, to iterate towards the fixed point $\lim_{n \rightarrow \infty}\textbf{f}^{(n)}[0:k]$. By Corollary \ref{corollary:unique_power}, we then have convergence to $\textbf{f}[0:k]$. \qed
\end{theorem}

\section{The approximate power method}
\label{sec:approx_power}
To speed up computation, we may employ the iterative analogue of the approximate direct solution in Section \ref{sec:approx_direct}. Doing so, we obtain the iterative sequence
$$
\textbf{x}^{(n+1)} =(A\circ\textbf{C}_m)\textbf{x}^{(n)} + \textbf{d}_m
$$
where the sizes of terms and definitions remain identical (see Section \ref{sec:approx_direct}). Note however that $m$ is used to identify the zero padding length used in \eqref{eqn:pi_widehat}. We do not use $n$ as notation since we are already using it for the iteration count.

\paragraph{}This power method is termed \textit{approximate} since $\lim_{n \rightarrow \infty} \mathcal{D}^{-1}\{\textbf{x}^{(n+1)}\} \approx \textbf{f}$. With the same reasoning, we term the previous power method as \textit{exact}.  For a given padding length $m$, the solution from the LU approximate method (see Section  \ref{sec:approx_direct}) equals the limit of the approximate power method $\lim_{n \rightarrow \infty}\mathcal{D}^{-1}\{\textbf{x}^{(n)}\}$. Therefore, as with the LU approximation method, when $m \rightarrow \infty$, then $\lim_{n \rightarrow \infty} \mathcal{D}^{-1}\{\textbf{x}^{(n+1)}\} = \textbf{f}$.

\paragraph{} In a similar fashion as with the exact power method, we can simplify the above by calculating $\bm{\mathcal{G}}_m \triangleq (A\circ\textbf{C}_m)$, and $\bm{\eta}_m \triangleq \textbf{d}_m$, with the  subscript $m$ highlighting the zero-padding length used. We can then use the iteration
\begin{flalign}
\label{eqn:approx_power_method}
\textbf{x}^{(n+1)}(\tau) =\bm{\mathcal{G}}_m(\tau)\textbf{x}^{(n)}(\tau) + \bm{\eta}_m(\tau)
\end{flalign}
with $\textbf{x}^{(0)}(\tau)= \textbf{0}$ as the initial hypervector.

\paragraph{} Since $\textbf{x}$ is a complex number, we can define one convergence criteria as
\begin{flalign}
\label{eqn:approx_method_power_convergence_criteria}
\Big(max_{s,\tau}|Real(\textbf{x}^{(n+1)}) - Real(\textbf{x}^{(n)})| \leq \epsilon\Big) \cap \Big(max_{s,\tau}|Imag.(\textbf{x}^{(n+1)}) - Imag.(\textbf{x}^{(n)})| \leq \epsilon\Big)
\end{flalign}
where the threshold is applied to the real and imaginary parts of $\bm{\phi}_{d,m}$ separately. Another convergence criteria is simply  $max_{s,\tau}|\textbf{x}^{(n+1)} - \textbf{x}^{(n)}| \leq \epsilon$.

\subsubsection{Algorithmic complexity} For each iteration, asymptotically in the worst-case we have a complexity of $O(|S_?|^2(m+k))$ attributed to the multiplications and additions. There is still a one-time cost of $O(|S_?|^2(m+k)log_2(m+k))$ for the preparation-phase of the solution, attributed to the computation of $(A\circ\textbf{C}_m)$.

\begin{theorem}[Convergence of the approximate power method]
	\label{theorem:convergence_power_method_approx}
	The approximate power method \eqref{eqn:approx_power_method} converges to a unique solution \textbf{x} of the set of equations
	$$(I - \bm{\mathcal{G}}_m(\tau))\textbf{x}(\tau) = \bm{\eta}_m(\tau)$$
	\paragraph{} \textit{Proof: } Firstly, the proof that a unique solution exists for the system above is due to Theorem \ref{theorem:unique_approx_DFT}. Then, let $\textbf{e}^{(n+1)} = \textbf{x}^{(n+1)}(\tau) - \textbf{x}(\tau)$. Substituting the terms on the RHS, we obtain  $\textbf{e}^{(n+1)} = \bm{\mathcal{G}}_m(\tau)\textbf{e}^{(n)}$. The absolute value of the error gives $|\textbf{e}^{(n+1)}| \leq  |\bm{\mathcal{G}}_m(\tau)||\textbf{e}^{(n)}|$. Since the matrix $|\bm{\mathcal{G}}_m(\tau)|$ is \textit{substochastic}, that is each row sums up to at most one, with at least one row summing up to strictly less than one, then \cite[Theorem 6.2.27]{horn_johnson_1985} proves this is a necessary condition for $|\bm{\mathcal{G}}_m(\tau)|$ to be convergent (to the zero matrix). Since $|\textbf{e}^{(n)}| \leq  |\bm{\mathcal{G}}_m(\tau)|^{(n)}|\textbf{e}^{(0)}|$, then $\lim_{n \rightarrow \infty }|\textbf{e}^{(n)}| = 0$. \qed
\end{theorem}

\section{The Jacobi and Gauss-Seidel methods} We proceed to derive analogues of other common iterative methods - the Jacobi and Gauss-Seidel methods. These algorithms belong to a larger family of iterative algorithms known as \textit{coordinate descent} methods.

\subsection{The Jacobi method} The \textit{Jacobi method} can be derived from the original equations \eqref{eqn:lineq_matrix_trans_pdf} as follows. For each $s \in S_?$ we have
\begin{flalign}
\label{eqn:jacobi_or_gauss_seidel_form}
&f_s = \sum_{t \in S_?}(A_{s,t}G_{s,t} *  f_t) + h_s \notag\\
&f_s - (A_{s,s}G_{s,s} * f_s) = \sum_{t \in S_?  \slash s}(A_{s,t}G_{s,t}*  f_t) + h_s \notag\\
&f_s * (\delta - A_{s,s}G_{s,s}) = \sum_{t \in S_?  \slash s}(A_{s,t}G_{s,t} *  f_t) + h_s \notag\\
&f_s =  \sum_{t \in S_?  \slash s}((A_{s,t}G_{s,t} \deconv (\delta - A_{s,s}G_{s,s})) *  f_t)  + (h_s \deconv (\delta - A_{s,s}G_{s,s}))
\end{flalign}
Then the Jacobi method is just the above but in iterative form, i.e.
\begin{flalign}
\label{eqn:jacobi_infinite_form}
f_s^{(n+1)} =  \sum_{t \in S_?  \slash s}((A_{s,t}G_{s,t} \deconv (\delta - A_{s,s}G_{s,s})) * f_t^{(n)}) + (h_s \deconv (\delta - A_{s,s}G_{s,s}))
\end{flalign}
or alternatively, if we are solving for the property $\mrmpdf$, for $r= 0,1,\cdots,k-1$, then we have the practical form
\begin{flalign*}
f_s^{(n+1)} =  \sum_{t \in S_?  \slash s}\mathtt{conv}_k(\mathtt{deconv}_k(A_{s,t}G_{s,t}, \delta - A_{s,s}G_{s,s}),\  f_t^{(n)})  + \mathtt{deconv}_k(h_s,\ \delta - A_{s,s}G_{s,s})    
\end{flalign*}   
where again we have $\lim_{n \rightarrow \infty}\textbf{f}^{(n+1)} $ converging to the solution for \eqref{eqn:lineq_matrix_trans_pdf}. The proof for this is presented later. In the derivation above, it is helpful to know that $*$ and $\pentagon$ behave like multiplication and division respectively.
\paragraph{} The algorithm for the Jacobi method is the same as the power method, except that the update rule above is used instead. Before analysing the time complexity of the algorithm, let us make the following improvements. Let us define a hypervector $\bm{\kappa}$ where $\kappa_s \triangleq \mathtt{deconv}_k(h_s,\ \delta - A_{s,s}G_{s,s})$, and a hypermatrix $\mathbf{H}$ where $H_{s,t} \triangleq \mathtt{deconv}_k(A_{s,t}G_{s,t}, \delta - A_{s,s}G_{s,s})$, but $H_{s,s} = \textbf{0}$, for all $s,t \in S_?$. Then the equations simplify to
\begin{flalign}
\label{eqn:jacobi_method_pdf}
f_s^{(n+1)} = \sum_{t \in S_?  \slash s}\mathtt{conv}_k(H_{s,t},\  f_t^{(n)}) + \kappa_s
\end{flalign}
which gives us the same equational form as the power method, of which we can simplify again by defining a hypermatrix $\bm{\mathcal{H}}$ where $\mathcal{H}_{i,j}[\tau] \triangleq \mathcal{D}\{\widetilde{H}_{i,j}\}[\tau]$, and therefore the Jacobi method in matrix notation reduces to
$$
\textbf{f}^{(n+1)} =  \mathcal{D}^{-1}\{\bm{\mathcal{H}}(\mathcal{D}\{\widetilde{\textbf{f}}^{(n)}\})\} + \bm{\kappa}
$$

\subsubsection{Algorithmic complexity}
\paragraph{} The complete time complexity is almost identical to the power method. In fact, per iteration the time complexity is still the same; asymptotically in the worst-case we still have $O(|S_?|(2k-1)log_2(2k-1)) + O(|S_?|^2(2k-1))$. However, there are now additional costs to transform the original equations into the form above. Firstly, $\bm{\kappa}$ has complexity roughly $O(|S_?|klog_2k)$, \textbf{H} has complexity roughly $O(|S_?|^2klog_2k)$, and finally $\bm{\mathcal{H}}$ has complexity $O(|S_?|^2(2k-1)log_2(2k-1))$.

\paragraph{} If no fast deconvolution library is available, and one has to repeatedly deconvolve multiple terms by some constant divisor, then one can make use of the equality: $$\mathtt{deconv}_k(a,b) = \mathtt{conv}_k(a,\mathtt{deconv}_k(\delta,b)) $$    
where $\mathtt{deconv}_k(\delta,b)$ can be saved and used when required. Thus deconvolution is only needed $|S_?|$ times at most, once for each row.

\subsection{The Gauss-Seidel method}
\label{subsec:Gauss_Seidel_method}
The Gauss-Seidel method is similar to the Jacobi method. Firstly we change the indices of our hypervector \textbf{f} from $f_s[r], f_t[r], f_u[r], \cdots$ to an ordered set of indices $f_0[r], f_1[r], f_2[r],\cdots,f_{|S_?|-1}$. Then for example, we can write the Jacobi method as
$$
f_i^{(n+1)} = \sum_{j=0, j \neq i}^{|S_?|-1}\mathtt{conv}_k(H_{i,j},\  f_j^{(n)}) + \kappa_i
$$
\paragraph{} The Gauss-Seidel method is simply
\begin{flalign}
\label{eqn:gs_practical_form}
f_i^{(n+1)} = \sum_{j= i+1}^{|S?|-1}\mathtt{conv}_k(H_{i,j},\  f_j^{(n)})
+\sum_{j=1}^{i}\mathtt{conv}_k(H_{i,j},\  f_j^{(n+1)}) + \kappa_i
\end{flalign}
and notice that to compute $f_i^{(n+1)}$, we require terms $f_j^{(n+1)}$ for all $j \leq i$. Hence, a sequential computation is required, i.e. we compute $f_i^{(n+1)}$ in the order $i= 0,1,\cdots,|S_?|-1$.

\paragraph{} To improve time complexity, let $\bm{\mathcal{H}}$ be defined as previously. Then, we can partition the hypermatrix into the strictly lower and strictly upper hypermatrix via $\bm{\mathcal{H}}[\tau] = \bm{\mathcal{H}}_{L}[\tau]  + \bm{\mathcal{H}}_{U}[\tau]$. Note that the diagonal entries $({\mathcal{H}_{s,s}})_{s \in S_?}$ are all zero.

\paragraph{} Then, the Gauss-Seidel method can be written in matrix form as
$$
\textbf{f}^{(n+1)} =  \mathcal{D}^{-1}\{\bm{\mathcal{H}}_{U}(\mathcal{D}\{\widetilde{\textbf{f}}^{(n)}\})+ \bm{\mathcal{H}}_{L}(\mathcal{D}\{\widetilde{\textbf{f}}^{(n+1)}\}) \} + \bm{\kappa}
$$

\subsubsection{Algorithmic complexity}
The time complexity is the same as the Jacobi method. The difference however is that the solution vector of the Jacobi method can be computed in parallel, i.e. computing $f_0[r],f_1[r],\cdots,f_{|S_?|-1}[r]$ in parallel, but the naive Gauss-Seidel method requires a sequential computation. However, the Gauss-Seidel method offers a better convergence rate (see Section \ref{sec:convergence_rate}).

\subsubsection{Approximate variants}
\label{subsec:approx_jacobi_gauss}
Just like the LU approximation method (Sec. \ref{sec:approx_direct}) and the approximate power method (Sec. \ref{sec:approx_power}), we can derive an approximate Jacobi and Gauss-Seidel methods. We have chosen not to elaborate on them however. The difference in algorithm is mainly two: 1) The inverse DFT is not applied to the solution hypervector after each iteration, but rather only after convergence is met. 2) A parameter $n$ (or $m$) is introduced which denotes the zero-padding size that, if increased, reduces the effects of time-aliasing, as discussed previously.

\subsection{Convergence of the Jacobi \& Gauss-Seidel methods}
The Jacobi method in the Fourier domain (i.e. taking the Fourier transform of \eqref{eqn:jacobi_infinite_form}) is
\begin{flalign}
\label{eqn:jacobi_DFT_form}
\textbf{x}^{(n+1)}(\tau) =  \bm{\mathcal{Q}}(\tau)\textbf{x}^{(n)}(\tau) + \textbf{z}(\tau)
\end{flalign}
and the Gauss-Seidel method is
\begin{flalign}
\label{eqn:gs_DFT_form}
\textbf{x}^{(n+1)}(\tau) =  \bm{\mathcal{Q}}_L(\tau)\textbf{x}^{(n)}(\tau) + \bm{\mathcal{Q}}_U(\tau)\textbf{x}^{(n+1)}(\tau) +  \textbf{z}(\tau)
\end{flalign}  
where $z_s(\tau) \triangleq  \mathcal{D}\{h_s \deconv \ \delta - A_{s,s}G_{s,s}\}$, and $\mathcal{Q}_{s,t} \triangleq \mathcal{D}\{ A_{s,t}G_{s,t} \deconv \delta - A_{s,s}G_{s,s}\}$.

\begin{theorem}[Convergence of the Jacobi and Gauss-Seidel methods]
	\label{theorem:convergence_jacobi}  Both the Jacobi \eqref{eqn:gs_DFT_form} and Gauss-Seidel \eqref{eqn:jacobi_DFT_form} methods converge to the true solution of the system defined in  \eqref{eqn:lineqn_matrix_trans}, i.e. $\lim_{n \rightarrow \infty}\textbf{x}^{(n)} = \textbf{x}$.
	
	\paragraph{} \textit{Proof: }  From Proposition \ref{prop:wirr_diagonal_dom},  $(I - (A \circ \mathbf{C})(\tau))$ is weak and irreducibly diagonally dominant. Then, \cite{bagnara1995unified} proves that this condition is sufficient to ensure that the Jacobi and Gauss-Seidel methods above converge for each $\tau$.  The methods also converge to a unique fixed-point due to the matrix $(I - (A \circ \mathbf{C})(\tau))$ being non-singular \cite{bagnara1995unified}. By taking the inverse Fourier transform of the solution, we obtain $\lim_{n \rightarrow \infty }\textbf{f}^{(n)} = \textbf{f}$, with $\textbf{f}$ being the solution to \eqref{eqn:lineq_matrix_trans_pdf}, and $\textbf{f}^{(n)}$ is the $n^{th}$ iteration of the Jacobi or Gauss-Seidel method in the time domain, e.g. via equation \eqref{eqn:jacobi_infinite_form} for the Jacobi method.  \qed
\end{theorem}     

\section{Convergence rate analysis}
\label{sec:convergence_rate} The convergence rate of the exact power, Jacobi, and Gauss-Seidel methods for convolution systems can be studied in the Fourier domain. The system of equations to be solved is \eqref{eqn:eqn_matrix_FT_re_arranged}, i.e.
$$ (I - (A \circ \mathbf{C})(\tau))\textbf{x} = \textbf{d}$$
Let us denote $(I - (A \circ \mathbf{C})(\tau))$ as $\pmb{\mathscr{A}}(\tau)$.
Then, both the Jacobi \eqref{eqn:jacobi_DFT_form} and Gauss-Seidel  \eqref{eqn:gs_DFT_form} methods can be written in the form
$$    \textbf{x}^{(n+1)}(\tau) =\textbf{B}(\tau)\textbf{x}^{(n)}(\tau)  + \textbf{c}(\tau)  $$
\paragraph{} From \cite[p. 190]{dahlquist1974numerical} , we can decompose $\pmb{\mathscr{A}}(\tau)$ into a lower triangular, diagonal and upper triangular matrix, in the manner $\pmb{\mathscr{A}}(\tau) = \textbf{D}(\tau)(\textbf{L}(\tau) + I + \textbf{U}(\tau))$ where $\textbf{D}(\tau) = diag(\mathscr{A}_{i,i}(\tau))$. Then, for both methods we have $\textbf{c}(\tau) =\textbf{D}(\tau)^{-1}\textbf{d}(\tau)$, but for the Jacobi method
$\textbf{B}_J(\tau) = - (\textbf{L}(\tau)+ \textbf{U}(\tau))$
whilst for the Gauss-Seidel,
$\textbf{B}_{GS}(\tau) = -(I + \textbf{L}(\tau))^{-1}\textbf{U}(\tau)$.

\paragraph{} From \cite[p. 191]{dahlquist1974numerical}, the asymptotic convergence rate of both methods is
$$ R(\tau) = -log_{10}(\rho(\textbf{B}(\tau)))$$
where $\rho(\textbf{B}(\tau))$ is the spectral radius of $\textbf{B}(\tau)$. It is expected that the convergence rate of Gauss-Seidel is faster than the Jacobi since $x_j^{(n+1)}(\tau)$ can be computed using the updated values $x_i^{(n+1)}(\tau)$ for all $i < j$. Whereas for the Jacobi method, $x_j^{(n+1)}(\tau)$ will only use the values $x_i^{(n)}(\tau)$. However, for general problems, \cite[pg. 291]{phillips1996theory} states that the Gauss-Seidel method may diverge to some problems where the Jacobi method converges.

\paragraph{} Before detailing the convergence rate of the power method, it is helpful to rewrite the system of convolution of equations in a different form. From Theorem \ref{theorem:derivation_of_system}, we saw that
\begin{flalign*}
\phi_s(\tau)  = \sum_{t \in S}\textbf{P}(s,t)\phi_{rew(s,t)}(\tau)\phi_{t}(\tau)
\end{flalign*}
which equates to the system $\textbf{x}(\tau) = \textbf{A}(\tau)\textbf{x}(\tau) $ where $(A_{s,t})_{s,t \in S^2} =  (\textbf{P}(s,t)\phi_{rew(s,t)}(\tau))_{s,t \in S^2}$ and $\textbf{x}(\tau) = (\phi_s(\tau))_{s \in S}$. Then, the (normalization-free) power method \cite{zapreev2006safe} is written as
$$\textbf{x}^{(n+1)}(\tau) = \textbf{A}(\tau)\textbf{x}^{(n)}(\tau) $$
which is equivalent to the Fourier transform of \eqref{eqn:lineq_power_method_pdf}, and has a convergence rate
$$ R(\tau) = \frac{\lambda_2(\tau)}{\lambda_1(\tau)}$$
where $\lambda_1(\tau),\lambda_2(\tau)$ are the dominant and sub-dominant eigenvalues respectively, of the matrix $\textbf{A}(\tau)$. Baier \& Katoen \cite[p. 754]{principles_book} states that the convergence of the power method is often less efficient than the Jacobi and Gauss-Seidel methods. Hence, one may expect that the Gauss-Seidel is the fastest, followed by the Jacobi and then the power method.

\paragraph{} Note that the rates above are for each $\tau$. If we require the convergence rate for all $\tau$, then it is determined by the smallest rate, i.e. $min_{\tau}(R(\tau))$. However, if we assume that $R_{pow}(\tau) \leq R_{J}(\tau) \leq R_{GS}(\tau)$ for all $\tau$, then we know that the Gauss-Seidel method is fastest for all $\tau$. And therefore, in the time-domain, the Gauss-Seidel is expected to be fastest there too (for convolution systems).

\paragraph{} Note that the same results hold for the approximate variants of the iterative methods when applied to sMRMs with discrete-lattice reward. However, the discrete Fourier transform is to be used in Theorem \ref{theorem:convergence_jacobi} when deriving the system $\pmb{\mathscr{A}} \textbf{x} = \textbf{d}$.

\section{Experiments}
\label{sec:exp_MC_pmf}
\begin{figure}[H]
	\centering
	\includegraphics[width=1\linewidth]{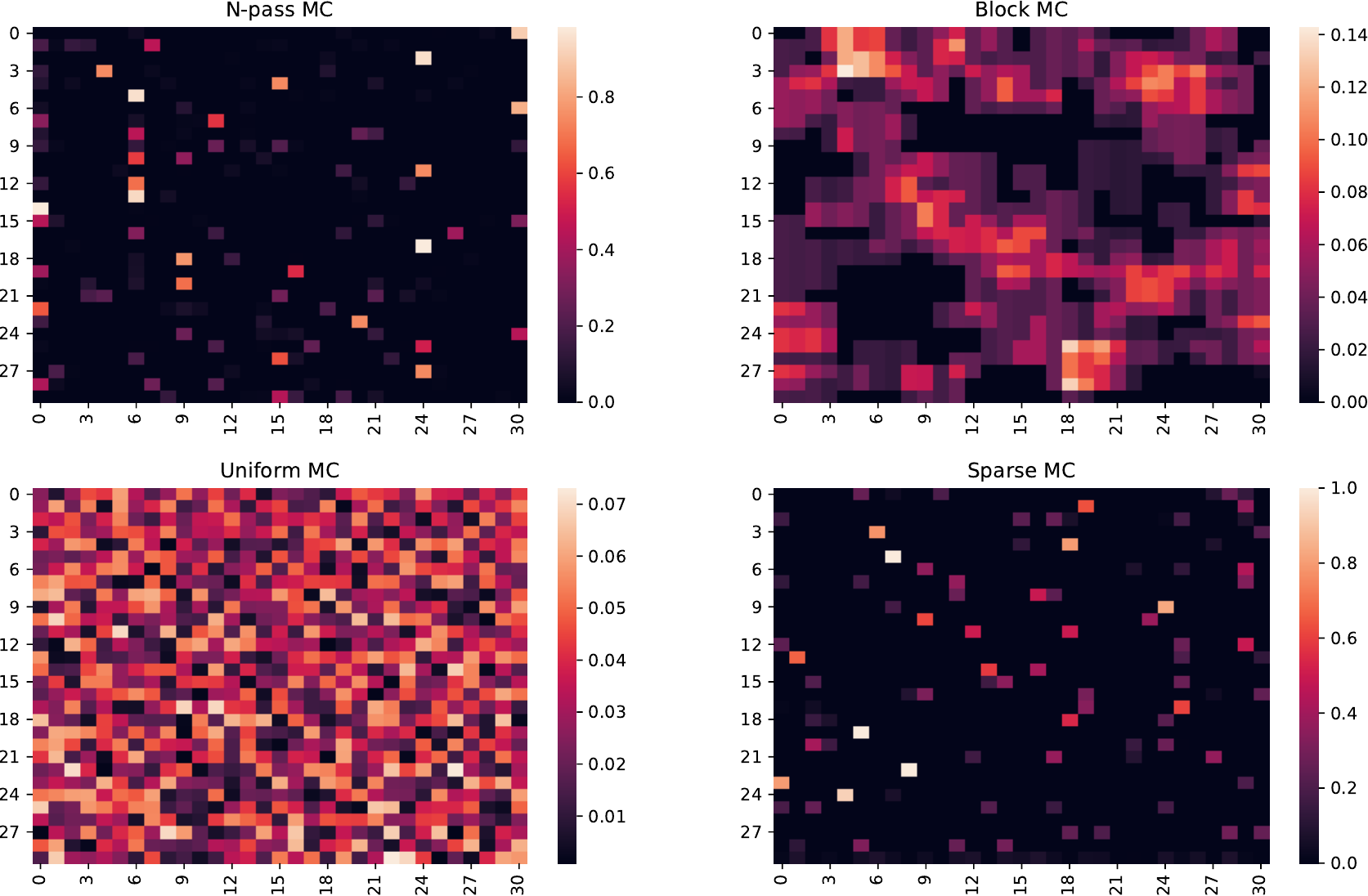}
	\caption{A  sample of each type of Markov chain used in the experiments, represented via a heatmap. Each rectangle above represents a probability matrix of a MC. The last column is used as the probability of entering the (set of) goal states $B$.}
	\label{fig:sampled_MCs}    
\end{figure}

\begin{figure}[H]
	\centering
	\includegraphics[width=1\linewidth]{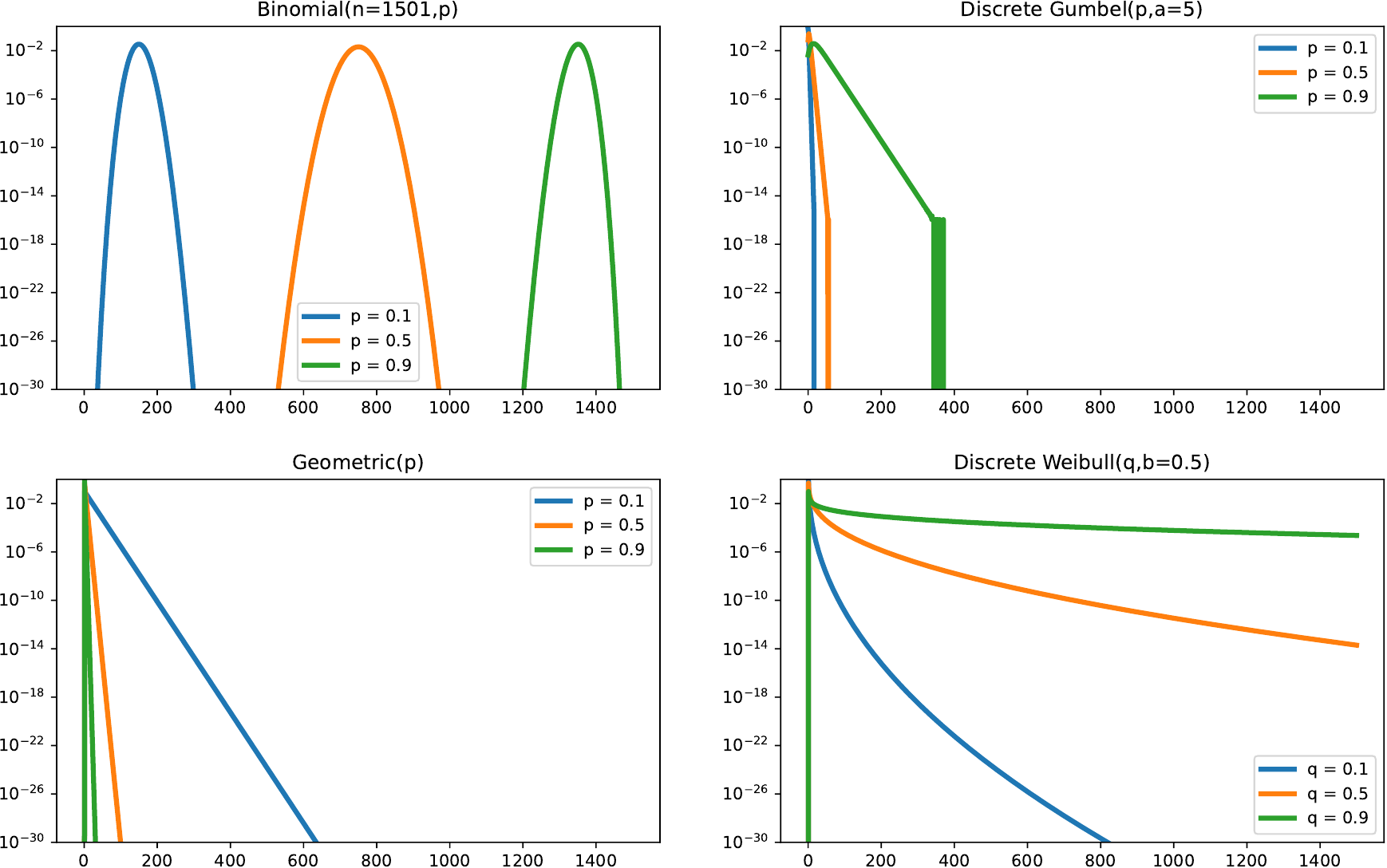}
	\caption{The different distributions used for the reward random variables. Each has one free parameter which we varied. The Weibull is used for testing sMRMs with long-tail distributions. The geometric is used for representing medium-tail distributions, whilst the Gumbel is used for light tailed distributions. The Binomial distribution is bounded in the interval $[0,N-1]$.}
	\label{fig:sampletails}
\end{figure}

We investigate here how different Markov chain structures and types of probability mass functions (pmfs) may affect the convergence rates of the power, Jacobi and Gauss-Seidel methods (with respect to sMRMs). We experiment with four types of Markov chains (Fig. \ref{fig:sampled_MCs}), and well as four types of pmfs (Fig. \ref{fig:sampletails}). We first construct sMRMs with different combinations of MC and pmf type, and then solve for the property $\mrmpdf$ over them with the aforementioned iterative methods. The main metric for determining convergence rate is the number of iterations taken to solve for the property. The time taken to solve is also measured, however the metric is not as useful due to the implementation of the iterative methods not being equally optimized. For each method, the recorded time is both the solving time, i.e. the time until convergence (solving phase), and the time required to build the relevant system of equations (preparation phase) (e.g. $\textbf{G}, \textbf{h}$). As for computer and software specs., the first two experiments (Sec. \ref{subsec:3exact} and Sec. \ref{subsec:powers}) are done on \textit{Computer 1}, and the remainder with \textit{Computer 2}. See Section \ref{sec:computer details} for computer and software details.

\paragraph{} The Markov chains are of four types: N-pass, Block, Uniform and Sparse. The Uniform model is used to simulate dense MCs, the Block for MCs with somewhat strongly connected components, the N-pass is a sparse model that is generated differently to Sparse. The probability matrix of each model is generated as follows. The Uniform and Sparse models are generated by creating matrices of which elements are sampled uniformly between zero and one. Then each row is normalized to ensure that each row sums to one forming the probability matrix \textbf{P}. For the Sparse model, roughly $10\%$ of \textbf{P} will have values. The Block model is created by repeatedly adding one to a random block of indices to a matrix with values all initially set to zero. Then, the resulting matrix is row-normalized. For our experiments we used 200 blocks to generate these matrices. Lastly, the $N$-pass model is a model that involves no normalization. Firstly, a set of states $Q$ are sampled of which their probabilities of reaching the goal state are generated uniformly in $[0,1)$. Then, for $N$ steps, at each step, the probability for each state entering a random state in $Q$ is generated, however in such a way that normalization is not required. For implementations of these algorithms, please refer to Appendix \ref{app_chap:mc_simulation}.

\paragraph{} The pmfs are of four types: Binomial, (discrete) Gumbel, Geometric, and discrete Weibull. Each pmf has different characteristics that may affect the convergence rate. In particular, we are interested in whether the tail properties of these pmfs may be a significant factor. The Weibull is considered to have the heaviest tail, followed by the Geometric, and lastly the Gumbel. The Binomial has a tail which does not extend infinitely, i.e. it's distribution has a bounded support. Each pmf has a free single parameter that will be varied for the experiments. Their definitions are as follows:

\begin{enumerate}
	\item $Binomial(n=1501,p)(t) = \binom{1501}{t}p^{t}(1-p)^{n-t}$
	\item $Gumbel(p,a=5)(t) = e^{-5p^{t+1}} - e^{-5p^t}$
	\item $Geometric(p)(t) = (1-p_1)^{t-1}p_1$
	\item $Weibull(q,b=0.5)(t) = f(t;q,0.5) = \begin{cases}
	q^{(t-1)^{0.5}} - q^{t^{0.5}} & t \in [1,2,3,\cdots] \\
	0 & otherwise
	\end{cases}$
\end{enumerate}
\paragraph{} The free parameter is $p$ for the first three pmf above, and $q$ for the Weibull. The fixed parameters above are chosen in a way to allow us to increase the tail strengths for the pmfs with tails, and to move the peak of the binomial distribution across the interval $[0,1,\cdots,1501]$. The discrete Weibull comes from \cite{warr2014numerical}, and the discrete Gumbel is from \cite{chakraborty2014discrete}.

\paragraph{Experiment set-up} The experiments below consist of randomly generating sMRMs, to solve $\mrmpdf$ via the system of equations \eqref{eqn:lineq_matrix_trans_pdf}. The sMRM is generated via a random sample of a selected MC type and pmf type. All the reward random variables for a given sMRM are of the same type. Their free parameters are sampled uniformly within a selected range shown later. For all experiments below, unless otherwise mentioned, we set $S_? = 30$, and we solve for $\mrmpdf$, for $r=0,1,\cdots,N-1$ with $N = 1501$. For each iterative method, we used the following termination criteria $$max_{s,r}|\textbf{f}^{(n+1)} - \textbf{f}^{(n)}| \leq \epsilon$$ where $\epsilon = 1e-7$. Additionally, a max iteration of 2000 was used to terminate the method if convergence was too slow.      

\subsection{Comparison of the exact iterative methods}
\label{subsec:3exact} 
\begin{figure}[h]
	\centerline{
		\includegraphics[width=1.25\linewidth]{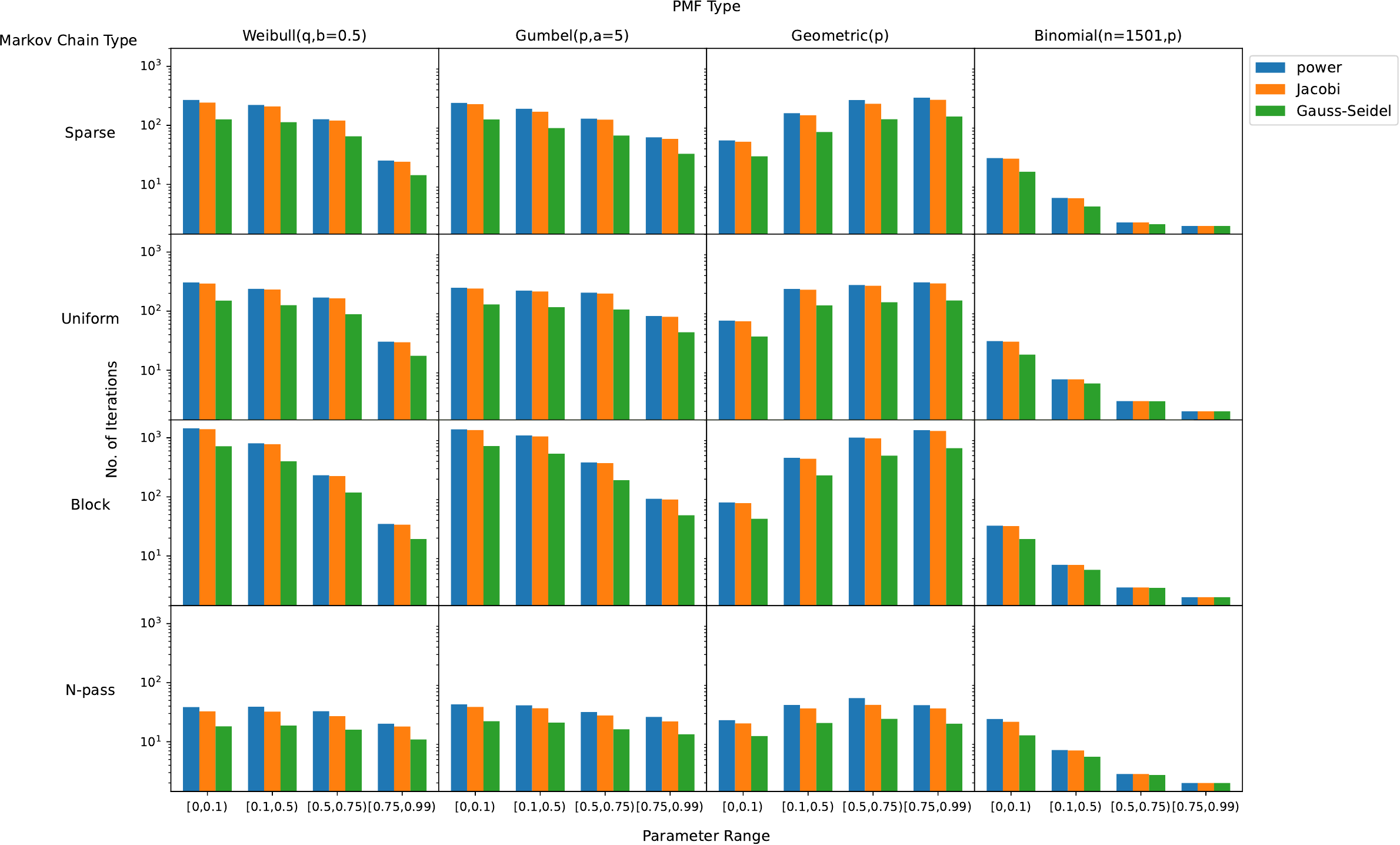}}
	\caption{Average no. of iterations with each method for a given sMRM type. Each triplet of bars (blue, orange, green) is an average of 50 unique samples. The free parameter ranges are labelled on the x-axis. The free parameter of each pmf of an sMRM is sampled within that range.}
	\label{fig:noiter3exact}
\end{figure}

\begin{figure}[h]
	\centerline{
		\includegraphics[width=1.25\linewidth]{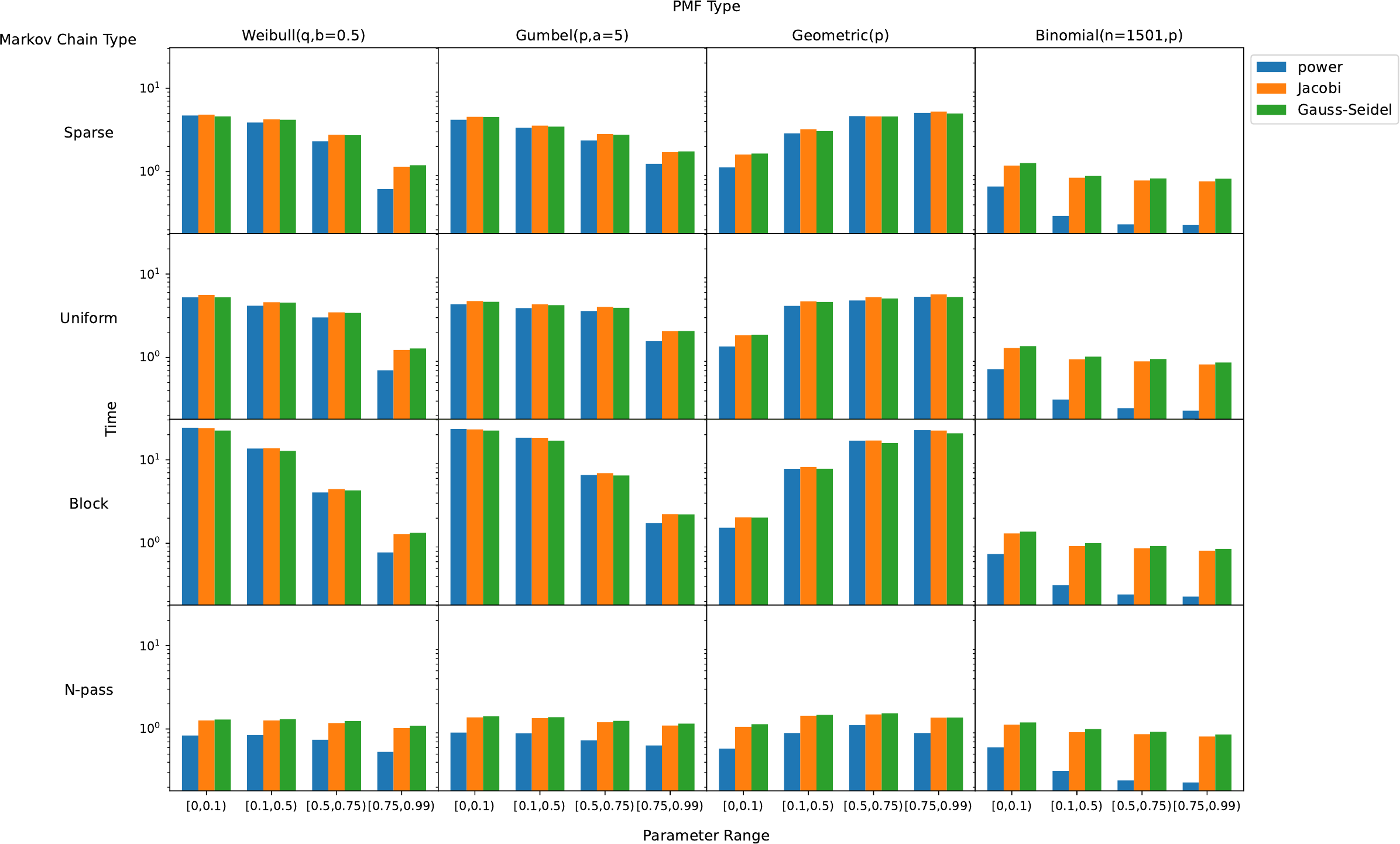}}
	\caption{Average time (secs.) taken for a given sMRM type. Whilst the Gauss-Seidel uses the fewest iteration, it's solving speed is slower as a result of how the method was implemented.}
	\label{fig:time3exact}
\end{figure}

\begin{figure}[h]
	\centerline{
		\includegraphics[width=1.25\linewidth]{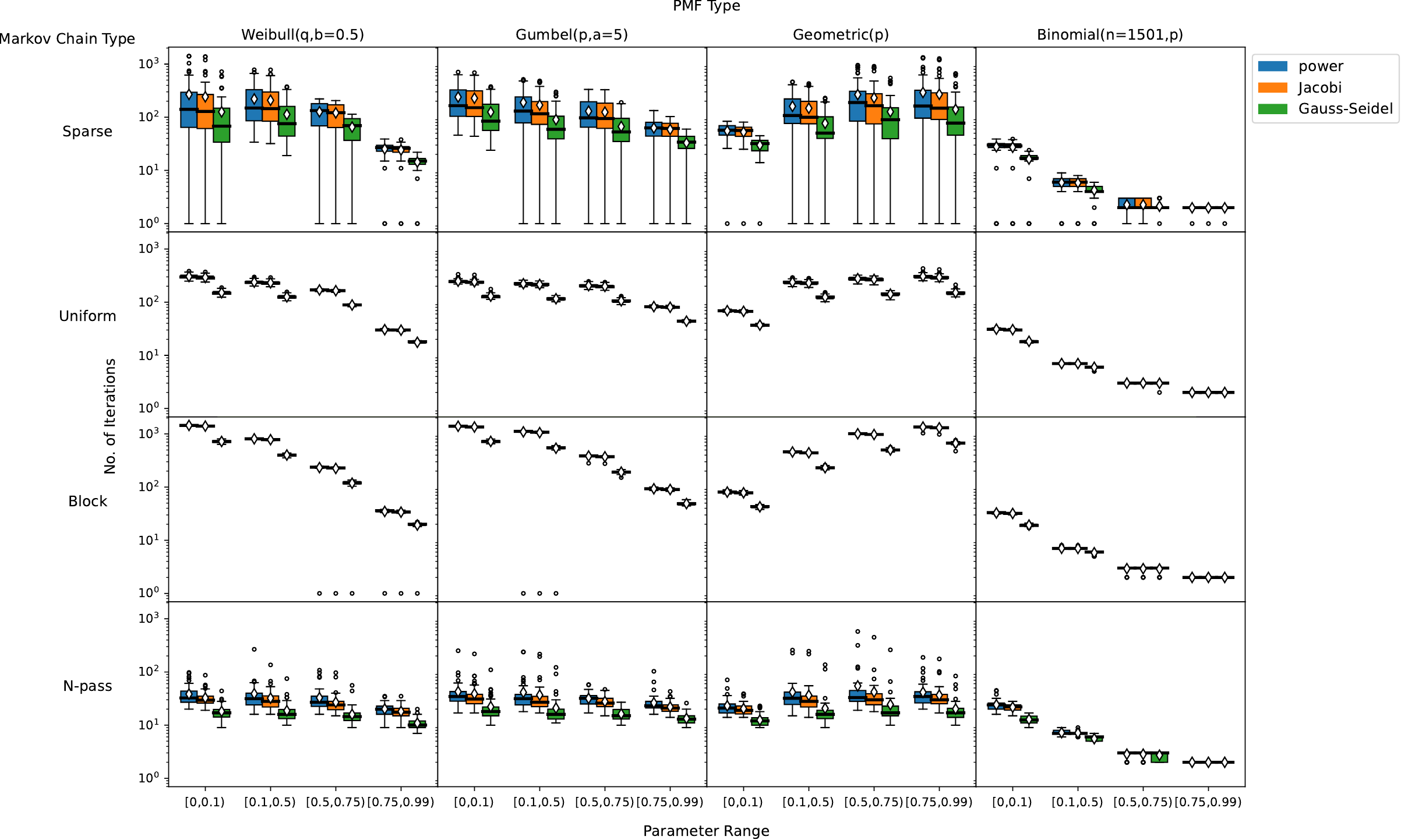}}
	\caption{Box-plot results for the number of iterations are shown for each (MC,pmf,pmf param.) combination, over 50 samples. The white diamond represents the mean, whilst the circle dots are outliers. The median is the black solid line between the interquartile range. }
	\label{fig:noiter3exactboxes}
\end{figure}

The average number of iterations taken by each of the exact iterative methods is presented in Fig. \ref{fig:noiter3exact}, with box-plots shown in Fig. \ref{fig:noiter3exactboxes}. The average time taken (secs.) is shown in Fig. \ref{fig:time3exact}. For each combination of (MC type, pmf type, pmf param.), we obtain 50 samples. In effect, each triplet of bars in Fig. \ref{fig:noiter3exact} is the average result of 50 unique samples. Hence, each cell is the result of 200 samples, with a total of 3200 samples for the whole plot.

\paragraph{} For both the Weibull and Gumbel pmf, a consistent trend is that as the parameter range increases, i.e. as the heaviness of their tails increases, then the iteration count decreases. However, for the Geometric pmf, we have the inverse effect. A reason could be because the Weibull and Gumbel distribution reduces its `peakiness' as the parameter range increases, whilst the Geometric pmf increases in 'peakiness'. Hence, a consistent hypothesis is that the peakiness of the distribution near zero leads to longer convergence rates. For the binomial pmf, the further away from zero its sampled mean, the lower the average number of iteration. This appears consistent with the previous hypothesis. However, an alternative hypothesis is that the further the mean is centred nearer $N$, the more likely the mean of the cumulated reward pmf is centred outside the interval $[0,1,\cdots,N]$. Thus, the values in this interval could be really small, meaning that convergence can be achieved faster.

\paragraph{} As for the effects of MC types, the most prominent result from the plot is that the $N$-pass MC took the fewest (average) iterations. The Block MC took the most, whilst the Sparse and Uniform MC had similar averages. From the box-plots (Fig. \ref{fig:noiter3exactboxes}),  the Sparse MC yields results with the largest spread, followed by $N$-pass. Some Sparse problems took longer than problems with the Uniform MC type, and others took near zero iterations. This latter case can happen if probability matrices that have been generated do not have any states reaching $B$. We found that the iterative methods applied to such MCs would terminate after one iteration, which is a correct result.

\paragraph{} Lastly, we find that the convergence rate (in terms of number of iterations) ordered by fastest to slowest is the Gauss-Seidel, followed by the Jacobi and then the power method. As for the times taken, then we find the power method being optimal due to implementation reasons.

\subsection{Comparison of the exact and approximate power methods}
\label{subsec:powers} 
\begin{figure}[h]
	\centerline{
		\includegraphics[width=1.25\linewidth]{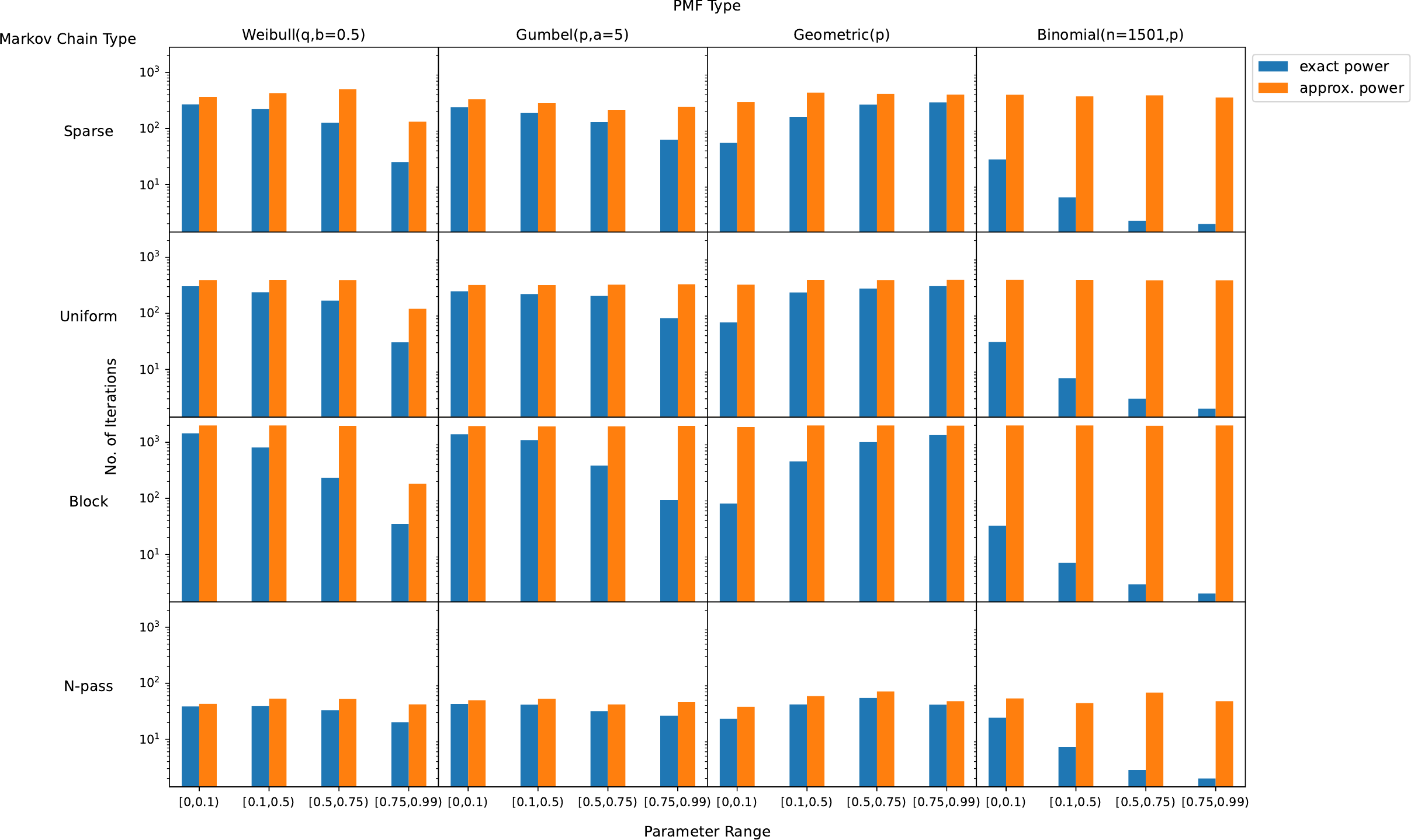}}
	\caption{Average no. of iterations for specific sMRMs. The no. of iterations themselves do not paint a full picture for solving times. Although, it is almost obvious that for the Binomial column, the exact method is faster.}
	\label{fig:noiterpowers}
\end{figure}

\begin{figure}[h]
	\centerline{
		\includegraphics[width=1.25\linewidth]{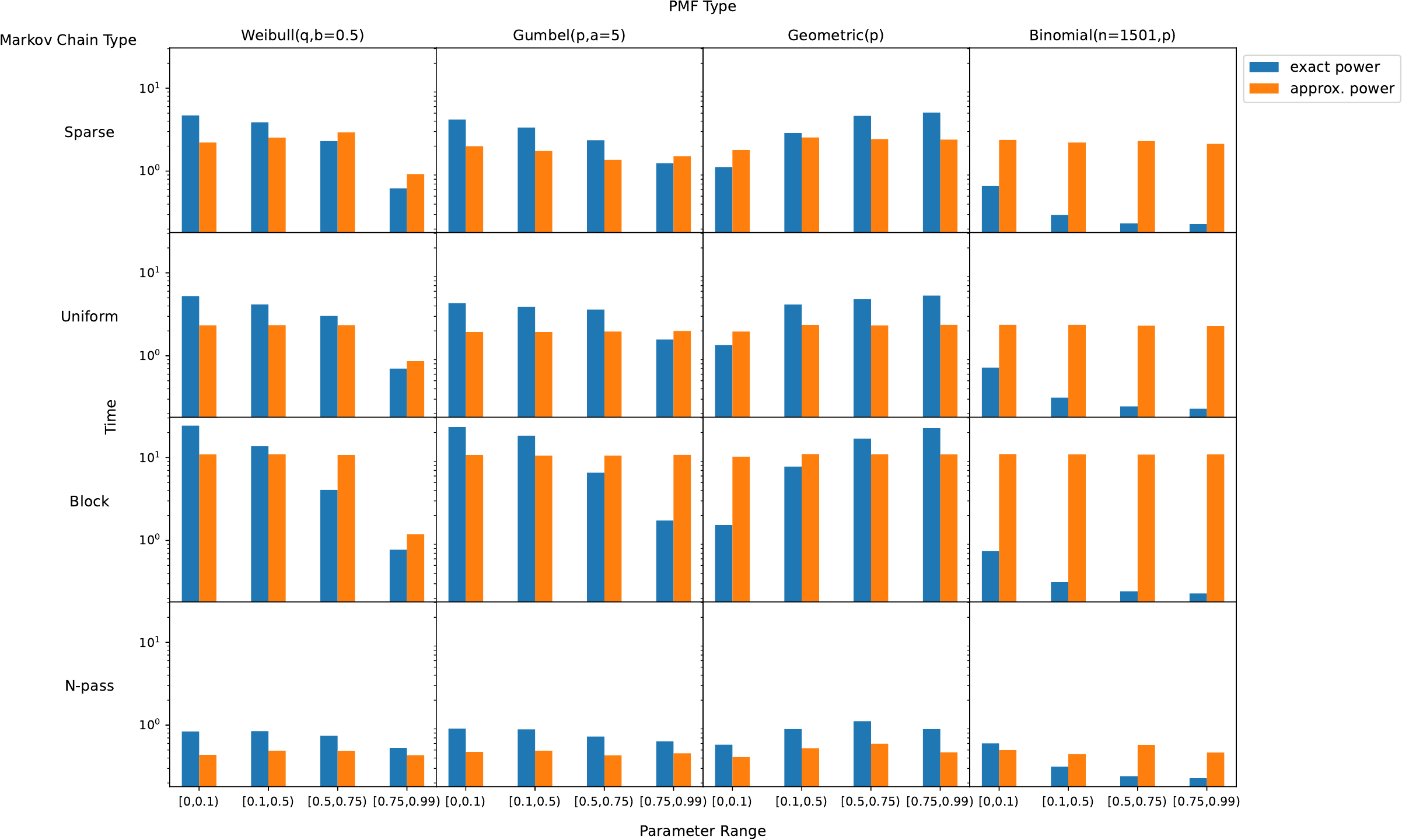}}
	\caption{Average time (secs.) taken for specific sMRMs.}
	\label{fig:timepowers}
\end{figure}

\begin{figure}[h]
	\centerline{
		\includegraphics[width=1.25\linewidth]{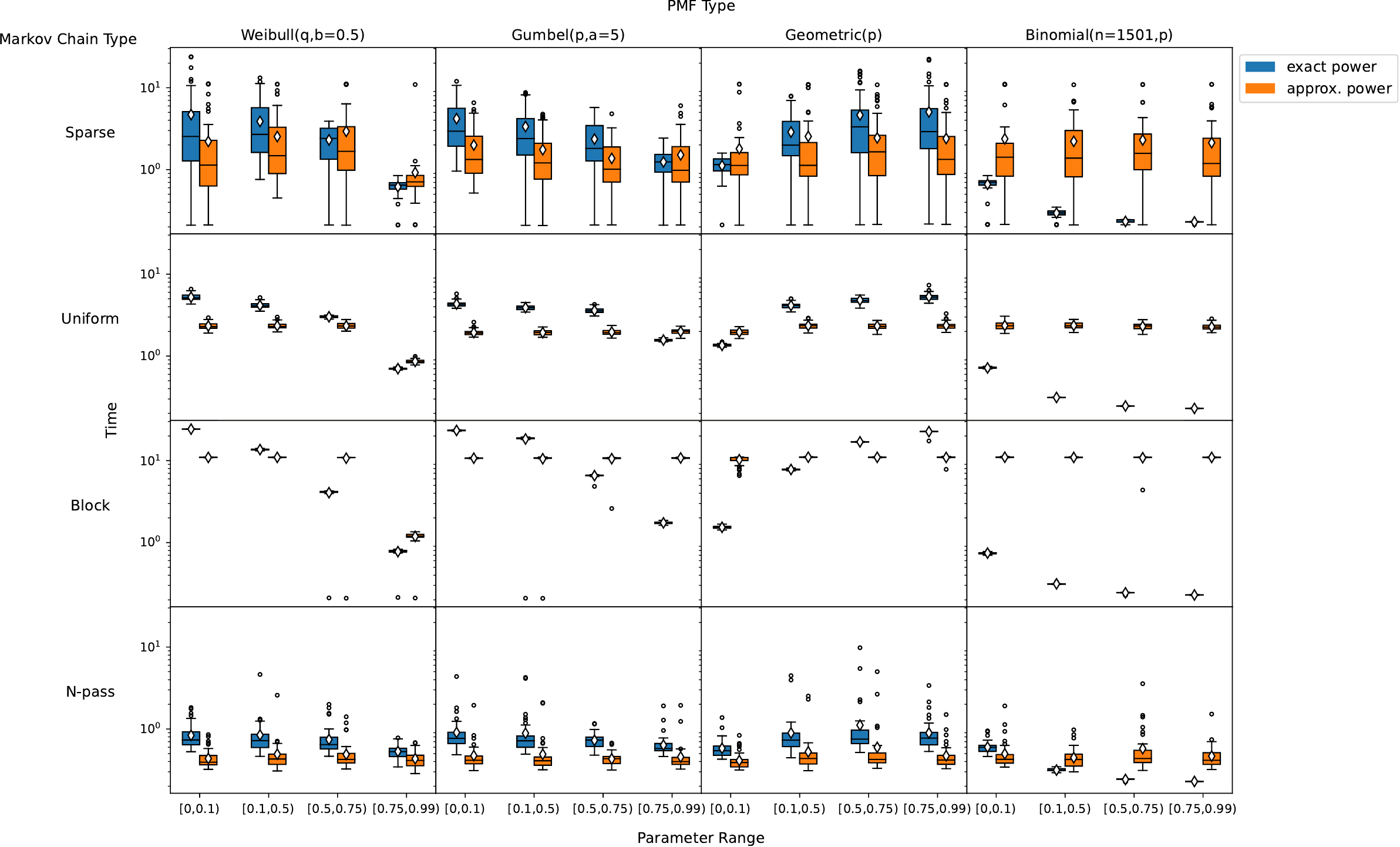}}
	\caption{Box-plot results for solving time. Definitions of symbols are as in Fig. \ref{fig:noiter3exactboxes}.}
	\label{fig:timepowersboxes}
\end{figure}

We now compare the exact and approximate method power using the same experiment as above. Hence, the samples of this experiment are identical to the previous one. The average number of iterations taken by each method is presented in Fig. \ref{fig:noiterpowers} with the average time taken shown in Fig. \ref{fig:timepowers}, and the box-plot for the time taken shown in Fig. \ref{fig:timepowersboxes}.

\paragraph{} We set the approximate method's padding length to $N-1$, which makes it identical to the exact power method in terms of space requirements. Therefore the main difference between the two methods is that the exact power method requires the DFT and the inverse DFT to be applied in each iteration, whereas the approximate method does not. For the exact power method, the same convergence threshold is used as before. For the approximate power method, we used the termination criteria defined in \eqref{eqn:approx_method_power_convergence_criteria} with the same threshold -  $\epsilon = 1e$-7. A max iteration of 2000 was also used to terminate the method if convergence was too slow, just like the power method.    

\paragraph{} A prominent trend is that the approximate power method requires more iterations on average than the exact method for any cell in Fig. \ref{fig:noiterpowers}. We find that the average no. of iteration for the approx. method is near equal for the first three rows. An exception is in the Weibull column for the parameter range $[0.75,0.99)$ where it decreases. This trend also highlights the Binomial column, where the average iteration count of the exact method reduces as the parameter range increases, but the average for the approximate method remains near constant. As for the time taken to solve (see Fig. \ref{fig:timepowers}), then for the Binomial column, it appears that the exact method is the optimal choice, whereas for the other columns, the results are mixed, although there appears to be some discernible albeit complex rule to determine roughly which method leads to faster times.

\subsection{Scalability of selected methods} We study here how solving speeds are affected by increasing $|S_?|$ and $N$. We test just three methods: exact power, approx. power and approx. LU methods. We also investigate the approximate maximum (abs.) error between the approximate methods and the exact power method.

\subsubsection{Varying $|S_?|$}

\begin{figure}[h]
	\centerline{
		\includegraphics[width=1.25\linewidth]{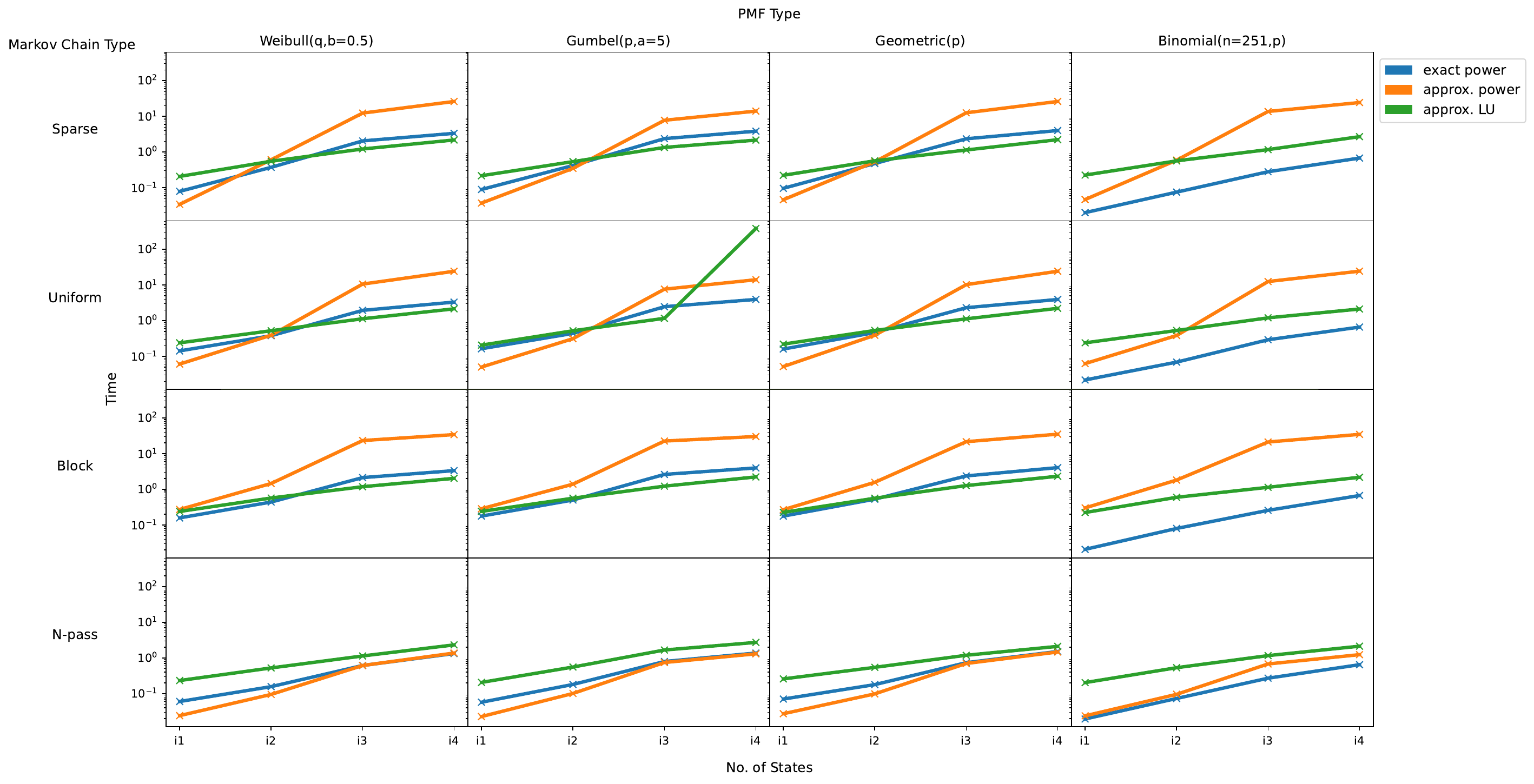}}
	\caption{Average time taken (secs.). For each cell, each interval involves computing 50 unique samples.}
	\label{fig:time3scale}
\end{figure}

\begin{figure}[h]
	\centerline{
		\includegraphics[width=1.25\linewidth]{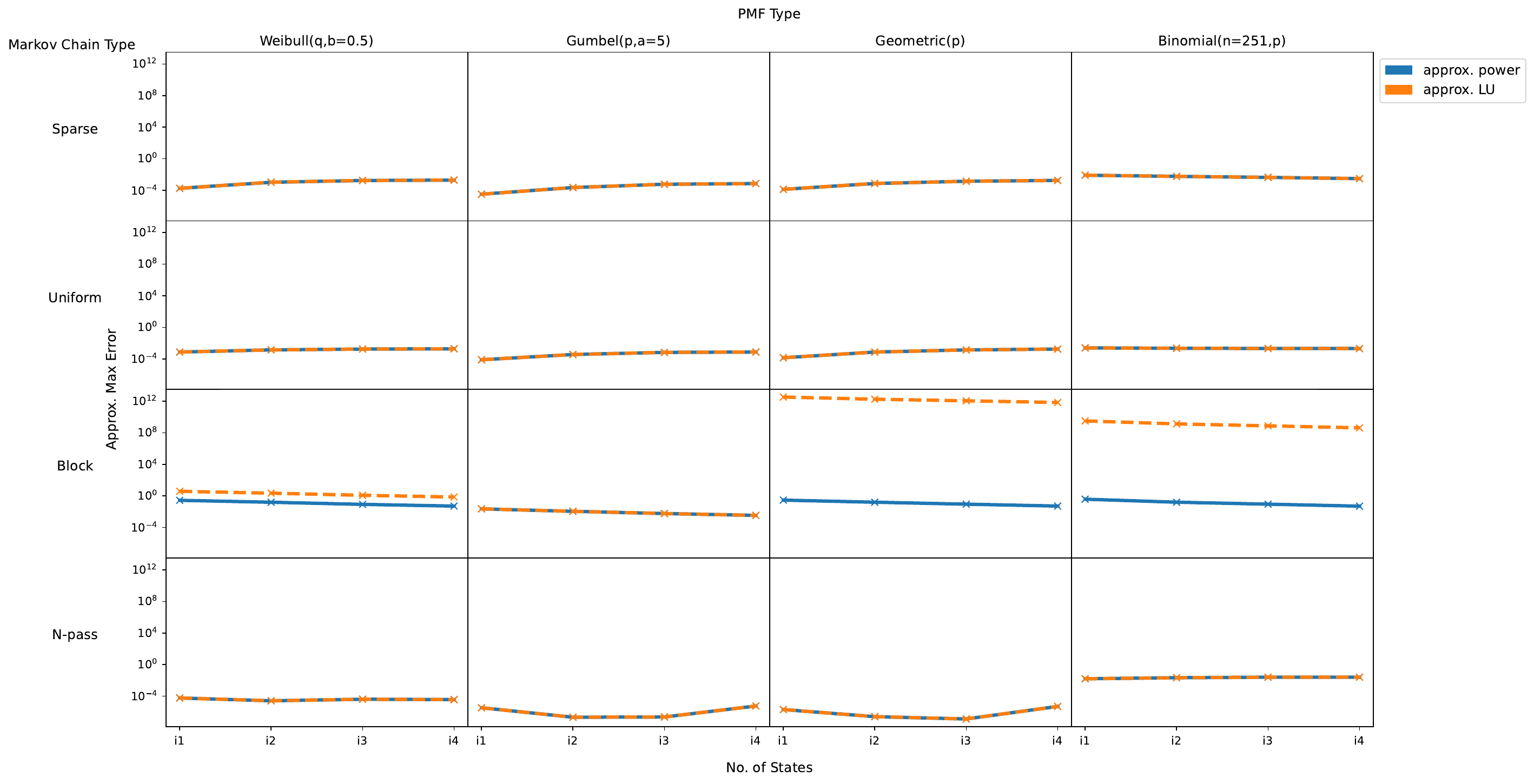}}
	\caption{Average max abs. error. Notice the huge gap in error for the Block column. }
	\label{fig:err2scale}
\end{figure}

\begin{figure}[H]
	\centerline{
		\includegraphics[width=1.25\linewidth]{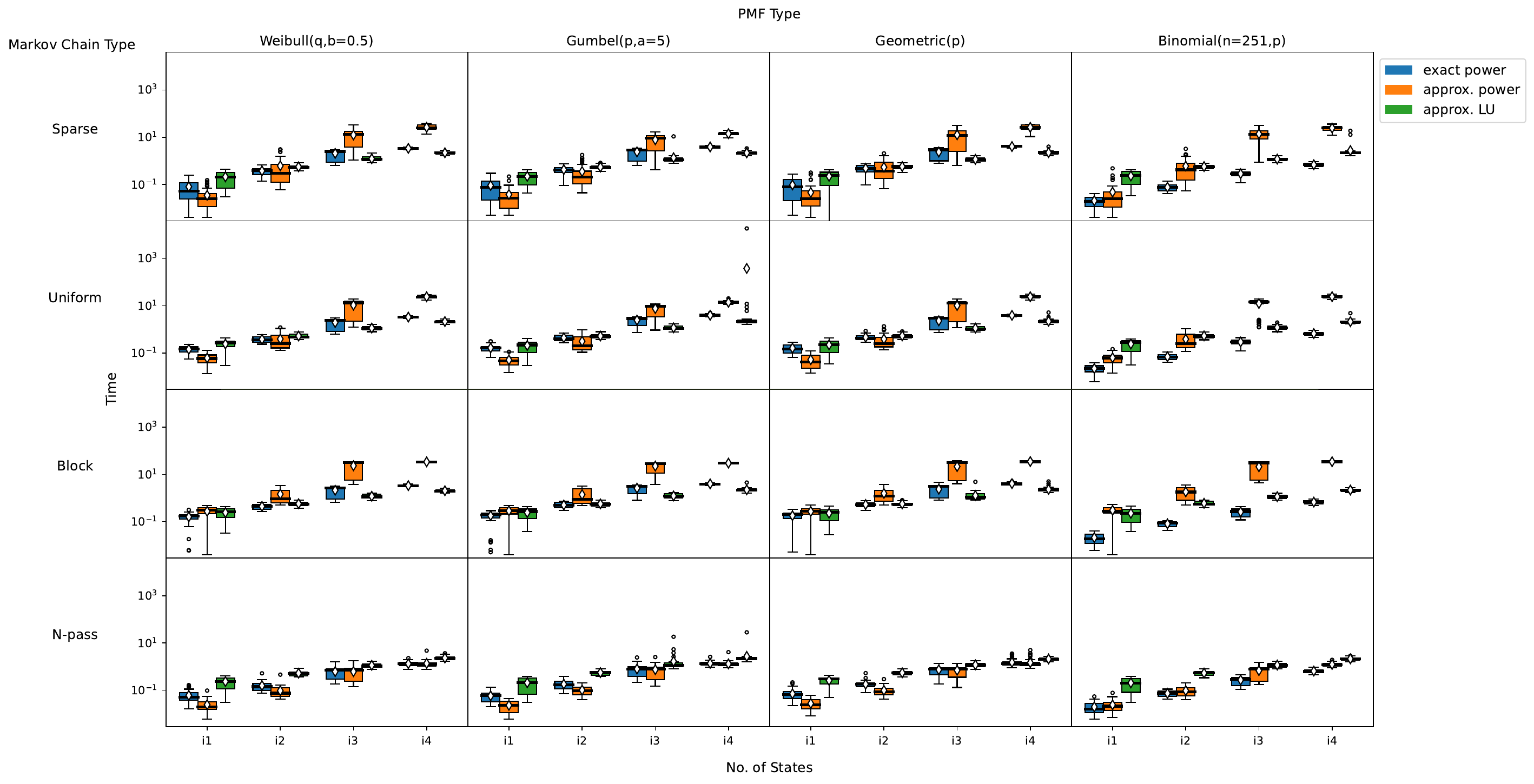}}
	\caption{Box-plot results for time taken.  Definitions of symbols are as in Fig. \ref{fig:noiter3exactboxes}.}
	\label{fig:time3scaleboxes}
\end{figure}

\begin{figure}[H]
	\centerline{
		\includegraphics[width=1.25\linewidth]{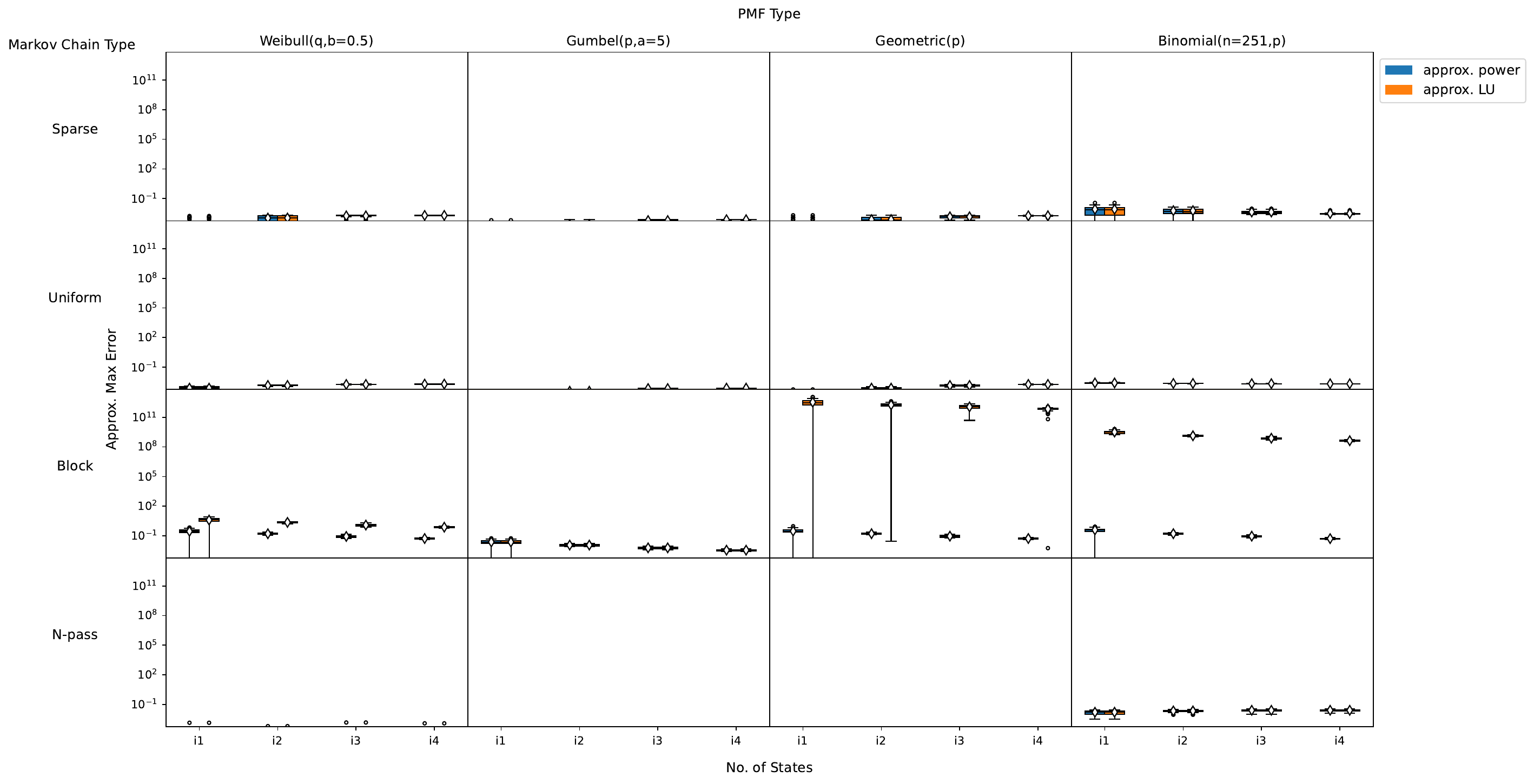}}
	\caption{Box-plot results for error of approximations.  Definitions of symbols are as in Fig. \ref{fig:noiter3exactboxes}. The results do not show outliers being the cause of the discrepancy between the errors in the Block row.}
	\label{fig:err2approxboxes}
\end{figure}

\paragraph{} The experiment setup is as before. For each combination of MC type and pmf type and parameter value, we generate 50 samples. Our first experiment involves generating samples for different values of $|S_?|$ whilst keeping $N$ fixed to $251$. For the Binomial pmf, its fixed parameter $n$ was set to $251$. The free-parameters of all pmf types are sampled uniformly in the range $[0.3,0.6]$. The average time taken (secs.) is shown in Fig. \ref{fig:time3scale} with the box-plots for time shown in Fig. \ref{fig:time3scaleboxes}. The average maximum error for the approximate methods are found in Fig. \ref{fig:err2scale}, with box-plots shown in Fig. \ref{fig:err2approxboxes}.

\paragraph{} In each plot,the x-axis labels (i1, i2, i3, i4) correspond to the intervals $[10,30),[30,50),[50,100)$ and $[100,150)$ respectively.  For each interval, for each sMRM sample generated, $|S_?|$ was randomly and uniformly selected between that interval. Note that for all results, for the Geometric column and the Sparse row, the interval (i1) uses only 49 samples since a sample resulted in failure when using the LU approx. method. This is due to a slice (matrix) of the hypermatrix used being singular, which can happen.

\paragraph{} In terms of average time, for the Binomial column, the exact power method is consistently best. For the $N$-pass row, the approx. LU method is consistently worst. For the remaining columns and rows, the approx. power method becomes worse than the exact power method as the number of states increases. As for the surprising upwards kink shown by the approx. LU method in the cell for the Gumbel pmf and Uniform MC, then the box-plot of Fig. \ref{fig:time3scaleboxes} shows that this is due to an outlier.

\paragraph{} Studying the error, we should find that the results are equal for each cell. However, surprisingly, the approx. LU method returns extremely wrong results in the cells of the Block MC row. A hypothesis may be that Block MCs are more susceptible to severe round-off errors. Another observation is that the Binomial pmf is hardest to obtain good accuracy, and that sMRMs with the $N$-pass MC is easiest to obtain good accuracy.


\subsubsection{Varying N}

\begin{figure}[h]
	\centerline{
		\includegraphics[width=1.25\linewidth]{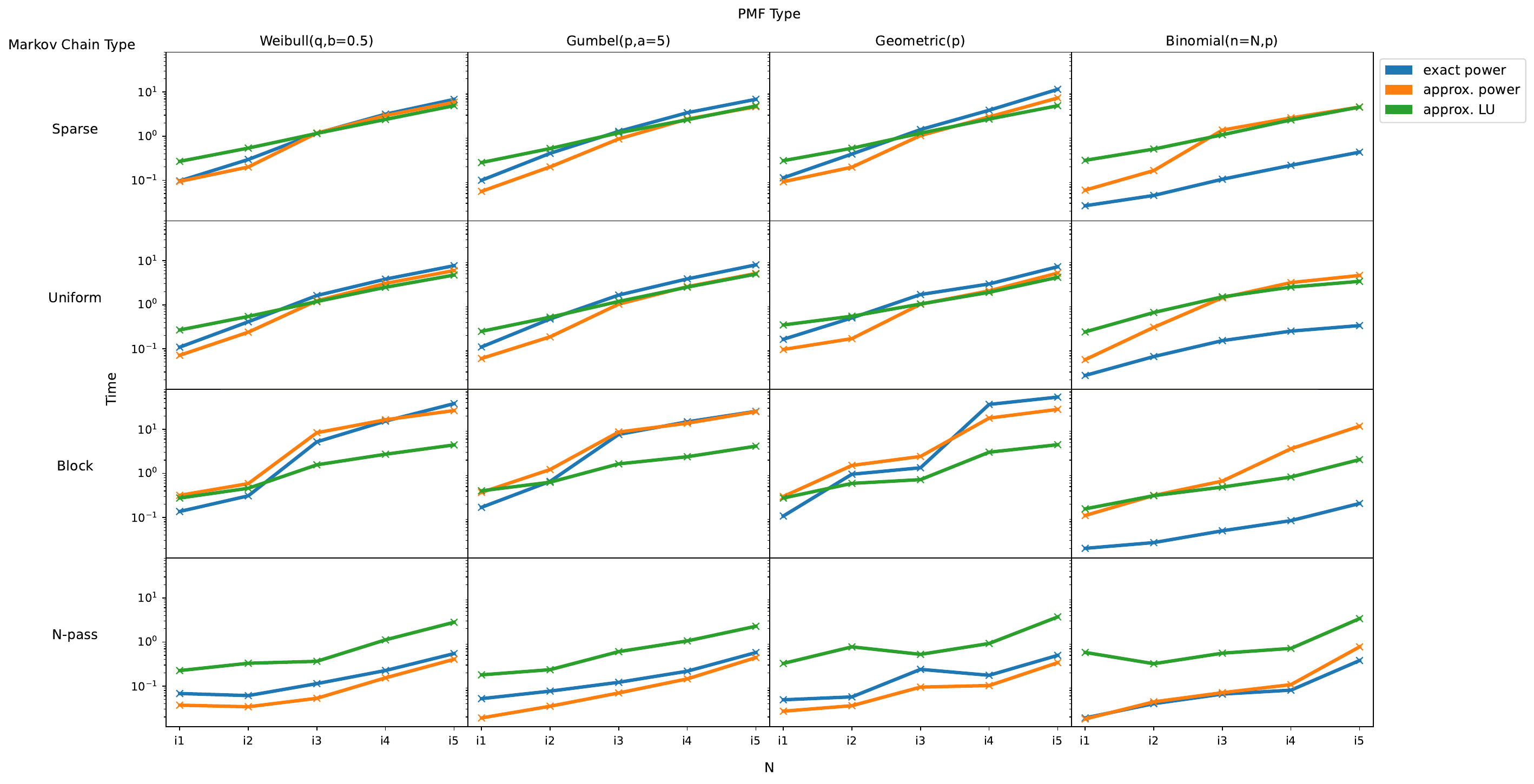}}
	\caption{Average time taken (secs.). In each plot,the x-axis labels (i1, i2, i3, i4,i5) correspond to the intervals $[1,251),[251,501),[501,1501), [1501,3001)$ and $[3001,6001)$ respectively. Comparing the speeds of the approximate methods here to their respective errors in Fig. \ref{fig:err2pmf}, neither approximate methods appear to be valuable.}
	\label{fig:time3pmf}
\end{figure}

\begin{figure}[h]
	\centerline{
		\includegraphics[width=1.25\linewidth]{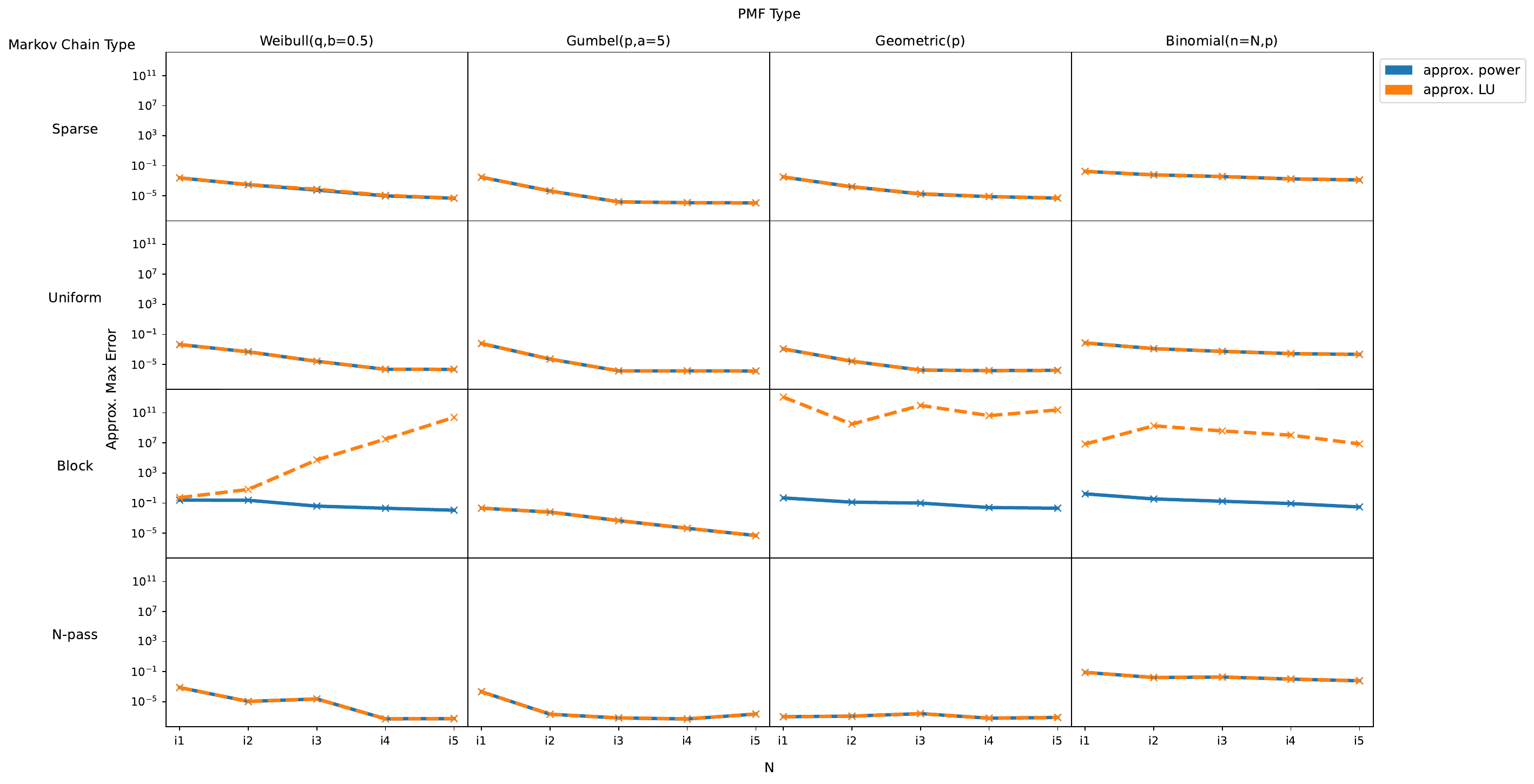}}
	\caption{Average max approximate abs. error. We find that increasing $N$, leads to a worser solution for the Block MC \& Weibull cell.}
	\label{fig:err2pmf}
\end{figure}

\paragraph{} Our second experiment involves generating samples for different values of $N$ whilst keeping $|S_?|$ fixed to $30$. The free parameters of the pmfs are once more sampled uniformly in the range $[0.3,0.6]$. The average time taken, and the average approximate max abs. errors are plotted in Fig.  \ref{fig:time3pmf} and \ref{fig:err2pmf} respectively. In each plot,the x-axis labels (i1, i2, i3, i4,i5) correspond to the intervals $[1,251),[251,501),[501,1501), [1501,3001)$ and $[3001,6001)$ respectively.  Note that for the Binomial pmf, the fixed parameter $n$ is set to $N$.

\paragraph{} For the first row of results (Sparse) of both Fig  \ref{fig:time3pmf} and \ref{fig:err2pmf}, five samples out of 1000 yielded singular matrices and therefore failure for the approximate LU method. One sample was from the Weibull column, with interval (i4). The remaining samples were from the Geometric column, with one from (i1), two from (i3) and one from (i4).


\paragraph{} The results show that the time complexity of each method is similar, which is expected. As for the error of the approximate methods, we find the same strange result where the approx. LU method fails drastically for sMRMs with the Block MC type.

\section{Summary and discussion}
\label{sec:summ_disc}
In this chapter, we presented iterative algorithms that employ the FFT for solving discrete-univariate non-negative reward sMRMs. For each of them we also presented approximate variants that could have been considered initially to be faster at solving. But it turns out when comparing the power method to its approximate variant, the results show that they may not be very useful. The iterative methods are optimizable making their solving easier, resulting in only the solution vector that requires transforming and inverting per iteration. We find that the power method is practical and can scale to small-moderate size problems without space optimizations.

\paragraph{} We also found that the algorithms here are more scalable relative to the existing symbolic approach as we could solve problems with 80 states and $k=1000$ on our laptop, whereas previously we could not solve a problem with six states (of a full graph) symbolically.

\subsection{Limitations of the iterative methods}
\label{subsec:limitations_iterative}

\paragraph{} A limitation of our solution is that the reward random variables must be defined only over the same lattice set (e.g. $\mathbb{N}$). A reason to avoid arbitrary spaces is that convolution with lattices does not require us to keep track of the abscissas. Otherwise, our solution (hyper-)vector grows per iteration, and the algorithm suffers from \textit{representation explosion}. Using lattices also allows use of the DFT, hence the FFT. Thus reward r.v.s not defined on lattices, should be approximated by lattices if possible. Even, with (non-stochastic) MRMs, the benefit of using lattices was that an algorithm for solving quantile probabilities can be done in pseudo-polynomial time, but generally exponential-time otherwise \cite{ummels2013computing}.

\paragraph{} The iterative methods presented here used for regular systems of linear equations have a deficiency in that they may prematurely converge to the wrong solution possibly due to extremely slow convergence. For example, given the iteration $\textbf{x}^{(n)} = A\textbf{x}^{(n-1)} + \textbf{b}$. Then $\textbf{x}^{(l)}$ is the solution after $l$ iterations, and $\textbf{x}^{(l+M)}$ is the solution after $l+M$ iterations where $M$ is a really large whole number. Assume that the maximum error between each iteration between $l$ and $l+M$ is less than the absolute error tolerance level $\epsilon$. Then, it is not necessarily the case that $\textbf{x}^{(l)} -\textbf{x}^{(l+M)} \leq \epsilon$.

\paragraph{} An example of where convergence to a wrong result for DTMCs is found in \cite{haddad2018interval}. To overcome this, they presented an algorithm called \textit{interval value iteration} which solves the problem twice, once from above giving an over-approximation and another from below giving an underapproximation. Both however converge to the true solution when allowed to iterate infinitely (theoretically). Note that value iteration algorithms are usable over DTMCs, whilst being defined for Markov decision processes (MDPs). And that even though the algorithm is written for systems of linear equations, this is adaptable for systems of convolution equations. Alternatively, \cite{quatmann2018sound} also provides a solution to the problem, which is stated to yield faster convergence for a particular set of problems. The solution by \cite{haddad2018interval} however requires values of the solution vector to be bounded between $[0,1]$, e.g. be probabilities, hence it cannot be used for expected rewards. An algorithm to resolve this was presented by \cite{baier2017ensuring}, and can be adapted for variance (and covariance) problems (see \cite{verhoeff2004reward}). Even with these methods, interval value iteration \cite{quatmann2018sound} does not guarantee fast convergence, i.e. convergence meeting the error tolerance level $\epsilon$ under reasonable time. See for example the criticism by \cite{mathur2020exact} and their resolution via exact model checking - an algorithm to sharpen approximate solutions derived via value iteration or interval value iteration.

\subsection{Direct convolution vs FFT convolution}
The benefit of using direct convolution versus the DFT and therefore the FFT is that the sparsity of reward random variables can be maintained. Take for example a discrete random variable defined on the points $[1,25,50]$. If we wanted to determine the probability of the accumulated reward being in $[0,99]$ (using the FFT) we will need to expand these three points into a vector of $2*100-1 = 199$ elements. Thus whilst we may gain speed-up using the FFT for dense problems, we worsen the space complexity. However, if we keep the sparse representation of the random variable, we are left with the direct convolution algorithm which is less tractable for large intervals ($O(n^2)$ vs $O(nlogn)$), and therefore we worsen time complexity. Note that whilst the hypervector $\textbf{h}$  of \eqref{eqn:power_method_early_conv} can be kept in sparse form, it is not the main contributor to space complexity growth, but rather the main contributor is the computation $(A\circ\textbf{G})\ \textcircled{$\ast$} {\textbf{f}}$.

\paragraph{} Realise that with the power method (at least), we can always partition convolutions (in $((A\circ\textbf{G})\ \textcircled{$\ast$} {\textbf{f}})$) into those more efficiently done with the FFT, and those done via direct convolution or otherwise. Therefore sparse regions can be isolated if necessary.

\paragraph{} Optimizations via overlap-add and overlap-save is possible, and appear to be  suitable areas to explore. For example, one may try to devise an algorithm to determine the optimal partition length for overlap add/save. This would be done before solving the system of equations once. Research into optimal kernel tiling \cite{pavel2013algorithms} and vector packing \cite{robertson1992computation} can be studied to reduce this blow up in space complexity without affecting time complexity drastically. Alternatively, hyper-graphs  may be used \cite{bradley2004hypergraph}.

    \chapter{Iterative Methods with Continuous Random Variables}
\section{Introduction} In this chapter we present numerical algorithms for solving first passage-reward densities (i.e. the property $\mrmpdf$) for continuous-reward sMRMs. The approach uses the DFT with numerical quadrature algorithms, as opposed to the continuous Fourier  or Laplace transform as found in literature for SMPs. We adapt some well-known quadrature rules like: The trapezoid rule, Simpson's rule and Romberg's method, for convolution. Then, the passage-reward densities can be solved for by combining these quadrature rules with the power method derived in the previous chapter for example. 

\paragraph{} With the iterative methods already defined, then if we want to extend our numerical solutions to probability density functions, i.e. to compute $\mrmpdf$ for the bounded continuous interval $[0,k]$ for $k \in \mathbb{R}_+$, all that is required is that we are able to perform at least two operations over pdfs;  summation/subtraction and convolution. This is sufficient for the power method. However, we are required to be able to perform deconvolution for the Jacobi/Gauss-Seidel iterative methods and the Gaussian elimination algorithm.

\paragraph{} We have mentioned previously in the literature review that there are existing algorithms for semi-Markov processes (SMP) that can be adapted for sMRMs. In fact, it is probably wise, for absolutely-continuous distributions to use these SMP techniques when possible, namely the Laplace-Euler technique \cite{bradley2004hypergraph}. However, it was stated that other techniques should also be considered for non-continuous pdfs \cite{warr2012introduction}, for example by the DFT. We present here a study of how the DFT can be used in such a case and investigates how well it succeeds with discontinuous pdfs. Such distributions can be derived from interpolated \textit{empirical cumulative distribution functions} (ecdf), or known distributions like the Uniform distribution, or compositions of continuous and discontinuous distributions. Note that discontinuous pdfs themselves can be approximated as smooth (absolutely continuous) distributions via a process known as \textit{mollification}, said to have been introduced by \cite{friedrichs1944identity}. The process involves smoothing the discontinuous distribution by convolving it with another distribution that is infinitely differentiable called a (Friedrich) \textit{mollifier}, named after the person who introduced it. Pacal is a (python) library that allows us to do this (by providing a mollifier distribution), and we see quite good results for the mollified Uniform distribution. However, upon smoothing a discontinuous pdf $f(x)$ (of a non-negative random variable) with a mollifier distribution, the resulting smoothed pdf $\hat{f}(x)$ will naturally have some error introduced. For example, it is possible that has  $\int_{\infty}^0\hat{f}(x)dx > 0 $, hence it is no longer a pdf of a strictly non-negative random variable.

\paragraph{} Another approach we could have taken is to use classes of functions (such as polynomials) to approximate our convolution integrals and convex combinations of pdfs. Then two issues arise in general when using analytical representations:

\begin{enumerate}
	\item The analytical representations will usually require a set of coefficients to represent. Then, convolutions of any two functions may yield a larger set of coefficients. In essence, for the finite case earlier, the re-zeroing of our DFT solution vector after each iteration (see Theorem \ref{theorem:length_of_DFTs}) (preventing either time-aliasing, or the growing of our DFT vector by subsequent zero-padding), will require an analogous counterpart in the continuous case. One simple solution for this is to re-approximate these resulting pdfs by a set of functions with a fixed finite number of coefficients, thus mitigating \textit{coefficient-growth}, or \textit{representation explosion} \cite{bradley2004hypergraph}.
	\item The second issue we have is that of convergence. To determine if the current approximation of the solution is good relative to the previous solution, we can use measures like the Kullback-Leibler divergence criteria, the Kolmogorov-Smirnov test, or Akaike's information criterion. However, we will still need to prove that convergence is guaranteed when enforcing a finite coefficient set for each function.
\end{enumerate}

\paragraph{} In this chapter, we experiment with perhaps a simpler approach, and experiment with variants of one method to deal with continuous random variables that are continuous over the non-negative real line (though it may be used with discontinuous pdfs). We will approximate convolutions via quadrature rules which will allow us to approximate the convolution of continuous random variables. This idea is already known, for example see \cite{qiang2010high}. As we shall see, this method allows us to rewrite our problem in the same form as the discrete case, and so we obtain similar time complexities as the algorithms for pmfs.

\paragraph{} In the remainder of this chapter, we first introduce numerical convolution and deconvolution techniques. We then present iterative methods for solving $\mrmpdf$ over sMRMs with continuous reward random variables, using these techniques. Finally, empirical results on the performance of the iterative methods are shown. Note that for the continuous reward random variables, we make the assumption in this chapter that \textit{no singularities exists} for them.

\section{Numerical Convolution}
\label{sec:num_conv}
In this section we develop rules for numerical convolution via old and simple ideas. We develop the Riemann sum approximations, the trapezoid rule, and the Romberg's method in that respective order. We will later show how these rules can be used to solve for the property $\mrmpdf$ within sMRMs. In this section, we assume that any functions defined below are \textit{Riemann integrable} (or continuous almost everywhere).

\subsection{Riemann sum approximation}
\label{subsection:riemann_sum}

\paragraph{} Let $X,Y$ be two non-negative continuous random variables distributed as $f(t),g(t)$ respectively and these pdfs are continuous (or continuous almost everywhere) also. Their summation gives us a new random variable  $Z = X+Y$ with pdf $h(t)$, which can be approximated via the \textit{right Riemann sum} with $N$ points as
\begin{flalign*}
h(t) = (f\ast g)(t) = \int_{0}^{t}f(t-x)g(x)dx \approx \frac{t}{N-1}\sum_{i=1}^{N-1}f(t-x_i)g(x_i) = h^{R}(t)
\end{flalign*}
or the \textit{left Riemann sum} with $N$ points as
\begin{flalign*}
h(t)   \approx \frac{t}{N-1}\sum_{i=0}^{N-2}f(t-x_i)g(x_i) = h^{L}(t)
\end{flalign*}
where $x_j = j\frac{t}{N-1}$ for $j=0,1,\cdots,N-1$.  It is known that as $N \rightarrow \infty$, then $h(t) =  h^{R}(t) = h^{L}(t)$.

\paragraph{} Then, the values ${(h^{R}(x_j))}_{0 \leq j < N}$ and ${(h^{L}(x_j))}_{0 \leq j < N}$ can be computed via a discrete convolution. Since $xs(t,N)$ is equidistant, we can employ the FFT to perform the convolution operation.

\paragraph{} Let $\textbf{x}=  \{x_j\}_{0 \leq j < N} =\{j\frac{k}{N-1}\}_{0 \leq j < N}$ be the set of $N$ equidistant points between $[0,k]$. For a given continuous function $l$, let $\bar{l}$ denote a vector of length $N$ where $(\bar{l}[j])_{0 \leq j < N} = (l(x_j))_{0 \leq j < N}$. Let $\bar{l}_0$ denote the same vector but where $\bar{l}_0[0] = 0$. Then, the right Riemann sum above can be written as
$$h^R(x_j) = \frac{k}{N-1}\mathtt{conv_N}(\bar{f},\bar{g}_0)[j]$$
and the left Riemann sum is
$$h^L(x_j) = \frac{k}{N-1}\mathtt{conv_N}(\bar{f}_0,\bar{g})[j] $$
for $j = 0,1,\cdots,N-1$.

\subsection{Trapezoid rule}
\label{subsection:trapezoid_rule}
\paragraph{} Previously, we presented rules for the Riemann approximation of $h(t)$. Now we present a way to compute the approximation via the \textit{trapezoid rule}.

\paragraph{} The trapezoid rule for $h(t)$ is simply
$$h(t) \approx \frac{1}{2}(h^{L}(t) + h^{R}(t)) = h^T(t)$$
that is, it can be computed as the average of the left and right Riemann sum. \paragraph{} Define $\textbf{x}=  \{x_j\}_{0 \leq j < N} =\{j\frac{k}{N-1}\}_{0 \leq j < N}$. Then, the trapezoid rule can be computed as
$$h^T(x_j) = \frac{k}{2(N-1)}(\mathtt{conv_N}(\bar{f},\bar{g}_0)[j] + \mathtt{conv_N}(\bar{f}_0,\bar{g})[j])$$
for $j = 0,1,\cdots,N-1$. Since the trapezoid rule is the average of the left \& right Riemann sums, it follows that as $N \rightarrow \infty$, then $h(t) =  h^{T}(t)$.

\subsection{Romberg's method}
\label{subsection:Romberg's method}
\textit{Romberg's method} \cite[p. 169]{phillips1996theory} named after W. Romberg who introduced it in  \cite{romberg1955vereinfachte},  allows us to achieve higher-order approximations and relies on the trapezoid rule. Romberg's method for $h(t)$ is defined recursively as
$$h(t) \approx h^{\mathcal{R}}_{l,N}(t) = \frac{4^{l-1}}{4^{l-1} - 1}h^{\mathcal{R}}_{l-1,(2N-1)}(t) + \frac{1}{4^{l-1}-1}h^{\mathcal{R}}_{l-1,N}(t)$$  
where $l$ is the \textit{level} or \textit{order} of the approximation with $l \geq 2$, $N$ is the number of points used, and $h^{\mathcal{R}}_{1,N}(t)$ is equal to the trapezoid rule approximation $h^T_N(t)$ with $N$ equidistant point used.

\paragraph{} To compute $h^{\mathcal{R}}_{l,N}(t)$, we first obtain the trapezoidal approximations $h^T_N(t)$, $h^T_{2N-1}(t)$, $h^T_{4N - 3}(t)$, $\cdots$, $h^T_{(2^{l-1}N) - (2^{l-1}-1)}(t)$. Then, from these terms, we compute the set of approximations $h^{\mathcal{R}}_{o,N}(t)$, $h^{\mathcal{R}}_{o,2N-1}(t)$, $\cdots $, $h^{\mathcal{R}}_{o,(2^{l-o}N) - (2^{l-o}-1)}(t)$, in the order $o=2, 3,\cdots,l$.

\paragraph{}For example, if we let $l=2$, then the Romberg method is equal to
$$h(t) \approx h^{\mathcal{R}}_{2,N}(t) = \frac{4}{3}h^{\mathcal{R}}_{1,(2N-1)}(t) + \frac{1}{3}h^{\mathcal{R}}_{1,N}(t) = \frac{4}{3}h^T_{2N-1}(t) + \frac{1}{3}h^T_{N}(t)$$
which is equal to a rule called Simpson's rule. Define $\textbf{x}=  \{x_j\}_{0 \leq j \leq {2N-1}} =\{j\frac{k}{2N-2}\}_{0 \leq j \leq 2N-1}$, and $\textbf{y}=  \{y_j\}_{0 \leq j < N} =\{j\frac{k}{N-1}\}_{0 \leq j < N}$. Then,
\begin{flalign*}
h^{\mathcal{R}}_{2,N}(y_j)
&= \frac{4}{3}h^T_{2N-1}(y_j) + \frac{1}{3}h^T_{N}(y_j) &\\
&= \frac{4}{3}\frac{k}{2(2N-2)}(\mathtt{conv_{2N - 1}}(\bar{f}_{\textbf{x}},\bar{g}_{\textbf{x}_0})[2j] + \mathtt{conv_{2N - 1}}(\bar{f}_{\textbf{x}_0},\bar{g}_{\textbf{x}})[2j]) &\\
&+\frac{1}{3}\frac{k}{2(N-1)}(\mathtt{conv_{N}}(\bar{f}_\textbf{y},\bar{g}_{\textbf{y}_0})[j] + \mathtt{conv_{N}}(\bar{f}_{\textbf{y}_0},\bar{g}_\textbf{y})[j])
\end{flalign*}
for $j = 0,1,\cdots,N-1$, and where for example $\bar{f}_\textbf{x}$ is equal to $\bar{f}$ derived from the set of points \textbf{x}. As $N \rightarrow \infty$, then $h(t) =  h^{\mathcal{R}}(t)$, furthermore as $l$ increases, the order of the accuracy increases. See Section \ref{section:error analysis} for details. An example of the Romberg method used for convolution is shown in Appendix \ref{sec:romberg_example}.

\section{Numerical Deconvolution}
\label{sec:numerical_deconv}
To use the iterative Jacobi/Gauss-Seidel algorithms, we will need to be able to numerically deconvolve.  Care has to be taken when performing numerical deconvolution since it cannot be naively integrated into the convolution equations \eqref{eqn:matrix_sol_pdf} as it involves dealing with functions that have singularities. Additionally, numerical deconvolution is generally unstable, for example see \cite{bini1986polynomial}.

\paragraph{} Note however, that the deconvolution operations within these iterative algorithms all take the same form - $pf_X \deconv (\delta_{x_0} - qf_Y)$ - where $f_X,f_Y$ are pdfs, $\delta_{x,0}$ is the Dirac delta, and $p,q$ are probabilities such that $0 \leq p+q \leq 1$, but $q \neq 1$. We exclude the case where $q = 0$, since this is trivial, as $pf_X \deconv (\delta_{x_0} - qf_Y) = pf_X$. See for example the Jacobi equation \eqref{eqn:jacobi_infinite_form}. Deriving a numerical algorithm for this deconvolution operation is made difficult due to the singularity introduced by the Dirac delta, combined with the fact that deconvolution is not distributive.

\paragraph{} However, note that this form is analogous to the limit of a geometric series. For example, by taking the FT of the deconvolution operation, we obtain
$$\frac{pF_X(\tau)}{1 - qF_Y(\tau)} = pF_X(\tau) + pF_X(\tau)qF_Y(\tau) + pF_X(\tau)q^2F_Y^2(\tau) +  \cdots = \sum_{k=0}^{\infty}pF_X(\tau)q^kF_Y^k(\tau)$$
where $|qF_Y(\tau)| < 1$. Taking the inverse FT of the above gives us
\begin{flalign}
\label{eqn:deconv_by_conv}
pf_X \deconv (\delta_{x_0} - qf_Y) = pf_X + pf_X*qf_Y + pf_X*(qf_Y*qf_Y) + pf_X*(qf_Y*qf_Y*qf_Y) + \cdots
\end{flalign}

\paragraph{} Then, the continuous convolutions above can be approximated via one of the quadrature rules previously described. For example, if using the right Riemann sum, let $\textbf{x}=  \{x_j\}_{0 \leq j < N} =\{j\frac{k}{N-1}\}_{0 \leq j < N}$. Denote $v_1 = pf_X$, and $v_2 =\delta -  qf_Y$. Then, we can introduce the deconvolution operation using convolutions via the right Riemann sum as
\begin{flalign}
\label{eqn:dvc}
\mathtt{dvc_{N,m}}^{R}(v_1,v_2) =\bar{v}_{1}  + \sum_{n=1}^{m}\mathtt{conv_N}(\bar{v}_{1},h_{N_0}^{(n)})
\end{flalign}
where $ h_{N}^{(n)} = \mathtt{conv_{N}}(\bar{h}_{N}^{(n-1)},\bar{v}_{2_0}) $, and $ h_{N}^{(1)} = \bar{v_2}$. As $n \rightarrow \infty$ and $N \rightarrow \infty$, the geometric series converges to the true solution.

\section{The continuous power method} We present a way to solve the set of equations \eqref{eqn:lineq_matrix_trans_pdf}, i.e.
$$\textbf{f} = (A\circ \textbf{G})\textcircled{$\ast$} \textbf{f} + \textbf{h}$$
for the property $\mrmpdf$, using the power method
$$
\textbf{f}^{(n+1)} = (A\circ \textbf{G})\textcircled{$\ast$} \textbf{f}^{(n)} + \textbf{h}
$$
and by incorporating the numerical convolution methods previously described. We \textbf{assume} that all pdfs in the system above are Riemann integrable (\textit{continuous almost everywhere}).

\subsection{The power method with Riemann sums}

Define $\textbf{x}=  \{x_j\}_{0 \leq j < N} =\{j\frac{k}{N-1}\}_{0 \leq j < N}$, the set of equidistant points between the interval $[0,k]$. Then, for each $n$, we can approximate $f_s^{(n+1)}$ with the right Riemann sum as
\begin{flalign*}
f_s^{(n+1)}(x_i)
&= \sum_{t \in S_?}((A_{s,t}G_{s,t}) \ast f_{t})(x_i) + h_s(x_i) &\\
&\approx \sum_{t \in S_?}\frac{k}{N-1}\mathtt{conv_N}(A_{s,t}\bar{G}_{s,t},\bar{f}_{t_0}^{(n)})[i] + h_s(x_i)
\end{flalign*}
for $i = 0,1,\cdots,N-1$. Thus, repeating this approximation for each $n$, the \textit{power method with the right Riemann sum} is
\begin{flalign}
\label{eqn:power_method_riemann_sum}
f_s^{R(n+1)}(x_i)
&=\sum_{t \in S_?}\frac{k}{N-1}\mathtt{conv_N}(A_{s,t}\bar{G}_{s,t},\bar{f}_{t_0}^{R(n)})[i] + h_s(x_i)
\end{flalign}
where $f_s^{R(0)}(x_i) = 0$. The \textit{power method with the left Riemann sum} is derived similarly. The power method with the left and right Riemann sums have the respective matrix forms
\begin{flalign}
\label{eqn:power_method_rieman_one_step}
\bar{\textbf{f}}^{L(n+1)}
= \frac{k}{N-1}(A\circ \bar{\textbf{G}}_0)\textcircled{$\ast$} \bar{\textbf{f}}^{L(n)} + \bar{\textbf{h}} \hspace{2cm}     \bar{\textbf{f}}^{R(n+1)}
= \frac{k}{N-1}(A\circ \bar{\textbf{G}})\textcircled{$\ast$} \bar{\textbf{f}}_0^{R(n)} + \bar{\textbf{h}}  
\end{flalign}


\subsection{The power method with the trapezoid rule} Let \textbf{x} be defined as above, then $f_s^{(n+1)}$ can be approximated with the trapezoid rule using
\begin{flalign*}
f_s^{(n+1)}(x_i)    &= \sum_{t \in S_?}((A_{s,t}G_{s,t}) \ast f_{t})(x_i) + h_s(x_i) &\\
&\approx \sum_{t \in S_?}\frac{k}{2(N-1)}(\mathtt{conv_N}(A_{s,t}\bar{G}_{s,t},\bar{f}_{t_0}^{(n)})[i] + \mathtt{conv_N}(A_{s,t}{{}\bar{G}_{{s,t}_0}},\bar{f}_{t}^{(n)})[i]) &\\
&\hspace{1cm} + h_s(x_i)
\end{flalign*}
for $i = 0,1,\cdots,N-1$. Then, repeating the approximation for each $n$, \textit{the power method with the trapezoid rule} is
\begin{flalign}
\label{eqn:power_method_trapezoid_rule}
f_s^{T(n+1)}(x_i)
&= \sum_{t \in S_?}\frac{k}{2(N-1)}(\mathtt{conv_N}(A_{s,t}\bar{G}_{s,t},\bar{f}_{t_0}^{T(n)})[i] + \mathtt{conv_N}(A_{s,t}{{}\bar{G}_{{s,t}_0}},\bar{f}_{t}^{T(n)})[i]) \notag &\\
&\hspace{1cm} + h_s(x_i)
\end{flalign}
The matrix form for this method is
\begin{flalign*}
\bar{\textbf{f}}^{T(n+1)}
&= \frac{k}{2(N-1)}((A\circ \bar{\textbf{G}})\textcircled{$\ast$} \bar{\textbf{f}}_0^{T(n)}
+ (A\circ \bar{\textbf{G}}_0)\textcircled{$\ast$} \bar{\textbf{f}}^{T(n)}
)+ \bar{\textbf{h}}
\end{flalign*}

\subsection{The power method with Romberg's method} Define $\textbf{x}=  \{x_j\}_{0 \leq j \leq {2N-1}} =\{j\frac{k}{2N-2}\}_{0 \leq j \leq 2N-1}$, and $\textbf{y}=  \{y_j\}_{0 \leq j < N} =\{j\frac{k}{N-1}\}_{0 \leq j < N}$. Then the \textit{power method with Romberg's method} with level $l$ is
\begin{flalign*}
f_s^{(n+1)}(y_i) &\approx f_{s_{(l,N)}}^{\mathcal{R}(n+1)}(y_i) &\\
&= \frac{4^{l-1}}{4^{l-1} - 1}f_{s_{(l-1,2N-1)}}^{\mathcal{R}(n+1)}(y_i) + \frac{1}{4^{l-1}-1}f_{s_{(l-1,N)}}^{\mathcal{R}(n+1)}(y_i)
\end{flalign*}
for $i = 0,1,\cdots,N-1$. To compute the above, we begin from the order $o=1,2,\cdots,l$. For example, let $l=2$, then
\begin{flalign*}
f_s^{(n+1)}(y_i)   &\approx f_{s_{(2,N)}}^{\mathcal{R}(n+1)}(y_i) &\\
&= \frac{4}{3}f_{s_{(1,2N-1)}}^{\mathcal{R}(n+1)}(y_i) + \frac{1}{3}f_{s_{(1,N)}}^{\mathcal{R}(n+1)}(y_i)  
= \frac{4}{3}f_{s_{(2N-1)}}^{T(n+1)}(y_i) + \frac{1}{3}f_{s_{(N)}}^{T(n+1)}(y_i)
\end{flalign*}
where we iterate using the power method with the trapezoid rule (until convergence) twice: once with the set of points \textbf{x} to obtain the approximation $f_{s_{(2N-1)}}^{T}$, and a second time with \textbf{y} to obtain $f_{s_{(N)}}^{T}$. Then, their weighted averages are taken as above.

\section{The continuous Jacobi method} This time, we use the Jacobi iterative sequence \eqref{eqn:jacobi_infinite_form}, $$\textbf{f}^{(n+1)} =  \sum_{t \in S_?  \slash s}\textbf{H} \textcircled{$\ast$} \textbf{f}^{(n)} + \bm{\kappa} $$
where $H_{s,t} = (A_{s,t}G_{s,t} \deconv (\delta - A_{s,s}G_{s,s}))$ and $\kappa_s= h_s \deconv (\delta - A_{s,s}G_{s,s})$.

\paragraph{} We elaborate here only the Jacobi method with the Riemann sum approximations. The method with other numerical integration rules are derived in the same fashion as with the power method. Replacing Jacobi with Gauss-Seidel is by analogy of Section \ref{subsec:Gauss_Seidel_method}.

\subsection{The Jacobi method with Riemann sums}

Define $\textbf{x}=  \{x_j\}_{0 \leq j < N} =\{j\frac{k}{N-1}\}_{0 \leq j < N}$, the set of equidistant points between the interval $[0,k]$. Then, for each $n$, we can approximate $f_s^{(n+1)}$ with the right Riemann sum as
\begin{flalign*}
f_s^{(n+1)}(x_i) &=  \sum_{t \in S_?  \slash s}((A_{s,t}G_{s,t} \deconv (\delta - A_{s,s}G_{s,s})) * f_t^{(n)})(x_i) + (h_s \deconv (\delta - A_{s,s}G_{s,s}))(x_i)&\\
&\approx \sum_{t \in S_?}\frac{k}{N-1}\mathtt{conv_N}(\bar{H}^{R}_{s,t},\bar{f}_{t_0}^{(n)})[i] + \kappa^R_s(x_i)
\end{flalign*}
where $\bar{H}^{R}_{s,t}(x_i) = \mathtt{dvc_{N,m}}^{R}(A_{s,t}G_{s,t},\delta - A_{s,s}G_{s,s})[i]$ and $\kappa^R_s(x_i) =  \mathtt{dvc_{N,m}}^{R}(h_s,\delta - A_{s,s}G_{s,s})[i]$, for $i = 0,1,\cdots,N-1$. Thus, repeating this approximation for each $n$, the \textit{Jacobi method with the right Riemann sum} is
\begin{flalign}
\label{eqn:jacobi_method_riemann_sum}
f_s^{R(n+1)}(x_i)
&=\sum_{t \in S_?}\frac{k}{N-1}\mathtt{conv_N}(\bar{H}^{R}_{s,t},\bar{f}_{t_0}^{R(n)})[i] + \kappa^R_s(x_i)
\end{flalign}
where $f_s^{R(0)}(x_i) = 0$. The Jacobi method with the left and right Riemann sums have the respective matrix forms
\begin{flalign}
\label{eqn:jacobi_method_rieman_one_step}
\bar{\textbf{f}}^{L(n+1)}
= \frac{k}{N-1}(A\circ \bar{\textbf{H}}^L_0)\textcircled{$\ast$} \bar{\textbf{f}}^{L(n)} + \bar{\bm{\kappa}}^L   \hspace{1.5cm}     \bar{\textbf{f}}^{R(n+1)}
= \frac{k}{N-1}(A\circ \bar{\textbf{H}}^R)\textcircled{$\ast$} \bar{\textbf{f}}_0^{R(n)} + \bar{\bm{\kappa}}^R
\end{flalign}

\section{Convergence analysis} For any variant of the power method described above, as the number of points used to approximate the convolution integrals increases, i.e. as $N \rightarrow \infty$, then the summations converge to the convolution integrals. Then, by Theorem \ref{theorem:unique sol}, a unique solution exists. And by Theorem \ref{theorem:fixed_point}, the power method converges to the unique solution as the number of iteration $n$ tends to infinity.

\paragraph{} If using the Jacobi/Gauss-Seidel methods, as $N \rightarrow \infty$, by Theorem \ref{theorem:convergence_jacobi}, they both converge to a unique solution. Additionally, the numerical deconvolution procedure described in Section \ref{sec:numerical_deconv} also converges to the analytical deconvolution operation as $N \rightarrow \infty$, and as more terms of the geometric series are computed.

\section{Error analysis}
\label{section:error analysis}
Let $\textbf{f}$ be the solution to the system \eqref{eqn:lineq_matrix_trans_pdf} and $\hat{\textbf{f}}_N, \hat{\textbf{f}}_M$ be approximate solutions derived by one of the methods described in this chapter, with $\hat{\textbf{f}}_N$ for example, using $N$ points for quadrature. Let $M > N$, then a general and practical way of approximating the \textit{absolute truncation error} (or \textit{absolute true error}) $|\textbf{E}| = |\hat{\textbf{f}}_N -\textbf{f}|$, is by using the \textit{absolute approximate error}, i.e. $|\textbf{E}| \approx |\hat{\textbf{f}}_N - \hat{\textbf{f}}_M|$. Hence, we can use the accuracy criteria
$$max_{s,r}|\hat{f_s}_N[r] - \hat{f_s}_M[r]| \leq \epsilon$$
to gauge if the accuracy is good enough (of either approximations). We know that as $N \rightarrow \infty$, and with $M > N$, then $|\textbf{E}| \rightarrow 0$. We may also use multiple approximations with the number of points $N < M < \cdots < Q$, and use the criteria above repeatedly for each consecutive pair to obtain a more secure estimate of the error.

\paragraph{} For some systems we may be able to determine an \textit{error bound} $\bm{\mathcal{E}}_N$, such that $|\hat{\textbf{f}}_N(t) -\textbf{f}(t)| = |\textbf{E}_N(t)| \leq |\bm{\mathcal{E}}_N(t)|$. Under certain conditions, one might to use the truncation error $\textbf{E}_N(t)$ for the quadrature rules (Riemann, trapezoid, Simpson \& Romberg) to determine  $\bm{\mathcal{E}}_N$.  We now proceed to derive the error bounds.

\paragraph{} Firstly, we define the truncation error for the quadrature rules above for one-dimensional convolution integrals. Define $p,q$ to be pdfs of non-negative random variables, then let $z(t) = (p \ast q)(t) = \int_{a}^{b}p(x)q(t-x)dx = \int_{a}^{b}f(x,t)dx$. Define \textbf{x} as the set of $N$ equidistant points between the interval $[a,b]$ that will be used for quadrature and let $h= \frac{(b-a)}{N}$ denote the step size between these points. Each quadrature rule can be written in the form
\begin{flalign}
\label{eqn:quadrature_form}
\int_{a}^{b}f(x,t)dx = \sum_{i=0}^{N}w_if(x_i,t) + E_N(t)
\end{flalign}
where $E_N(t)$ is the truncation error. For simplicity, we assume that $f$ is \textit{infinitely differentiable}. Then, the errors for the various quadrature rules (taken from \cite{phillips1996theory}) are as follows.
\begin{enumerate}
	\item For both the left and right Riemann sums, $E_N(t) = (b-a)\frac{h}{2}f^{'}(\xi,t)$, for some $\xi \in [a,b]$, and where  $f^{'}(x,t)$ is the first derivative with respect to $x$. The truncation error is of order $O(h)$.
	\item For the trapezoid rule, $E_N(t) = -(b-a)\frac{h^2}{12}f^{''}(\xi,t) \sim O(h^2)$.
	\item For Simpson's rule, $E_N(t) = -(b-a)\frac{h^4}{180}f^{''''}(\xi,t) \sim O(h^4)$.
	\item Romberg's method has a truncation error of order $O(h^{2l+2})$ where $l$ is the level for the method. From \cite[pg. 170]{phillips1996theory}, we obtain  $$E_N(t) = \alpha_0h^{2(l+1)}\mathscr{E}_{2(l+1)}(t) + \alpha_1h^{4(l+1)}\mathscr{E}_{4(l+1)}(t) + \alpha_2h^{6(l+1)}\mathscr{E}_{6(l+1)}(t) \cdots$$
	for some constants $a_0,a_1,a_2,\cdots$, where
	$\mathscr{E}_{2(r)} = (f_{(2r-1)x}(b,t) - f_{(2r-1)x}(a,t))$, and $f_{(2r-1)x}$ denotes the partial derivative $\partial^{(2r-1)} f/ \partial x^{(2r-1)}$. A more precise formula is given in the reference, that enables the determination of the constants  $a_0,a_1,a_2,\cdots$.
\end{enumerate}
\paragraph{} The truncation errors above can be written in a more general form. Let $f(x,\vec{y})$ be a multivariable function being integrated over $x$ as above. Then, the error bounds above take the form $$(b-a)g(h)f_{cx}(\xi_{\vec{y}},\vec{y})$$
where $f_{cx}$ denotes the partial derivative $\partial^c f/ \partial x^c$,  $g$ is a function of $h$, and $\xi_{\vec{y}}$ is value in $[a,b]$ that depends on the values $\vec{y}$. For the derivations later, note that for the degenerate function $f(x,\vec{y}) \equiv 1$, its partial derivatives are all strictly zero (as are all further partial derivatives). Hence, the truncation error is zero. Therefore the quadrature rules above can all integrate degenerate functions over any interval $[a,b]$ without error.

\paragraph{} Secondly, when approximating multiple integrals by using the same quadrature rule for each dimension, we find that the truncation error is of the same order. Philips \& Taylor \cite{phillips1996theory} showed how the error bounds for Simpson's rule can be generalized for double integrals. We further generalize this to multiple integrals with any of the rules above. The error bounds are derived as follows. Define the multiple integral
\begin{flalign}
\label{eqn:main_multi_integral}
z(t) = \idotsint_a^{b} f(\vec{x},t) \,d\vec{x}
\end{flalign}
where $\vec{x} = (x_0,x_2,\cdots,x_{\omega-1})$. Let $I = \{0,1,\cdots,N\}$. For $i = (i_0,i_1,i_2,\cdots,i_{k-1}) \in I^k$, we write $\textbf{x}^k_i = (x_{i_0},x_{i_1},\cdots,x_{i_{k-1}}) \in \textbf{x}^k$  to denote a $k$-dimensional point from the set of equidistant points  $\textbf{x}^k$ defined on the grid $[a,b]^k$. Let the set of $N$ weights for the quadrature rule being used be denoted as \textbf{w}. Then $\textbf{w}^k_i = (w_{i_0},w_{i_1},\cdots,w_{i_{k-1}}) \in \textbf{w}^k$ is the weight associated to the $k$-dimensional point $\textbf{x}^k_i$.

\paragraph{} Now, the approximation of $z(t)$ can be written as
$$z(t) \approx \sum_{i \in I^{\omega}} w_{i_0}\cdots w_{i_{\omega-1}} f(x_{i_0},\cdots,x_{i_{\omega-1}},t) = \sum_{i \in I^{\omega}}^{N} (\prod_j\textbf{w}^{\omega}_{i,j} ) f(\textbf{x}^{\omega}_i,t) $$    
\paragraph{} As for determining the error bounds, applying a quadrature rule to the innermost integral of equation \eqref{eqn:main_multi_integral} yields
$$\int_{a}^{b}f(x_0,x_1,\cdots,x_{\omega-1},t)dx_0 =  \sum_{i \in I}w_{i_0}f(x_{i_0},x_1,\cdots,x_{\omega-1},t) + E_{x_0}(t) $$
where $E_{x_0}(t) = $$(b-a)g(h)f_{cx_0}(\xi,x_1,\cdots,x_{\omega-1},t)$. Then, applying the quadrature rule to the next innermost integral for each $i_0$ leads to
\begin{align*}
&\int_{a}^{b}\Big(\sum_{i \in I}w_{i_0}f(x_{i_0},x_1,\cdots,x_{\omega-1},t) + E_{x_0}(t)\Big)dx_1 &\\
&= \sum_{i \in I^2}w_{i_1}w_{i_0}f(x_{i_0},x_{i_1},\cdots,x_{\omega-1},t) + w_{i_0}E_{x_{1}}(t,i_0) + \int_a^bE_{x_0}(t) dx_1
\end{align*}
where $E_{x_{1}}(t,i_0) = $$(b-a)g(h)f_{cx_1}(x_{i_0},\xi,x_2,\cdots,x_{\omega-1},t)$. In this fashion, after applying the quadrature rule to all integrals, we find for example that the error from the $1^{st}$ quadrature accumulates to $ \hat{E}_{x_0}(t)$ defined as
\begin{align*}
\hat{E}_{x_0}(t) &= \idotsint_a^{b}E_{x_0}(t)\,dx_1 \cdots dx_{\omega-1}
\end{align*}
\paragraph{} We find generally that the total error accumulated from the $\gamma^{th}$ quadrature, where $\gamma=1,\cdots,\omega$, is
\begin{flalign}
\label{eqn:error_x_general_form}
\hat{E}_{x_{\gamma-1}}(t) = \sum_{i \in I^{\gamma-1}}w_{i_0}\cdots w_{i_{\gamma-2}}\idotsint_a^{b}E_{x_{\gamma-1}}(t,i_0,\cdots,i_{\gamma-2})\,dx_{\gamma} \cdots dx_{\omega-1}
\end{flalign}
where $E_{x_{\gamma-1}}(t,i_0,\cdots,i_{\gamma-2}) = $$(b-a)g(h)f_{c\gamma-1}(x_{i_0},\cdots,x_{i_{\gamma-2}},\xi,x_{\gamma}
,\cdots,x_{\omega-1},t)$. When $\gamma = 0$, there are no summation variables, i.e. since $I^0 = \{\}$. When $\gamma = \omega$, there are no integration variables, since $d_{x_{\omega}} > d_{x_{\omega-1}}$. The total truncation error is then
\begin{flalign}
\label{eqn:total_truncation_error}
\sum_{j=0}^{\omega-1}\hat{E}_{x_j}(t)
\end{flalign}

\paragraph{} Note however, that if the modulus of $f_{c{x_j}}(\cdots,t)$ is bounded above by $M_x(t)$, for each of $j=0,1,\cdots,\omega-1$. Then,
\begin{flalign*}
|\hat{E}_{x_j}(t)|
&\leq (b-a)|g(h)|M_{x_{j}}(t)\sum_{i \in I^{j}}w_{i_0}\cdots w_{i_{j-1}}\idotsint_a^{b}1\,dx_{j+1} \cdots dx_{\omega-1} &\\
&= (b-a)|g(h)|M_{x_{j}}(t)\sum_{i \in I^{j}}w_{i_0}\cdots w_{i_{j-1}}(b-a)^{(\omega-1) -j} &\\
&= (b-a)|g(h)|M_{x_{\gamma-1}}(t)(b-a)^{j}(b-a)^{(\omega-1) -j} &\\
&= (b-a)^{\omega}|g(h)|M_{x_{\gamma-1}}
\end{flalign*}
\par{} The integrals above are applied over the constant 1, and therefore integrate to $(b-a)^{(\omega-1) -j}$. The summations with weights $w_{i_0}\cdots w_{i_{j-1}}$ can be interpreted as applying a quadrature approximation to an integral of a degenerate function $g(x_0,\cdots,x_{j-1}) \equiv 1$, which therefore leads to no truncation error and we obtain the result $(b-a)^{j}$.

\paragraph{} Hence, the total truncation error \eqref{eqn:total_truncation_error} is bounded above as
$$|\sum_{j=0}^{\omega-1}\hat{E}_{x_j}(t)| \leq (b-a)^{\omega}|g(h)|\sum_{j=0}^{\omega-1}M_{x_{j}}(t) = (b-a)^{\omega}|g(h)|\sum\vec{M}_{f}(t)$$

\paragraph{} We can now address the truncation error and its order for the methods described in this chapter for solving systems of convolution equations. The pdf $f_s$ has the form
$$f_s(x) = \sum_{\hat{\pi} \in \Pi.s}Pr(\hat{\pi})f_{Rew(\hat{\pi})}(x)$$
where each $f_{Rew(\hat{\pi})}(x_i)$  is equal to a series of convolutions between the pdfs on each transition in $\hat{\pi}$. The series of convolution equates to an $(|\hat{\pi}|-1)$-dimensional integral, i.e.
\begin{align*}
f_{Rew(\hat{\pi})}(t) &= (f_{rew(s_0,s_1)}* f_{rew(s_1,s_2)} * \cdots
*f_{rew(s_{|\hat{\pi}|-2},s_{|\hat{\pi}|-1})})(t) &\\
&= \idotsint_0^{t} f(x_0,x_1,\cdots,x_{|\hat{\pi}|-2},t)  \,dx_0 \cdots dx_{|\hat{\pi}|-2}    
\end{align*}
for some function $f$. Let $\hat{f}_{Rew(\hat{\pi})}(t)$ be the approximation of the integral via a quadrature rule. Denote the length $|\hat{\pi}|- 1= p$, then
$$|\hat{f}_{Rew(\hat{\pi})}(t) - f_{Rew(\hat{\pi})}(t)| \leq t^{p}|g(h)| \sum \vec{M}_{f_{Rew(\hat{\pi})}}(t)$$
hence,
\begin{flalign}
\label{eqn:conv_equations_error_bound}
|\hat{f}_s(t) - f_s(t)| \leq  \sum_{\hat{\pi} \in \Pi.s}Pr(\hat{\pi})t^{p}|g(h)|\sum \vec{M}_{f_{Rew(\hat{\pi})}}(t)
\end{flalign}
or in terms of the truncation error order,
\begin{flalign}
\label{eqn:conv_equations_order_bound}
|\hat{f}_s(t) - f_s(t)| \leq  \sum_{\hat{\pi} \in \Pi.s}Pr(\hat{\pi})O(|g(h)|) = O(|g(h)|)
\end{flalign}

\paragraph{} The order proves that the various quadrature rules converge at rates identical to the one-dimensional convolution case. As for generally using the error bounds in practice, then it does not appear to be tractable due to the requirement of computing $\vec{M}_{f_{Rew(\hat{\pi})}}(t)$ for example. This involves convolving all pdfs of the path analytically,  being able to derive its partial derivatives, and then the ability to determine the maximum value of their moduli.     

\section{Experiments with sMRMs} We move forward to experimenting with continuous-reward sMRMs. We will use the power method as the iterative algorithm of choice. 

\subsection{Example 1a: A toy problem with discontinuities}
\label{problem:5statesMRM_continuous}
In this problem, each reward random variable of a given sMRM is the same - a mixture of uniform distributions; a discontinuous distribution. We seek to determine the first-passage reward density between the interval $[0,100]$, of reaching $B$ from a state $s_0$. We also investigate if mollification of our reward pdfs - the smoothing of our discontinuous reward pdfs, can help us attain better accuracy. We used the power method to derive the solution and the absolute error tolerance level was set to $1e$-8.

\paragraph{} Let us use the following problem to evaluate our results. We will re-use the matrix and vector $A,b$ respectively from the previous problem in Section \ref{problem:5statesMRM}:

\[ A =
\begin{blockarray}{ccccc}
& s_0 & s_1 & s_2 & s_3 &  \\
\begin{block}{c(cccc)}
s_0 & 0.1288838  &0.38242891& 0.12495781 &0.13139189\\
s_1 & 0.27758284 &0.09654253& 0.15592425 &0.24690511 \\
s_2 & 0.10418887 &0.18054794& 0.1492027  &0.32815053 \\
s_3 & 0.33540355 &0.31410283& 0.16746947 &0.1316041 \\
\end{block}
\end{blockarray}
\]
\[    
b = \begin{pmatrix}
0.23233759 &0.22304527 &0.23790995 &0.05142005
\end{pmatrix}
\]
\paragraph{} The transition probabilities from $s_i$ to the three other states corresponds to the $i^{th}$ row of the matrix $A$. Likewise, its probability of entering $B$ immediately is the $i^{th}$ element of $b$.

\begin{figure}[H]
	\centering
	\includegraphics[width=0.8\linewidth]{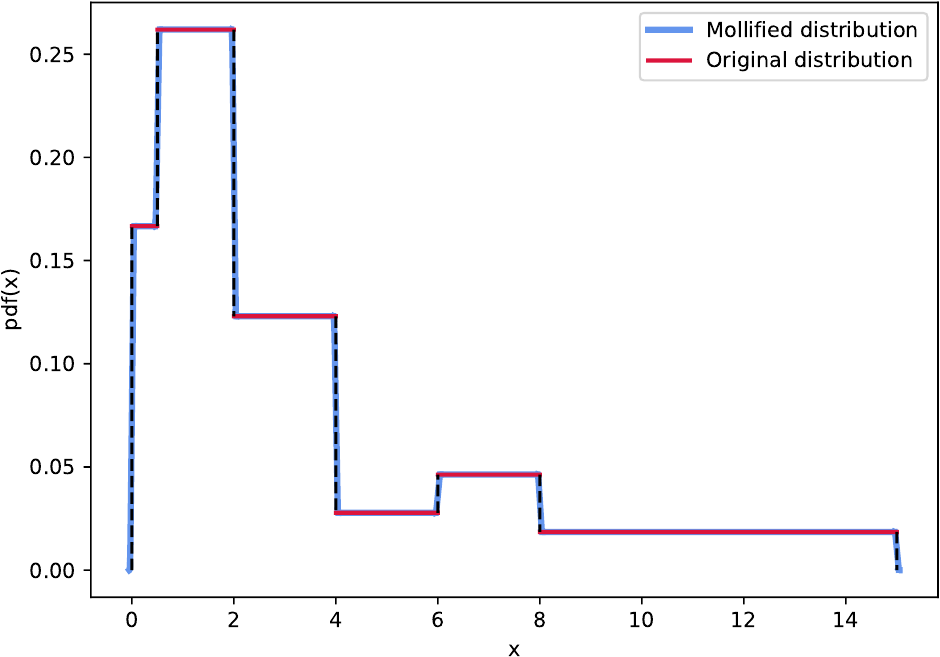}
	\caption{Pdf of a mixture of uniform distributions.}
	\label{fig:mollifieddistribution-crop}
\end{figure}

\paragraph{} As for the rewards, each element is the mixture of uniform distributions:
$$f_{rew(s,t)}(x) = \frac{1}{3}f_{U(0,2)}(x) +  \frac{1}{3}f_{U(0.5,4)}(x) + \frac{1}{6}f_{U(2,8)}(x) + \frac{1}{6}f_{U(6,15)}(x)$$
where $f_U(a,b)$ is the pdf of the uniform distribution defined in the interval $[a,b]$. We also mollify the distribution using the Pacal library, via a mollifier distribution. The mollifier distribution is added to the reward random variable (i.e. their pdfs are convolved via Pacal) to smooth it out. It has a parameter to determine how much smoothing occurs, which we set to $0.05$. The mixture distribution, unadulterated and mollified, is shown in Figure \ref{fig:mollifieddistribution-crop}. Note that the mollification procedure causes the distribution to have a slight measure for values $x < 0$.

\begin{figure}[H]
	\centerline{
		\includegraphics[width=1\linewidth]{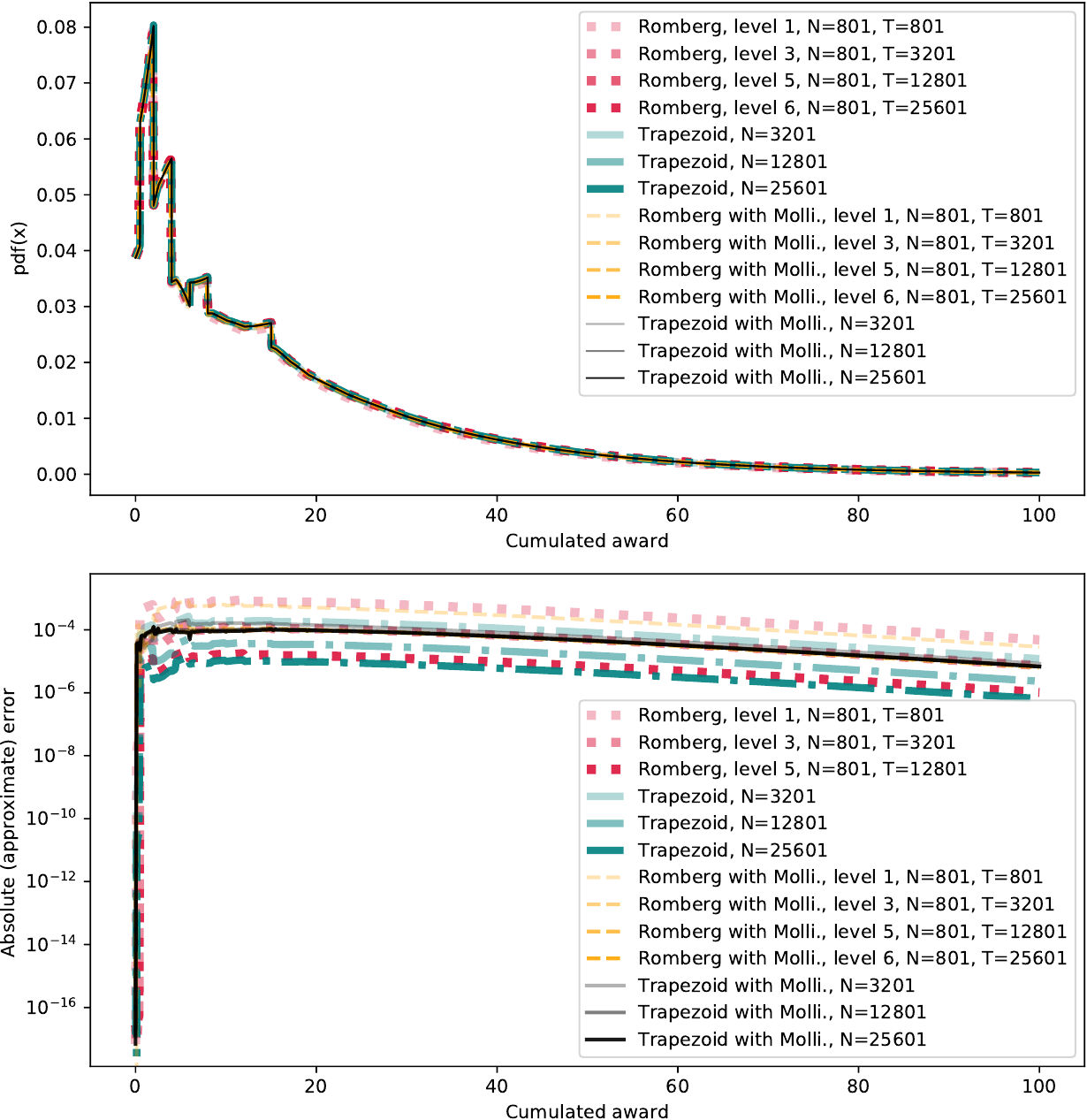}}
	\caption{Top: Approximations of the pdf $f_{s_0}(r) = Pr(r \cap s_0 \vDash \Diamond B)$ via various quadrature rules and number of points $N$. Bottom: Error of approximate pdfs relative to the $6^{th}$ level Romberg approximation. $T$ is the maximum number of points needed, i.e. the solution via the power method with the trapezoid rule at $N=T$ is required.}
	\label{fig:contsmrmdiscontinuous-crop}
\end{figure}

\paragraph{} The density for $Pr(r \cap s_0 \vDash \Diamond B)$ is plotted in Figure \ref{fig:contsmrmdiscontinuous-crop}. Romberg's method of various levels and the trapezoid rule with various numbers of samples were used. The error plotted is the absolute error relative to Romberg's method, at the $6^{th}$ level (or order). The error is also shown for the same set of points of each solution. This is possible with the use of equidistant points, as well as the number of points used.

\paragraph{} We find that a considerable number of points are required to attain an absolute error around four decimal places. With an increase of points, mollification  appears to either converge more slowly or to the wrong solution. Since mollification introduces error, we should expect convergence to a slightly wrong solution. The approximate pdf has values strictly less than two decimal places, therefore true accuracy is better measured beyond two decimal places.

\begin{figure}[H]
	\centerline{
		\includegraphics[width=1\linewidth]{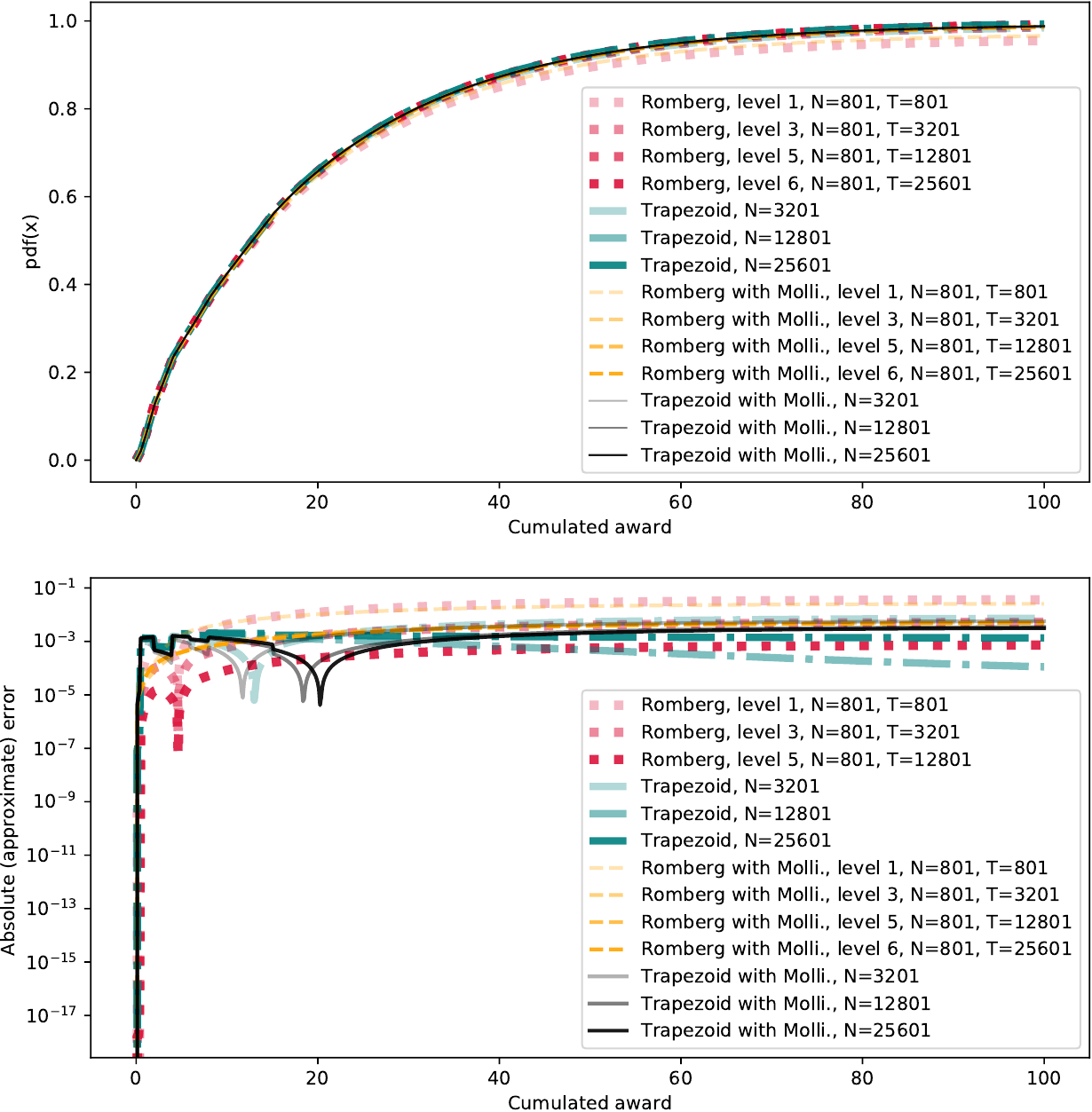}}
	\caption{Top: Approximations of the cdf $F_{s_0}(r) = Pr(Rew \leq r \cap s_0 \vDash \Diamond B)$ via various quadrature rules and number of points $N$. Bottom: Error of cdfs relative to the $6^{th}$ level Romberg approximation.  }
	\label{fig:contsmrmcdfdiscontinuous-crop}
\end{figure}

\paragraph{} The respective cdfs and their absolute errors are also plotted in Figure \ref{fig:contsmrmcdfdiscontinuous-crop}. The cdfs are derived using the cumulative trapezoid rule over the pdfs. We find some rules, like the trapezoid rule with mollification, gain a significant digit in error. No rule obtains four decimal places of error. Note that when determining probabilities with cdfs or computing quantile queries, we may have to interpolate them which potentially introduces further error.

\subsection{Example 1b: A toy problem without discontinuities}  This problem is identical to the previous problem $\textbf{1a}$, except that each reward distribution is replaced with the Weibull distribution with parameters $k=3,\lambda=1$, hence each reward pdf is continuous. The power method was used again, but the absolute error tolerance level $\epsilon$ was set to $1e$-16.

\begin{figure}[H]
	\centerline{
		\includegraphics[width=1\linewidth]{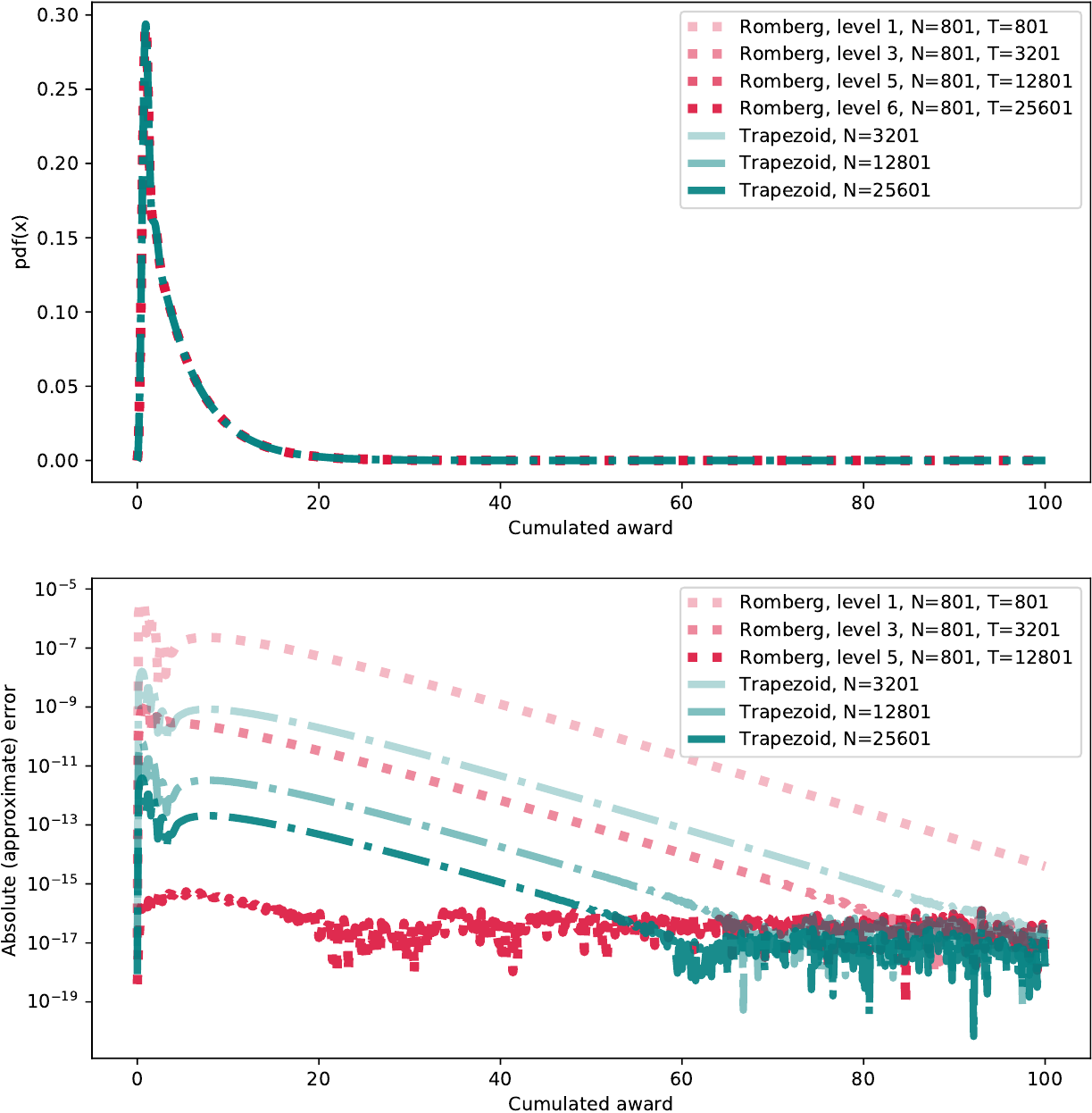}}
	\caption{Approximations of the pdf $f_{s_0}(r) = Pr(r \cap s_0 \vDash \Diamond B)$ via various quadrature rules and number of points $N$. Bottom: Error of approximate pdfs relative to the $6^{th}$ level Romberg approximation.}
	\label{fig:contsmrmcontinuous-crop}
\end{figure}

\paragraph{} We find in this case greater accuracy (or low absolute error relative to Romberg's method at level 6) can be achieved with far fewer points. The pdfs and their respective errors are plotted in Figure \ref{fig:contsmrmcontinuous-crop}.

\subsection{Example 2: Movement of coronary patients}
\label{problem:coronary_patients}
We obtain a problem from \cite{warr2012introduction} with a model originally from \cite{kao1974modeling}, that captures the movement of myocardial infarction positive patients in a hospital. The model approximates the likelihood of patients transiting between 9 states: coronary care unit (CCU), post-coronary care unit (PCCU), intensive-care unit (ICU), medical unit (MED), surgery (SERG), ambulatory care (AMB), extended care facility (ECF), HOME and DIED. The model is an sMRM, represented by the probability matrix $\textbf{P}$ and the reward matrix $rew$ below. In this instance the reward represents time, and thus the sMRM is an SMP. Such a model may be useful for operational research purposes, i.e. to aid in planning or decision making.

\paragraph{} The following is the probability matrix of the sMRM:
\begin{figure}[h]
	\centerline{
		$
		\textbf{P} =
		\begin{blockarray}{cccccccccc}
		\ & CCU& PCCU& ICU& MED& SURG& AMB& ECF& HOME& DIED\\
		\begin{block}{c(ccccccccc)}
		CCU     &0.0000 &0.7447& 0.0084 &0.1339& 0.0042& 0.0063& 0.0000 &0.0063& 0.0962\\
		PCCU    &0.0192 &0.0000& 0.0137 &0.0247& 0.0027& 0.0027& 0.0577 &0.8298& 0.0495 \\
		ICU        &0.0000 &0.5833& 0.0000 &0.1667& 0.0833& 0.0000& 0.0000 &0.0000& 0.1667 \\
		MED        &0.0000 &0.0135& 0.0405 &0.0000& 0.0135& 0.0270& 0.0811 &0.7028& 0.1216 \\
		SURG    &0.0000 &0.0000& 0.0000 &0.0000& 0.0000& 0.0000& 0.0000 &1.0000& 0.0000 \\
		AMB        &0.0000 &0.0000& 0.0000 &0.0000& 0.0000& 0.0000& 0.0000 &1.0000& 0.0000 \\
		ECF        &0.0000 &0.0000& 0.0000 &0.0000& 0.0000& 0.0000& 1.0000 &0.0000& 0.0000 \\
		HOME    &0.0000 &0.0000& 0.0000 &0.0000& 0.0000& 0.0000& 0.0000 &1.0000& 0.0000 \\
		DIED    &0.0000 &0.0000& 0.0000 &0.0000& 0.0000& 0.0000& 0.0000 &0.0000& 1.0000 \\
		\end{block}
		\end{blockarray}
		$}
\end{figure}
\par{} The reward matrix consists of only five unique Weibull($\gamma,\theta$) distributions:
\begin{flalign*}
&f_1 \sim Weibull(4.738025, 4566277818.13) &\\
&f_2 \sim Weibull(2.207438, 14541.6089) &\\
&f_3 \sim Weibull(0.766338, 16.6991) &\\
&f_4 \sim Weibull(2.303331, 1017649.5158) &\\
&f_6 \sim Weibull(1.623492, 4707.3132)
\end{flalign*}
and where the Weibull distribution is defined as
$$
f(k;\gamma,\theta) = \frac{\theta}{\gamma}k^{\gamma-1}e^{-(k^{\gamma}/\theta)}
$$    
with $k>0, \gamma>0, \theta>0$.     The random cost in this problem is the amount of time elapsed. The reward/cost (pdf) matrix is defined as
$$
rew =
\begin{pmatrix}
&    f_1    &    f_1 &    f_1    &    f_1    &     f_2    &          &    f_2    &     f_3\\
f_4 &         &     f_1 &    f_4    &    f_1    &     f_1    &     f_4 &    f_4    &     f_6\\
&     f_4 &         &    f_1 &     f_1 &          &         &          &     f_3\\
&    f_4 &    f_4 &        &    f_4    &    f_4    &    f_4    &    f_4    &     f_6\\
&        &        &        &        &        &        &    f_4    &      \\
&        &        &        &        &        &        &    f_4    &      \\
&        &        &        &        &        &        &        &      \\
&        &        &        &        &        &        &        &      \\
&        &        &        &        &        &        &        &      \\
\end{pmatrix}$$
where empty spaces are placeholders for the (zero) Dirac delta distribution $\delta_{k,0}$. Now, we can use both \textbf{P} and $rew$ to obtain the set of terms: $A,\textbf{G}$ and $\textbf{h}$ for each problem we are required to solve below. Note that there appears to be a discrepancy between the reward matrix of \cite{warr2012introduction} and \cite{kao1974modeling}, in the second row of the second-to-last column on page 690 (of matrix $H^1$). The latter had the element as blank, but the former used $f_4$ instead. It may be that there is a typographical error in the latter, and so we will plot our experiments relative to \cite{warr2012introduction}. Perhaps the validity of this is that in doing so, our experiments align more closely with the results of both works.

\paragraph{} Formally, we are solving for the reward-bounded reachability probability $\mrmpdfcostbounded$, described in Section \ref{sec:reachability_problems}. This as we have mentioned is considered a cumulative distribution function with the variable $r$. The set of absorbing states is $Abs = \{\text{ECF, HOME, DIED}\}$. In this problem we are interested in the cdf $Pr(\text{CCU} \vDash_{\leq r} B)$ where $B \in Abs$. We solve this by computing the passage-reward density $Pr(r\ \cap\ \text{CCU} \vDash \Diamond B)$ or $f_{\text{CCU} \vDash \Diamond B}$, and then integrating it to obtain the cdf. Hence, we solved three systems of convolution equations, once for each state in $Abs$, each system with the form $\textbf{f} = A\circ\textbf{G} \textcircled{$\ast$}\textbf{f} + \textbf{h}$,
defined in Section \ref{subsubsec:definitions_of_terms}. 

\paragraph{} The power method with the trapezoid rule was used with $N = 4001$ (the number of points for quadrature). After obtaining an approximation for \textbf{f}, each $f_s$  was integrated via the cumulative trapezoid rule. Note that since $Pr(s \vDash \Diamond B) \leq 1$ (because $B$ is a proper subset of $Abs$), then each \textbf{f} is a hypervector of partial pdfs (see Section \ref{subsection:partial_densities}). Ultimately, we obtain three partial cdfs of the form $F_{\text{CCU} \vDash \Diamond B}(r)$ for each $B \in Abs$.

\begin{figure}[h]
	\centerline{
		\includegraphics[width=1\linewidth]{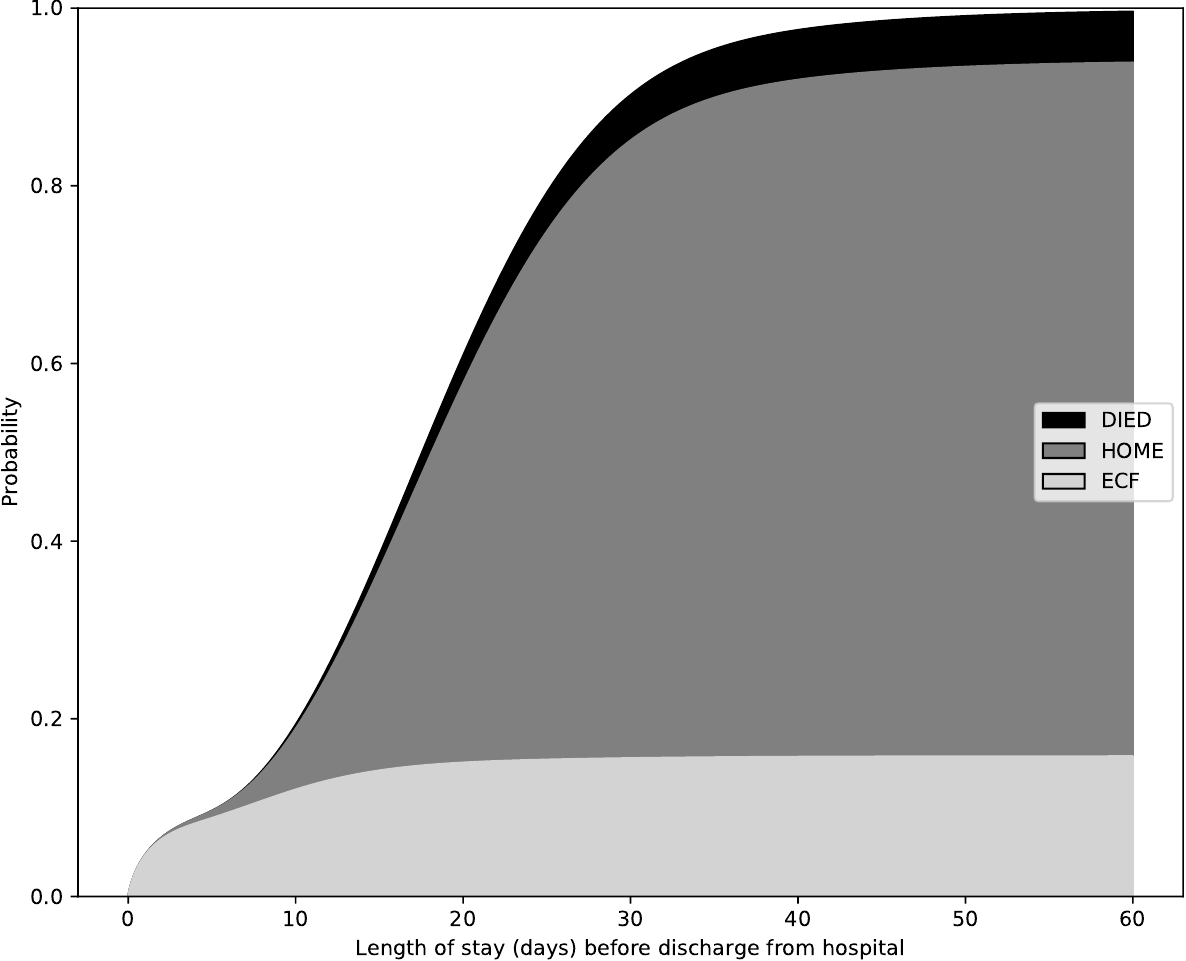}}
	\caption{The partial cdfs of being in state DIED, HOME and ECF  are plotted against time. The process started in the state CCU. The cdfs are stacked on top of each other (by addition) to create a full cdf that converges to 1 as time tends to infinity. Each colour encapsulates the regions of a particular cdf.}
	\label{fig:coronarypatientscdfs-crop}
\end{figure}

\paragraph{} We plot the (partial) cdf of reaching each absorbing state in Fig. \ref{fig:coronarypatientscdfs-crop} starting from the state CCU. The result appears similar to that found in \cite{kao1974modeling} who used a different discretization technique instead. The plot is read by drawing a vertical line from the x-axis, and the proportion of colour intersecting the line is the probability of being in the state with that respective color, at the chosen time. The cdf $F_{\text{CCU} \vDash \Diamond Abs}(r)$ can be computed as
$$F_{\text{CCU} \vDash \Diamond Abs}(r) \equiv \sum_{B \in Abs}F_{\text{CCU} \vDash \Diamond B}(r)$$ and is represented by the uppermost edge, i.e. where the white density intersects the black density. If the survival functions are computed instead (i.e. flipping the graph below upside down), then the results are similar to those of \cite{warr2012introduction} for their plot of time vs. the probability of being in a particular state.

\subsection{Empirical evaluation of convergence rate, time and error} 
\begin{figure}[H]
	\centering
	\includegraphics[width=1\linewidth]{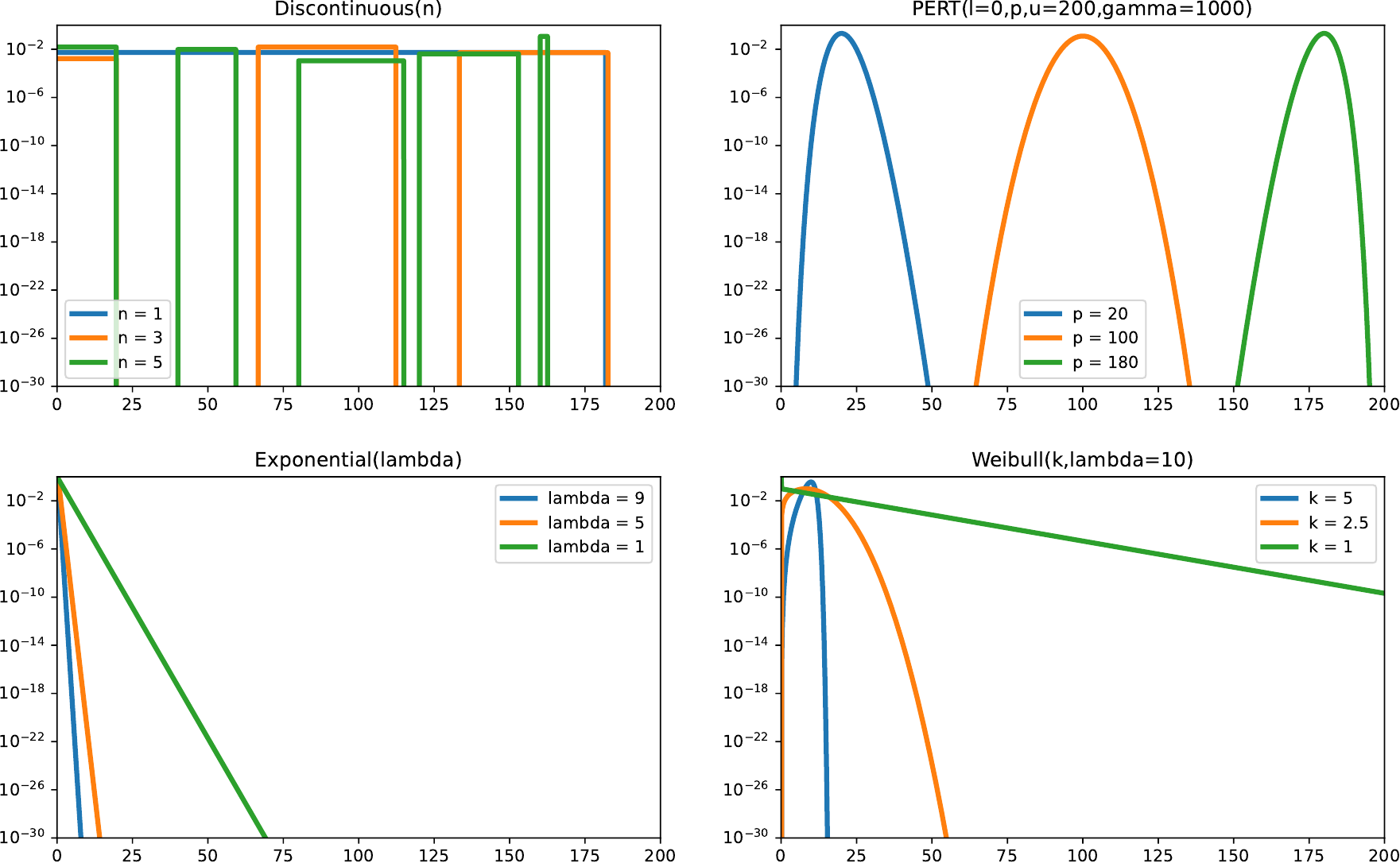}
	\caption{The shapes of the various distributions used in the experiment. Weibull is considered a heavy tail-distribution whilst the Exponential is a medium (weight)-tail distribution. The PERT and Discontinuous (piecewise uniform) distribution have pdfs supported on the bounded interval $[0,200]$.}
	\label{fig:sampleconttails}
\end{figure}

We investigate here the effectiveness of the Romberg method for solving specific classes of sMRMs. Specifically, this experiment is similar to the experiment found in Section \ref{sec:exp_MC_pmf}. Thus, it should be read first. The difference now is that we are working with four continuous reward random variables. Additionally, the iterative methods of choice are: Romberg's method of levels 1, 2 and 3, with the power method. This corresponds to the power method with the trapezoid rule, Simpson's rule and \textit{Boole's rule} respectively. Once more, we use the four MC types to generate the probability matrices (see Fig. \ref{fig:sampled_MCs}), but now we have four types of pdfs (see Fig. \ref{fig:sampleconttails}). The computer used for this entire section is \textit{Computer 2}. See Section \ref{sec:computer details} for details of the computer and software used.

\paragraph{} The pdfs are of four types: A mixture of uniform distributions (we name this as Discontinuous), (Modified) PERT, Exponential, and Weibull. We test two tail distributions: Weibull (heavy-tail) and Exponential (medium-tail). We use two distributions with bounded support: PERT and Discontinuous, with the latter providing us results for discontinuous distributions. Their definitions are as follows:

\begin{enumerate}
	\item $Discontinuous(n)(t) = \sum_{i=1}^{n}p_i f_{U[x_i,x_i + hq_i]}(t)$, where $f_{U[c,d]}$ is the pdf of a uniform distribution supported on the interval $[c,d]$. Then, $p_i,q_i \in [0,1]^2$, but we also have $\sum_i p_i = 1$. Lastly, $\{x_j\}_{0 \leq j < n} =\{j\frac{b}{n-1}\}_{0 \leq j < n}$ is the set of $n$ equidistant points between $[0,b]$, with $h$ being the step size.
	\item $PERT(l=0,p,u=200,\gamma=1000)(t)$. This is the modified variant of the PERT distribution, extending it with a $\gamma$ parameter that controls the shape of the tail values of the distribution. This is stated to have been proposed by Vose \cite{vose2000risk}. A definition of the pdf can be found in \cite{vosesite}. In it, the symbols $l,p,u$ are replaced by $min,mode,max$ respectively.
	\item $Exponential(\lambda)(t)$. This is defined as per usual.
	\item $Weibull(k,\lambda=10)(t)$. This is defined as per usual.
\end{enumerate}
\paragraph{} The free parameters of these distributions are $n,p,\lambda$ and $k$ respectively. The fixed parameters above are chosen in a way to allow us to increase the tail strengths for the pdf with unbounded tails, and to move the peak of the PERT distribution across the interval $[0,200]$. The other PERT parameters are chosen such that the shape is similar to that of the discrete binomial distribution (used for the experiments on discrete-reward sMRMs (see Fig. \ref{fig:sampletails})). For each \textit{Discontinuous} variable, the values $(p_i)_{i \leq n}$, $(q_i)_{i \leq n}$ are sampled uniformly between $[0,1]$, and the $p_i$'s are then normalized to sum to one.

\paragraph{Experiment set-up} The experiments below consist of randomly generating sMRMs, to solve $\mrmpdf$ via the system of equations \eqref{eqn:lineq_matrix_trans_pdf}. The sMRM is generated via a random sample of a selected MC type and pmf type. All the reward random variables for a given sMRM are of the same type. Their free parameters are sampled uniformly within a selected range shown later. For all experiments below, unless otherwise mentioned, we set $S_? = 30$, and we solve for $\mrmpdf$ (for all $s \in S_?$), for each $r \in \textbf{x} \subset [0,b]$ with $b = 200$, and where $\textbf{x}=  \{x_j\}_{0 \leq j < N} =\{j\frac{b}{N-1}\}_{0 \leq j < N}$ is the set of $N$ equidistant points between $[0,b]$. We set $N = 801$. For each iterative method, we used the following termination criteria $$max_{s,r}|\textbf{f}^{(n+1)} - \textbf{f}^{(n)}| \leq \epsilon$$ where $\epsilon = 1e-7$.
Additionally, a max iteration of 2000 was used to terminate the method if convergence was too slow.  For each combination of MC type, pdf type and parameter range, 50 unique sMRMs were sampled.

\begin{figure}[h]
	\centerline{
		\includegraphics[width=1.25\linewidth]{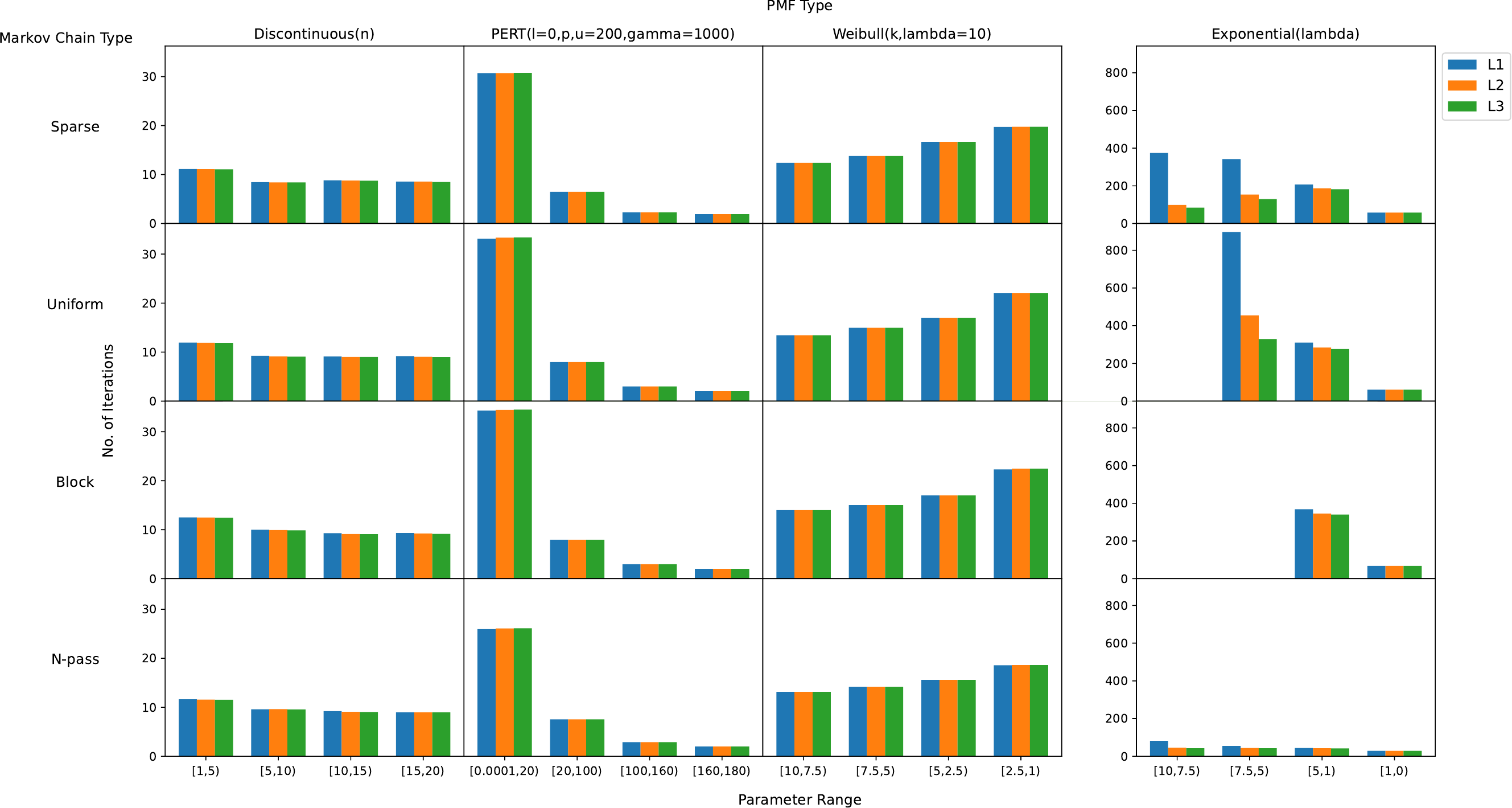}}
	\caption{Average no. of iterations. L1, L2, L3 refers to the Trapezoid rule with $N,2N-1,  4N-2$ points respectively (between [0,200]). The Exponential column is separated due to some of its values being significantly larger than the remainder of the plots. Additionally, some entries are blank due to lack of convergence. See Table \ref{table:plots_expon_column}. For each triplet of bars, 50 unique sMRMs were sampled.}
	\label{fig:noiter3cont}
\end{figure}

\begin{figure}[h]
	\centerline{
		\includegraphics[width=1.25\linewidth]{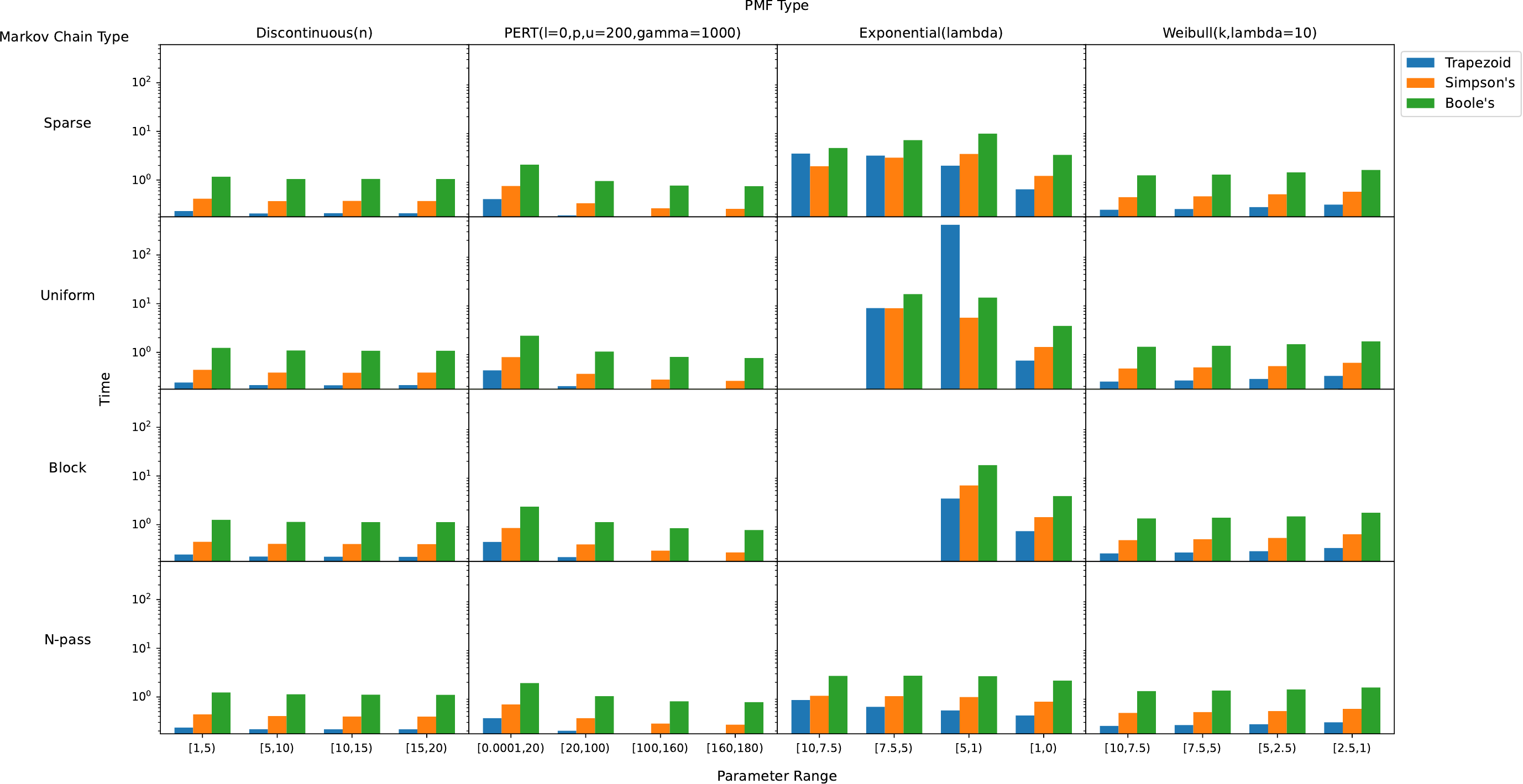}}
	\caption{Avg. time taken (secs.). The time excludes the combining of the trapezoid rules to form the Romberg approximations. This latter procedure is easy and does not take much time. The average time taken is not surprising and is reflective of the no. of iterations (see Fig. \ref{fig:noiter3cont}). However, it does help gauge the speed of the methods.}
	\label{fig:time3cont}
\end{figure}

\begin{figure}[h]
	\centerline{
		\includegraphics[width=1.25\linewidth]{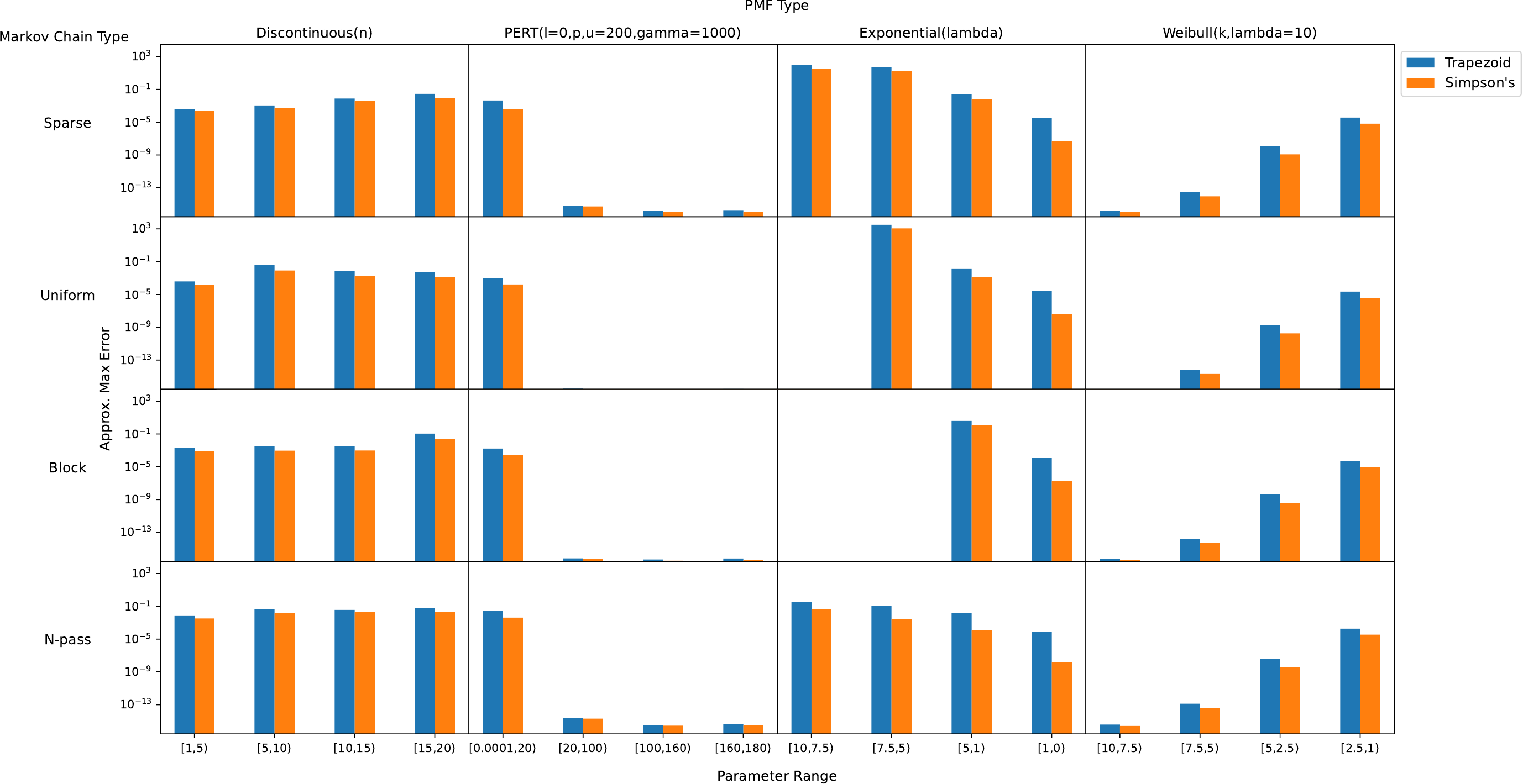}}
	\caption{Average (approximate) max abs. error (compared to Boole's). We find for the Discontinuous distribution, the max error is quite large for all MC types. Note that some entries above are blank due to lack of convergence (see \ref{table:plots_expon_column}).}
	\label{fig:err3cont}
\end{figure}

\begin{figure}[h]
	\centerline{
		\includegraphics[width=1.25\linewidth]{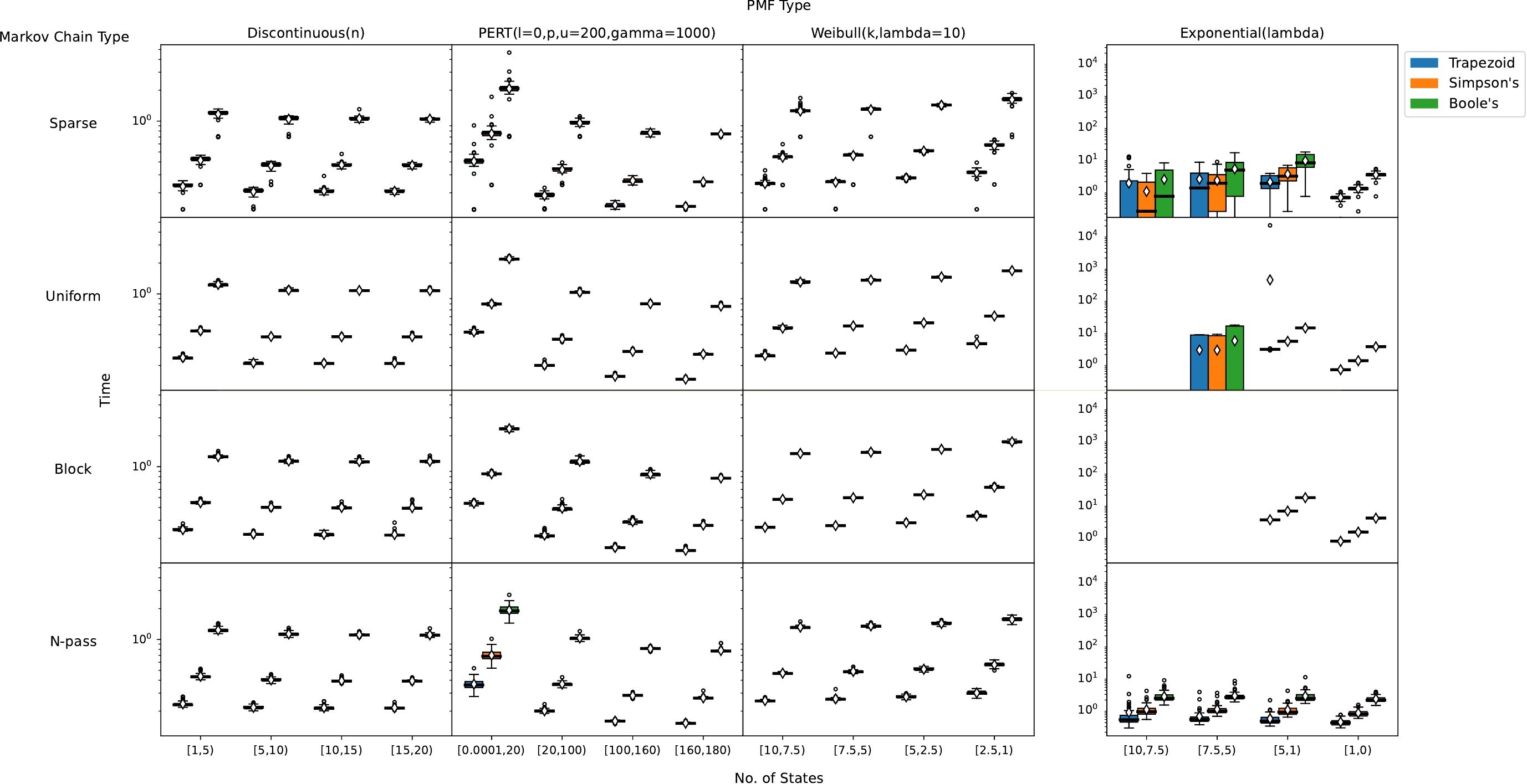}}
	\caption{Box-plots for time taken (secs.). Apart from the Exponential column, the results appear to be relatively close to the mean. For sMRMs of the Exponential type and Uniform MC type, notice the outlier (for the [5,1) param. range) that has pulled up the mean significantly. It is not clear what has caused that result, which took more than 2 hours to complete. Definitions of symbols in the plot are given in Fig. \ref{fig:noiter3exactboxes}.}
	\label{fig:time3cont_boxes}
\end{figure}

\begin{figure}[h]
	\centerline{
		\includegraphics[width=1.25\linewidth]{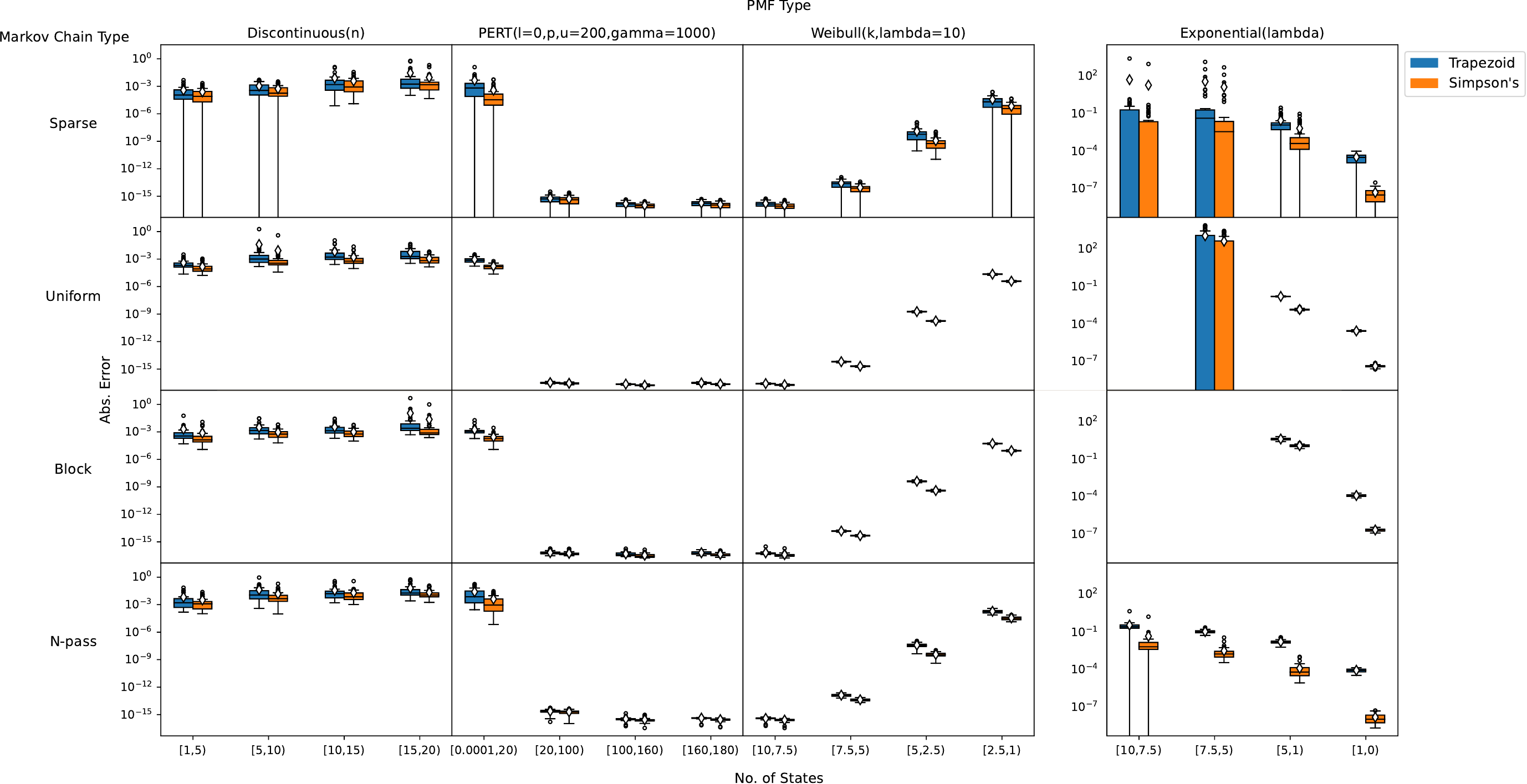}}
	\caption{Box-plots for the (approximate) max abs. error (compared to Boole's). The long (lower) whiskers for the Sparse row is potentially due to Sparse MC's being sampled of which have states that can reach $B$. In this case, the result will just be the zero vector. This issue was discussed previously in Section \ref{sec:exp_MC_pmf}. It is not clear what is causing there to be large interquantile ranges for the Exponential column.  Definitions of symbols in the plot are given in Fig. \ref{fig:noiter3exactboxes}.}
	\label{fig:err3cont_boxes}
\end{figure}

\paragraph{} The average no. of iterations required for the power method with the trapezoid rule to converge, is plotted in Fig. \ref{fig:noiter3cont}. The terms $L1,L2,L3$ in the legend refers to the trapezoid rule with $N,2N-1,4N-2$ points. Then, these rules are combined to form the $N$ point Romberg approximation of level 2 and level 3, i.e. Simpson's and Boole's approximation. See Section \eqref{subsection:Romberg's method} for details. In Fig. \ref{fig:time3cont}, we plot the average time the methods take to solve the problems and the box-plots for time are in Fig. \ref{fig:time3cont_boxes}. The time taken includes the time to build the relevant hypermatrices, but excludes the time it takes to combine the trapezoid rules to form the higher level approximations. The average max absolute approximate error is shown in Fig. \ref{fig:err3cont} with box-plots for the error in Fig. \ref{fig:err3cont_boxes}. Note that for all our results, the Exponential pdf prevented convergence in particular situations, hence some samples were excluded (and we did not re-sample). The no. of samples included are shown in Table \ref{table:plots_expon_column}. We assumed divergence if the max error between two iterations is larger than $1e2$, i.e. if $$max_{s,r}|\textbf{f}^{(n+1)} - \textbf{f}^{(n)}| >  1e2 $$

\begin{table}
	\begin{center}
		\begin{tabular}{ |c|c|c|c|c| }
			\hline
			\multicolumn{5}{|c|}{Exponential Column} \\
			\hline
			& \multicolumn{4}{|c|}{Parameter Range} \\
			\hline
			MC type & [10,7.5) &     [7.5,5) &     [5,1) & [1,0)\\
			\hline
			Sparse & 26 & 38 & 50 & 50\\
			Uniform & 0 & 17 & 50 & 50\\
			Block & 0 & 0 & 50 & 50\\
			$N$-pass & 48 & 50 & 50 & 50\\
			\hline
		\end{tabular}
	\end{center}
	\caption{\label{table:plots_expon_column} No. of samples out of 50, where the iterative method converges, for a given parameter range and MC type, of the Exponential column.}
\end{table}

\paragraph{} Studying the no. of iterations, a consistent hypothesis is that convergence is fastest when either concentration of the pdf's density is close to the zeroth point (of the abscissa). In fact, having a high concentration near zero may be responsible for the lack of convergence seen in the Exponential column (see Table \ref{table:plots_expon_column}). Note that increasing the parameter range introduces more density to the Weibull near zero, but reduces that of PERT and Exponential types. It would appear that setting the zeroth value of each pdf to zero in a sMRM can help ensure convergence. Although, it is false to state that it guarantees convergence. A simple counterexample consisting of mixtures of extremely tall uniformly distributed pdfs (defined strictly on the lattice of points for quadrature) can be used to demonstrate this.

\paragraph{} Next, the error for the Exponential column improves as the parameter range increases, whilst the Weibull increases. A hypothesis for this may be tied to the fact that as its parameter range increases, the Weibull increases in density around zero, whilst for the Exponential it decreases. However, the theory is quite weak in our eyes. Another observation is that the PERT distribution decreases in error as the mean shifts towards $b$. A hypothesis for this may be that having the mean of each reward pdf further away from $b$ implies that the mean of the accumulated reward pdf is more likely to lie outside of the interval $[0,N]$, thus the density is not likely centred in $[0,N]$ suggesting that the values of the pdf there is small, hence easier to obtain accuracy for.

\section{Summary and discussion}
\label{sec:summary_continuous}
In this chapter we presented a new quadrature methodology for convolutions of functions defined over $[0,k]$. The method was based on equidistant points and allowed us to use the FFT once more. We also found a manner in which we could adapt Romberg's method for convolution, and displayed it's capability of improving accuracy further. Additionally, we derived some proofs concerning convergence of particular iterative algorithms, but more work is needed to be done here, e.g. to find a general condition for convergence with the trapezoid rule.

\paragraph{} As stated in the literature review, the benefit of our DFT approach is that it removes the need of algebraically transforming pdfs, or numerically integrating them to derive the transforms, hence simplifying the overall problem. Additionally we stated that the Laplace-Euler algorithm was shown not accurate for discontinuous distributions. For this, perhaps our algorithm can be recommended in this case. We cannot however recommend our algorithm generally (for continuous pdfs for example) without having it first experimentally compared to \cite{bradley2004hypergraph} which used the Laplace transform, with the power method, and the Laplace-Euler inversion algorithm.

\paragraph{} For future work, one may want to consider experimenting with the continuous Fourier transform instead of the Laplace. And using the inverse Fourier transform to derive the pdfs once the system has been solved in the Fourier domain. One reason being that there appears to be some work towards numerically inverting characteristic functions; Fourier transforms of random variables defined also on the negative-real line. The foundation of one direction of work was by Gil-Palaez in \cite{gil1951note}. A unified and multivariate approach is presented in \cite{shephard1991characteristic}. However, more recently the inversion can now be done with FFTs as shown by \cite{witkovsky2016numerical}, based on the earlier work of \cite{hurlimann2013improved}.

\paragraph{} Aside from using point-based methods, there are still other approaches left for experimentation and investigation. For example, one may attempt instead to derive a theory for sMRMs using Fourier extensions and the convolution algorithm presented by \cite{xu2017fast} such that it leads to convergent and unique results. Alternatively, faster function approximation and convolution algorithms relative to the Pacal library \cite{korzen2014pacal} may be attempted instead, along with proving convergence and uniqueness. Alternatively, the complexity within the Pacal library could be relaxed, i.e. the handling of operations including both discrete and continuous random variables together can be removed, and the handling of singularities also.

\paragraph{} Lastly, we are aware of \textit{convolution quadrature} techniques \cite{lubich2004convolution}. Without having delved deep into the theory, these formulas here appear unrelated to our formulas of convolution quadrature. Additionally, the trapezoid rule for convolution quadrature \cite{trefethen2014exponentially} is not necessarily the same one as the trapezoid rule presented here.

\subsection{Drawbacks of the proposed algorithms}
Firstly, we have previously mentioned in Section \ref{subsec:limitations_iterative}, that iterative methods may prematurely converge in particular scenarios. The solution by \cite{haddad2018interval} requires the values of the solutions to be bounded in $[0,1]$. However, pmfs and not pdfs satisfy this property. Whilst the passage-reward densities are not bounded in this interval, the cdf of these densities are. Therefore one could try adapting the algorithm for the cdfs instead to determine convergence more soundly. Alternatively, rather than computing the passage-reward densities, we could instead compute their cdfs directly. Let $Pr(Rew \leq r\ \cap\ \vDash \Diamond B) = F_s(r) = \int_{0}^{r}f_s(x)dx$. Then the system of convolution equations become
$$\textbf{F} = (A \circ \textbf{G})\ \textcircled{$\ast$}\ \textbf{F} + \textbf{H}$$
where \textbf{G} is still a the matrix of reward pdfs, i.e. we have $G_{s,t} = f_{rew(s,t)}$. However $H_s = \sum_{u \in B}P_{s,u}F_{rew(s,u)}$ is a vector of reward (partial) cdfs. Solving this instead may be preferred as it avoids having to numerically integrate the passage-reward densities.

\paragraph{} Secondly, \cite{bradley2004hypergraph}  presented a comprehensive algorithm that can be used for resolving $\mrmpdf$ over an sMRM. The algorithm involves taking the Laplace transform of the system of convolution equations \eqref{eqn:lineq_matrix_trans_pdf}. This leads to a system identical in form to \eqref{eqn:lineqn_matrix_trans}, i.e.
\begin{flalign}
\label{eqn:system_via_laplace}
\textbf{x}(s) =(A\circ \textbf{C}(s))\textbf{x}(s) + \textbf{d}(s)
\end{flalign}
where $s$ is a complex number. The definitions of the terms above are now using the Laplace transform and not the FT, i.e.
\begin{align*}
\textbf{C} &\triangleq  (C_{s,t})_{s,t \in S_?} \triangleq (\mathcal{L}\{f_{(rew(s,t))}\})_{s,t \in S_?}&\\
\textbf{d} &\triangleq  (d_s)_{s \in S_?} \triangleq (\sum_{u \in B}\textbf{P}(s,u)\mathcal{L}\{f_{rew(s,u)}\})_{s \in S_?}
\end{align*}
with the Laplace operator $\mathcal{L}$ applied to a function $f$ being defined as $$\mathcal{L}\{f\}(\sigma + i\omega) = \int_{0}^{\infty}f(t)e^{(\sigma + i\omega)t}dt$$

\par{} Then, the power method is used to resolve each linear system  in \eqref{eqn:system_via_laplace} for a set of $N$ sampled points $(s_i)_{0<i\leq N-1}$. Once the power method converges (for each $s_i$), either the Laplace-Euler or the Laplace-Laguerre inversion algorithms is applied to $\textbf{x}(s)$ (described in \cite{bradley2004hypergraph}) to obtain the passage-reward density $\textbf{f}(r)$ at $M$ points, i.e. for each $r \in Q \subset [a,b]$, where $|Q| = M$. For the Laplace-Euler inversion algorithm (their general strategy), then $M=N/c$ for a choice of $c\in\{15,16,\dots,50\}$, which determines the accuracy of the inversion. It is expected that as $c$ increases, the accuracy increases. For the Laguerre inversion algorithm, which is less applicable, $N$ is fixed to 400, and $M$ is independent of $N$.

\paragraph{} This technique has perhaps two main advantages over our algorithms. Using their algorithm, the property $\mrmpdf$ can be resolved for the set of points $r \in Q \subset [a,b]$. Our algorithm enforces that $Q$ be a set of equidistant points between $[0,b]$ only. Secondly, for strictly continuous pdfs, their algorithm appears to be quite efficient in obtaining high accuracy for the $M$ points for the pdfs, relative to $N$, the number points sampled from the Laplace transforms. However, a disadvantage of their method could be that it requires the analytical Laplace transforms of pdfs to be known, whereas ours do not. Nevertheless, these transforms may be approximated numerically or obtained via a CAS (computer algebra system). If the approximations are poor, then like our algorithm, the power method may fail to converge. Disregarding the accuracy of solutions, the worst-case time complexities for both ours and their power method are similar.

    \chapter{Conclusion}     
    
    We began our thesis presenting the theoretical derivations for reachability problems over a new class of models (in the model-checking literature) - stochastic Markov reward models. We were also able to find practical algorithms for model-checking reachability properties via the DFT for the univariate reward problem. This solution extends to the case where our rewards are random vectors, with mutually independent components.
    
    \paragraph{} We found in the discrete-reward setting that the best algorithm (relative to our experiments) turns out to just be the naive exact-power method. The direct algorithms are quite slow, and our experiments show that the power method outperforms its approximate counterpart, and also the approximate direct method (i.e. the LU approximation method). The other iterative algorithms were slower due to our deconvolution algorithm not being optimized (via the FFT). Additionally, it was harder to vectorize the computations of the Gauss-Seidel algorithm to take effective use of the $numpy$ library. Note however, we did not experiment with approximate variants of the other iterative methods. But they suffer from the same problem as their exact counterparts. Nevertheless, at this time we recommend the exact-power method for use.
    
    \paragraph{} We stated earlier that some optimizations possible would be: 1) That zero-padding can be removed almost entirely by the work of \cite{robertson1992computation}, and 2) that kernel tiling algorithms can be used for improving speed \cite{pavel2013algorithms}.  However, it is worth considering whether algebraic decision diagrams (ADD) \cite{bahar1997algebric} can be used in the discrete-reward setting for reducing space complexity. Current model-checking tools such as \cite{kwiatkowska2002prism} employ decision diagrams such as the multi-terminal binary decision diagrams (MTBDD) to help circumvent state space explosion. The paper \cite{bahar1997algebric} presents an algorithm for matrix multiplication when matrices are represented as ADDs, and the exact-power method for example requires not much more than that in terms of unique operations.
    
    \paragraph{} In the continuous-reward setting, we showed that our algorithm via the DFT can handle discontinuous distributions. Unfortunately, many samples were required to ensure sufficient accuracy. If improving on our approach, one may consider replacing Romberg's method with another algorithm of a similar form, e.g. \cite{youngberg2012alternative}.  Additionally, we have stated that there exists different equations for deriving the first-passage reward densities \cite{kao1974modeling,warr2012introduction}. These approaches may yield better convergence properties.
   
    \paragraph{} Whilst our solution works for continuous pdfs and of which required fewer samples for sufficient accuracy, we have not compared our results with the algorithm of \cite{bradley2004hypergraph} using the Laplace transform and the Laplace-Euler inversion algorithm, with the power method yet due to lack of time. We have, however, discussed their effectiveness in the literature review. Thus, we cannot recommend our algorithm over theirs for problems with continuous-pdfs yet.
 
    \section{Future work}

    \paragraph{Negative rewards} The sMRM affords random variables being negative, or having positive measure over the negative portion of the real line. It is not clear how the DFT can be used since it leads to representation growth of the solution vector if applied naively. Alternatively if the continuous Fourier Transform is employed, no representation growth occurs. Then, future work is in determining both highly precise numerical methods for the transform itself (hence not relying on algebraic techniques) and its inversion. Numerical algorithms for the transform include those that use the FFT, for example \cite{witkovsky2016numerical}, based on the work of \cite{hurlimann2013improved}.
    
    \paragraph{Non-determinism} It is possible to extend sMRMs with \textit{actions} creating stochastic Markov reward Decision Processes. That is, at any state, the process makes a non-deterministic choice between a set of actions, each of which has a transition probability distribution over the set of states of the process.  The goal here would be to design efficient algorithms for problems over this class of models. Such problems may involve determining optimal policies or determining the maximum/minimum probability of reaching a set of states $B$.

    \bibliographystyle{abbrv}
    \bibliography{sums_rv_dense_corrections.bib}
    
        \appendix
        \appendixpage
        \chapter{Supplementary Proofs}
    	
    	\section{Proof of Theorem \ref{theorem:gauss_solution_set}}
        \begin{proposition}  	\label{theorem:gauss_solution_set_full} Given a system of convolution equations $\textbf{E}\textcircled{$\ast$}\textbf{f} = \textbf{h}$, with equations $E_0,E_2,\cdots E_N$, and the unique solution $\textbf{x}$ (a hypervector), where each $E_i$ represents the equation $$\sum_j E_{i,j} \ast f_j = h_i$$
        	If the system is transformed via these two updating/substitution operations:
        	
        	\begin{enumerate}
        		\item replacing an equation by a convolution of itself with a non-zero vector.
        		\item replacing an equation with the equation itself added to another equation that has been convolved with a non-zero vector.
        	\end{enumerate}
        	or these two pivoting operations:
        	\begin{enumerate}
        		\item swapping one equation with another (i.e. swapping rows of \textbf{E}).
        		\item changing the order of terms in each and every equation (i.e. swapping columns of \textbf{E}).
        	\end{enumerate}    	
        	then the transformed system will still have the same solution set \textbf{x}, with the exception of the case where columns are swapped. In this latter situation, the transformed system will have a solution set that is a reordering of the terms of \textbf{x}.

        	\paragraph{} \textit{Proof: } Let $E_i(\textbf{x})$ be the equation $E_i$ where the unknown terms \textbf{f} are replaced with the solution \textbf{x}. Then the conjunction 
        	$E_0(\textbf{x}) \cap E_1(\textbf{x}) \cap \cdots \cap E_N(\textbf{x})$ is true. Concerning the substitution rules, note that convolving any equation $E_i$ both sides by any non-zero vector $g$ will also maintain the truth of the conjunction, e.g. $E_0(\textbf{x}) \cap \cdots \cap\ g\ast E_i(\textbf{x})\ \cap \cdots \cap E_N(\textbf{x})$ holds. And secondly, replacing any equation by adding the equation itself to another equation convolved with a non-zero vector $g$ will also maintain the conjunction, e.g. $E_0(\textbf{x}) \cap \cdots \cap\ (g\ast E_j(\textbf{x})) + E_i(\textbf{x})\ \cap \cdots \cap E_N(\textbf{x})$.
        	
        	\paragraph{} Concerning pivoting rules, we can re-arrange the order of the equations and the conjunction will remain true, e.g. swapping the first two rows, we have $E_1(\textbf{x}) \cap E_0(\textbf{x}) \cap \cdots \cap E_N(\textbf{x})$. Lastly, if we swap the order of the terms in \textbf{f}, e.g. the term $f_j$ with $f_k$, then let $E'_i$ denote equation $i$ after the terms $E_{i,j}$ and $E_{i,k}$ are swapped. Let $\textbf{x}$ have its indices swapped in the same fashion as \textbf{f} and denote this $\textbf{x}'$. Then the conjunction $E_1'(\textbf{x}') \cap E_1'(\textbf{x}') \cap \cdots \cap E_{N}'(\textbf{x}')$ holds. \qed

        \end{proposition} 
        
         \section{Proof of Theorem \ref{theorem:fixed_point}}
        \begin{theorem}[Fixed point characterization (for pmfs)] 
        	\label{theorem:fixed_point_full}
        	The cumulated reward hypervector $\textbf{f} = (\textbf{f}[r])_{r \in \mathbb{N}} = (f_{s \vDash \Diamond B}[r])_{r \in \mathbb{N},s \in S_?} =  ( Pr(r \cap s \vDash \Diamond B) )_{r \in \mathbb{N},s \in S_?}$ is the (unique) fixed point of the operator $\Upsilon: [0,1]^{\mathbb{N} \times  S_?\times 1} \rightarrow [0,1]^{\mathbb{N} \times S_? \times 1}$. This operator is defined as:
        	\[\Upsilon(\textbf{f}) = ((A\circ\textbf{G})\ \textcircled{$\ast$} \ {\textbf{f}}) + \textbf{h} \]
        	
        	\paragraph{} Additionally, let $\textbf{f}^{(0)} = \textbf{0}$, and $\textbf{f}^{(n+1)} = \Upsilon(\textbf{f}^{(n)})$ where $n \geq 0$. Then, for any ${r \in \mathbb{N}}$ and  $s \in S_?$ the following three statements hold:
        	
        	\begin{enumerate}
        		\item $f^{(n)}_s[r] = Pr(\ r\  \cap [s \vDash \Diamond^{\leq n} B])$, for all $n \geq 0$.
        		\item $ \lim\limits_{n \rightarrow \infty} f_s^{(n)}[r] = f_s[r]$.
        		\item $f_s^{(0)}[r] \leq  f_s^{(1)}[r] \leq  f_s^{(2)}[r]  \leq \cdots \leq  f_s[r]$.
        	\end{enumerate}
        	where statement 2. states that the solution converges to a fixed point, and 3. states that the convergence is monotonic.
        	
        	\paragraph{} \textit{Proof of 1.:} 
        	\label{proof:time_bounded_reachability_sMRMs}
        	The proof involves induction on $n$. Firstly with the base case $n=0$, $f_s^{(0)}[r] = \textbf{0} \equiv Pr(\ r\  \cap [s \vDash \Diamond^{\leq 0} B]) = Pr(\ r\  \cap [s \vDash B]) = \textbf{0}$ for all $s \in S_?$ (since $S_? \subseteq S - B$). This proves the base case. Now to prove the inductive step. Firstly, note that the event $[s \vDash \Diamond^{(\leq n+1)} B]$ is equal to
        	$$[s \vDash \bigcirc B] \cup
        	[s \vDash \bigcirc S_?  \ \cap \ s \vDash \Diamond^{(\leq n+1)} B ]
        	$$
        	i.e. the event of reaching $B$ under or equal to $n+1$ steps, is equal to the event of reaching $B$ in the next step OR the event of reaching $B$ under $n+1$ steps on the condition of passing through $S_?$ at the next step. Note that $[s \vDash \bigcirc B]$ is disjoint to $   	 [s \vDash \bigcirc S_? ]$, hence
        	$$Pr([s \vDash \Diamond^{(\leq n+1)} B]) = Pr([s \vDash \bigcirc B]) + Pr(   	 [s \vDash \bigcirc S_?  \ \cap \ s \vDash \Diamond^{(\leq n+1)} B ])$$	
        	\paragraph{} Now, starting the proof with a top-down approach, we have   	 
        	\begin{align*}
        	f_s^{(n+1)}[r]
        	&=     Pr( r \cap   [s \vDash \Diamond^{(\leq n+1)} B]) & \\
        	&= Pr(r \cap [s \vDash \bigcirc B]) +  Pr(r \cap    	 [s \vDash \bigcirc S_?  \ \cap \ s \vDash \Diamond^{(\leq n+1)} B ])
        	\end{align*}  	 
        	where $Pr(r \cap [s \vDash \bigcirc B]) =  \sum_{u \in B}\textbf{P}(s,u)f_{rew(s,u)}[r]
        	= h_s[r]$.
        	Therefore this corresponds to the $\textbf{h}$ portion of the $\Upsilon()$ operator. Next, by letting  $\Pi.s$ denote the set of all paths beginning in $s$, satisfying $[s \vDash \bigcirc S_?  \ \cap \ s \vDash \Diamond^{(\leq n+1)} B ]$, we have    
        	\begin{flalign*}
        	&Pr(r \cap    	 [s \vDash \bigcirc S_?  \ \cap \ s \vDash \Diamond^{(\leq n+1)} B ]) &\\
        	&= \sum_{\hat{\pi} \in \Pi.s}Pr(\hat{\pi})f_{Rew(\hat{\pi})}[r] & \\
        	&= \sum_{t \in S_?}\textbf{P}(s,t)(f_{rew(s,t)} \ast (\sum_{\hat{\phi} \in \Pi.t}
        	(Pr(\hat{\phi})f_{rew(\hat{\phi})})))[r] & \\
        	&= \sum_{t \in S_?} A_{s,t}(G_{s,t}	\ast f_t^{(n)})[r]  
        	\end{flalign*}  	 
        	where $f_t^{(n)}[k] = Pr( k \cap  [s \vDash \Diamond^{\leq n} B])$; the induction hypothesis which we have assumed correct. Note that the result corresponds to the $A\circ\textbf{G}$ portion of the $\Upsilon()$ operator. Thus we have that    		 
        	\begin{flalign*}
        	f_s^{(n+1)}[r]
        	&= Pr( r \cap  [s \vDash \Diamond^{\leq n+1} B]) &\\
        	&=  Pr(r \cap [s \vDash \bigcirc S_?  \ \cap \ s \vDash \Diamond^{(\leq n+1)} B ]) + Pr(r \cap [s \vDash \bigcirc B])  & \\    
        	&= \sum_{t \in S_?}((A_{s,t}G_{s,t}) \ast {f}_{t}^{(n)})[r] + h_s[r]
        	\end{flalign*}	
        	which proves the equivalence. 
        	
        	\paragraph{} \textit{Proof of 2.:} We know that the event $[s \vDash \Diamond B]$ is equal to the event $[\lim_{n \rightarrow \infty} s \vDash \Diamond^{\leq n} B]$. Thus
        	\begin{flalign*}
        	f_s[r] = Pr(r \cap s \vDash \Diamond B) 
        	= \lim_{n \rightarrow \infty} Pr(r \cap s  \vDash \Diamond^{\leq n} B) 
        	= \lim\limits_{n \rightarrow \infty} f_s^{(n)}[r]
        	\end{flalign*}
        	\paragraph{} \textit{Proof of 3.:} It is true that $[s \vDash \Diamond^{\leq n} B] \subseteq [s \vDash \Diamond^{\leq (n+1)} B]$. Then for any $r \in \mathbb{N}$,  
        	$$[Rew = r\  \cap [s \vDash \Diamond^{\leq n} B]] \subseteq [Rew = r\  \cap [s \vDash \Diamond^{\leq n+1} B]]$$
        	\paragraph{} From statement 1. we know that $f^{(n)}_s[r] = Pr(\ r\  \cap [s \vDash \Diamond^{\leq n} B])$. Thus, 
        	$$f^{(n)}_s[r] = Pr(\ r\  \cap [s \vDash \Diamond^{\leq n} B]) \leq Pr(\ r\  \cap [s \vDash \Diamond^{\leq n+1} B]) = f^{(n+1)}_s[r]$$ 
        	
        \end{theorem} \qed
    
    	\chapter{Supplementary Results }
    	
    	\paragraph{} We present here additional experiments that we have performed, but deemed to not be as useful as the results in the main thesis. For the empirical results in this chapter, unless stated otherwise, only the \textit{uniform} Markov chain was used (see Fig. \ref{fig:sampled_MCs}). That is, given a fixed size empty matrix, we sampled values uniformly between zero and one for each and normalized the matrix to ensure the rows would sum to one. This resulting matrix would then be a valid probability matrix. Uniform pmfs were sampled in a similar fashion.
        
        \section{Experiments with exact and approximate power method}
        	
        \subsection{Previous toy problem}
        We solve the same problem earlier from Section \ref{problem:5statesMRM}, this time including solutions from the power method - both exact and approximate. For the power method, we use an absolute tolerance level of $1e-16$ for convergence (\eqref{eqn:stopping_criteria}). The property being solved for is $Pr(r\ \cap\ s \vDash \Diamond s_4)$, for all $s \in S_? = \{s_0, s_1, s_2, s_3\}$ and $r = 0,1,2,\cdots,N-1$, where $N= 150$. The Gaussian elimination algorithm, and the exact power method computes $Pr(r\ \cap\ s \vDash \Diamond s_4)$ for $r = 0,1,2,\cdots,N-1$, whilst the approximate LU method, and approximate power method computes the property for $r= 0,1,\cdots T$, where $T = N + n$, and $n$ is the zero padding length (needed for the approximate methods).
        
\begin{figure}
        	\centerline{
        		\includegraphics[width=1\linewidth]{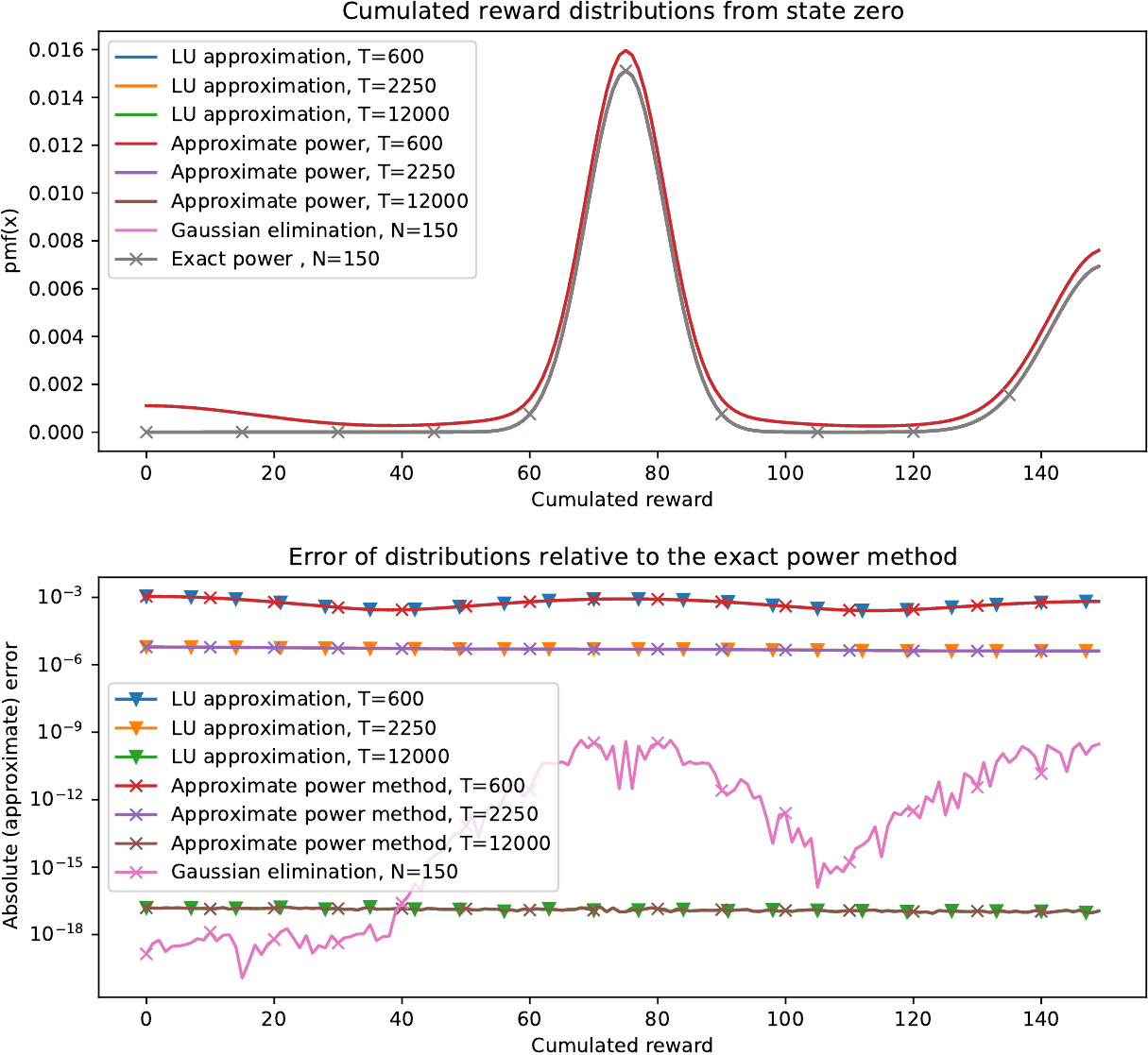}}
        	\caption{Similar to Fig. \ref{fig:gevslu}. We also include results for the power method. Note in the bottom plot that the approximate (direct and iterative) methods for a given $T$ overlap one another in error. For the LU approximations, the different padding lengths can be determined by $T-N$ where $T$ stands for Total length.}
        	\label{fig:powergeapproxexact}
        \end{figure}

        \paragraph{} The results are shown in Fig. \ref{fig:powergeapproxexact}. We find that the approximate power methods align with the LU approximate method of equal $T$. There appears to be a discrepancy between them when $T=12000$ possibly due to numerical errors. In fact, the approximate power method could not converge under 1000 iterations - a hard-limit for the number of iterations we enforced. Increasing the limit to 10,000 iterations, the algorithm still did not converge (see Table \ref{fig:powergeapproxexacttimes-crop}).
        
        \paragraph{} The plot on the bottom of Fig. \ref{fig:powergeapproxexact} is the absolute (approximate) error relative to the \textit{exact power method}, and not Gaussian elimination as previously. We find that Gaussian elimination is not as accurate as the power method. Supporting this is the fact that the approximate methods' accuracy nears that of the power method for large $T$. R This may be attributed to rounding errors, and methods to mitigate such problems (for linear systems) are presented in \cite{dahlquist1974numerical}. More specifically, it could be due to the numerical instability of deconvolution, which is less known. Nevertheless, if accuracy is important, in practice the results of Gaussian elimination can be used as an initial guess for an iterative method like the power method, since solutions can be made more precise via successive iterations.
        
        \begin{table}[H]
        	\centering
        	\includegraphics[width=1\linewidth]{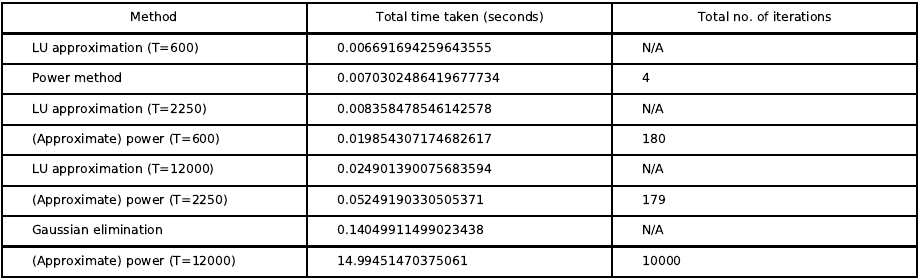}
        	\caption{ The approximate power method (with T=12000) did not converge (reaching the maximum iteration limit of 10000), probably due to the really small threshold. The LU approximation (with T=12000) method is significantly faster than Gaussian elimination and also achieved greater accuracy rel. to the power method (shown in Fig. \ref{fig:powergeapproxexact}). Although, it did require significantly more space ($\times \frac{12000}{150}$).}
        	\label{fig:powergeapproxexacttimes-crop}
        \end{table}
        
        \paragraph{} We see in Table \ref{fig:powergeapproxexacttimes-crop} that the exact power method required fewer iterations relative to the approximate power method to reach convergence. For $T=12000$, whilst the approximate power method failed to converge under 10000 iterations, Fig. \ref{fig:powergeapproxexact} suggests (at least for state $s_0$) that the hypervector has converged seeing that the absolute error relative to the power method is near 1e-16.
        
        \paragraph{} In this one experiment, in comparison to the exact power method, the approximate power method was slower. Gaussian elimination was non-competitive for either speed or accuracy. The LU approximation would not be considered competitive if high speed and accuracy are required. However from a more general perspective, if we know the characteristics of an sMRM problem beforehand, we may select the direct algorithms if for example it is known that the iterative algorithms will fail to converge quickly, and then apply the iterative methods on top of their results to improve their accuracies. However if $|S_?|$ is large, both the LU approximation method and GE have time complexities which are cubic in $|S_?|$. Therefore, from a scalability perspective the iterative algorithms may be considered more scalable with respect to $|S_?|$, and therefore preferred. However, whilst we cannot offer a general recommendation for which algorithm to use (without a problem context), our current preferred choice would be to use the iterative algorithms as a first.
        
        \paragraph{}  We now proceed to evaluate the usefulness of the approximate power method relative to the exact power method. The following series of tests show that it is not necessarily faster at convergence.
        
        \subsection{Convergence performance of exact vs approximate power method} 
        \paragraph{} We investigated the performance of the \textit{approximate} power method and the \textit{exact} power method over random systems, and timed how long they took to solve for the property $\mrmpdf$. These experiments should help us determine whether the exact power method's requirement for transforms and inversions of the solution hypervector (for each iteration) contributes significantly to solving time. We found three factors that determine the solving time: 1) $|S_?|$, the size of the set of states that can eventually reach $B$ our goal states. 2) $k$, such that $\mrmpdf$ is computed for $r = 0,1,\cdots,k-1$. 3) $\lambda$, which defines the points $r= 0,1,\cdots,\lfloor\lambda k\rfloor - 1$ for which $f_{rew(s,t)}[r]$ is allowed to be non-zero in, for all $s,t \in S^2$. Note that $\sum_{x=0}^{\lfloor\lambda k\rfloor - 1}f_{rew(s,t)}[x] = 1$, and $f_{rew(s,t)}[r] = 0$ for $r >\lfloor\lambda k\rfloor - 1$. Therefore each $f_{rew(s,t)}$ is still a pmf.
        
        \paragraph{} We set the approximate method's padding length to $k-1$, which makes it identical to the exact power method in terms of space requirements. Therefore the main difference between the two methods is that the exact power method requires the DFT and the inverse DFT to be applied in each iteration, whereas the approximate method does not.
        
        \paragraph{} For the exact power method, we used the following termination criteria $$max_{s,r}|\textbf{f}^{(n+1)} - \textbf{f}^{(n)}| \leq \epsilon$$ where $\epsilon = 1e-7$. For the approximate power method (Sec.  \ref{sec:approx_power}), the criteria is the same, $$max_{s,r}|\textbf{x}^{(n+1)} - \textbf{x}^{(n)}| \leq \epsilon$$
        where it is to be noted that $\textbf{x}^{(n)}$ is a complex number.
        
        \paragraph{Experiment set-up} For our experiment, we varied each of the factors above, whilst keeping others fixed. Firstly, we will set as fixed $k = 1000,|S_?|=  30$, and $\lambda = 0.5$. Then we will vary each parameter separately to get a measure of their contribution to time. For any valuation of the parameters ($k,|S_?|,\lambda$), 200 experiments are performed, their average and worst times recorded and plotted. For each of these 200 experiments, since we are solving for the property $\mrmpdf$, we randomly and uniformly generated the terms in the system
        $$\textbf{f} = (A\circ \textbf{G})\textcircled{$\ast$} \textbf{f} + \textbf{h}$$
        i.e. $A,\textbf{G},\textbf{h}$ are randomly generated. The pmfs $(G_{s,t}[r])_{s,t \in S_?}$ and $(h_s[r])_{s \in S_?}$ are generated to satisfy $\lambda$.
        
        \paragraph{} For each experiment, we only recorded the solving time, i.e. time until convergence (solving phase). We did not record the time required to build the relevant system of equations (preparation phase) (e.g. $\textbf{G}, \textbf{h}$) since it is the same for both methods.
        
\begin{figure}[H]
        	\centerline{
        		\includegraphics[width=1\linewidth]{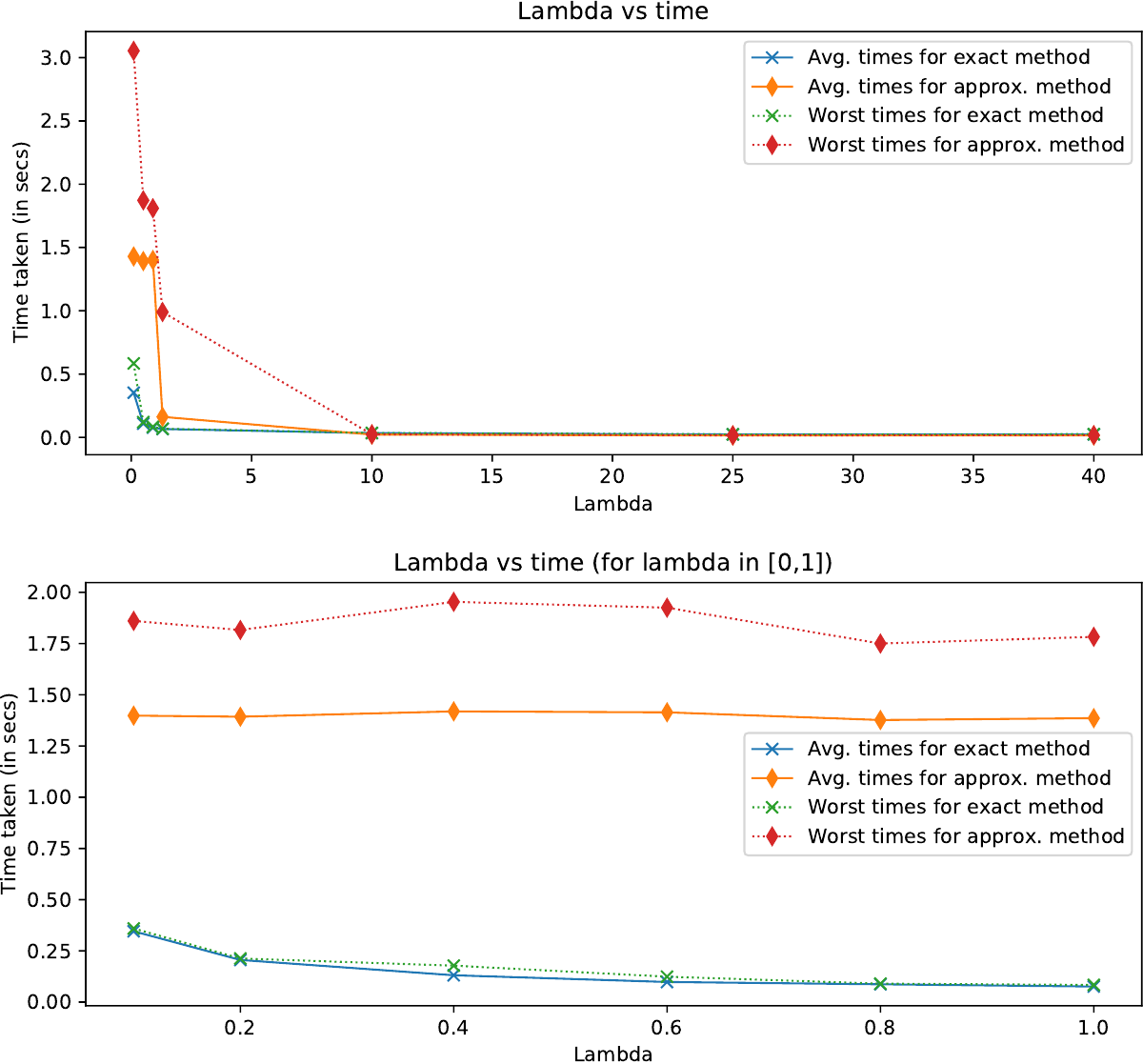}}
        	\caption{Experiments performed with $ k = 1000, |S_?| = 30$ and $\lambda$ was varied. The top graph was achieved using $7\times200$ experiments, whilst the bottom took $6\times200$ experiments. Each experiment consisted of a different randomly generated system. For the top graph, the values of lambda were $[0.1, 0.5, 0.9, 1.3, 10, 25, 40]$. For the bottom, the values were $[0.1, 0.2, 0.4, 0.6, 0.8, 1.0]$.}
        	\label{fig:convergence_rate_len_vs_j-crop}
        \end{figure}

        \paragraph{} Our first experiment (Fig. \ref{fig:convergence_rate_len_vs_j-crop}) is to vary $\lambda$. We find generally that for $\lambda \in [0,1]$, the exact power method is faster at terminating. However for $\lambda > 1$, the approximate method increases in competitiveness.
        
\begin{figure}[H]
        	\centering
        	\includegraphics[width=1\linewidth]{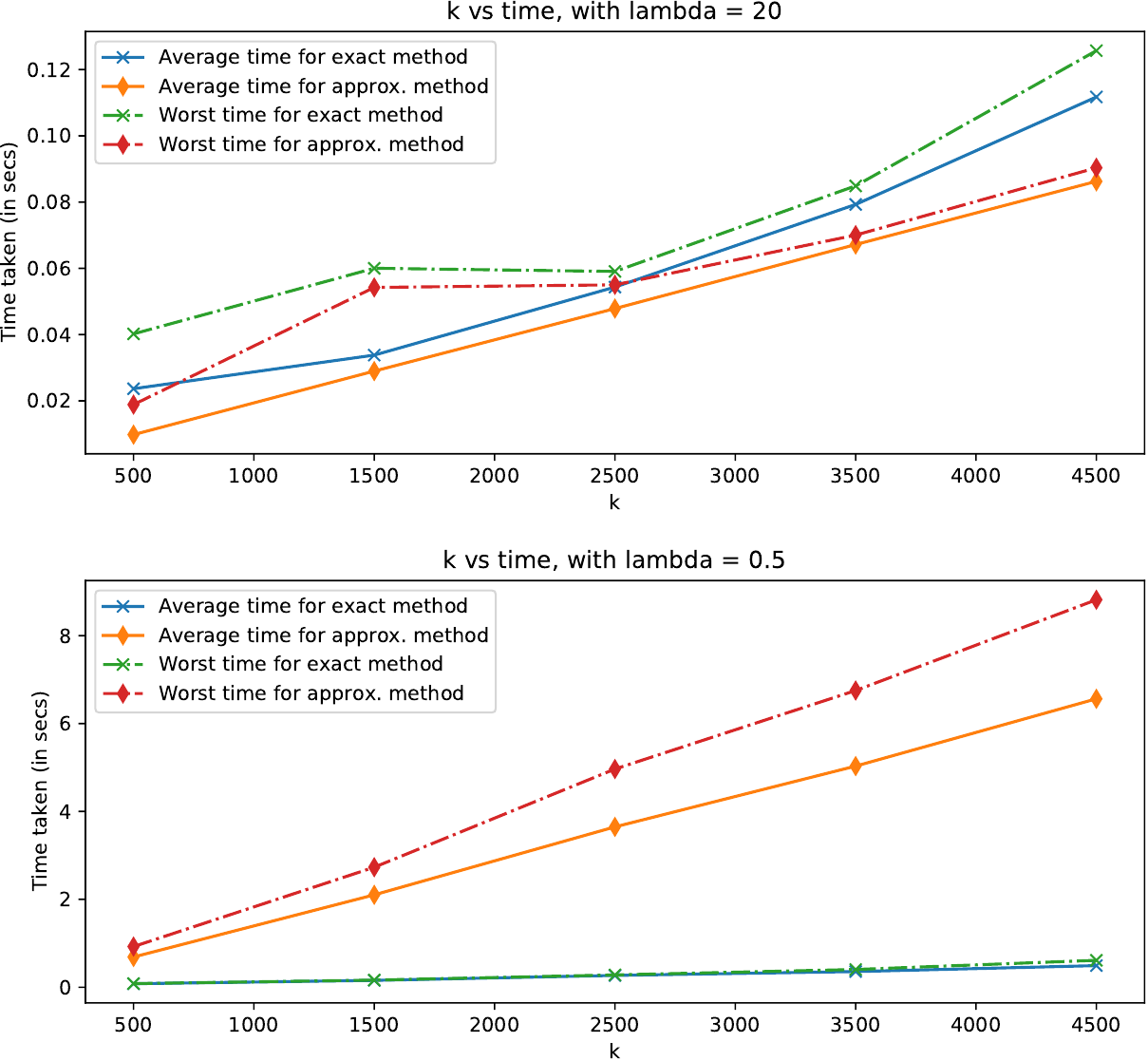}
        	\caption{Experiments performed with $|S_?| = 30$ and $k$ was varied. The top and bottom graphs were achieved with 1000 experiments each.}
        	\label{fig:convergenceratekvstime-crop}
        \end{figure}

        \paragraph{} Next, we experimented with $k$. We see in Figure \ref{fig:convergenceratekvstime-crop} that when $\lambda = 20$, the time taken for each method was generally small, which corresponds with our previous results shown in Fig. \ref{fig:convergence_rate_len_vs_j-crop}. With $\lambda = 0.5$, we find there to be a significant growth in time requirements for the approximate method.

\begin{figure}[H]
        	\centering
        	\includegraphics[width=1\linewidth]{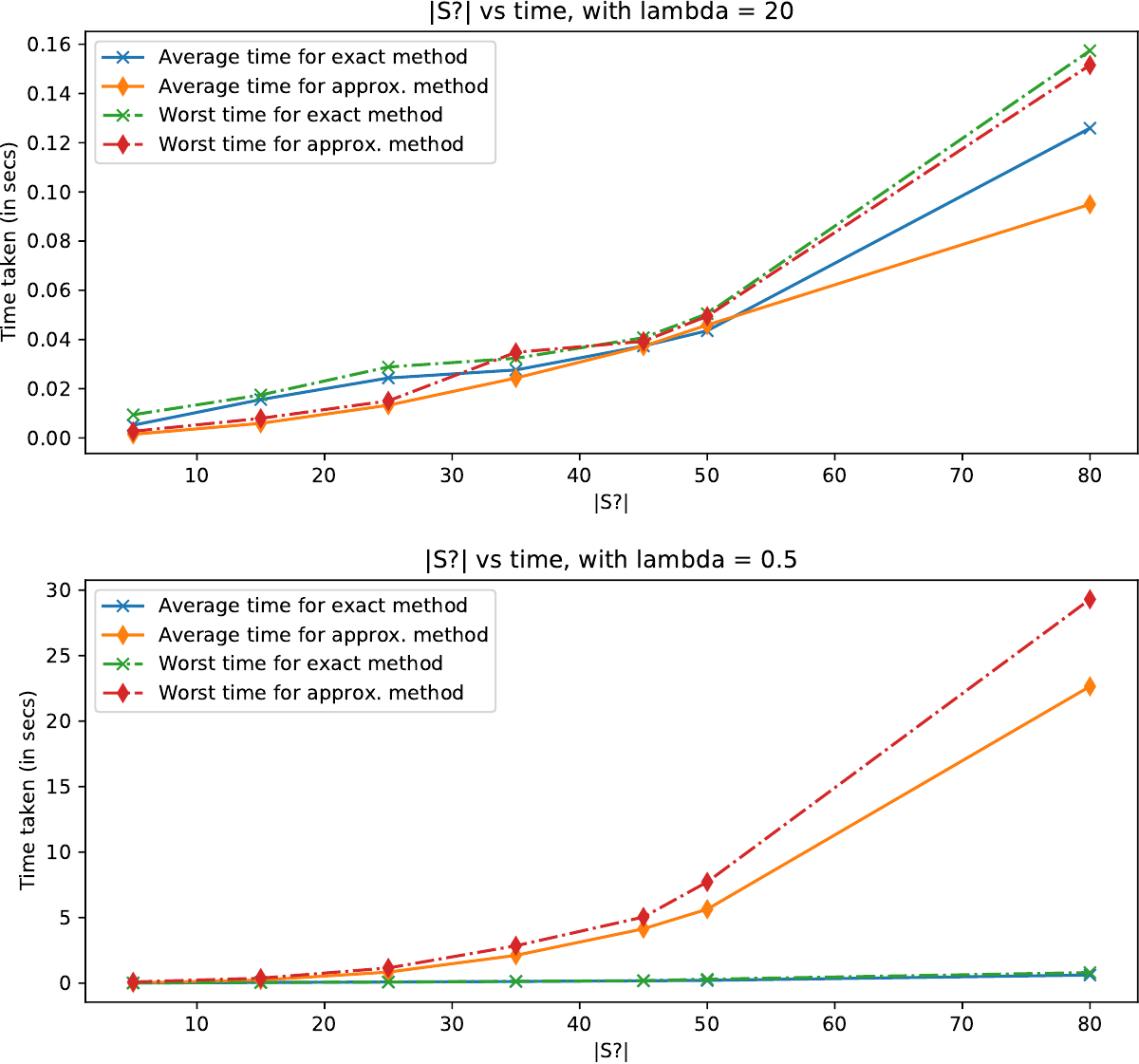}
        	\caption{Experiments performed with $k = 1000$ and $S_?$ was varied. The top and bottom graphs were achieved with 1400 experiments each.}
        	\label{fig:convergenceratessizevstime-crop}
        \end{figure}

        \paragraph{}Finally, we investigated the complexity of the system when increasing $S_?$  (Fig. \ref{fig:convergenceratessizevstime-crop}). In this setting the exact power method seems to prevail for $\lambda=20$, but we find these methods perform similarly when $\lambda=0.5$.
        
        \paragraph{} To conclude our findings, there does not currently appear to be much benefit in using the approximate power method for time savings alone. If further testing the effectiveness of the methods, the pmfs $(f_{rew(s,t)}[r])_{s,t \in S^2}$ could be sampled instead from known tailed random variables that vary in their tails' strength: light, medium and heavy. Our pmfs were randomly generated, where we sampled values between $[0,1]$ for each $f_{rew(s,t)}[r]$ and then normalized  to sum to one. 
        
        \section{Experiments with exact iterative methods} 
        
        \paragraph{Implementation details} For the following experiments, any conclusions derived for time should take into consideration that the different methods used do not have fully optimized implementations. For example, any deconvolutions performed was computed via an algorithm with $O(k^2)$ time complexity (where $k$ is the length of the resulting vector) whilst there are $O(klog_2k)$ algorithms. Another reason would be that the \textit{numpy} package in python was used to implement a portion of the iterative methods. However, its usage was not equal between the three exact methods: power, Jacobi, and Gauss-Seidel. The package would enable us to speed up the methods significantly. Hence, the fairest indicator of performance for the iterative methods would be the number of iterations required to solve a problem.

        \subsection{Scalability of exact iterative methods}
        
        \paragraph{} Our next experiment is almost identical to that described in Section \ref{subsection:scalability_GE_vs_LU}. It differs in that it includes results from the exact Jacobi and Gauss-Seidel methods.
        
        \paragraph{} As before, $|S_?|$ and $k$ are varied independently. Now, performance is characterized by time and number of iterations. When varying $|S_?|$, $k$ was fixed to $1501$. And when varying $k$, $|S_?|$ was fixed to ten. For each valuation of the parameters ($|S_?|,k$), 200 experiments are performed. For each experiment, the system of convolution equations to be solved is sampled as before, their average and worst solving times (for each method) is recorded and plotted. Additionally, we now also record the number of iterations taken by the iterative methods to solve the problem. The convergence threshold set once more to $1e$-16. 
        
\begin{figure}[H]
        	\centering
        	\includegraphics[width=1\linewidth]{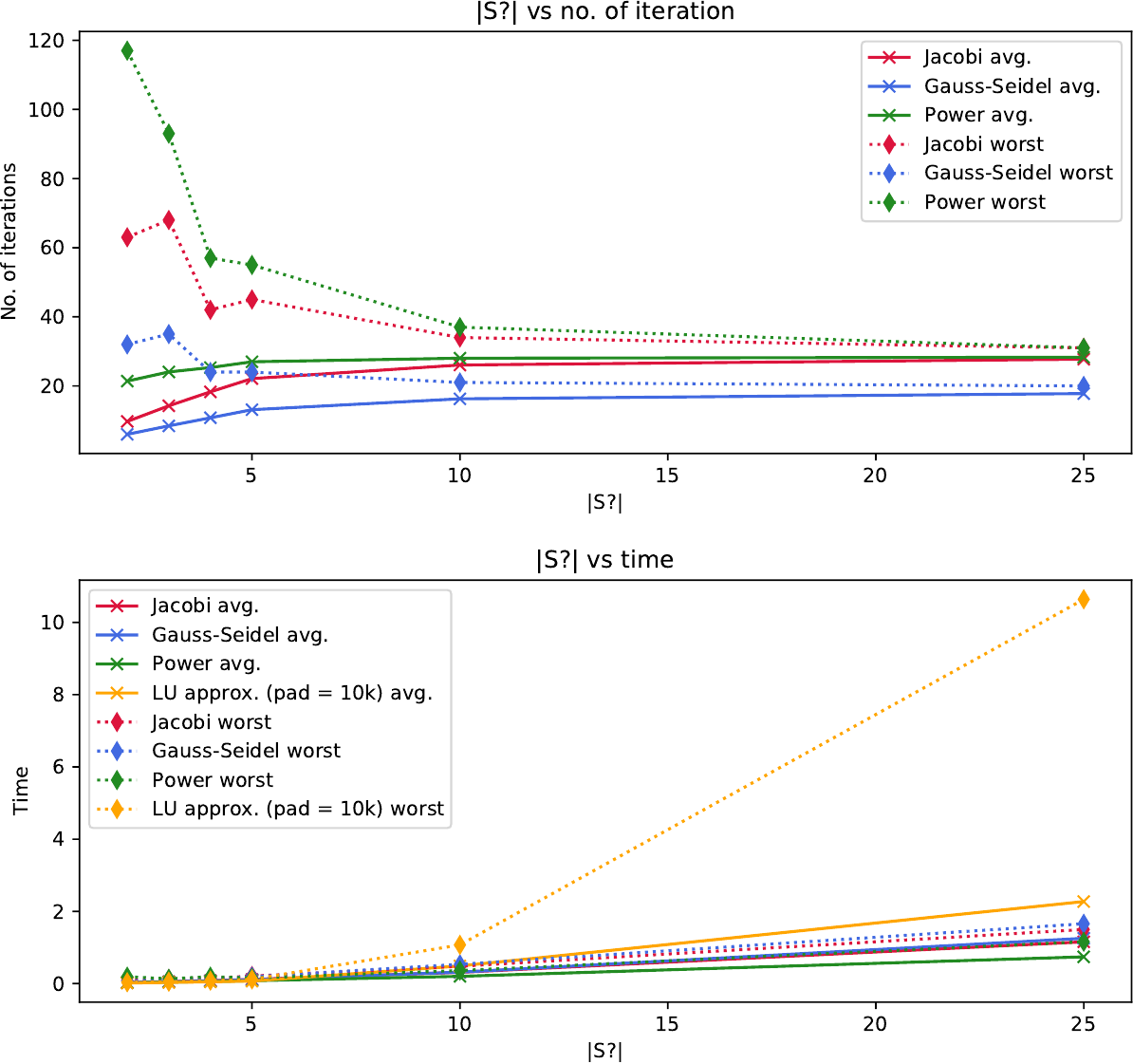}
        	\caption{The LU approximation algorithm requires no iterations, therefore it does not have a plot in the top graph. As $|S_?|$ increases, there is a slight growth in the average no. of iterations for each method. Additionally, the discrepancy between the worst and average no. of iterations decreases. In the legend, (pad = 10k) means that the LU approx. method uses a padding length (see \eqref{eqn:pi_widehat}) of $10k$. Each graph (top and bottom) was achieved after performing ($200 \times 6$) experiments.}
        	\label{fig:scalabilityiterativemethodss-crop}
        \end{figure}

\begin{figure}[H]
        	\centering
        	\includegraphics[width=1\linewidth]{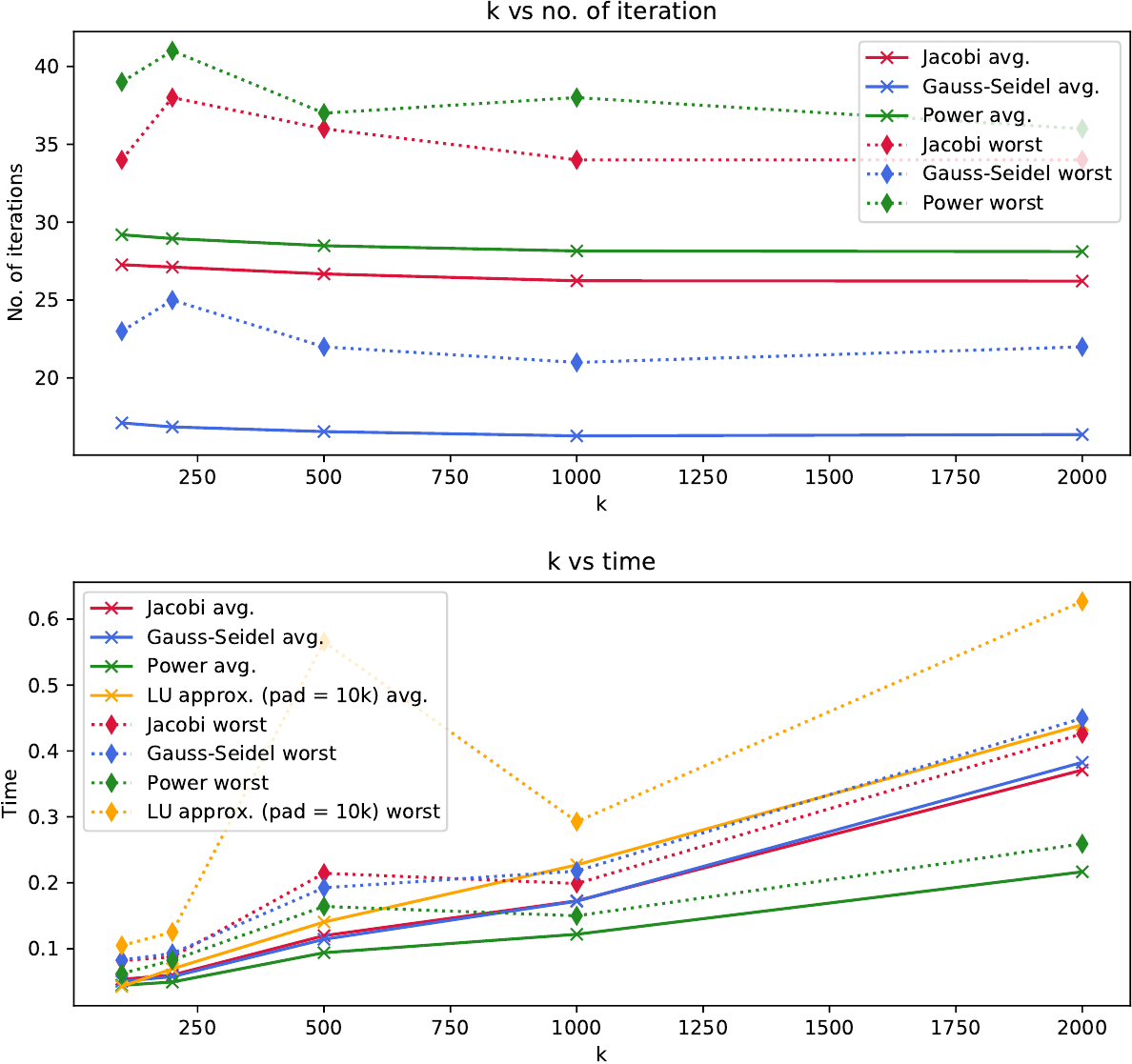}
        	\caption{The LU approx. method is the worst performer here in terms of (average and worst) time. Each graph (top and bottom) was achieved after performing ($200 \times 5$) experiments.}
        	\label{fig:scalabilityiterativemethodsn-crop}
        \end{figure}

        \paragraph{} In Figure \ref{fig:scalabilityiterativemethodss-crop} we plot the results for when $|S_?|$ is varied. As expected, the Gauss-Seidel algorithm requires the fewest iterations on average. However, due to the time complexity of our deconvolution algorithm and the lack of an optimized implementation, the libraries used to implement these algorithms lead to the power method being quickest on average. The LU approximation algorithm was found to be the least scalable. When we varied $k$ instead (see Fig. \ref{fig:scalabilityiterativemethodsn-crop}), similar results are obtained.

    \section{An example with the Romberg method}
    \label{sec:romberg_example}
\begin{figure}[H]
     	 \centerline{
     		 \includegraphics[width=1\linewidth]{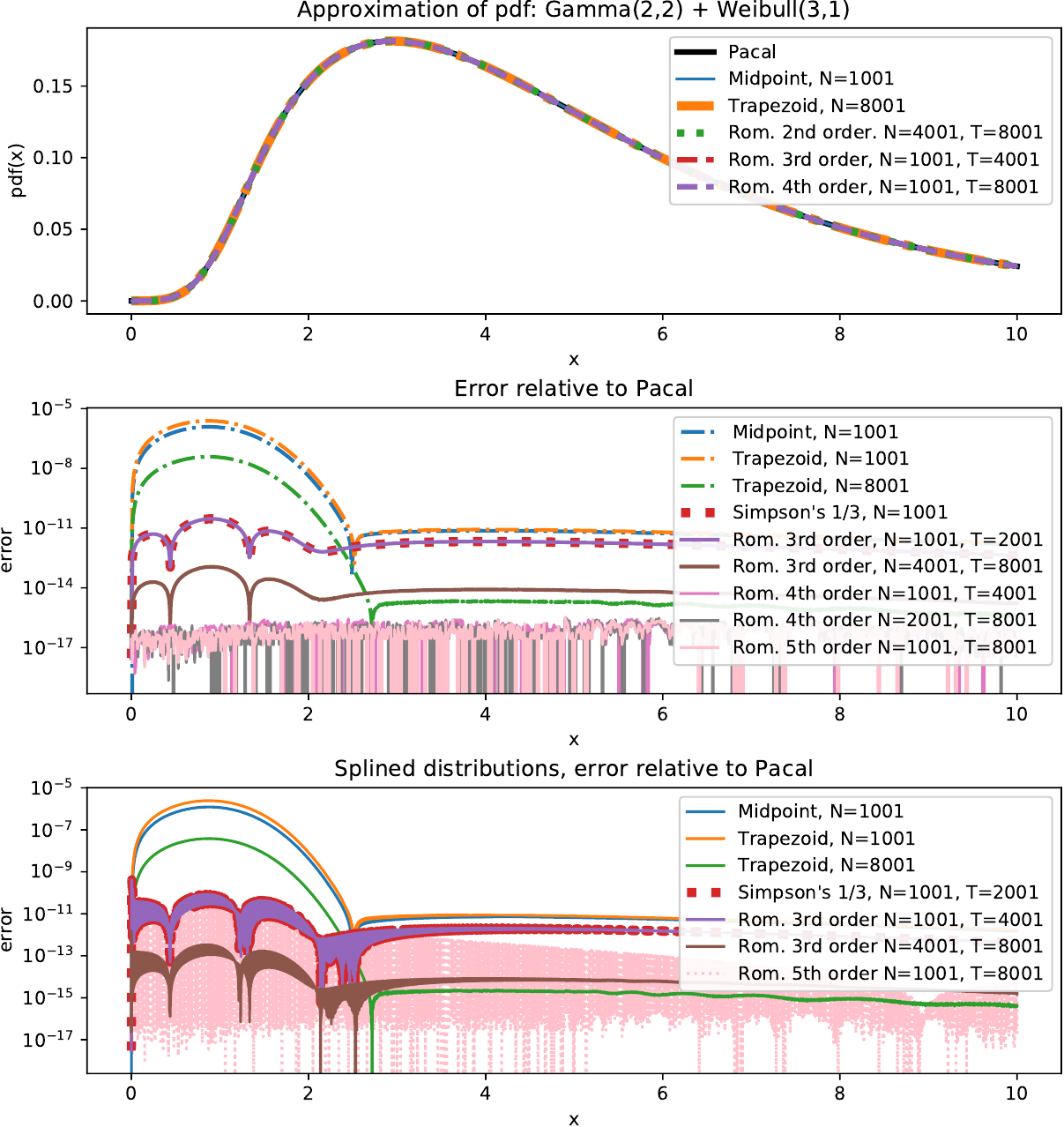}}
     	 \caption{Top: The approximation of the convolution via some quadrature rules and also Pacal. Middle: The absolute error of certain rules (not specifically those from Top) with respect to the Romberg method. Bottom: The absolute error of these rules when interpolated via a cubic spline relative to Pacal. Selected rules were chosen again for each plot.}
     	 \label{fig:rombergexperiments-crop}
    \end{figure}

     In Figure \ref{fig:rombergexperiments-crop}, we apply the Romberg method of different levels to approximate the convolution integral of a Gamma distribution with parameters ($k=2, \theta=2$), and a Weibull distribution with parameters ($k=3, \lambda=1$). We find here that Romberg's method at level three shows an equivalence to Simpson's method derived using the midpoint rule (not presented in this thesis). Additionally, the trapezoid rule with 8001 points yields weaker accuracy relative to Simpson's rule with 1001 points (but using a total of around 2002 points for computing the Midpoint and trapezoid rule).
    
    \paragraph{} The last sub-figure of Fig. \ref{fig:rombergexperiments-crop} shows what happens when we interpolate our pdf approximations. We find there to be a significant loss of accuracy with Romberg at level five and less so at level 3. The trapezoid rule and midpoint rule do not appear to be affected much.
        
    \chapter{Markov Chain Simulation}
    \label{app_chap:mc_simulation}
      We present here the four algorithms used to simulate the respective Markov chains in the experiment section of Chapter 6. They are written in python and tested to run on python 3.8. Firstly, the necessary imports are:
        
    \begin{python}
    import matplotlib.pyplot as plt
    import numpy as np
    import sparse
    import seaborn as sns
	\end{python}

	The code to generate uniform MCs is directly below.
           
    \begin{python} 
    def generate_random_MC_uniform(num_states):
    P_matrix = np.random.random((num_states + 1, num_states))  
    P_matrix += 0.01  #prevents columns summing to zero (highly improbable)
    P_matrix = (P_matrix / P_matrix.sum(0)).T
       
    b_vector = P_matrix[:, -1]  
    A_matrix = P_matrix[:, :-1] 
    
    return A_matrix, b_vector, P_matrix
   	\end{python}
	Next is the code to generate block MCs. For our experiments, we used the default parameters.
	\begin{python}                         
    def generate_random_MC_block(num_states, num_pass=200, block_scale=5):
    P_matrix = np.zeros((num_states, num_states + 1))
    full_state_idxs = np.arange(num_states)
    
    # for each state, create the probability that it can reach B
    reach_idxs = np.random.choice(full_state_idxs, int(num_states * np.random.rand()), replace=False)
    reach_probs = np.random.rand(len(reach_idxs))
    P_matrix[reach_idxs, -1] = reach_probs
    
    # perform n passes
    for _ in range(num_pass):
    block_size = int((np.random.rand() * num_states) / (2 * block_scale))
    initial_idx = np.maximum((np.random.rand(2) * num_states).astype(int) - block_size, [0, 0])
    P_matrix[initial_idx[0]:initial_idx[0] + 2 * block_size,
    initial_idx[1]:initial_idx[1] + 2 * block_size] += 1
    
    row_sums = P_matrix.sum(0)  # for each col, sum up all the rows
    non_zero_idxs = row_sums.nonzero()
    P_matrix[:, non_zero_idxs] /= row_sums[non_zero_idxs]
    
    b_vector = P_matrix[:, -1] 
    A_matrix = P_matrix[:, :-1] 
    
    return A_matrix, b_vector, P_matrix
	\end{python}
     Then, we have the algorithm to generate $N$-pass MCs. Again, the default parameter is used.
    \begin{python}              
    def generate_random_MC_npass(num_states, num_pass=1000):
    P_matrix = np.zeros((num_states, num_states + 1))
    full_state_idxs = np.arange(num_states)
    
    # for each state, create probabilty that it can reach B
    reach_idxs = np.random.choice(full_state_idxs, np.maximum(1, int(num_states * np.random.rand())), replace=False)
    reach_probs = np.random.rand(len(reach_idxs))
    P_matrix[reach_idxs, -1] = reach_probs
    
    # perform n passes
    for _ in range(num_pass):
    sel_states = np.random.permutation(full_state_idxs)
    choices = np.random.choice(reach_idxs, num_states, replace=True)
    temp = np.random.uniform(0, 1 - P_matrix[sel_states, :].sum(1), len(sel_states)) 
    P_matrix[sel_states, choices] += temp
    reach_idxs = sel_states
    
    b_vector = P_matrix[:, -1]  
    A_matrix = P_matrix[:, :-1]  
    
    return A_matrix, b_vector, P_matrix
   	\end{python}
	Lastly, we have the algorithm to sample sparse MCs. Once more, the default param. is used.
	\begin{python}
    def generate_random_MC_sparse(num_states, density=0.1):
    P_matrix = sparse.random((num_states + 1, num_states), density) 
    P_matrix = P_matrix.todense()
    full_state_idxs = np.arange(num_states)
    
    # for each state, create probabilty that it can reach B
    reach_idxs = np.random.choice(full_state_idxs, np.maximum(1, int(num_states * np.random.rand())), replace=False)
    reach_probs = np.random.rand(len(reach_idxs))
    P_matrix[reach_idxs, -1] = reach_probs    
    
    # normalize
    row_sums = P_matrix.sum(0)  
    non_zero_idxs = row_sums.nonzero()
    P_matrix[:, non_zero_idxs] /= row_sums[non_zero_idxs]
    
    P_matrix = P_matrix.T
    b_vector = P_matrix[:, -1] 
    A_matrix = P_matrix[:, :-1] 
    
    return A_matrix, b_vector, P_matrix      
    \end{python}
        
\end{document}